# SOME NEUTROSOPHIC ALGEBRAIC STRUCTURES AND NEUTROSOPHIC N-ALGEBRAIC STRUCTURES

**W. B. Vasantha Kandasamy**
e-mail: vasantha@iitm.ac.in
web: http://mat.iitm.ac.in/~wbv

**Florentin Smarandache**
e-mail: smarand@unm.edu

# CONTENTS

# PREFACE

In this book, for the first time we introduce the notion of neutrosophic algebraic structures for groups, loops, semigroups and groupoids and also their neutrosophic N-algebraic structures.

One is fully aware of the fact that many classical theorems like Lagrange, Sylow and Cauchy have been studied only in the context of finite groups. Here we try to shift the paradigm by studying and introducing these theorems to neutrosophic semigroups, neutrosophic groupoids, and neutrosophic loops. We have intentionally not given several theorems for semigroups and groupoid but have given several results with proof mainly in the case of neutrosophic loops, biloops and N-loops.

One of the reasons for this is the fact that loops are generalizations of groups and groupoids. Another feature of this book is that only meager definitions and results are given about groupoids. But over 25 problems are suggested as exercise in the last chapter. For groupoids are generalizations of both semigroups and loops.

This book has seven chapters. Chapter one provides several basic notions to make this book self-contained. Chapter two introduces neutrosophic groups and neutrosophic N-groups and gives several examples. The third chapter deals with neutrosophic semigroups and neutrosophic N-semigroups, giving several interesting results. Chapter four introduces neutrosophic loops and neutrosophic N-loops. We introduce several new, related definitions. In fact we construct a new class of neutrosophic loops using modulo integer $Z_n$, $n > 3$, where n is odd. Several properties of these structures are proved using number theoretic techniques. Chapter five just introduces the concept of neutrosophic groupoids and neutrosophic N-groupoids. Sixth chapter innovatively gives mixed neutrosophic structures and their duals. The final chapter gives problems for the interested reader to solve. Our main motivation is to attract more researchers towards algebra and its various applications.

We express our sincere thanks to Kama Kandasamy for her help in the layout and Meena for cover-design of the book. The authors express their whole-hearted gratefulness to Dr.K.Kandasamy whose invaluable support and help, and patient proofreading contributed to a great extent to the making of this book.

W.B.VASANTHA KANDASAMY
FLORENTIN SMARANDACHE
2006

## Chapter One

# INTRODUCTION

In this chapter we introduce certain basic concepts to make this book a self contained one. This chapter has 5 sections. In section one the notion of groups and N-groups are introduced. Section two just mentions about semigroups and N-semigroups. In section 3 loops and N-loops are recalled. Section 4 gives a brief description about groupoids and their properties. Section 5 recalls the mixed N algebraic structure.

## 1.1 Groups, N-group and their basic Properties

It is a well-known fact that groups are the only algebraic structures with a single binary operation that is mathematically so perfect that an introduction of a richer structure within it is impossible. Now we proceed on to define a group.

**DEFINITION 1.1.1:** *A non empty set of elements G is said to form a group if in G there is defined a binary operation, called the product and denoted by '•' such that*

i. *a, b $\in$ G implies that a • b $\in$ G (closed).*
ii. *a, b, c $\in$ G implies a • (b • c) = (a • b) • c (associative law).*
iii. *There exists an element e $\in$ G such that a • e = e • a = a for all a $\in$ G (the existence of identity element in G).*
iv. *For every a $\in$ G there exists an element $a^{-1} \in$ G such that a • $a^{-1}$ = $a^{-1}$ • a = e (the existence of inverse in G).*

**DEFINITION 1.1.2:** *A subgroup N of a group G is said to be a normal subgroup of G if for every g $\in$ G and n $\in$ N, g n $g^{-1}$ $\in$ N.*

Equivalently by $gNg^{-1}$ we mean the set of all $gng^{-1}$, n $\in$ N then N is a normal subgroup of G if and only if $gNg^{-1} \subset N$ for every g $\in$ G.

**THEOREM 1.1.1:** N is a normal subgroup of G if and only if $gNg^{-1} = N$ for every g $\in$ G.

**DEFINITION 1.1.3:** *Let G be a group. Z(G) = {x $\in$ G | gx = xg for all g $\in$ G}. Then Z(G) is called the center of the group G.*

**DEFINITION 1.1.4:** *Let G be a group, A, B be subgroups of G. If x, y $\in$ G define x ~ y if y = axb for some a $\in$ A and b $\in$ B. We call the set AxB = {axb / a $\in$ A, b $\in$ B} a double coset of A, B in G.*

**DEFINITION 1.1.5:** *Let G be a group. A and B subgroups of G, we say A and B are conjugate with each other if for some g $\in$ G, $A = gBg^{-1}$.*

*Clearly if A and B are conjugate subgroups of G then o(A) = o(B).*

**THEOREM: (LAGRANGE).** If G is a finite group and H is a subgroup of G then o(H) is a divisor of o(G).

**COROLLARY 1.1.1:** *If G is a finite group and a $\in$ G, then o(a) / o(G).*

**COROLLARY 1.1.2:** *If G is a finite group and a $\in$ G, then $a^{o(G)} = e$.*

In this section we give the two Cauchy's theorems one for abelian groups and the other for non-abelian groups. The main result on finite groups is that if the order of the group is n ($n < \infty$) if p is a prime dividing n by Cauchy's theorem we will always be able to pick up an element $a \in G$ such that $a^p = e$. In fact we can say Sylow's theorem is a partial extension of Cauchy's theorem for he says this finite group G has a subgroup of order $p^\alpha$($\alpha \geq 1$, p, a prime).

**THEOREM: (CAUCHY'S THEOREM FOR ABELIAN GROUPS).** *Suppose G is a finite abelian group and p / o(G), where p is a prime number. Then there is an element $a \neq e \in G$ such that $a^p = e$.*

**THEOREM: (CAUCHY):** *If p is a prime number and p | o(G), then G has an element of order p.*

Though one may marvel at the number of groups of varying types carrying many different properties, except for Cayley's we would not have seen them to be imbedded in the class of groups this was done by Cayley's in his famous theorem. Smarandache semigroups also has a beautiful analog for Cayley's theorem which is given by A(S) we mean the set of all one to one maps of the set S into itself. Clearly A(S) is a group having n! elements if $o(S) = n < \propto$, if S is an infinite set, A(S) has infinitely many elements.

**THEOREM: (CAYLEY)** *Every group is isomorphic to a subgroup of A(S) for some appropriate S.*

The Norwegian mathematician Peter Ludvig Mejdell Sylow was the contributor of Sylow's theorems. Sylow's theorems serve double purpose. One hand they form partial answers to the converse of Lagrange's theorem and on the other hand they are the complete extension of Cauchy's Theorem. Thus Sylow's work interlinks the works of two great mathematicians Lagrange and Cauchy. The following theorem is one, which makes use of Cauchy's theorem. It gives a nice partial converse to Lagrange's theorem and is easily understood.

**THEOREM: (SYLOW'S THEOREM FOR ABELIAN GROUPS)** *If G is an abelian group of order o(G), and if p is a prime number,*

*such that $p^{\alpha} \mid o(G)$, $p^{\alpha+1} \nmid o(G)$, then G has a subgroup of order $p^{\alpha}$.*

**COROLLARY 1.1.3:** *If G is an abelian group of finite order and $p^{\alpha} \mid o(G)$, $p^{\alpha+1} \nmid o(G)$, then there is a unique subgroup of G of order $p^{\alpha}$.*

**DEFINITION 1.1.6:** *Let G be a finite group. A subgroup G of order $p^{\alpha}$, where $p^{\alpha} / o(G)$ but $p^{\alpha} \nmid o(G)$, is called a p-Sylow subgroup of G. Thus we see that for any finite group G if p is any prime which divides o(G); then G has a p-Sylow subgroup.*

**THEOREM (FIRST PART OF SYLOW'S THEOREM):** *If p is a prime number and $p^{\alpha}/ o(G)$ and $p^{\alpha+1} \nmid o(G)$, G is a finite group, then G has a subgroup of order $p^{\alpha}$.*

**THEOREM: (SECOND PART OF SYLOW'S THEOREM):** *If G is a finite group, p a prime and $p^{n} \mid o(G)$ but $p^{n+1} \nmid o(G)$, then any two subgroup of G of order $p^{n}$ are conjugate.*

**THEOREM: (THIRD PART OF SYLOW'S THEOREM):** *The number of p-Sylow subgroups in G, for a given prime, is of the form 1 + kp.*

**DEFINITION 1.1.7:** *Let $\{G, *_1, \ldots, *_N\}$ be a non empty set with N binary operations. $\{G, *_1, \ldots, *_N\}$ is called a N-group if there exists N proper subsets $G_1, \ldots, G_N$ of G such that*

i. $G = G_1 \cup G_2 \ldots \cup G_N$.
ii. $(G_i, *_i)$ *is a group for i = 1, 2, …, N.*

*We say proper subset of G if $G_i \not\subseteq G_j$ and $G_j \not\subseteq G_i$ if $i \neq j$ for $1 \leq i, j \leq N$. When N = 2 this definition reduces to the definition of bigroup.*

**DEFINITION 1.1.8:** *Let $\{G, *_1, \ldots, *_N\}$ be a N-group. A subset H $(\neq \phi)$ of a N-group $(G, *_1, \ldots, *_N)$ is called a sub N-group if H* *itself is a N-group under $*_1, *_2, \ldots, *_N$, binary operations defined on G.*

**THEOREM 1.1.2:** *Let $(G, *_1, \ldots, *_N)$ be a N-group. The subset $H \neq \phi$ of a N-group G is a sub N-group then $(H, *_i)$ in general are not groups for i =1, 2, ..., N.*

**DEFINITION 1.1.9:** *Let $(G, *_1, \ldots, *_N)$ be a N-group where $G = G_1 \cup G_2 \cup \ldots \cup G_N$. Let $(H, *_1, \ldots, *_N)$ be a sub N-group of $(G, *_1, \ldots, *_N)$ where $H = H_1 \cup H_2 \cup \ldots \cup H_N$ we say $(H, *_1, \ldots, *_N)$ is a normal sub N-group of $(G, *_1, \ldots, *_N)$ if each $H_i$ is a normal subgroup of $G_i$ for i = 1, 2,..., N.*

*Even if one of the subgroups $H_i$ happens to be non normal subgroup of $G_i$ still we do not call H a normal sub-N-group of the N-group G.*

**DEFINITION 1.1.10:** *Let $(G = G_1 \cup G_2 \cup \ldots \cup G_N, *_1, *_2, \ldots, *_N)$ and $(K = K_1 \cup K_2 \cup \ldots \cup K_N, *_1, \ldots, *_N)$ be any two N- groups. We say a map $\phi : G \to K$ to be a N-group homomorphism if $\phi \mid G_i$ is a group homomorphism from $G_i$ to $K_i$ for i = 1, 2,..., N. i.e. $\phi \big|_{G_i} : G_i \to K_i$ is a group homomorphism of the group $G_i$ to the group $K_i$; for i = 1, 2, ..., N.*

## 1.2 Semigroups and N-semigroups

In this section we just recall the notion of semigroup, bisemigroup and N-semigroups. Also the notion of symmetric semigroups. For more refer [49-50].

**DEFINITION 1.2.1:** *Let (S, o) be a non empty set S with a closed, associative binary operation 'o' on S. (S, o) is called the semigroup i.e., for a, b $\in$ S, a o b $\in$ S.*

**DEFINITION 1.2.2:** *Let S(n) denote the set of all mappings of (1, 2, ..., n) to itself S(n) under the composition of mappings is a semigroup. We call S(n) the symmetric semigroup of order $n^n$.*

**DEFINITION 1.2.3:** *Let* $(S = S_1 \cup S_2, *, o)$ *be a non empty set with two binary operations * and o S is a bisemigroup if*

i. $S = S_1 \cup S_2$, *$S_1$ and $S_2$ are proper subsets of S.*
ii. $(S_1, *)$ *is a semigroup.*
iii. $(S_2, o)$ *is a semigroup.*

More about bisemigroups can be had from [48-50]. Now we proceed onto define N-semigroups.

**DEFINITION 1.2.4:** *Let* $S = (S_1 \cup S_2 \cup \ldots \cup S_N, *_1, *_2, \ldots, *_N)$ *be a non empty set with N binary operations. S is a N-semigroup if the following conditions are true.*

i. $S = S_1 \cup S_2 \cup \ldots \cup S_N$ *is such that each $S_i$ is a proper subset of S.*
ii. $(S_i, *_i)$ *is a semigroup for 1, 2, …, N.*

We just give an example.

***Example 1.2.1:*** Let $S = \{S_1 \cup S_2 \cup S_3 \cup S_4, *_1, *_2, *_3, *_4\}$ where

$S_1$ = $Z_{12}$, semigroup under multiplication modulo 12,
$S_2$ = S(4), symmetric semigroup,
$S_3$ = Z semigroup under multiplication and
$S_4$ = $\left\{ \begin{pmatrix} a & b \\ c & d \end{pmatrix} \middle| a, b, c, d \in Z_{10} \right\}$ under matrix multiplication.

S is a 4-semigroup.

## 1.3 Loops and N-loops

We at this juncture like to express that books solely on loops are meager or absent as, R.H.Bruck deals with loops on his book "*A Survey of Binary Systems*", that too published as early as 1958, [3]. Other two books are on "*Quasigroups and Loops*" one by H.O. Pflugfelder, 1990 which is introductory and the other book co-edited by Orin Chein, H.O. Pflugfelder and J.D. Smith in 1990 [25]. So we felt it important to recall almost all the properties and definitions related with loops [3, 47]. We just recall a few of the properties about loops which will make this book a self contained one.

**DEFINITION 1.3.1:** *A non-empty set L is said to form a loop, if on L is defined a binary operation called the product denoted by '$\bullet$' such that*

i. *For all a, b $\in$ L we have a $\bullet$ b $\in$ L (closure property).*
ii. *There exists an element e $\in$ L such that a $\bullet$ e = e $\bullet$ a = a for all a $\in$ L (e is called the identity element of L).*
iii. *For every ordered pair (a, b) $\in$ L $\times$ L there exists a unique pair (x, y) in L such that ax = b and ya = b.*

**DEFINITION 1.3.2:** *Let L be a loop. A non-empty subset H of L is called a subloop of L if H itself is a loop under the operation of L.*

**DEFINITION 1.3.3:** *Let L be a loop. A subloop H of L is said to be a normal subloop of L, if*

i. *xH = Hx.*
ii. *(Hx)y = H(xy).*
iii. *y(xH) = (yx)H*

*for all x, y $\in$ L.*

**DEFINITION 1.3.4:** *A loop L is said to be a simple loop if it does not contain any non-trivial normal subloop.*

**DEFINITION 1.3.5:** *The commutator subloop of a loop L denoted by L' is the subloop generated by all of its commutators, that is, $\langle\{x \in L / x = (y, z)$ for some $y, z \in L\}\rangle$ where for $A \subseteq L$, $\langle A\rangle$ denotes the subloop generated by A.*

**DEFINITION 1.3.6:** *If x, y and z are elements of a loop L an associator (x, y, z) is defined by, (xy)z = (x(yz)) (x, y, z).*

**DEFINITION 1.3.7:** *The associator subloop of a loop L (denoted by A(L)) is the subloop generated by all of its associators, that is* $\langle\{x \in L / x = (a, b, c)$ *for some* $a, b, c \in L\}\rangle$.

**DEFINITION 1.3.8:** *The centre Z(L) of a loop L is the intersection of the nucleus and the Moufang centre, that is* $Z(L) = C(L) \cap N(L)$.

**DEFINITION [35]**: *A normal subloop of a loop L is any subloop of L which is the kernel of some homomorphism from L to a loop.*

Further Pflugfelder [25] has proved the central subgroup Z(L) of a loop L is normal in L.

**DEFINITION [35]**: *Let L be a loop. The centrally derived subloop (or normal commutator- associator subloop) of L is defined to be the smallest normal subloop* $L' \subset L$ *such that* $L / L'$ *is an abelian group.*

*Similarly nuclearly derived subloop (or normal associator subloop) of L is defined to be the smallest normal subloop* $L_1 \subset L$ *such that* $L / L_1$ *is a group. Bruck proves* $L'$ *and* $L_1$ *are well defined.*

**DEFINITION [35]**: *The Frattini subloop* $\phi(L)$ *of a loop L is defined to be the set of all non-generators of L, that is the set of all* $x \in L$ *such that for any subset S of L,* $L = \langle x, S\rangle$ *implies* $L = \langle S\rangle$. *Bruck has proved as stated by Tim Hsu* $\phi(L) \subset L$ *and* $L / \phi(L)$ *is isomorphic to a subgroup of the direct product of groups of prime order.*

**DEFINITION [22]**: *Let L be a loop. The commutant of L is the set* $(L) = \{a \in L / ax = xa \ \forall x \in L\}$. *The centre of L is the set of all* $a \in C(L)$ *such that* $a \bullet xy = ax \bullet y = x \bullet ay = xa \bullet y$ *and* $xy \bullet a = x \bullet ya$ *for all* $x, y \in L$. *The centre is a normal subloop. The commutant is also known as Moufang Centre in literature.*

**DEFINITION [23]:** *A left loop (B, •) is a set B together with a binary operation '•' such that (i) for each a $\in$ B, the mapping x $\rightarrow$ a • x is a bijection and (ii) there exists a two sided identity 1 $\in$ B satisfying 1 • x = x • 1 = x for every x $\in$ B. A right loop is defined similarly. A loop is both a right loop and a left loop.*

**DEFINITION [11]** : *A loop L is said to have the weak Lagrange property if, for each subloop K of L, |K| divides |L|. It has the strong Lagrange property if every subloop K of L has the weak Lagrange property.*

**DEFINITION 1.3.9:** *A loop L is said to be power-associative in the sense that every element of L generates an abelian group.*

**DEFINITION 1.3.10:** *A loop L is diassociative loop if every pair of elements of L generates a subgroup.*

**DEFINITION 1.3.11:** *A loop L is said to be a Moufang loop if it satisfies any one of the following identities:*

i. *(xy) (zx) = (x(yz))x*
ii. *((xy)z)y = x(y(zy))*
iii. *x(y(xz) = ((xy)x)z*

*for all x, y, z $\in$ L.*

**DEFINITION 1.3.12:** *Let L be a loop, L is called a Bruck loop if x(yx)z = x(y(xz)) and $(xy)^{-1} = x^{-1}y^{-1}$ for all x, y, z $\in$ L.*

**DEFINITION 1.3.13:** *A loop (L, •) is called a Bol loop if ((xy)z)y = x((yz)y) for all x, y, z $\in$ L.*

**DEFINITION 1.3.14:** *A loop L is said to be right alternative if (xy)y = x(yy) for all x, y $\in$ L and L is left alternative if (xx)y = x(xy) for all x, y $\in$ L. L is said to be an alternative loop if it is both a right and left alternative loop.*

**DEFINITION 1.3.15:** *A loop (L, •) is called a weak inverse property loop (WIP-loop) if (xy)z = e imply x(yz) = e for all x, y, z $\in$ L.*

**DEFINITION 1.3.16:** *A loop L is said to be semi alternative if (x, y, z) = (y, z, x) for all x, y, z ∈ L, where (x, y, z) denotes the associator of elements x, y, z ∈ L.*

**THEOREM (MOUFANG'S THEOREM):** *Every Moufang loop G is diassociative more generally, if a, b, c are elements in a Moufang loop G such that (ab)c = a(bc) then a, b, c generate an associative loop.*

The proof is left for the reader; for assistance refer Bruck R.H. [3].

**DEFINITION 1.3.17:** *Let L be a loop, L is said to be a two unique product loop (t.u.p) if given any two non-empty finite subsets A and B of L with |A| + |B| > 2 there exist at least two distinct elements x and y of L that have unique representation in the from x = ab and y = cd with a, c ∈ A and b, d ∈ B.*

*A loop L is called a unique product (u.p) loop if, given A and B two non-empty finite subsets of L, then there always exists at least one x ∈ L which has a unique representation in the from x = ab, with a ∈ A and b ∈ B.*

**DEFINITION 1.3.18:** *Let (L, •) be a loop. The principal isotope (L, ∗) of (L, •) with respect to any predetermined a, b ∈ L is defined by x ∗ y = XY, for all x, y ∈ L, where Xa = x and bY = y for some X, Y ∈ L.*

**DEFINITION 1.3.19:** *Let L be a loop, L is said to be a G-loop if it is isomorphic to all of its principal isotopes.*

The main objective of this section is the introduction of a new class of loops with a natural simple operation. As to introduce loops several functions or maps are defined satisfying some desired conditions we felt that it would be nice if we can get a natural class of loops built using integers.

Here we define the new class of loops of any even order, they are distinctly different from the loops constructed by other researchers. Here we enumerate several of the properties enjoyed by these loops.

**DEFINITION [41]:** *Let $L_n(m) = \{e, 1, 2, ..., n\}$ be a set where $n > 3$, n is odd and m is a positive integer such that $(m, n) = 1$ and $(m - 1, n) = 1$ with $m < n$.*

*Define on $L_n(m)$ a binary operation '$\bullet$' as follows:*

i. $e \bullet i = i \bullet e = i$ *for all* $i \in L_n(m)$
ii. $i^2 = i \bullet i = e$ *for all* $i \in L_n(m)$
iii. $i \bullet j = t$ *where* $t = (mj - (m - 1)i) \pmod n$

*for all $i, j \in L_n(m)$; $i \neq j$, $i \neq e$ and $j \neq e$, then $L_n(m)$ is a loop under the binary operation '$\bullet$'.*

***Example 1.3.1:*** Consider the loop $L_5(2) = \{e, 1, 2, 3, 4, 5\}$. The composition table for $L_5(2)$ is given below:

| $\bullet$ | e | 1 | 2 | 3 | 4 | 5 |
|---|---|---|---|---|---|---|
| e | e | 1 | 2 | 3 | 4 | 5 |
| 1 | 1 | e | 3 | 5 | 2 | 4 |
| 2 | 2 | 5 | e | 4 | 1 | 3 |
| 3 | 3 | 4 | 1 | e | 5 | 2 |
| 4 | 4 | 3 | 5 | 2 | e | 1 |
| 5 | 5 | 2 | 4 | 1 | 3 | e |

This loop is of order 6 which is both non-associative and non-commutative.

*Physical interpretation of the operation in the loop $L_n(m)$*:

We give a physical interpretation of this class of loops as follows: Let $L_n(m)= \{e, 1, 2, \ldots , n\}$ be a loop in this identity element of the loop are equidistantly placed on a circle with e as its centre.

We assume the elements to move always in the clockwise direction.

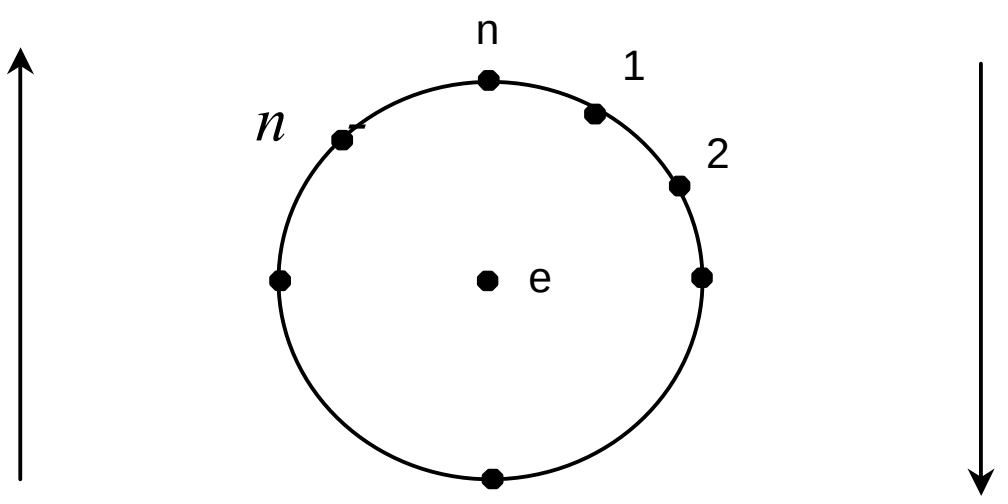

Let i, j ∈ $L_n(m)$ (i ≠ j, i ≠ e, j ≠ e). If j is the $r^{th}$ element from i counting in the clockwise direction the i • j will be the $t^{th}$ element from j in the clockwise direction where t = (m –1)r. We see that in general i • j need not be equal to j • i. When i = j we define $i^2$ = e and i • e = e • i = i for all i ∈ $L_n(m)$ and e acts as the identity in $L_n(m)$.

***Example 1.3.2***: Now the loop $L_7(4)$ is given by the following table:

| • | e | 1 | 2 | 3 | 4 | 5 | 6 | 7 |
|---|---|---|---|---|---|---|---|---|
| e | e | 1 | 2 | 3 | 4 | 5 | 6 | 7 |
| 1 | 1 | e | 5 | 2 | 6 | 3 | 7 | 4 |
| 2 | 2 | 5 | e | 6 | 3 | 7 | 4 | 1 |
| 3 | 3 | 2 | 6 | e | 7 | 4 | 1 | 5 |
| 4 | 4 | 6 | 3 | 7 | e | 1 | 5 | 2 |
| 5 | 5 | 3 | 7 | 4 | 1 | e | 2 | 6 |
| 6 | 6 | 7 | 4 | 1 | 5 | 2 | e | 3 |
| 7 | 7 | 4 | 1 | 5 | 2 | 6 | 3 | e |

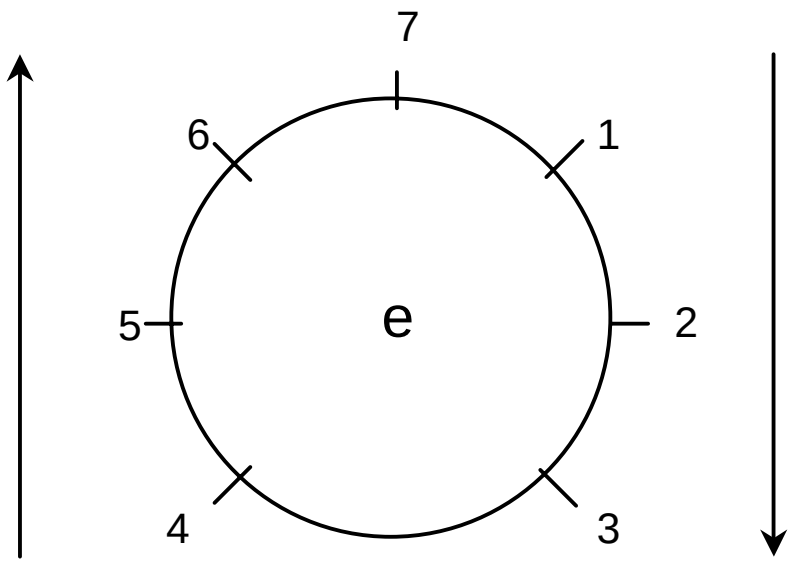

Let 2, 4 ∈ $L_7(4)$. 4 is the $2^{nd}$ element from 2 in the clockwise direction. So 2.4 will be (4 – 1)2 that is the $6^{th}$ element from 4 in the clockwise direction which is 3. Hence 2.4 = 3.

**Notation**: Let $L_n$ denote the class of loops. $L_n(m)$ for fixed n and various m's satisfying the conditions m < n, (m, n) = 1 and (m – 1, n) = 1, that is $L_n = \{L_n(m) \mid n > 3, n \text{ odd}, m < n, (m, n) = 1 \text{ and } (m-1, n) = 1\}$.

***Example 1.3.3:*** Let n = 5. The class $L_5$ contains three loops; viz. $L_5(2)$, $L_5(3)$ and $L_5(4)$ given by the following tables:

$L_5(2)$

| • | e | 1 | 2 | 3 | 4 | 5 |
|---|---|---|---|---|---|---|
| e | e | 1 | 2 | 3 | 4 | 5 |
| 1 | 1 | e | 3 | 5 | 2 | 4 |
| 2 | 2 | 5 | e | 4 | 1 | 3 |
| 3 | 3 | 4 | 1 | e | 5 | 2 |
| 4 | 4 | 3 | 5 | 2 | e | 1 |
| 5 | 5 | 2 | 4 | 1 | 3 | e |

$L_5(3)$

| • | e | 1 | 2 | 3 | 4 | 5 |
|---|---|---|---|---|---|---|
| e | e | 1 | 2 | 3 | 4 | 5 |
| 1 | 1 | e | 4 | 2 | 5 | 3 |
| 2 | 2 | 4 | e | 5 | 3 | 1 |
| 3 | 3 | 2 | 5 | e | 1 | 4 |
| 4 | 4 | 5 | 3 | 1 | e | 2 |
| 5 | 5 | 3 | 1 | 4 | 2 | e |

$L_5(4)$

| • | e | 1 | 2 | 3 | 4 | 5 |
|---|---|---|---|---|---|---|
| e | e | 1 | 2 | 3 | 4 | 5 |
| 1 | 1 | e | 5 | 4 | 3 | 2 |
| 2 | 2 | 3 | e | 1 | 5 | 4 |
| 3 | 3 | 5 | 4 | e | 2 | 1 |
| 4 | 4 | 2 | 1 | 5 | e | 3 |
| 5 | 5 | 4 | 3 | 2 | 1 | e |

**THEOREM [27]:** *Let $L_n$ be the class of loops for any $n > 3$, if $n = p_1^{\alpha_1} p_2^{\alpha_2} \ldots p_k^{\alpha_k}$ ($\alpha_i > 1$, for $i = 1, 2, \ldots, k$), then $|L_n| = \prod_{i=1}^{k} (p_i - 2)\, p_i^{\alpha_i - 1}$ where $|L_n|$ denotes the number of loops in $L_n$.*

The proof is left for the reader as an exercise.

**THEOREM [27]**: *$L_n$ contains one and only one commutative loop. This happens when $m = (n + 1) / 2$. Clearly for this m, we have $(m, n) = 1$ and $(m – 1, n) = 1$.*

It can be easily verified by using simple number theoretic techniques.

**THEOREM [27]**: *Let $L_n$ be the class of loops. If $n = p_1^{\alpha_1} p_2^{\alpha_2} \ldots p_k^{\alpha_k}$, then $L_n$ contains exactly $F_n$ loops which are strictly non-commutative where $F_n = \prod_{i=1}^{k} (p_i - 3)\, p_i^{\alpha_i - 1}$.*

The proof is left for the reader as an exercise.

<u>Note</u>: If n = p where p is a prime greater than or equal to 5 then in $L_n$ a loop is either commutative or strictly non-commutative. Further it is interesting to note if n = 3t then the class $L_n$ does not contain any strictly non-commutative loop.

**THEOREM [32]**: *The class of loops $L_n$ contains exactly one left alternative loop and one right alternative loop but does not contain any alternative loop.*

*Proof*: We see $L_n(2)$ is the only right alternative loop that is when m = 2 (Left for the reader to prove using simple number theoretic techniques). When m = n –1 that is $L_n(n – 1)$ is the only left alternative loop in the class of loops $L_n$.

From this it is impossible to find a loop in $L_n$, which is simultaneously right alternative and left alternative. Further it is clear from earlier result both the right alternative loop and the left alternative loop is not commutative.

**THEOREM [27]:** *Let $L_n$ be the class of loops:*

i. *$L_n$ does not contain any Moufang loop.*
ii. *$L_n$ does not contain any Bol loop.*
iii. *$L_n$ does not contain any Bruck loop.*

The reader is requested to prove these results using number theoretic techniques.

**THEOREM [41]**: *Let $L_n(m) \in L_n$. Then $L_n(m)$ is a weak inverse property (WIP) loop if and only if $(m^2 - m + 1) \equiv 0 (\mathrm{mod}\ n)$.*

*Proof*: It is easily checked that for a loop to be a WIP loop we have "if (xy)z = e then x(yz) = e where x, y, z ∈ L." Both way conditions can be derived using the defining operation on the loop $L_n(m)$.

***Example 1.3.4:*** L be the loop $L_7(3)$ = {e, 1, 2, 3, 4, 5, 6, 7} be in $L_7$ given by the following table:

| • | e | 1 | 2 | 3 | 4 | 5 | 6 | 7 |
|---|---|---|---|---|---|---|---|---|
| e | e | 1 | 2 | 3 | 4 | 5 | 6 | 7 |
| 1 | 1 | e | 4 | 7 | 3 | 6 | 2 | 5 |
| 2 | 2 | 6 | e | 5 | 1 | 4 | 7 | 3 |
| 3 | 3 | 4 | 7 | e | 6 | 2 | 5 | 1 |
| 4 | 4 | 2 | 5 | 1 | e | 7 | 3 | 6 |
| 5 | 5 | 7 | 3 | 6 | 2 | e | 1 | 4 |
| 6 | 6 | 5 | 1 | 4 | 7 | 3 | e | 2 |
| 7 | 7 | 3 | 6 | 2 | 5 | 1 | 4 | e |

It is easily verified $L_7(3)$ is a WIP loop. One way is easy for $(m^2 - m + 1) \equiv 0 \ (\mathrm{mod}\ 7)$ that is $9 - 3 + 1 = 9 + 4 + 1 \equiv 0 (\mathrm{mod}\ 7)$. It is interesting to note that no loop in the class $L_n$ contain any associative loop.

**THEOREM [27]**: *Let $L_n$ be the class of loops. The number of strictly non-right (left) alternative loops is $P_n$ where* $P_n = \prod_{i=1}^{k} (p_i - 3) p_i^{\alpha_i - 1}$ *and* $n = \prod_{i=1}^{k} p_i^{\alpha_i}$.

The proof is left for the reader to verify.

Now we proceed on to study the associator and the commutator of the loops in $L_n$.

**THEOREM [27]**: *Let $L_n(m) \in L_n$. The associator $A(L_n(m)) = L_n(m)$.*

For more literature about the new class of loops refer [41, 47].

**DEFINITION 1.3.20:** *Let $(L, *_1, \ldots, *_N)$ be a non empty set with N binary operations $*_i$. L is said to be a N loop if L satisfies the following conditions:*

1. *$L = L_1 \cup L_2 \cup \ldots \cup L_N$ where each $L_i$ is a proper subset of L; i.e., $L_i \not\subseteq L_j$ or $L_j \not\subseteq L_i$ if $i \neq j$ for $1 \leq i, j \leq N$.*
2. *$(L_i, *_i)$ is a loop for some i, $1 \leq i \leq N$.*
3. *$(L_j, *_j)$ is a loop or a group for some j, $1 \leq j \leq N$.*

*For a N-loop we demand atleast one $(L_j, *_j)$ to be a loop.*

**DEFINITION 1.3.21:** *Let $(L = L_1 \cup L_2 \cup \ldots \cup L_N, *_1, \ldots, *_N)$ be a N-loop. L is said to be a commutative N-loop if each $(L_i, *_i)$ is commutative, i = 1, 2, …, N. We say L is inner commutative if each of its proper subset which is N-loop under the binary operations of L are commutative.*

**DEFINITION 1.3.22:** *Let $L = \{L_1 \cup L_2 \cup \ldots \cup L_N, *_1, *_2, \ldots, *_N\}$ be a N-loop. We say L is a Moufang N-loop if all the loops $(L_i, *_i)$ satisfy the following identities.*

1. *$(xy)(zx) = (x(yz))x$*
2. *$((xy)z)y = x(y(zy))$*

*iii.* $x(y(xz)) = ((xy)x)z$

*for all* $x, y, z \in L_i$, $i = 1, 2, \ldots, N$.

Now we proceed on to define a Bruck N-loop.

**DEFINITION 1.3.23:** *Let* $L = (L_1 \cup L_2 \cup \ldots \cup L_N, *_1, \ldots, *_N)$ *be a N-loop. We call L a Bruck N-loop if all the loops* $(L_i, *_i)$ *satisfy the identities*

*i.* $(x(yz))z = x(y(xz))$
*ii.* $(xy)^{-1} = x^{-1}y^{-1}$

*for all* $x, y \in L_i$, $i = 1, 2, \ldots, N$.

**DEFINITION 1.3.24:** *Let* $L = (L_1 \cup L_2 \cup \ldots \cup L_N, *_1, \ldots, *_N)$ *be a N-loop. A non empty subset P of L is said to be a sub N-loop, if P is a N-loop under the operations of L i.e.,* $P = \{P_1 \cup P_2 \cup P_3 \cup \ldots \cup P_N, *_1, \ldots, *_N\}$ *with each* $\{P_i, *_i\}$ *is a loop or a group.*

**DEFINITION 1.3.25:** *Let* $L = \{L_1 \cup L_2 \cup \ldots \cup L_N, *_1, \ldots, *_N\}$ *be a N-loop. A proper subset P of L is said to be a normal sub N-loop of L if*

*i.* *P is a sub N-loop of L.*
*ii.* $x_i P_i = P_i x_i$ *for all* $x_i \in L_i$.
*iii.* $y_i (x_i P_i) = (y_i x_i) P_i$ *for all* $x_i y_i \in L_i$.

*A N-loop is said to be a simple N-loop if L has no proper normal sub N-loop.*

Now we proceed on to define the notion of Moufang center.

**DEFINITION 1.3.26:** *Let* $L = \{L_1 \cup L_2 \cup \ldots \cup L_N, *_1, \ldots, *_N\}$ *be a N-loop. We say* $C_N(L)$ *is the Moufang N-centre of this N-loop if* $C_N(L) = C_1(L_1) \cup C_2(L_2) \cup \ldots \cup C_N(L_N)$ *where* $C_i(L_i) = \{x_i \in L_i / x_i y_i = y_i x_i \text{ for all } y_i \in L_i\}$, $i = 1, 2, \ldots, N$.

**DEFINITION 1.3.27:** *Let L and P to two N-loops i.e.* $L = \{L_1 \cup L_2 \cup \ldots \cup L_N, *_1, \ldots, *_N\}$ *and* $P = \{P_1 \cup P_2 \cup \ldots \cup P_N, o_1, \ldots, o_N\}$. *We say a map* $\theta : L \to P$ *is a N-loop homomorphism if* $\theta = \theta_1 \cup \theta_2 \cup \ldots \cup \theta_N$ *'*$\cup$*' is just a symbol and* $\theta_i$ *is a loop homomorphism from* $L_i$ *to* $P_i$ *for each* $i = 1, 2, \ldots, N$.

**DEFINITION 1.3.28:** *Let* $L = \{L_1 \cup L_2 \cup \ldots \cup L_N, *_1, \ldots, *_N\}$ *be a N-loop. We say L is weak Moufang N-loop if there exists atleast a loop* $(L_i, *_i)$ *such that* $L_i$ *is a Moufang loop.*

*Note:* $L_i$ should not be a group it should be only a loop.

**DEFINITION 1.3.29:** *Let* $L = \{L_1 \cup L_2 \cup \ldots \cup L_N, *_1, \ldots, *_N\}$ *be a N-loop. If x and y* $\in L$ *are elements of* $L_i$ *the N-commutator (x, y) is defined as* $xy = (yx)(x, y)$, $1 \leq i \leq N$.

**DEFINITION 1.3.30:** *Let* $L = \{L_1 \cup L_2 \cup \ldots \cup L_N, *_1, \ldots, *_N\}$ *be a N-loop. If x, y, z are elements of the N-loop L, an associator (x, y, z) is defined only if* $x, y, z \in L_i$ *for some i* $(1 \leq i \leq N)$ *and is defined to be* $(xy)z = (x(yz))(x, y, z)$.

**DEFINITION 1.3.31:** *Let* $L = \{L_1 \cup L_2 \cup \ldots \cup L_N, *_1, *_2, \ldots, *_N\}$ *be a N-loop of finite order. For* $\alpha_i \in L_i$ *define* $R_{\alpha_i}$ *as a permutation of the loop* $L_i$, $R_{\alpha_i} : x_i \to x_i\alpha_i$. *This is true for i = 1, 2,…, N we define the set*

$$\left\{R_{\alpha_1} \cup R_{\alpha_2} \cup \ldots \cup R_{\alpha_N} \,\middle|\, \alpha_i \in L_i; i = 1, 2, \ldots, N\right\}$$

*as the right regular N-representation of the N loop L.*

**DEFINITION 1.3.32:** *Let* $L = \{L_1 \cup L_2 \cup \ldots \cup L_N, *_1, \ldots, *_N\}$ *be a N-loop. For any pre determined pair* $a_i, b_i \in L_i$, $i \in \{1, 2, \ldots, N\}$ *a principal isotope* $(L, o_1, \ldots, o_N)$, *of the N loop L is defined by* $x_i \, o_i \, y_i = X_i *_i Y_i$ *where* $X_i + a_i = x_i$ *and* $b_i + Y_i = y_i$, $i = 1, 2, \ldots, N$. *L is called G-N-loop if it is isomorphic to all of its principal isotopes.*

## 1.4 Groupoids and N-groupoids

In this section we just recall the notion of groupoids. We also give some new classes of groupoids constructed using the set of modulo integers. This book uses in several examples the groupoids from these new classes of groupoids For more about groupoids please refer [45].

**DEFINITION 1.4.1:** *Given an arbitrary set P a mapping of $P \times P$ into P is called a binary operation on P. Given such a mapping $\sigma : P \times P \to P$ we use it to define a product $*$ in P by declaring $a * b = c$ if $\sigma(a, b) = c$.*

**DEFINITION 1.4.2:** *A non empty set of elements G is said to form a groupoid if in G is defined a binary operation called the product denoted by $*$ such that $a * b \in G$ for all a, b $\in$ G.*

**DEFINITION 1.4.3:** *A groupoid G is said to be a commutative groupoid if for every a, b $\in$ G we have $a * b = b * a$.*

**DEFINITION 1.4.4:** *A groupoid G is said to have an identity element e in G if $a * e = e * a = a$ for all a $\in$ G.*

**DEFINITION 1.4.5:** *Let (G, $*$) be a groupoid a proper subset H $\subset$ G is a subgroupoid if (H, $*$) is itself a groupoid.*

**DEFINITION 1.4.6:** *A groupoid G is said to be a Moufang groupoid if it satisfies the Moufang identity $(xy)(zx) = (x(yz))x$ for all x, y, z in G.*

**DEFINITION 1.4.7:** *A groupoid G is said to be a Bol groupoid if G satisfies the Bol identity $((xy)z)y = x((yz)y)$ for all x, y, z in G.*

**DEFINITION 1.4.8:** *A groupoid G is said to be a P-groupoid if $(xy)x = x(yx)$ for all x, y $\in$ G.*

**DEFINITION 1.4.9:** *A groupoid G is said to be right alternative if it satisfies the identity $(xy)y = x(yy)$ for all x, y $\in$ G.*

*Similarly we define G to be left alternative if (xx) y = x (xy) for all x, y ∈ G.*

**DEFINITION 1.4.10:** *A groupoid G is alternative if it is both right and left alternative, simultaneously.*

**DEFINITION 1.4.11:** *Let (G, ∗) be a groupoid. A proper subset H of G is said to be a subgroupoid of G if (H, ∗) is itself a groupoid.*

**DEFINITION 1.4.12:** *A groupoid G is said to be an idempotent groupoid if* $x^2 = x$ *for all x ∈ G.*

**DEFINITION 1.4.13:** *Let G be a groupoid. P a non empty proper subset of G, P is said to be a left ideal of the groupoid G if 1) P is a subgroupoid of G and 2) For all x ∈ G and a ∈ P, xa ∈ P. One can similarly define right ideal of the groupoid G. P is called an ideal if P is simultaneously a left and a right ideal of the groupoid G.*

**DEFINITION 1.4.14:** *Let G be a groupoid A subgroupoid V of G is said to be a normal subgroupoid of G if*

i. *aV = Va*
ii. *(Vx)y = V(xy)*
iii. *y(xV) = (yx)V*

*for all x, y, a ∈ V.*

**DEFINITION 1.4.15:** *A groupoid G is said to be simple if it has no non trivial normal subgroupoids.*

**DEFINITION 1.4.16:** *A groupoid G is normal if*

i. *xG = Gx*
ii. *G(xy) = (Gx)y*
iii. *y(xG) = (yx)G for all x, y ∈ G.*

**DEFINITION 1.4.17:** *Let G be a groupoid H and K be two proper subgroupoids of G, with $H \cap K = \phi$. We say H is conjugate with K if there exists a $x \in H$ such that $H = x K$ or $Kx$ ('or' in the mutually exclusive sense).*

**DEFINITION 1.4.18:** *Let $(G_1, \theta_1)$, $(G_2, \theta_2)$, ... , $(G_n, \theta_n)$ be n groupoids with $\theta_i$ binary operations defined on each $G_i$, i = 1, 2, 3, ... , n. The direct product of $G_1$, ... , $G_n$ denoted by $G = G_1 \times \ldots \times G_n = \{(g_1, \ldots , g_n) \mid g_i \in G_i\}$ by component wise multiplication on G, G becomes a groupoid.*

*For if $g = (g_1, \ldots , g_n)$ and $h = (h_1, \ldots , h_n)$ then $g \bullet h = \{(g_1\theta_1h_1, g_2\theta_2h_2, \ldots , g_n\theta_nh_n)\}$. Clearly, $gh \in G$. Hence G is a groupoid.*

**DEFINITION 1.4.19:** *Let G be a groupoid we say an element $e \in G$ is a left identity if $ea = a$ for all $a \in G$. Similarly we can define right identity of the groupoid G, if $e \in G$ happens to be simultaneously both right and left identity we say the groupoid G has an identity.*

**DEFINITION 1.4.20:** *Let G be a groupoid. We say a in G has right zero divisor if $a * b = 0$ for some $b \neq 0$ in G and a in G has left zero divisor if $b * a = 0$. We say G has zero divisors if $a \bullet b = 0$ and $b * a = 0$ for $a, b \in G \setminus \{0\}$ A zero divisor in G can be left or right divisor.*

**DEFINITION 1.4.21:** *Let G be a groupoid. The center of the groupoid $C(G) = \{x \in G \mid ax = xa \text{ for all } a \in G\}$.*

**DEFINITION 1.4.22:** *Let G be a groupoid. We say $a, b \in G$ is a conjugate pair if $a = bx$ (or $xa$ for some $x \in G$) and $b = ay$ (or $ya$ for some $y \in G$).*

**DEFINITION 1.4.23:** *Let G be a groupoid of order n. An element a in G is said to be right conjugate with b in G if we can find $x, y \in G$ such that $a \bullet x = b$ and $b \bullet y = a$ ($x * a = b$ and $y * b = a$). Similarly, we define left conjugate.*

**DEFINITION 1.4.24:** *Let $Z^+$ be the set of integers. Define an operation $*$ on $Z^+$ by $x * y = mx + ny$ where $m, n \in Z^+$, $m < \propto$ and $n < \propto$ $(m, n) = 1$ and $m \neq n$. Clearly $\{Z^+, *, (m, n)\}$ is a groupoid denoted by $Z^+ (m, n)$. We have for varying m and n get infinite number of groupoids of infinite order denoted by $\mathbf{Z}^+$.*

Here we define a new class of groupoids denoted by Z(n) using $Z_n$ and study their various properties.

**DEFINITION 1.4.25:** *Let $Z_n = \{0, 1, 2, \ldots , n - 1\}$ $n \geq 3$. For a, b $\in Z_n \setminus \{0\}$ define a binary operation $*$ on $Z_n$ as follows. $a * b = ta + ub \pmod n$ where t, u are 2 distinct elements in $Z_n \setminus \{0\}$ and $(t, u) = 1$ here ' + ' is the usual addition of 2 integers and ' ta ' means the product of the two integers t and a. We denote this groupoid by $\{Z_n, (t, u), *\}$ or in short by $Z_n (t, u)$.*

It is interesting to note that for varying $t, u \in Z_n \setminus \{0\}$ with $(t, u) = 1$ we get a collection of groupoids for a fixed integer n. This collection of groupoids is denoted by Z(n) that is $Z(n) = \{Z_n, (t, u), * \mid$ for integers $t, u \in Z_n \setminus \{0\}$ such that $(t, u) = 1\}$. Clearly every groupoid in this class is of order n.

***Example 1.4.1:*** Using $Z_3 = \{0, 1, 2\}$. The groupoid $\{Z_3, (1, 2), *\} = (Z_3 (1, 2))$ is given by the following table:

| $*$ | 0 | 1 | 2 |
|---|---|---|---|
| 0 | 0 | 2 | 1 |
| 1 | 1 | 0 | 2 |
| 2 | 2 | 1 | 0 |

Clearly this groupoid is non associative and non commutative and its order is 3.

**THEOREM 1.4.1:** *Let $Z_n = \{0, 1, 2, \ldots , n\}$. A groupoid in Z (n) is a semigroup if and only if $t^2 \equiv t \pmod n$ and $u^2 \equiv u \pmod n$ for $t, u \in Z_n \setminus \{0\}$ and $(t, u) = 1$.*

**THEOREM 1.4.2:** *The groupoid $Z_n$ (t, u) is an idempotent groupoid if and only if $t + u \equiv 1 \pmod n$.*

**THEOREM 1.4.3:** *No groupoid in Z (n) has {0} as an ideal.*

**THEOREM 1.4.4:** *P is a left ideal of $Z_n$ (t, u) if and only if P is a right ideal of $Z_n$ (u, t).*

**THEOREM 1.4.5:** *Let $Z_n$ (t, u) be a groupoid. If n = t + u where both t and u are primes then $Z_n$ (t, u) is simple.*

**DEFINITION 1.4.26:** *Let $Z_n = \{0, 1, 2, \ldots , n-1\}$ $n \geq 3$, $n < \infty$. Define $*$ a closed binary operation on $Z_n$ as follows. For any a, $b \in Z_n$ define $a * b = at + bu \pmod n$ where (t, u) need not always be relatively prime but $t \neq u$ and $t, u \in Z_n \setminus \{0\}$.*

**THEOREM 1.4.6:** *The number of groupoids in $Z^*(n)$ is (n – 1) (n – 2).*

**THEOREM 1.4.7:** *The number of groupoids in the class Z (n) is bounded by (n – 1) (n – 2).*

**THEOREM 1.4.8:** *Let $Z_n$ (t, u) be a groupoid in $Z^*(n)$ such that (t, u) = t, n = 2m, t / 2m and t + u = 2m. Then $Z_n$ (t, u) has subgroupoids of order 2m / t or n / t.*

*Proof:* Given n is even and t + u = n so that u = n – t.
Thus $Z_n$ (t, u) = $Z_n$ (t, n – t). Now using the fact $t \cdot Z_n = \left\{0, t, 2t, 3t, \ldots, \left(\frac{n}{t} - 1\right)t\right\}$ that is t . $Z_n$ has only n / t elements and these n / t elements from a subgroupoid. Hence $Z_n$ (t, n – t) where (t, n – t) = t has only subgroupoids of order n / t.

**DEFINITION 1.4.27:** *Let $Z_n = \{0, 1, 2, \ldots , n-1\}$ $n \geq 3$, $n < \infty$. Define $*$ on $Z_n$ as $a * b = ta + ub \pmod n$ where t and $u \in Z_n \setminus \{0\}$ and t can also equal u. For a fixed n and for varying t and u* *we get a class of groupoids of order n which we denote by $Z^{**}(n)$.*

**DEFINITION 1.4.28:** *Let $Z_n = \{0, 1, 2, \ldots , n - 1\}$ $n \geq 3$, $n < \propto$. Define $*$ on $Z_n$ as follows. $a * b = ta + ub \pmod n$ where $t, u \in Z_n$. Here t or u can also be zero.*

**DEFINITION 1.4.29:** *Let $G = (G_1 \cup G_2 \cup \ldots \cup G_N;\ *_1, \ldots, *_N)$ be a non empty set with N-binary operations. G is called the N-groupoid if some of the $G_i$'s are groupoids and some of the $G_j$'s are semigroups, $i \neq j$ and $G = G_1 \cup G_2 \cup \ldots \cup G_N$ is the union of proper subsets of G.*

**DEFINITION 1.4.30:** *Let $G = (G_1 \cup G_2 \cup \ldots \cup G_N;\ *_1, *_2, \ldots, *_N)$ be a N-groupoid. The order of the N-groupoid G is the number of distinct elements in G. If the number of elements in G is finite we call G a finite N-groupoid. If the number of distinct elements in G is infinite we call G an infinite N-groupoid.*

**DEFINITION 1.4.31:** *Let $G = \{G_1 \cup G_2 \cup \ldots \cup G_N, *_1, *_2, \ldots, *_N\}$ be a N-groupoid. We say a proper subset $H = \{H_1 \cup H_2 \cup \ldots \cup H_N, *_1, *_2, \ldots, *_N\}$ of G is said to be a sub N-groupoid of G if H itself is a N-groupoid of G.*

**DEFINITION 1.4.32:** *Let $G = (G_1 \cup G_2 \cup \ldots \cup G_N, *_1, *_2, \ldots, *_N)$ be a finite N-groupoid. If the order of every sub N-groupoid H divides the order of the N-groupoid, then we say G is a Lagrange N-groupoid.*

It is very important to mention here that in general every N-groupoid need not be a Lagrange N-groupoid.
Now we define still a weaker notion of a Lagrange N-groupoid.

**DEFINITION 1.4.33:** *Let $G = (G_1 \cup G_2 \cup \ldots \cup G_N, *_1, \ldots, *_N)$ be a finite N-groupoid. If G has atleast one nontrivial sub N-groupoid H such that $o(H) / o(G)$ then we call G a weakly Lagrange N-groupoid.*

**DEFINITION 1.4.34:** *Let* $G = (G_1 \cup G_2 \cup \ldots \cup G_N, *_1, *_2, \ldots, *_N)$ *be a N groupoid. A sub N groupoid* $V = V_1 \cup V_2 \cup \ldots \cup V_N$ *of G is said to be a normal sub N-groupoid of G; if*

i. $a V_i = V_i a$, $i = 1, 2, \ldots, N$ *whenever* $a \in V_i$
ii. $V_i (x y) = (V_i x) y$, $i = 1, 2, \ldots, N$ *for* $x, y \in V_i$
iii. $y (x V_i) = (xy) V_i$, $i = 1, 2, \ldots, N$ *for* $x, y \in V_i$.

*Now we say a N-groupoid G is simple if G has no nontrivial normal sub N-groupoids.*

Now we proceed on to define the notion of N-conjugate groupoids.

**DEFINITION 1.4.35:** *Let* $G = (G_1 \cup G_2 \cup \ldots \cup G_N, *_1, \ldots, *_N)$ *be a N-groupoid. Let* $H = \{H_1 \cup \ldots \cup H_N; *_1, \ldots, *_N\}$ *and* $K = \{K_1 \cup K_2 \cup \ldots \cup K_N, *_1, \ldots, *_N\}$ *be sub N-groupoids of* $G = G_1 \cup G_2 \cup \ldots \cup G_N$*; where* $H_i$, $K_i$ *are subgroupoids of* $G_i$ $(i = 1, 2, \ldots, N)$.

*Let* $K \cap H = \phi$*. We say H is N-conjugate with K if there exists* $x_i \in H_i$ *such that* $x_i K_i = H_i$ *(or* $K_i x_i = H_i$*) for* $i = 1, 2, \ldots, N$ *'or' in the mutually exclusive sense.*

**DEFINITION 1.4.36:** *Let* $G = (G_1 \cup G_2 \cup \ldots \cup G_N, *_1, *_2, \ldots, *_N)$ *be a N-groupoid. An element x in G is said to be a zero divisor if their exists a y in G such that* $x *_i y = y *_i x = 0$ *for some i in* $\{1, 2, \ldots, N\}$.

We define N-centre of a N-groupoid G.

**DEFINITION 1.4.37:** *Let* $G = (G_1 \cup G_2 \cup \ldots \cup G_N, *_1, *_2, \ldots, *_N)$ *be a N-groupoid. The N-centre of* $(G_1 \cup G_2 \cup \ldots \cup G_N, *_1, \ldots, *_N)$ *denoted by* $NC(G) = \{x_1 \in G_1 \mid x_1 a_1 = a_1 x_1 \text{ for all } a_1 \in G_1\} \cup \{x_2 \in G_2 \mid x_2 a_2 = a_2 x_2 \text{ for all } a_2 \in G_2\} \cup \ldots \cup \{x_N \in G_N \mid x_N a_N = a_N x_N \text{ for all } a_N \in G_N\} = \bigcup_{i=1}^{N} \{x_i \in G \mid x_i a_i = a_i x_i \text{ for all } a_i \in G_i\} = NC(G)$.

## 1.5 Mixed N-algebraic Structures

In this section we proceed onto define the notion of N-group-semigroup algebraic structure and other mixed substructures and enumerate some of its properties.

**DEFINITION 1.5.1:** *Let $G = \{G_1 \cup G_2 \cup \ldots \cup G_N, *_1, \ldots, *_N\}$ where some of the $G_i$'s are groups and the rest of them are semigroups. $G_i$'s are such that $G_i \not\subseteq G_j$ or $G_j \not\subseteq G_I$ if $i \neq j$, i.e. $G_i$'s are proper subsets of G. $*_1, \ldots, *_N$ are N binary operations which are obviously are associative then we call G a N-group semigroup.*

We can also redefine this concept as follows:

**DEFINITION 1.5.2:** *Let G be a non empty set with N-binary operations $*_1, \ldots, *_N$. We call G a N-group semigroup if G satisfies the following conditions:*

1. *$G = G_1 \cup G_2 \cup \ldots \cup G_N$ such that each $G_i$ is a proper subset of G (By proper subset $G_i$ of G we mean $G_i \not\subseteq G_j$ or $G_j \not\subseteq G_i$ if $i \neq j$. This does not necessarily imply $G_i \cap G_j = \phi$).*
2. *$(G_i, *_i)$ is either a group or a semigroup, $i = 1, 2, \ldots, N$.*
3. *At least one of the $(G_i, *_i)$ is a group.*
4. *At least one of the $(G_j, *_j)$ is semigroup $i \neq j$.*

*Then we call $G = \{G_1 \cup G_2 \cup \ldots \cup G_N, *_1, \ldots, *_N\}$ to be a N-group semigroup ( $1 \leq i, j \leq N$).*

**DEFINITION 1.5.3:** *Let $G = \{G_1 \cup G_2 \cup \ldots \cup G_N, *_1, \ldots, *_N\}$ be a N-group-semigroup. We say G is a commutative N-group semigroup if each $(G_i, *_i)$ is a commutative structure, $i = 1, 2, \ldots, N$.*

**DEFINITION 1.5.4:** *Let $G = \{G_1 \cup G_2 \cup ... \cup G_N, *_1, ..., *_N\}$ be a N-group. A proper subset P of G where $(P_1 \cup P_2 \cup ... \cup P_N, *_1, ..., *_N)$ is said to be a N-subgroup of the N-group semigroup G if each $(P_i, *_i)$ is a subgroup of $(G_i, *_i)$; $i = 1, 2, ..., N$.*

**DEFINITION 1.5.5:** *Let $G = \{G_1 \cup G_2 \cup ... \cup G_N, *_1, ..., *_N\}$ be a N-group semigroup where some of the $(G_i, *_i)$ are groups and rest are $(G_j, *_j)$ are semigroups, $1 \leq i, j \leq N$. A proper subset P of G is said to be a N-subsemigroup if $P = \{P_1 \cup P_2 \cup ... \cup P_N, *_1, ..., *_N\}$ where each $(P_i, *_i)$ is only a semigroup under the operation $*_i$.*

Now we proceed on to define the notion of N-subgroup semigroup of a N-group semigroup.

**DEFINITION 1.5.6:** *Let $G = \{G_1 \cup G_2 \cup ... \cup G_N, *_1, ..., *_N\}$ be a N-group semigroup. Let P be a proper subset of G. We say P is a N-subgroup semigroup of G if $P = \{P_1 \cup P_2 \cup ... \cup P_N, *_1, ..., *_N\}$ and each $(P_i, *_i)$ is a group or a semigroup.*

**DEFINITION 1.5.7**: *Let $G = \{G_1 \cup G_2 \cup ... \cup G_N, *_1, ..., *_N\}$ be a N-group semigroup. We call a proper subset P of G where $P = \{P_1 \cup P_2 \cup ... \cup P_N, *_1, ..., *_N\}$ to be a normal N-subgroup semigroup of G if $(G_i, *_i)$ is a group then $(P_i, *_i)$ is a normal subgroup of $G_i$ and if $(G_j, *_j)$ is a semigroup then $(P_j, *_j)$ is an ideal of the semigroup $G_j$. If G has no normal N-subgroup semigroup then we say N-group semigroup is simple.*

**DEFINITION 1.5.8:** *Let $L = \{L_1 \cup L_2 \cup ... \cup L_N, *_1, ..., *_N\}$ be a non empty set with N-binary operations defined on it. We call L a N-loop groupoid if the following conditions are satisfied:*

*i.* $L = L_1 \cup L_2 \cup ... \cup L_N$ *where each* $L_i$ *is a proper subset of L i.e.* $L_i \not\subseteq L_j$ *or* $L_j \not\subseteq L_i$ *if* $i \neq j$, *for* $1 \leq i, j \leq N$.

*ii.* $(L_i, *_i)$ *is either a loop or a groupoid.*

iii. *There are some loops and some groupoids in the collection* $\{L_1, \ldots, L_N\}$.

*Clearly L is a non associative mixed N-structure.*

**DEFINITION 1.5.9:** *Let* $L = \{L_1 \cup L_2 \cup \ldots \cup L_N, *_1, \ldots, *_N\}$ *be a N-loop groupoid. L is said to be a commutative N-loop groupoid if each of* $\{L_i, *_i\}$ *is commutative.*

Now we give an example of a commutative N-loop groupoid.

**DEFINITION 1.5.10:** *Let* $L = \{L_1 \cup L_2 \cup \ldots \cup L_N, *_1, \ldots, *_N\}$ *be a N-loop groupoid. A proper subset P of L is said to be a sub N-loop groupoid if* $P = \{P_1 \cup P_2 \cup \ldots \cup P_N, *_1, \ldots, *_N\}$ *be a N-loop groupoid. A proper subset* $P = \{P_1 \cup P_2 \cup \ldots \cup P_N, *_1, \ldots, *_N\}$ *is such that if P itself is a N-loop groupoid then we call P the sub N-loop groupoid of L.*

**DEFINITION 1.5.11:** *Let* $L = \{L_1 \cup L_2 \cup \ldots \cup L_N, *_1, \ldots, *_N\}$ *be a N-loop groupoid. A proper subset* $G = \{G_1 \cup G_2 \cup \ldots \cup G_N, *_1, \ldots, *_N\}$ *is called a sub N-group if each* $(G_i, *_i)$ *is a group.*

**DEFINITION 1.5.12:** *Let* $L = \{L_1 \cup L_2 \cup \ldots \cup L_N, *_1, \ldots, *_N\}$ *be a N-loop groupoid. A proper subset* $T = \{T_1 \cup T_2 \cup \ldots \cup T_N, *_1, \ldots, *_N\}$ *is said to be a sub N-groupoid of the N-loop groupoid if each* $(T_i, *_i)$ *is a groupoid.*

**DEFINITION 1.5.13:** *Let* $L = \{L_1 \cup L_2 \cup \ldots \cup L_N, *_1, \ldots, *_N\}$ *be a N-loop groupoid. A non empty subset* $S = \{S_1 \cup S_2 \cup \ldots \cup S_N, *_1, \ldots, *_N\}$ *is said to be a sub N-loop if each* $\{S_i, *_i\}$ *is a loop.*

**DEFINITION 1.5.14:** *Let* $L = \{L_1 \cup L_2 \cup \ldots \cup L_N, *_1, \ldots, *_N\}$ *be a N-loop groupoid. A non empty subset* $W = \{W_1 \cup W_2 \cup \ldots \cup W_N, *_1, \ldots, *_N\}$ *of L said to be a sub N-semigroup if each* $\{W_i, *_i\}$ *is a semigroup.*

**DEFINITION 1.5.15:** *Let $L = \{L_1 \cup L_2 \cup \ldots \cup L_N, *_1, \ldots, *_N\}$ be a N-loop groupoid. Let $R = \{R_1 \cup R_2 \cup \ldots \cup R_N, *_1, \ldots, *_N\}$be a proper subset of L. We call R a sub N-group groupoid of the N-loop groupoid L if each $\{R_i, *_i\}$ is either a group or a groupoid.*

**DEFINITION 1.5.16:** *Let $L = \{L_1 \cup L_2 \cup \ldots \cup L_N, *_1, \ldots, *_N\}$ be a N-loop groupoid of finite order. $K = \{K_1 \cup K_2 \cup K_3, *_1, \ldots, *_N\}$ be a sub N-loop groupoid of L. If every sub N-loop groupoid divides the order of the N-loop groupoid L then we call L a Lagrange N-loop groupoid. If no sub N-loop groupoid of L divides the order of L then we say L is a Lagrange free N-loop groupoid.*

**DEFINITION 1.5.17**: *Let $L = \{L_1 \cup L_2 \cup \ldots \cup L_N, *_1, \ldots, *_N\}$ be a N-loop groupoid. We call L a Moufang N-loop groupoid if each $(L_i, *_i)$ satisfies the following identities:*

- *i.* $(xy)(zx) = (x(yz))x$.
- *ii.* $((xy)z)y = x(y(2y))$.
- *iii.* $x(y(xz)) = (xy)x)z$ *for* $x, y, z \in L_i$, $1 \le i \le N$.

*Thus for a N-loop groupoid to be Moufang both the loops and the groupoids must satisfy the Moufang identity.*

**DEFINITION 1.5.18:** *Let $L = \{L_1 \cup L_2 \cup \ldots \cup L_N, *_1, \ldots, *_N\}$ be a N-loop groupoid. A proper subset P ($P = P_1 \cup P_2 \cup \ldots \cup P_N, *_1, \ldots, *_N$) of L is a normal sub N-loop groupoid of L if*

- *i.* *If P is a sub N-loop groupoid of L.*
- *ii.* $x_i P_i = P_i x_i$ *(where* $P_i = P \cap L_i$*).*
- *iii.* $y_i (x_i P_i) = (y_i x_i) P_i$ *for all* $x_i, y_i \in L_i$.

*This is true for each $P_i$, i.e., for i = 1, 2, …, N.*

**DEFINITION 1.5.19:** *Let $L = \{L_1 \cup L_2 \cup \ldots \cup L_N, *_1, \ldots, *_N\}$ and $K = \{K_1 \cup K_2 \cup \ldots \cup K_N, *_1, \ldots, *_N\}$ be two N-loop groupoids such that if $(L_i, *_i)$ is a groupoid then $\{K_i, *_i\}$ is also a groupoid. Likewise if $(L_j, *_j)$ is a loop then $(K_j, *_j)$ is also a* *loop true for $1 \leq i, j \leq N$. A map $\theta = \theta_1 \cup \theta_2 \cup \ldots \cup \theta_N$ from L to K is a N-loop groupoid homomorphism if each $\theta_i : L_i \to K_i$ is a groupoid homomorphism and $\theta_j : L_j \to K_j$ is a loop homomorphism $1 \leq i, j \leq N$.*

**DEFINITION 1.5.20:** *Let $L = \{L_1 \cup L_2 \cup \ldots \cup L_N, *_1, \ldots, *_N\}$ be a N-loop groupoid and $K = \{K_1 \cup K_2 \cup \ldots \cup K_M, *_1, \ldots, *_M\}$ be a M-loop groupoid. A map $\phi = \phi_1 \cup \phi_2 \cup \ldots \cup \phi_{N'}$ from L to K is called a pseudo N-M-loop groupoid homomorphism if each $\phi_i$: $L_t \to K_s$ is either a loop homomorphism or a groupoid homomorphism, $1 \leq t \leq N$ and $1 \leq s \leq M$, according as $L_t$ and $K_s$ are loops or groupoids respectively (we demand $N \leq M$ for if $M > N$ we have to map two or more $L_i$ onto a single $K_j$ which can not be achieved easily).*

**DEFINITION 1.5.21:** *Let $L = \{L_1 \cup L_2 \cup \ldots \cup L_N, *_1, \ldots, *_N\}$be a N-loop groupoid. We call L a Smarandache loop N-loop groupoid (S-loop N-loop groupoid) if L has a proper subset $P = \{P_1 \cup P_2 \cup \ldots \cup P_N, *_1, \ldots, *_N\}$ such that each $P_i$ is a loop i.e. P is a N-loop.*

Now we proceed on to define the mixed N-algebraic structures, which include both associative and non associative structures. Here we define them and give their substructures and a few of their properties.

**DEFINITION 1.5.22:** *Let A be a non empty set on which is defined N-binary closed operations $*_1, \ldots, *_N$. A is called as the N-group-loop-semigroup-groupoid (N-glsg) if the following conditions, hold good.*

i. *$A = A_1 \cup A_2 \cup \ldots \cup A_N$ where each $A_i$ is a proper subset of A (i.e. $A_i \not\subseteq A_j \not\subseteq$ or $A_j \not\subseteq A_i$ if $(i \neq j)$.*
ii. *$(A_i, *_i)$ is a group or a loop or a groupoid or a semigroup (or used not in the mutually exclusive sense) $1 \leq i \leq N$. A is a N –glsg only if the collection $\{A_1, \ldots, A_N\}$ contains groups, loops, semigroups and groupoids.*

**DEFINITION 1.5.23:** *Let $A = \{A_1 \cup \ldots \cup A_N, *_1, \ldots, *_N\}$ where $A_i$ are groups, loops, semigroups and groupoids. We call a non empty subset $P = \{P_1 \cup P_2 \cup \ldots \cup P_N, *_1, \ldots, *_N\}$ of A, where $P_i = P \cap A_i$ is a group or loop or semigroup or groupoid according as $A_i$ is a group or loop or semigroup or groupoid. Then we call P to be a sub N-glsg.*

**DEFINITION 1.5.24**: *Let $A = \{A_1 \cup A_2 \cup \ldots \cup A_N; *_1, \ldots, *_N\}$ be a N-glsg. A proper subset $T = \{T_{i_1} \cup \ldots \cup T_{i_K}, *_{i_1}, \ldots, *_{i_K}\}$ of A is called the sub K-group of N-glsg if each $T_{i_t}$ is a group from $A_r$ where $A_r$ can be a group or a loop or a semigroup of a groupoid but has a proper subset which is a group.*

**DEFINITION 1.5.25:** *Let $A = \{A_1 \cup A_2 \cup \ldots \cup A_N; *_1, \ldots, *_N\}$ be a N-glsg. A proper subset $T = \{T_{i_1} \cup \ldots \cup T_{i_r}\}$ is said to be sub r-loop of A if each $T_{i_j}$ is a loop and $T_{i_j}$ is a proper subset of some $A_p$. As in case of sub K-group r need not be the maximum number of loops in the collection $A_1, \ldots, A_N$.*

**DEFINITION 1.5.26:** *Let $A = \{A_1 \cup A_2 \cup \ldots \cup A_N; *_1, \ldots, *_N\}$ be a N-glsg. Let $P = \{P_1 \cup P_2 \cup \ldots \cup P_N\}$ be a proper subset of A where each $P_i$ is a semigroup then we call P the sub u-semigroup of the N-glsg.*

**DEFINITION 1.5.27:** *Let $A = \{A_1 \cup A_2 \cup \ldots \cup A_N; *_1, \ldots, *_N\}$ be a N-glsg. A proper subset $C = \{C_1 \cup C_2 \cup \ldots \cup C_t\}$ of A is said to be a sub-t-groupoid of A if each $C_i$ is a groupoid.*

**DEFINITION 1.5.28:** *Let $A = \{A_1 \cup A_2 \cup \ldots \cup A_N; *_1, \ldots, *_N\}$ be a N-glsg. Suppose A contains a subset $P = P_{L_1} \cup \ldots \cup P_{L_k}$ of A such that P is a sub K-group of A. If every P-sub K-group of A is commutative we call A to be a sub-K-group commutative N-glsg.*

*If atleast one of the sub-K-group P is commutative we call A to be a weakly sub K-group-commutative N-glsg. If no sub K-group of A is commutative we call A to be a non commutative sub-K-group of N-glsg.*

For more about these notions please refer [50].

**Chapter Two**

# NEUTROSOPHIC GROUPS AND NEUTROSOPHIC N-GROUPS

This chapter has three sections; first section deals with neutrosophic groups and their properties. In section two neutrosophic bigroups are introduced for the first time and analyzed. Section three defines the notion of neutrosophic N-groups.

The notion of neutrosophic structures are defined for the first time. As in general case we, define the neutrosophic N-structure. Also when we define neutrosophic algebraic structure it need not be having all the properties. For we define the neutrosophic element as I where I is an indeterminate and I is such that $I^2 = I$.

This equation $I^2 = I$ does not imply I (I – 1) = 0 or any such relations. It is just like saying $1^2 = I$. So it is not an easy task to talk of inverse for we add or multiply I by a scalar c and call it as c + I or as cI which is a neutrosophic element.

For instance when we say, let $G = \langle Z_2 \cup I\rangle$ generate a neutrosophic group under '+', we have N(G) = {I, 1, 1 + I, 0} where 1.I = I.1 = I also (I + I) = 0 for I + I = 2. I = 0I, as $2 \equiv 0$ (mod 2). We call {0, 1, I, 1 + I} to be the neutrosophic group generated by $Z_2 \cup I$. Here N(G) is a group under '+'.

*Note:* We do not demand a neutrosophic group to be a group. But clearly it contains a group. For if $\langle Z_2 \setminus \{0\} \cup I\rangle$ generates the neutrosophic group under multiplication modulo $2\{\langle(Z_2 \setminus \{0\} \cup I\rangle, \times\}$ is not a group under multiplication modulo 2.

## 2.1 Neutrosophic Groups and their Properties

In this section for the first time the notion of neutrosophic groups are introduced, neutrosophic groups in general do not have group structure. We also define yet another notion called pseudo neutrosophic groups which have group structure. As neutrosophic groups do not have group structure the classical theorems viz. Sylow, Lagrange or Cauchy are not true in general which forces us to define notions like Lagrange neutrosophic groups, Sylow neutrosophic groups and Cauchy elements. Examples are given for the understanding of these new concepts.

**DEFINITION 2.1.1:** *Let (G, *) be any group, the neutrosophic group is generated by I and G under * denoted by* $N(G) = \{\langle G \cup I\rangle, *\}$.

***Example 2.1.1:*** Let $Z_7 = \{0, 1, 2, \ldots, 6\}$ be a group under addition modulo 7. $N(G) = \{\langle Z_7 \cup I\rangle$, '+' modulo 7$\}$ is a neutrosophic group which is in fact a group. For $N(G) = \{a + bI / a, b \in Z_7\}$ is a group under '+' modulo 7. Thus this neutrosophic group is also a group.

***Example 2.1.2:*** Consider the set $G = Z_5 \setminus \{0\}$, G is a group under multiplication modulo 5. Consider $N(G) = \{\langle G \cup I\rangle$, multiplication modulo 5$\}$. N(G) is called the neutrosophic group generated by $G \cup I$. Clearly N(G) is not a group for $I^2 = I$ and I is not the identity but only an indeterminate, but N(G) is defined as the neutrosophic group.

Thus based on this we have the following theorem:

**THEOREM 2.1.1:** *Let (G, *) be a group,* $N(G) = \{\langle G \cup I\rangle, *\}$ *be the neutrosophic group.*

i. *N(G) in general is not a group.*
ii. *N(G) always contains a group.*

*Proof:* To prove N(G) in general is not a group it is sufficient we give an example; consider $\langle Z_5 \setminus \{0\} \cup I\rangle$ = G = {1, 2, 4, 3, I, 2 I, 4 I, 3 I}; G is not a group under multiplication modulo 5. In fact {1, 2, 3, 4} is a group under multiplication modulo 5.

N(G) the neutrosophic group will always contain a group because we generate the neutrosophic group N(G) using G and I. So $G \underset{\neq}{\subset} N(G)$ hence N(G) will always contain a group.

Now we proceed onto define the notion of neutrosophic subgroup of a neutrosophic group.

**DEFINITION 2.1.2:** *Let N(G) = $\langle G \cup I\rangle$ be a neutrosophic group generated by G and I. A proper subset P(G) is said to be a neutrosophic subgroup if P(G) is a neutrosophic group i.e. P(G) must contain a (sub) group.*

***Example 2.1.3:*** Let $N(Z_2) = \langle Z_2 \cup I\rangle$ be a neutrosophic group under addition. $N(Z_2)$ = {0, 1, I, 1 + I}. Now we see {0, I} is a group under + in fact a neutrosophic group {0, 1 + I} is a group under '+' but we call {0, I} or {0, 1 + I} only as pseudo neutrosophic groups for they do not have a proper subset which is a group. So {0, I} and {0, 1 + I} will be only called as pseudo neutrosophic groups (subgroups).

We can thus define a pseudo neutrosophic group as a neutrosophic group, which does not contain a proper subset which is a group. Pseudo neutrosophic subgroups can be found as a substructure of neutrosophic groups. Thus a Pseudo neutrosophic group though has a group structure is not a neutrosophic group and a neutrosophic group cannot be a pseudo neutrosophic group. Both the concepts are different.

Now we see a neutrosophic group can have substructures which are pseudo neutrosophic groups which is evident from the example.

***Example 2.1.4:*** Let $N(Z_4) = \langle Z_4 \cup I\rangle$ be a neutrosophic group under addition modulo 4. $\langle Z_4 \cup I\rangle$ = {0, 1, 2, 3, I, 1 + I, 2I, 3I, 1 + 2I, 1 + 3I, 2 + I, 2 + 2I, 2 + 3I, 3 + I, 3 + 2I, 3 + 3I}. $o(\langle Z_4 \cup I\rangle) = 4^2$.

Thus neutrosophic group has both neutrosophic subgroups and pseudo neutrosophic subgroups. For T = {0, 2, 2 + 2I, 2I} is a neutrosophic subgroup as {0 2} is a subgroup of $Z_4$ under addition modulo 4. P = {0, 2I} is a pseudo neutrosophic group under '+' modulo 4.

Now we are not sure that general properties which are true of groups are true in case of neutrosophic groups for neutrosophic groups are not in general groups. We see that in case of finite neutrosophic groups the order of both neutrosophic subgroups and pseudo neutrosophic subgroups do not divide the order of the neutrosophic group. Thus we give some problems in the chapter 7.

**THEOREM 2.1.2:** *Neutrosophic groups can have non trivial idempotents.*

*Proof:* For I ∈ N(G) and $I^2 = I$ .

*Note:* We cannot claim from this that N(G) can have zero divisors because of the idempotent as our neutrosophic structures are algebraic structures with only one binary operation multiplication in case $I^2 = I$ .

We illustrate these by examples.

***Example 2.1.5:*** Let N(G) = {1, 2, I, 2I} a neutrosophic group under multiplication modulo three. We see $(2I)^2 \equiv I \pmod 3$, $I^2 = I$. (2I) I = 2I, $2^2 \equiv 1 \pmod 3$. So P = {1, I, 2I} is a pseudo neutrosophic subgroup. Also $o(P) \nmid o(N(G))$.

Thus we see order of a pseudo neutrosophic group need not divide the order of the neutrosophic group.

We give yet another example which will help us to see that Lagrange's theorem of finite groups in case of finite neutrosophic groups is not true.

***Example 2.1.6:*** Let N(G) = {1, 2, 3, 4, I, 2I, 3I, 4I} be a neutrosophic group under multiplication modulo 5. Now consider P = {1, 4, I, 2I, 3I, 4I} $\subset$ N(G). P is a neutrosophic subgroup. o(N(G)) = 8 but o(P) = 6, $6 \nmid 8$. So clearly neutrosophic groups in general do not satisfy the Lagrange theorem for finite groups.

So we define or characterize those neutrosophic groups, which satisfy Lagrange theorem as follows:

**DEFINITION 2.1.3:** *Let N(G) be a neutrosophic group. The number of distinct elements in N(G) is called the order of N(G). If the number of elements in N(G) is finite we call N(G) a finite neutrosophic group; otherwise we call N(G) an infinite neutrosophic group, we denote the order of N(G) by o(N(G)) or |N(G)|.*

**DEFINITION 2.1.4:** *Let N(G) be a finite neutrosophic group. Let P be a proper subset of N(G) which under the operations of N(G) is a neutrosophic group. If o(P) / o(N(G)) then we call P to be a Lagrange neutrosophic subgroup. If in a finite neutrosophic group all its neutrosophic subgroups are Lagrange then we call N(G) to be a Lagrange neutrosophic group.*

*If N(G) has atleast one Lagrange neutrosophic subgroup then we call N(G) to be a weakly Lagrange neutrosophic group. If N(G) has no Lagrange neutrosophic subgroup then we call N(G) to be a Lagrange free neutrosophic group.*

We have already given examples of these. Now we proceed on to define the notion called pseudo Lagrange neutrosophic group.

**DEFINITION 2.1.5:** *Let N(G) be a finite neutrosophic group. Suppose L is a pseudo neutrosophic subgroup of N(G) and if o(L) / o(N(G)) then we call L to be a pseudo Lagrange neutrosophic subgroup. If all pseudo neutrosophic subgroups of N(G) are pseudo Lagrange neutrosophic subgroups then we call N(G) to be a pseudo Lagrange neutrosophic group.*

*If N(G) has atleast one pseudo Lagrange neutrosophic subgroup then we call N(G) to be a weakly pseudo Lagrange neutrosophic group. If N(G) has no pseudo Lagrange*

*neutrosophic subgroup then we call N(G) to be pseudo Lagrange free neutrosophic group.*

Now we illustrate by some example some more properties of neutrosophic groups, which paves way for more definitions. We have heard about torsion elements and torsion free elements of a group.

We in this book define neutrosophic element and neutrosophic free element of a neutrosophic group.

**DEFINITION 2.1.6:** *Let N(G) be a neutrosophic group. An element $x \in N(G)$ is said to be a neutrosophic element if there exists a positive integer n such that $x^n = I$, if for any y a neutrosophic element no such n exists then we call y to be a neutrosophic free element.*

We illustrate these by the following examples.

***Example 2.1.7:*** Let N(G) = {1, 2, 3, 4, 5, 6, I, 2I, 3I, 4I, 5I, 6I} be a neutrosophic group under multiplication modulo 7. We have $(3I)^6 = I$, $(4I)^3 = I$ $(6I)^2 = I$, $I^2 = I$, $(2I)^6 = I$, $(5I)^6 = I$. In this neutrosophic group all elements are either torsion elements or neutrosophic elements.

***Example 2.1.8:*** Let us now consider the set {1, 2, 3, 4, I, 2I, 3I, 4I, 1 + I, 2 + I, 3 + I, 4 + I, 1 + 2I, 1 + 3I, 1 + 4I, 2 + 2I, 2 + 3I, 2 + 4I, 3 + 2I, 3 + 3I, 3 + 4I, 4 + 2I, 4 + 3I, 4 + 4I}. This is a neutrosophic group under multiplication modulo 5. For {1, 2, 3, 4} = $Z_5 \setminus \{0\}$ is group under multiplication modulo 5. $(1 + I)^4 = 1$ (mod 5) we ask "Is it a neutrosophic element of N(G)?" $(2 + I)^4 = 1$ (mod 5). $(1 + 4I)^2 = 1 + 4I$ this will be called as neutrosophic idempotent and $(1 + 3I)^2 = 1$ (mod 5) neutrosophic unit.

In view of the above example we define the following:

**DEFINITION 2.1.7:** *Let N(G) be a neutrosophic group. Let $x \in N(G)$ be a neutrosophic element such that $x^m = 1$ then x is called the pseudo neutrosophic torsion element of N(G).*

In the above example we have given several pseudo neutrosophic torsion elements of N(G).

Now we proceed on to define Cauchy neutrosophic elements of a neutrosophic group N(G).

**DEFINITION 2.1.8:** *Let N(G) be a finite neutrosophic group. Let $x \in N(G)$, if x is a torsion element say $x^m = 1$ and if m/o N(G)) we call x a Cauchy element of N(G); if x is a neutrosophic element and $x^t = I$ with t / o(N(G)), we call x a Cauchy neutrosophic element of N(G). If all torsion elements of N(G) are Cauchy we call N(G) as a Cauchy neutrosophic group. If every neutrosophic element is a neutrosophic Cauchy element then we call the neutrosophic group to be a Cauchy neutrosophic, neutrosophic group.*

We now illustrate these concepts by the following examples:

***Example 2.1.9:*** Let N(G) = {0, 1, 2, 3, 4, I, 2I, 3I, 4I} be a neutrosophic group under multiplication modulo 5. {1, 2, 3, 4} is a group under multiplication modulo 5. Now we see o (N (G)) = 9, $4^2 \equiv 1 \pmod 5$ $2 \nmid 9$ similarly $(3I)^4 = I$ but $4 \nmid 9$. Thus none of these elements are Cauchy elements or Cauchy neutrosophic Cauchy elements of N(G).

Now we give yet another example.

***Example 2.1.10:*** Let N(G) be a neutrosophic group of finite order 4 where N(G) = {1, 2, I, 2I} group under multiplication modulo 3. Clearly every element in N(G) is either a Cauchy neutrosophic element or a Cauchy element.

Thus we give yet another definition.

**DEFINITION 2.1.9:** *Let N(G) be a neutrosophic group. If every element in N(G) is either a Cauchy neutrosophic element of N(G) or a Cauchy element of N(G) then we call N(G) a strong Cauchy neutrosophic group.*

The above example is an instance of a strong Cauchy neutrosophic group. Now we proceed on to define the notion of p-Sylow neutrosophic subgroup, Sylow neutrosophic group, weak Sylow neutrosophic group and Sylow free neutrosophic group.

**DEFINITION 2.1.10:** *Let N(G) be a finite neutrosophic group. If for a prime* $p^{\alpha}$ / *o(N(G)) and* $p^{\alpha+1} \nmid$ *o(N(G)), N(G) has a neutrosophic subgroup P of order* $p^{\alpha}$ *then we call P a p-Sylow neutrosophic subgroup of N(G).*

*Now if for every prime p such that* $p^{\alpha}$ / *o(N(G)) and* $p^{\alpha+1} \nmid$ *o(N(G)) we have an associated p-Sylow neutrosophic subgroup then we call N(G) a Sylow neutrosophic group.*

*If N(G) has atleast one p-Sylow neutrosophic subgroup then we call N(G) a weakly Sylow neutrosophic group. If N(G) has no p-Sylow neutrosophic subgroup then we call N(G) a Sylow free neutrosophic group.*

Now unlike in groups we have to speak about Sylow notion associated with pseudo neutrosophic groups.

**DEFINITION 2.1.11:** *Let N(G) be a finite neutrosophic group. Let P be a pseudo neutrosophic subgroup of N (G) such that o (P) =* $p^{\alpha}$ *where* $p^{\alpha}$ / *o (N(G)) and* $p^{\alpha+1} \nmid$ *o(N(G)), for p a prime, then we call P to be a p-Sylow pseudo neutrosophic subgroup of N(G).*

*If for a prime p we have a pseudo neutrosophic subgroup P such that o(P) =* $p^{\alpha}$ *where* $p^{\alpha}$ / *o(N(G)) and* $p^{\alpha+1} \nmid$ *o(N(G)), then we call P to be p-Sylow pseudo neutrosophic subgroup of N(G). If for every prime p such that* $p^{\alpha}$ / *o(N(G)) and* $p^{\alpha+1} \nmid$ *o(N(G)), we have a p-Sylow pseudo neutrosophic subgroup then we call N(G) a Sylow pseudo neutrosophic group.*

*If on the other hand N(G) has atleast one p-Sylow pseudo neutrosophic subgroup then we call N(G) a weak Sylow pseudo neutrosophic group. If N(G) has no p-Sylow pseudo neutrosophic subgroup then we call N(G) a free Sylow pseudo neutrosophic group.*

Now we proceed on to define neutrosophic normal subgroup.

**DEFINITION 2.1.12:** *Let N(G) be a neutrosophic group. Let P and K be any two neutrosophic subgroups of N(G). We say P and K are neutrosophic conjugate if we can find x, y $\in$ N(G) with x P = K y.*

We illustrate this by the following example:

***Example 2.1.11:*** Let N(G) = {0, 1, 2, 3, 4, 5, I, 2I, 3I, 4I, 5I, 1 + I, 2 + I, 3 + I, …, 5 + 5I} be a neutrosophic group under addition modulo 6. P = {0, 3, 3I, 3+3I} is a neutrosophic subgroup of N(G). K = {0, 2, 4, 2 + 2I, 4 + 4I, 2I, 4I} is a neutrosophic subgroup of N(G). For 2, 3 in N(G) we have 2P = 3K = {0}. So P and K are neutrosophic conjugate.

Thus in case of neutrosophic conjugate subgroups K and P we do not demand o(K) = o(P).

Now we proceed on to define neutrosophic normal subgroup.

**DEFINITION 2.1.13:** *Let N(G) be a neutrosophic group. We say a neutrosophic subgroup H of N(G) is normal if we can find x and y in N(G) such that H =xHy for all x, y $\in$ N (G) (Note x = y or y = $x^{-1}$ can also occur).*

***Example 2.1.12:*** Let N (G) be a neutrosophic group given by N (G) = {0, 1, 2, I, 2I, 1 + I, 2 + I, 2I + 1, 2I + 2} under multiplication modulo 3.

H = {1, 2, I, 2I} is a neutrosophic subgroup such that for no element in N (G) \ {0}; xHy = H so H is not normal. Take K = {1, 2, 1 + I, 2 + 2I} is a neutrosophic subgroup. K is not normal.

**DEFINITION 2.1.14:** *A neutrosophic group N(G) which has no nontrivial neutrosophic normal subgroup is called a simple neutrosophic group.*

Now we define pseudo simple neutrosophic groups.

**DEFINITION 2.1.15:** *Let N(G) be a neutrosophic group. A proper pseudo neutrosophic subgroup P of N(G) is said to be*

*normal if we have P = xPy for all x, y ∈ N(G). A neutrosophic group is said to be pseudo simple neutrosophic group if N(G) has no nontrivial pseudo normal subgroups.*

We do not know whether there exists any relation between pseudo simple neutrosophic groups and simple neutrosophic groups.

Now we proceed on to define the notion of right (left) coset for both the types of subgroups.

**DEFINITION 2.1.16:** *Let L (G) be a neutrosophic group. H be a neutrosophic subgroup of N(G) for n ∈ N(G), then H n = {hn / h ∈ H} is called a right coset of H in G.*

Similarly we can define left coset of the neutrosophic subgroup H in G.

It is important to note that as in case of groups we cannot speak of the properties of neutrosophic groups as we cannot find inverse for every x ∈ N (G).

So we make some modification before which we illustrate these concepts by the following examples.

***Example 2.1.13:*** Let N(G) = {1, 2, 3, 4, I, 2I, 3I, 4I} be a neutrosophic group under multiplication modulo 5. Let H = {1, 4, I, 4 I} be a neutrosophic subgroup of N(G). The right cosets of H are as follows:

| | | | | |
|---|---|---|---|---|
| H.2 | = | {2, 3, 2I, 3I} | | |
| H.3 | = | {3, 2, 3I, 2I} | | |
| H.1 | = | H4 | = | {1, 4, I, 4I} |
| H. I | = | {I 4I} | = | H. 4I |
| H.2 I | = | {2I, 3I} | = | H 3I = {3I, 2I}. |

Therefore the classes are

| | | | | |
|---|---|---|---|---|
| [2] | = | [3] | = | {2, 3, 2I, 3I} |
| [1] | = | [4] | = | H = {1, 4, I 4I} |
| [I] | = | [4I] | = | R {I, 4I} |
| [2I] | = | [3I] | = | {3I, 2I}. |

Now we are yet to know whether they will partition N(G) for we see here the cosets do not partition the neutrosophic group.

That is why we had problems with Lagrange theorem so only we defined the notion of Lagrange neutrosophic group.

We give yet another example before which we define the concept of commutative neutrosophic group.

**DEFINITION 2.1.17:** *Let N(G) be a neutrosophic group. We say N(G) is a commutative neutrosophic group if for every pair a, b ∈ N(G), a b = b a.*

We have seen several examples of commutative neutrosophic groups. So now we give an example of a non-commutative neutrosophic group.

***Example 2.1.14:*** Let

$$N(G) = \left\{ \begin{pmatrix} a & b \\ c & d \end{pmatrix} \mid a, b, c, d, \in \{0, 1, 2, I, 2I\} \right\}.$$

N(G) under matrix multiplication modulo 3 is a neutrosophic group which is non commutative.
We now give yet another example of cosets in neutrosophic groups.

***Example 2.1.15:*** Let N(G) = {0, 1, 2, I, 2I, 1 + I, 1 + 2I, 2 + I, 2 + 2I} be a neutrosophic group under multiplication modulo 3.

Consider P = {1, 2, I, 2I}. ⊂ N(G); P is a neutrosophic subgroup.

| | | |
|---|---|---|
| P. 0 | = | {0} |
| P. 1 | = | {1, 2, I, 2I} |
| | = | P2 |
| P. I | = | {I 2I} |
| | = | P. 2I |
| P (1 + I) | = | {1 + I, 2 + 2I, 2I, I} |

$P(2 + I) = \{2 + I, 1 + 2I, 0\}$
$P(1 + 2I) = \{1 + 2I, 2 + I, 0\}$
$= P(2 + I)$
$P(2 + 2I) = \{2 + 2I, 1 + I, I, 2I\}$
$= P(1 + I)$.

We see the coset does not partition the neutrosophic group. Now using the concept of pseudo neutrosophic subgroup we define pseudo coset.

**DEFINITION 2.1.18:** *Let N(G) be a neutrosophic group. K be a pseudo neutrosophic subgroup of N(G). Then for a ∈ N(G), Ka = {ka | k ∈ K} is called the pseudo right coset of K in N(G).*

On similar lines we define the notion of pseudo left coset of a pseudo neutrosophic subgroup K of N (G). We illustrate this by the following example.

***Example 2.1.16:*** Let N (G) = {0, 1, I, 2, 2I, 1 + I, 1 + 2I, 2 + I, 2 + 2I} be a neutrosophic group under multiplication modulo 3. Take K = {1, 1 + I}, a pseudo neutrosophic group.

Now we will study the cosets of K . K . 0 = {0}.

$K\,1 = \{1, 1 + I\}$
$K(1 + I) = \{1 + I\}$
$K\,2 = \{2, 2 + 2I\}$
$K.\,I = \{I, 2I\}$
$K\,2I = \{2I, I\}$
$= K.\,I.$
$K(1 + 2I) = \{1 + 2I\}$
$K(2 + I) = \{2 + I\}$
$K(2 + 2I) = \{2 + 2I, 2\}$
$= K.2.$

We see even the pseudo neutrosophic subgroups do not in general partition the neutrosophic group which is evident from the example.

Now we proceed on to define the concept of center of a neutrosophic group.

**DEFINITION 2.1.19:** *Let N(G) be a neutrosophic group, the center of N(G) denoted by C(N(G)) = {x ∈ N(G) | ax = xa for all a ∈ N(G)}.*

*Note:* Clearly C(N(G)) ≠ ϕ for the identity of the neutrosophic group belongs to C(N(G)). Also if N(G) is a commutative neutrosophic group then C(N(G)) = N(G). As in case of groups we can define in case of neutrosophic groups also direct product of neutrosophic groups N(G).

**DEFINITION 2.1.20**: *Let $N(G_1)$, $N(G_2)$, ..., $N(G_n)$ be n neutrosophic groups the direct product of the n-neutrosophic groups is denoted by $N(G) = N(G_1) \times ... \times N(G_n) = \{(g_1, g_2, ..., g_n) \mid g_i \in N(G_i); i = 1, 2, ..., n\}$.*

N(G) is a neutrosophic group for the binary operation defined is component wise; for if $*_1, *_2, \ldots, *_n$ are the binary operations on $N(G_1), \ldots, N(G_n)$ respectively then for $X = (x_1, \ldots, x_n)$ and $Y = (y_1, y_2, \ldots, y_n)$ in N(G), $X * Y = (x_1, \ldots, x_n) * (y_1, \ldots, y_n) = (x_1*y_1, \ldots, x_n*y_n) = (t_1, \ldots, t_n) \in N(G)$ thus closure axiom is satisfied. We see if $e = (e_1, \ldots, e_n)$ is identity element where each $e_i$ is the identity element of $N(G_i)$; $1 \leq i \leq n$ then $X * e = e * X = X$.

It is left as a matter of routine for the reader to check that N (G) is a neutrosophic group. Thus we see that the concept of direct product of neutrosophic group helps us in obtaining more and more neutrosophic groups.

*Note:* It is important and interesting to note that if we take in $N(G_i)$, $1 \leq i \leq n$. some $N(G_i)$ to be just groups still we continue to obtain neutrosophic groups.

We now give some examples as illustrations.

***Example 2.1.17:*** Let $N(G_1) = \{0, 1, I, 1 + I\}$ and $G_2 = \{g \mid g^3 = 1\}$. $N(G) = N(G_1) \times G_2 = \{(0, g) (0, 1) (0, g^2) (1, g) (1, 1) (1, g^2) (I, g) (1\ g^2) (I, 1) (1 + I, 1) (1 + I, g) ( 1 + I, g^2)\}$ is a

neutrosophic group of order 12. Clearly $\{(1, 1), (1, g), (1, g^2)\}$ is the group in N(G).

Now several other properties which we have left out can be defined appropriately.

*Note:* We can also define independently the notion of pseudo neutrosophic group as a neutrosophic group which has no proper subset which is a group but the pseudo neutrosophic group itself is a group. We can give examples of them, the main difference between a pseudo neutrosophic group and a neutrosophic group is that a pseudo neutrosophic group is a group but a neutrosophic group is not a group in general but only contains a proper subset which is a group. Now we give an example of a pseudo neutrosophic group.

***Example 2.1.18:*** Consider the set N(G) = {1, 1 + I} under the operation multiplication modulo 3. {1, 1 + I} is a group called the pseudo neutrosophic group for this is evident from the table.

| * | 1 | 1 + I |
|---|---|---|
| 1 | 1 | 1 + I |
| 1 + I | 1 + I | 1 |

Clearly {1, 1 + I} is group but has no proper subset which is a group. Also this pseudo neutrosophic group can be realized as a cyclic group of order 2.

## 2.2 Neutrosophic Bigroups and their Properties

Now we proceed onto define the notion of neutrosophic bigroups. However the notion of bigroups have been defined and dealt in [48]. The neutrosophic bigroups also enjoy special properties and do not satisfy most of the classical results. So substructures like neutrosophic subbigroups, Lagrange neutrosophic subbigroups, p-Sylow neutrosophic subbigroups are defined, leading to the definition of Lagrange neutrosophic bigroups, Sylow neutrosophic bigroups and super Sylow neutrosophic bigroups.

For more about bigroups refer [48].

**DEFINITION 2.2.1:** *Let $B_N(G) = \{B(G_1) \cup B(G_2), *_1, *_2\}$ be a non empty subset with two binary operation on $B_N(G)$ satisfying the following conditions:*

- *i. $B_N(G) = \{B(G_1) \cup B(G_2)\}$ where $B(G_1)$ and $B(G_2)$ are proper subsets of $B_N(G)$.*
- *ii. $(B(G_1), *_1)$ is a neutrosophic group.*
- *iii. $(B(G_2), *_2)$ is a group.*

*Then we define $(B_N(G), *_1, *_2)$ to be a neutrosophic bigroup. If both $B(G_1)$ and $B(G_2)$ are neutrosophic groups we say $B_N(G)$ is a strong neutrosophic bigroup. If both the groups are not neutrosophic group we see $B_N(G)$ is just a bigroup.*

We first illustrate this with some examples before we proceed on to analyze their properties.

***Example 2.2.1:*** Let $B_N(G) = \{B(G_1) \cup B(G_2)\}$ where $B(G_1) = \{g \mid g^9 = 1\}$ be a cyclic group of order 9 and $B(G_2) = \{1, 2, I, 2I\}$ neutrosophic group under multiplication modulo 3. We call $B_N(G)$ a neutrosophic bigroup.

***Example 2.2.2:*** Let $B_N(G) = \{B(G_1) \cup B(G_2)\}$ where $B(G_1) = \{1, 2, 3, 4, I, 4I, 3I, 2I\}$ a neutrosophic group under multiplication modulo 5. $B(G_2) = \{0, 1, 2, I, 2I, 1 + I, 2 + I, 1 + 2I, 2 + 2I\}$ is a neutrosophic group under multiplication modulo 3. Clearly $B_N(G)$ is a strong neutrosophic bigroup.

We now define the notion of finite neutrosophic bigroup.

**DEFINITION 2.2.2:** *Let $B_N(G) = \{B(G_1) \cup B(G_2), *_1, *_2\}$ be a neutrosophic bigroup. The number of distinct elements in $B_N(G)$ gives the order of the neutrosophic bigroup. If the number of elements in $B_N(G)$ is finite we call $B_N(G)$ a finite neutrosophic bigroup. If it has infinite number of elements then we call*

*$B_N$ (G) an infinite neutrosophic bigroup. We denote the order of $B_N$ (G) by o ($B_N$ (G)).*

We now proceed on to define the notion of neutrosophic bisubgroup / subbigroup (we can use both to mean the same structure).

**DEFINITION 2.2.3:** *Let $B_N(G) = \{B(G_1) \cup B(G_2), *_1, *_2\}$ be a neutrosophic bigroup. A proper subset $P = \{P_1 \cup P_2\ *_1, *_2\}$ is a neutrosophic subbigroup of $B_N(G)$ if the following conditions are satisfied $P = \{P_1 \cup P_2, *_1, *_2\}$ is a neutrosophic bigroup under the operations $*_1$, $*_2$ i.e. ($P_1$, $*_1$) is a neutrosophic subgroup of ($B_1$, $*_1$) and ($P_2$, $*_2$) is a subgroup of ($B_2$, $*_2$). $P_1 = P \cap B_1$ and $P_2 = P \cap B_2$ are subgroups of $B_1$ and $B_2$ respectively. If both of $P_1$ and $P_2$ are not neutrosophic then we call $P = P_1 \cup P_2$ to be just a bigroup.*

We illustrate this by an example.

***Example 2.2.3:*** Let B(G) = $\{B(G_1) \cup B(G_2), *_1, *_2\}$ be a neutrosophic bigroup, where

$B(G_1)$ = {0, 1, 2, 3, 4, I, 4I, 2I, 3I} is a neutrosophic group under multiplication modulo 5.
$B(G_2)$ = $\{g \mid g^{12} = 1\}$ is a cyclic group of order 12.

Let $P(G) = \{P(G_1) \cup P(G_2), *_1, *_2\}$ where

$P(G_1)$ = $\{1, I, 4, 4I\} \subset B(G_1)$ is a neutrosophic group.
$P(G_2)$ = $\{g^2, g^4, g^6, g^8, g^{10}, 1\} \subset B(G_2)$.

$P(G_1) \cup P(G_2) = P(G)$ is a neutrosophic subbigroup of the neutrosophic bigroup.

Let us consider $M = M_1 \cup M_2$ where $M_1 = \{1, 4\} \subset B(G_1)$ and $M_2 = \{1, g^6\}$. M is just a subbigroup of the neutrosophic bigroup $B_N(G) = B(G_1) \cup B(G_2)$.

*Note:* If both $B(G_1)$ and $B(G_2)$ are commutative groups then we call $B_N(G) = \{B(G_1) \cup B(G_2)$ to be a commutative bigroup. If both $B(G_1)$ and $B(G_2)$ are cyclic, we call $B_N(G)$ a cyclic bigroup. We wish to state the notion of normal bigroup.

**DEFINITION 2.2.4:** *Let $B_N(G) = \{B(G_1) \cup B(G_2), *_1, *_2\}$ be a neutrosophic bigroup. $P(G) = \{P(G_1) \cup P(G_2), *_1, *_2\}$ be a neutrosophic bigroup. $P(G) = \{P(G_1) \cup P(G_2), *_1, *_2\}$ is said to be a neutrosophic normal subbigroup of $B_N(G)$ if $P(G)$ is a neutrosophic subbigroup and both $P(G_1)$ and $P(G_2)$ are normal subgroups of $B(G_1)$ and $B(G_2)$ respectively.*

We just illustrate this by the following example.

***Example 2.2.4:*** $B_N(G) = \{B(G_1) \cup B(G_2), *_1, *_2\}$; where $B(G_1) = \{1, 4, 2, 3, I, 2I\ 3I, 4I\}$ and $B(G_2) = S_3$. $B_N(G)$ is a neutrosophic bigroup. This has no neutrosophic normal subbigroup.

Now we give some examples, which show that the order of neutrosophic subbigroup does not divide the order of the neutrosophic bigroup in general.

***Example 2.2.5:*** Let $B_N(G) = \{B(G_1) \cup B(G_2), *_1, *_2\}$; where

$B(G_1)$ = $\{1, 2, 3, 4, I, 2I, 3I, 4I\}$ a neutrosophic group under multiplication modulo 5 and
$B(G_2)$ = $\{g \mid g^9 = 1\}$, a cyclic group of order 9,

$o(B_N(G)) = 17$ a prime. But this neutrosophic bigroup has neutrosophic subbigroups.

Take $P = P(G_1) \cup P(G_2)$ where $P(G_1) = \{1, 4, I, 4I\}$ and $P(G_2) = \{1, g^3, g^6\}$. P is a neutrosophic subbigroup. $o(P) = 7$ and $(7, 17) = 1$.

In fact order of none of the neutrosophic subbigroups will divide the order of the neutrosophic bigroup as $o(B_N(G)) = 17$, a prime.

***Example 2.2.6:*** Let $B_N(G) = \{B(G_1) \cup B(G_2), *_1, *_2\}$ where

$B(G_1)$ = {0, 1, 2, 3, 4, I, 2I, 3I, 4I, 1 + I, 2 + I, 3 + I, 4 + I, 1 + 2I, 2 + 2I, 3 + 2I, 4 + 2I, 1 + 3I, 2 + 3I, 3 + 3I, 4 + 3I, 4 + 4I, 3 + 4I, 2 + 4I, 1 + 4I } be a neutrosophic group under multiplication modulo 5.

$B(G_2)$ = $\{g \mid g^{10} = 1\}$ a cyclic group of order 10.

$o(B_N(G)) = 35$. $P(G) = P(G_1) \cup P(G_2)$ where $P(G_1) = \{0, 1, I, 4I, 4\} \subset B(G_1)$ and $P(G_2) = \{g^2, g^4, g^6, g^8, 1\} \subset B(G_2)$. $o(P(G)) = 10$, $10 \nmid 35$.

Take $T = T(G_1) \cup T(G_2)$ where $T(G_1) = \{0, 1, I, 4I, 4\} \subset B(G_1)$, $T(G_2) = \{1, g^5\} \subset B(G_2)$, $o(T) = 7$. 7/35. So this neutrosophic subbigroup is such that the order divides the order of the neutrosophic bigroup.

Seeing these examples we venture to make the following definitions.

**DEFINITION 2.2.5:** *Let $B_N(G) = \{B(G_1) \cup B(G_2), *_1, *_2\}$ be a neutrosophic bigroup of finite order. Let $P(G) = \{P(G_1) \cup P(G_2), *_1, *_2\}$ be a neutrosophic subbigroup of $B_N(G)$. If $o(P(G)) / o(B_N(G))$ then we call P(G) a Lagrange neutrosophic subbigroup, if every neutrosophic subbigroup P is such that $o(P) / o(B_N(G))$ then we call $B_N(G)$ to be a Lagrange neutrosophic bigroup. That is if every proper neutrosophic subbigroup is Lagrange then we call $B_N(G)$ to be a Lagrange neutrosophic bigroup. If $B_N(G)$ has atleast one Lagrange neutrosophic subbigroup then we call $B_N(G)$ to be a weak Lagrange neutrosophic bigroup. If $B_N(G)$ has no Lagrange neutrosophic subbigroup then $B_N(G)$ is called Lagrange free neutrosophic bigroup.*

Now we proceed on to give some examples and results which guarantee the existence of Lagrange free neutrosophic bigroup.

**THEOREM 2.2.1:** *Let $B_N(G) = \{B(G_1) \cup B(G_2), *_1, *_2\}$ be a neutrosophic bigroup of prime order p, then $B_N(G)$ is a Lagrange free neutrosophic bigroup.*

*Proof:* Given $B_N(G)$ is a neutrosophic bigroup of order p, p a prime. So if $\{P_N(G_1) \cup P(G_2), *_1, *_2\}$ is any neutrosophic subbigroup of $B_N(G)$ clearly $(o(P_N(G), o(B_N(G))) = 1$. Hence the claim.

Now we proceed on to define the notion of Sylow property on the neutrosophic bigroup.

**DEFINITION 2.2.6:** *Let $B_N(G) = \{B(G_1) \cup B(G_2), *_1, *_2\}$ be a neutrosophic bigroup of finite order. Let p be a prime, such that $p^\alpha / o(B_N(G))$ and $p^{\alpha+1} \nmid o(B_N(G))$, if $B_N(G)$ has a neutrosophic subbigroup P of order $p^\alpha$ then we call P a p-Sylow neutrosophic subbigroup.*

*If for every prime p such that $p^\alpha / o(B_N(G))$ and $p^{\alpha+1} \nmid o(B_N(G))$ we have a p-Sylow neutrosophic subbigroup; then we call $B_N(G)$ a Sylow- neutrosophic bigroup. If $B_N(G)$ has atleast one p-Sylow neutrosophic subbigroup then we call $B_N(G)$ a weakly Sylow neutrosophic bigroup. If $B_N(G)$ has no p-Sylow neutrosophic subbigroup then we call $B_N(G)$ a free Sylow neutrosophic bigroup.*

We know the collection of all neutrosophic bigroups which are of order p, p a prime then we call $B_N(G)$ a Sylow free neutrosophic bigroup.

Now we proceed on to define the notion of Cauchy neutrosophic element and Cauchy element of a neutrosophic bigroup.

**DEFINITION 2.2.7:** *Let $B_N(G)$ be a neutrosophic bigroup of finite order, x in $B_N(G)$ is a Cauchy element if $x^m = 1$ and $m / o(B_N(G))$; y in $B_N(G)$ is a Cauchy neutrosophic element if $y^t = I$ and $t / o(B_N(G))$. If every element in $B_N(G)$ is either a Cauchy element or a Cauchy neutrosophic element then we call $B_N(G)$ to be a Cauchy neutrosophic bigroup. If $B_N(G)$ has atleast a Cauchy element or a Cauchy neutrosophic element then we call $B_N(G)$ a weakly Cauchy neutrosophic bigroup.*

*If no element in $B_N(G)$ is a Cauchy element or a Cauchy neutrosophic element then we call $B_N(G)$ a Cauchy free neutrosophic bigroup.*

Now we define the notion of conjugate neutrosophic subbigroup.

**DEFINITION 2.2.8:** *Let $B_N(G) = \{B(G_1) \cup B(G_2), *_1, *_2\}$ be a neutrosophic bigroup. Suppose $P = \{P(G_1) \cup P(G_2), *_1, *_2\}$ and $K = \{K(G_1) \cup K(G_2), *_1, *_2\}$ be any two neutrosophic subbigroups we say P and K are conjugate if each $P(G_i)$ is conjugate with $K(G_i)$, $i = 1, 2$, then we say P and K are neutrosophic conjugate subbigroups of $B_N(G)$.*

It is interesting to note that even if P and K are neutrosophic conjugate subbigroups o(P) need not be equal to o(K) which is a marked difference from the usual groups.

Now we proceed on to define the notion the neutrosophic bicentre of a neutrosophic bigroup.

**DEFINITION 2.2.9:** *Let $B_N(G) = \{B(G_1) \cup B(G_2), *_1, *_2\}$ be any neutrosophic bigroup. The neutrosophic bicentre of the bigroup $B_N(G)$ denoted $C_N(G) = C(G_1) \cup C(G_2)$ where $C(G_1)$ is the centre of $B(G_1)$ and $C(G_2)$ is the centre of $B(G_2)$. If the neutrosophic bigroup is commutative then $C_N(G) = B_N(G)$.*

It is important to note that $C_N(G)$ is non-empty for, atleast they have identity element in them.

***Example 2.2.7:*** Let $B_N(G) = \{B(G_1) \cup B(G_2) *_1, *_2\}$ where $B(G_1) = \{1, 2, I, 2I\}$ a neutrosophic group under multiplication modulo 3 and $B(G_2) = S_3$. The centre of $B_N(G)$ which is $C_N(G) = B(G_1) \cup \left\{\begin{pmatrix} 1 & 2 & 3 \\ 1 & 2 & 3 \end{pmatrix}\right\} \subset B_N(G)$. We see $C_N(G)$ is a neutrosophic bigroup of order 5 and $o(C_N(G)) / o(B_N(G))$.

Now we proceed on to define strong neutrosophic bigroups and enumerate some of its properties.

**DEFINITION 2.2.10:** *A set* $(\langle G \cup I\rangle +, o)$ *with two binary operations '+' and 'o' is called a strong neutrosophic bigroup if*

i. $\langle G \cup I\rangle = \langle G_1 \cup I\rangle \cup \langle G_2 \cup I\rangle$,
ii. $(\langle G_1 \cup I\rangle, +)$ *is a neutrosophic group and*
iii. $(\langle G_2 \cup I\rangle, o)$ *is a neutrosophic group.*

***Example 2.2.8:*** Let $\{\langle G \cup I\rangle, *_1, *_2\}$ be a neutrosophic strong bigroup where $\langle G \cup I\rangle = \langle Z \cup I\rangle \cup \{0, 1, 2, 3, 4, I, 2I, 3I, 4I\}$. $\langle Z \cup I \rangle$ under '+' is a neutrosophic group and {0, 1, 2, 3, 4, I, 2I, 3I, 4I} under multiplication modulo 5 is a neutrosophic group.

Now we proceed on to define neutrosophic subbigroup of a strong neutrosophic bigroup.

**DEFINITION 2.2.11:** *A subset* $H \neq \phi$ *of a strong neutrosophic bigroup* $(\langle G \cup I \rangle, *, o)$ *is called a strong neutrosophic subbigroup if H itself is a strong neutrosophic bigroup under '*' and 'o' operations defined on* $\langle G \cup I\rangle$.

We have a interesting theorem based on this definition.

**THEOREM 2.2.2:** *Let* $(\langle G \cup I\rangle, +, o)$ *be a strong neutrosophic bigroup. A subset* $H \neq \phi$ *of a strong neutrosophic bigroup* $\langle G \cup I\rangle$ *is a neutrosophic subbigroup then (H, +) and (H, o) in general are not neutrosophic groups.*

*Proof:* Given $(\langle G \cup I\rangle, +, o)$ is a strong neutrosophic bigroup and $H \neq \phi$ of G is a neutrosophic subbigroup of G to show (H, +) and (H, o) are not neutrosophic bigroups. We give an example, to prove this consider the strong neutrosophic bigroup; $\langle G \cup I\rangle = \langle Z \cup I\rangle \cup \{-1, 1, i, -i, I, -I, iI, -iI\}$ under the operations + and 'o'.

$(\langle Z \cup I\rangle, +)$ is a neutrosophic group under '+' and $\{-1, 1, i, -i, I, -I, iI, -iI\}$ is a neutrosophic group under multiplication 'o'. $H = \{0, 1, -1, \langle 2Z \cup 2I\rangle, i, -i\}$, H = proper subset of $\langle G \cup I\rangle$ and $H = H_1 \cup H_2 = \{1\ -1\ i\ -i\} \cup \{\langle 2Z \cup 2I \rangle\}$ is a neutrosophic subbigroup but H is not a group under '+' or 'o'.

Thus we give a nice characterization theorem about the strong neutrosophic subbigroup.

**THEOREM 2.2.3:** *Let $\{\langle G \cup I \rangle, +, o\}$ be a strong neutrosophic bigroup. Then the subset H ($\neq \phi$) is a strong neutrosophic subbigroup of $\langle G \cup I \rangle$ if and only if there exists two proper subsets $\langle G_1 \cup I \rangle$, $\langle G_2 \cup I \rangle$ of $\langle G \cup I \rangle$ such that*

i. *$\langle G \cup I \rangle = \langle G_1 \cup I \rangle \cup \langle G_2 \cup I \rangle$ with $(\langle G_1 \cup I \rangle, +)$ is a neutrosophic group and $(\langle G_2 \cup I \rangle, o)$ a neutrosophic group.*
ii. *$(H \cap \langle G_i \cup I \rangle, +)$ is a neutrosophic subbigroup of $\langle G_i \cup I \rangle$ for i = 1, 2.*

*Proof:* Let $H \neq \phi$ be a strong neutrosophic subbigroup of $\langle G \cup I \rangle$. Therefore there exists two proper subsets, $H_1$, $H_2$ of H such that

(1) $H = H_1 \cup H_2$.
(2) $(H_1, +)$ is a neutrosophic group.
(3) $(H_2, o)$ is a neutrosophic group.

Now we choose $H_1$ as $H \cap \langle G_1 \cup I \rangle$ then we see $H_1$ is a subset of $\langle G_1 \cup I \rangle$ and by (2) $(H_1, +)$ is a neutrosophic subgroup of $\langle G_1 \cup I \rangle$. Similarly choose $H_2 = H \cap \langle G_2 \cup I \rangle$ and we see $H_2$ as $H \cap \langle G_2 \cup I \rangle$ which is clearly a neutrosophic subgroup of $\langle G_2 \cup I \rangle$.

Conversely suppose (1) and (2) of the statements of the theorem be true. To prove, (H, +, o) is a strong neutrosophic bigroup it is enough to prove $(H \cap \langle G_1 \cup I \rangle) \cup (H \cap \langle G_2 \cup I \rangle) = H$.

By using set theoretic methods the relation is true. It is important to note that in the above theorem the condition (1) can be removed, we have included it only for easy working.

Now we proceed on to define the notion of strong neutrosophic commutative bigroup.

**DEFINITION 2.2.12:** *Let ($\langle G \cup I \rangle$, $*_1$, $*_2$) be a strong neutrosophic bigroup. $\langle G \cup I \rangle = \langle G_1 \cup I \rangle \cup \langle G_2 \cup I \rangle$ is said to be a commutative neutrosophic bigroup if both ($\langle G_1 \cup I \rangle$, $*_1$) and ($\langle G_2 \cup I \rangle$, $*_2$) are commutative.*

Now as in case of bigroups even in case of strong neutrosophic bigroups by the order of the strong neutrosophic bigroup we mean the number of distinct elements in it. If the number of distinct elements is finite we say the neutrosophic bigroup is of finite order, otherwise of infinite order.

**THEOREM 2.2.4:** *Let ($\langle G \cup I \rangle$, +, o) be a strong finite neutrosophic bigroup. Let $H \neq \phi$ be a proper neutrosophic subbigroups of $\langle G \cup I \rangle$. Then the order of H in general does not divide the order of $\langle G \cup I \rangle$.*

*Proof:* This is evident from the following example. Take $\langle G \cup I \rangle = (\langle G_1 \cup I \rangle \cup \langle G_2 \cup I \rangle)$ where

$\langle G_1 \cup I \rangle$ = {1, 2, 3, 4, I, 2I, 3I, 4I} a neutrosophic group under multiplication modulo 5.
$\langle G_2 \cup I \rangle$ = {0, 1, 2, I, 2I, 1 + I, 1 + 2 I, I + 2, 2 + 2I} multiplication modulo 3.

$o(\langle G \cup I \rangle) = 17$.

Take $\{1, I, 4I, 4\} \subset \{1, 2, 3, 4, I, 2I, 3I, 4I\}$ and $\{0, 1, I, 2, 2I\} \subset \langle G_2 \cup I \rangle$. $H = H_1 \cup H_2 = \{1, I, 4, 4I\} \cup \{0, 1, 2, I, 2I\}$ is a neutrosophic subbigroup. o(H) = 9. $9 \nmid 17$.

So it is not easy to derive Lagrange Theorem for neutrosophic bigroups. Now we proceed on to define neutrosophic normal subbigroup.

**DEFINITION 2.2.13:** *Let ($\langle G \cup I \rangle$ +, o) be a strong neutrosophic bigroup with $\langle G \cup I \rangle = \langle G_1 \cup I \rangle \cup \langle G_2 \cup I \rangle$. Let (H, + o) be a neutrosophic subbigroup where $H = H_1 \cup H_2$. We say H is a neutrosophic normal subbigroup of G if both $H_1$ and $H_2$ are* *neutrosophic normal subgroups of* $\langle G_1 \cup I\rangle$ *and* $\langle G_2 \cup I\rangle$ *respectively.*

***Example 2.2.9:*** Let ($\langle G \cup I\rangle$, +, o) be a neutrosophic bigroup where $\langle G \cup I\rangle = \langle G_1 \cup I\rangle \cup \langle G_2 \cup I\rangle$ with

$\langle G_1 \cup I\rangle$ = $\langle Z \cup I\rangle$ the neutrosophic group under addition and
$\langle G_2 \cup I\rangle$ = {0, 1, 2, 3, 4, I, 2I, 3I, 4I} a neutrosophic group under multiplication modulo 5.

Take $H = H_1 \cup H_2$ where

$H_1$ = $\{\langle 2Z \cup I\rangle, +\} \subset \{\langle Z \cup I\rangle, +\}$ is a neutrosophic subgroup and
$H_2$ = {0, 1, I, 4, 4I} is a neutrosophic subgroup.

Thus H is a strong neutrosophic normal subbigroup of $\langle G \cup I\rangle$.

As in case of strong neutrosophic subbigroup we see, the strong neutrosophic normal subbigroup of a finite neutrosophic bigroup does not in general divide the order of the neutrosophic bigroup.

Now we proceed on to define strong neutrosophic homomorphism of strong neutrosophic bigroups.

**DEFINITION 2.2.14:** *Let* ($\langle G \cup I\rangle$, +, o) *and* ($\langle K \cup I\rangle$, o', $\oplus$) *be any two strong neutrosophic bigroups where* $\langle G \cup I\rangle = \langle G_1 \cup I\rangle \cup \langle G_2 \cup I\rangle$ *and* $\langle K \cup I\rangle = \langle K_1 \cup I\rangle \cup \langle K_2 \cup I\rangle$. *We say a bimap* $\phi = \phi_1 \cup \phi_2$: $\langle G \cup I\rangle \to \langle K \cup I\rangle$ *(Here* $\phi_1(I) = I$ *and* $\phi_2(I) = I$*) is said to be a strong neutrosophic bigroup bihomomorphism if* $\phi_1 = \phi / \langle G_1 \cup I\rangle$ *and* $\phi_2 = \phi / \langle G_2 \cup I\rangle$ *where* $\phi_1$ *and* $\phi_2$ *are neutrosophic group homomorphism from* $\langle G_1 \cup I\rangle$ *to* $\langle K_1 \cup I\rangle$ *and* $\langle G_2 \cup I\rangle$ *to* $\langle K_2 \cup I\rangle$ *respectively.*

**THEOREM 2.2.5:** *Let* ($\langle G \cup I\rangle = \langle G_1 \cup I\rangle \cup \langle G_2 \cup I\rangle$, o, +) *be a strong neutrosophic bigroup of finite order n. If p/n then the neutrosophic bigroup may not in general have neutrosophic bigroup of order p.(p need not necessarily be a prime).*

*Proof:* We prove this by the following example. Let (⟨G ∪ I⟩, +, *) be a strong neutrosophic bigroup where ⟨G ∪ I⟩ = ⟨$G_1$ ∪ I⟩ ∪ ⟨$G_2$ ∪ I⟩ with

| | | |
|---|---|---|
| ⟨$G_1$ ∪ I⟩ | = | {0, 1, 2, 1 + I, I, 2I, 2 + I, 2I + 2, 1 + 2I} a neutrosophic group under multiplication modulo 3. |
| ⟨$G_2$ ∪ I⟩ | = | {0, 1, 2, 3, 4, I, 2I, 3I, 4I} a neutrosophic group under multiplication modulo 5. |

o(⟨G ∪ I⟩) = 18. Take H = $H_1$ ∪ $H_2$ where $H_1$ = {1, 2, I, 2I} and $H_2$ = {1, 4, I, 4I}. o(H) = 8, H is a neutrosophic subbigroup of ⟨G ∪ I ⟩. o(H) $\nmid$ o(⟨G ∪ I⟩) i.e. 8 $\nmid$ 18. Hence the claim.

One can develop the notion of biorder and the notion of pseudo divisor to strong neutrosophic bigroup from bigroups. Interested reader can refer [48, 50].

Here we define a new notion called Lagrange strong neutrosophic subbigroup and Lagrange strong neutrosophic bigroup.

**DEFINITION 2.2.15:** *Let (⟨G ∪ I⟩, *, o) be a strong neutrosophic bigroup of finite order. Let H ≠ $\phi$ be a strong neutrosophic subbigroup of (⟨G ∪ I⟩, *, o). If o (H) / o (⟨G ∪ I⟩) then we call H a Lagrange strong neutrosophic subbigroup of ⟨G ∪ I⟩. If every strong neutrosophic subbigroup of ⟨G ∪ I⟩ is a Lagrange strong neutrosophic subbigroup then we call ⟨G ∪ I⟩ a Lagrange strong neutrosophic bigroup.*

*If the strong neutrosophic bigroup has atleast one Lagrange strong neutrosophic subbigroup then we call ⟨G ∪ I⟩ a weakly Lagrange strong neutrosophic bigroup.*

*If ⟨G ∪ I⟩ has no Lagrange strong neutrosophic subbigroup then we call ⟨G ∪ I⟩ a Lagrange free strong neutrosophic bigroup.*

The following result is important for it gives us a class of Lagrange free strong neutrosophic bigroup.

**THEOREM 2.2.6:** *All strong neutrosophic bigroups of a prime order are Lagrange free strong neutrosophic bigroups.*

The proof of the above theorem is left as an exercise.

Now it may happen for a finite strong neutrosophic bigroup of order n, we may find a prime p such that $p^{\alpha}/n$ and $p^{\alpha+1} \nmid n$ and the strong neutrosophic bigroup having strong neutrosophic subbigroups of order $p^{\alpha}$; how to define such neutrosophic subbigroup?

**DEFINITION 2.2.16:** *Let ($\langle G \cup I \rangle$, o, *) be a strong neutrosophic bigroup of finite order n. If for a prime p such that $p^{\alpha} / o\langle G \cup I \rangle$ and $p^{\alpha+1} \nmid o\langle G \cup I \rangle$ we have strong neutrosophic subbigroup H of order $p^{\alpha}$ then we call H a p-Sylow strong neutrosophic subbigroup. If for each prime p, such that $p^{\alpha} / (\langle G \cup I \rangle)$ and $p^{\alpha+1} \nmid o(\langle G \cup I \rangle)$ we have a p-Sylow strong neutrosophic subbigroup then we call $\langle G \cup I \rangle$ to be a Sylow strong neutrosophic bigroup. If $\langle G \cup I \rangle$ has atleast one p-Sylow strong neutrosophic subbigroup then we call $\langle G \cup I \rangle$ a weakly Sylow strong neutrosophic bigroup. If $\langle G \cup I \rangle$ has no p-Sylow strong neutrosophic subbigroup then we say $\langle G \cup I \rangle$ is a Sylow free strong neutrosophic bigroup.*

Next we proceed on to define the notion of Cauchy element and Cauchy neutrosophic element.

**DEFINITION 2.2.17:** *Let ($\langle G \cup I \rangle$, *, o) be a strong neutrosophic bigroup of finite order. $x \in \langle G \cup I \rangle$ is a Cauchy element if $x^n = e$ (e the identity element of $G_i$) and $n / o(\langle G \cup I \rangle)$. An element $y \in \langle G \cup I \rangle$ is a Cauchy neutrosophic element if $y^m = I$ and $m / o(\langle G \cup I \rangle)$. If in a neutrosophic bigroup every element x is such that $x^n = 1$ is a Cauchy element or every y such that $y^m = I$ is a Cauchy neutrosophic element then we call the strong neutrosophic bigroup to be a Cauchy strong neutrosophic bigroup.*

*If ⟨G ∪ I⟩ has atleast a Cauchy element and a Cauchy neutrosophic element then we call ⟨G ∪ I⟩ a weakly Cauchy strong neutrosophic bigroup. If ⟨G ∪ I⟩ has no Cauchy element and no Cauchy strong neutrosophic element then we say ⟨G ∪ I⟩ is a Cauchy free strong neutrosophic bigroup. If the neutrosophic bigroup has only Cauchy elements or Cauchy neutrosophic elements then we call ⟨G ∪ I⟩ to be a semi Cauchy strong neutrosophic bigroup.('or' not used in the mutually exclusive sense).*

It is an easy task to verify all Cauchy strong neutrosophic bigroups are semi Cauchy neutrosophic bigroups.

We can also develop a new type of Sylow substructures of finite strong neutrosophic bigroups. Throughout this section we mean the neutrosophic bigroup is a strong neutrosophic bigroup.

**DEFINITION 2.2.18:** *Let (⟨G ∪ I⟩, o, *) be a neutrosophic bigroup with ⟨G ∪ I⟩ = ⟨$G_1$ ∪ I⟩ ∪ ⟨$G_2$ ∪ I⟩. Let H = $H_1$ ∪ $H_2$ be a neutrosophic subbigroup of ⟨G ∪ I⟩. We say H is a ($p_1$, $p_2$) Sylow strong neutrosophic subbigroup of ⟨G ∪ I⟩ if $H_1$ is a $p_1$-Sylow neutrosophic subgroup of ⟨$G_1$ ∪ I⟩ and $H_2$ is a $p_2$ -Sylow neutrosophic subgroup of ⟨$G_2$ ∪ I⟩.*

We have the following theorem:

**THEOREM 2.2.7:** *Let (⟨G ∪ I⟩, *, o) be a finite neutrosophic bigroup with ⟨G ∪ I⟩ = ⟨$G_1$ ∪ I⟩ ∪ ⟨$G_2$ ∪ I⟩. Let ($p_1$, $p_2$) be any pair of primes such that $p_1^{\alpha} \mid o\langle G_1 \cup I\rangle$ and $p_2^{\beta}$ / o⟨$G_2$ ∪ I⟩. Then (⟨G ∪ I⟩, o, *) has ($p_1$, $p_2$) Sylow strong neutrosophic subbigroup with biorder $p_1^{\alpha} + p_2^{\beta}$.*

*Proof:* Given (⟨G ∪ I⟩, o, *) = (⟨$G_1$ ∪ I⟩, o) ∪ (⟨$G_2$ ∪ I⟩, *) is a neutrosophic bigroup of finite order. Given ($p_1$, $p_2$) are a pair of primes such that $p_1^{\alpha} \mid (\langle G_1 \cup I\rangle)$ and $p_2^{\beta}$/ o (⟨$G_2$ ∪ I⟩) and $H_1$ is a $p_1$-Sylow neutrosophic subgroup of ⟨$G_1$ ∪ I⟩, and $H_2$ is a $p_2$-Sylow neutrosophic subgroup of ⟨$G_2$ ∪ I⟩. Thus H = $H_1$ ∪ $H_2$ is the required $(p_1, p_2)$- Sylow neutrosophic subbigroup of $(\langle G \cup I\rangle$, o, *). We see biorder is $p_1^{\alpha} + p_2^{\beta}$.

Now we proceed on to define conjugate neutrosophic sub-bigroups.

**DEFINITION 2.2.19:** *Let $G = \langle G_1 \cup I, *, \oplus\rangle$, be a neutrosophic bigroup. We say two neutrosophic strong subbigroups $H = H_1 \cup H_2$ and $K = K_1 \cup K_2$ are conjugate neutrosophic subbigroups of $\langle G \cup I\rangle = \langle G_1 \cup I\rangle \cup (\langle G_2 \cup I\rangle)$ if $H_1$ is conjugate to $K_1$ and $H_2$ is conjugate to $K_2$ as neutrosophic subgroups of $\langle G_1 \cup I\rangle$ and $\langle G_2 \cup I\rangle$ respectively.*

Now we proceed on to define normalizer of an element a of the neutrosophic bigroup $\langle G \cup I\rangle$.

**DEFINITION 2.2.20:** *Let $(\langle G \cup I\rangle, *, o) = (\langle G \cup I\rangle, *) \cup (\langle G \cup I\rangle, o)$ be a neutrosophic bigroup. The normalizer of a in $\langle G \cup I\rangle$ is the set $N(a) = \{x \in \langle G \cup I\rangle \mid x a = ax\} = N_1(a) \cup N_2(a) = \{x_1 \in \langle G_1 \cup I\rangle \mid x_1a = ax_1\} \cup \{x_2 \in \langle G_2 \cup I\rangle \mid x_2 a = ax_2\}$ if $a \in (G_1 \cup I) \cap (G_2 \cup I)$ if $a \in (\langle G_1 \cup I\rangle)$ and $a \notin (\langle G_2 \cup I\rangle)$, $N_2(a) = \phi$ like wise if $a \notin \langle G_1 \cup I\rangle$ and $a \in (\langle G_2 \cup I\rangle)$ then $N_1(a) = \phi$ and $N_2(a) = N(a)$ is a neutrosophic subbigroup of $\langle G \cup I\rangle$, clearly $N(a) \neq \phi$ for $I \in N(a)$.*

Now we proceed on to define the notion of right bicoset of a neutrosophic subbigroup.

**DEFINITION 2.2.21:** *Let $(\langle G \cup I\rangle, o, *) = (\langle G_1 \cup I\rangle, o) \cup (\langle G_2 \cup I\rangle, *)$ be a neutrosophic bigroup. Let $H = H_1 \cup H_2$ be a strong neutrosophic subbigroup of $\langle G \cup I\rangle$. The right bicoset of H in $\langle G \cup I\rangle$ for some a in $\langle G \cup I\rangle$ is defined to be $Ha = \{h_1 a \mid h_1 \in H_1$ and $a \in G_1 \cap G_2\} \cup \{h_2 a \mid h_2 \in H_2$ and $a \in G_1 \cap G_2\}$ if $a \in G_1$ and $a \notin G_2$ then $Ha = \{h_1a \mid h_1 \in H_1\} \cup H_2$ and $a \notin G_1$ then $Ha = \{h_2a \mid h_2 \in H_2\} \cup H_1$ and $a \in G_2$*

Thus we have bicosets depends mainly on the way we choose a. If $G_1 \cap G_2 = \phi$ then the bicoset is either $H_1a \cup H_2$ or $H_1 \cup H_2a$.

Similarly we define left bicoset of a neutrosophic subbigroup H of $\langle G \cup I\rangle$.

The next natural question would be if $H = H_1 \cup H_2$ and $K = K_1 \cup K_2$ be any two strong neutrosophic subbigroups of $\langle G \cup I\rangle = \langle G_1 \cup I\rangle \cup (\langle G_2 \cup I\rangle$ how to define HK.

We define $HK = \left\{ h_1k_1,\ h_2k_2 \,\middle|\, \begin{matrix} h_1 \in H_1 & h_2 \in H_2 \\ k_1 \in K_1 & \text{and } K_2 \in K_2 \end{matrix} \right\}$.

The following theorem is left as an exercise for the reader to prove.

**THEOREM 2.2.8:** *Let $\langle G \cup I\rangle = \langle G_1 \cup I\rangle \cup \langle G_2 \cup I\rangle$ be a neutrosophic bigroup $H = H_1 \cup H_2$ and $K = K_1 \cup K_2$ by any two neutrosophic subbigroups of $\langle G \cup I\rangle$. HK will be a neutrosophic subbigroup if and only if $H_1K_1 = K_1H_1$ and $H_2K_2 = K_2H_2$ are neutrosophic subgroups of $\langle G_1 \cup I\rangle$ and $\langle G_2 \cup I\rangle$ respectively.*

We define strong neutrosophic quotient bigroup.

**DEFINITION 2.2.22:** *Let $\langle G \cup I\rangle = \langle G_1 \cup I\rangle \cup \langle G_2 \cup I\rangle$ be a neutrosophic bigroup we say $N = N_1 \cup N_2$ is a neutrosophic normal subbigroup if and only if $N_1$ is a neutrosophic normal subgroup of $\langle G_1 \cup I\rangle$ and $N_2$ is a neutrosophic normal subgroup of $\langle G_2 \cup I\rangle$.*

*We define the strong neutrosophic quotient bigroup $\frac{\langle G \cup I\rangle}{N}$ as $\left[\frac{\langle G_1 \cup I\rangle}{N_1} \cup \frac{\langle G_2 \cup I\rangle}{N_2}\right]$ which is also a neutrosophic bigroup.*

It is important to note only strong neutrosophic bigroups of $\langle G \cup I\rangle$ can have just neutrosophic subbigroups of $\langle G \cup I\rangle$. All properties and definitions can be easily carried out for these neutrosophic subbigroups also. Infact we can have for strong neutrosophic bigroups just subbigroups also. Infact we can have for strong neutrosophic bigroups just subbigroups which are neither neutrosophic nor strong neutrosophic.

## 2.3 Neutrosophic N-groups and their properties

Now we proceed on to define the notion of neutrosophic N-group and give some of its basic properties, such that when they satisfy classical theorem like Lagrange theorem and Sylow theorem.

**DEFINITION 2.3.1:** *Let* $(\langle G \cup I\rangle, *_1,\ldots, *_N)$ *be a nonempty set with N-binary operations defined on it. We say* $\langle G \cup I\rangle$ *is a strong neutrosophic N-group if the following conditions are true.*

i. $\langle G \cup I\rangle = \langle G_1 \cup I\rangle \cup \langle G_2 \cup I\rangle \cup \ldots \cup \langle G_N \cup I\rangle$ *where* $\langle G_i \cup I\rangle$ *are proper subsets of* $\langle G \cup I\rangle$.
ii. $(\langle G_i \cup I\rangle, *_i)$ *is a neutrosophic group,* $i = 1, 2,\ldots, N$. *If in the above definition we have*
    a. $\langle G \cup I\rangle = G_1 \cup \langle G_2 \cup I\rangle \cup \langle G_3 \cup I\rangle \ldots \cup G_K \cup G_{K+1} \cup \ldots \cup G_N$.
    b. $(G_i, *_i)$ *is a group for some i or*
iii. $(\langle G_j \cup I\rangle, *_j)$ *is a neutrosophic group for some j.*

*then we call* $\langle G \cup I\rangle$ *to be a neutrosophic N-group.*

We now illustrate this by examples.

***Example 2.3.1:*** Let $(\langle G \cup I\rangle = (\langle G_1 \cup I\rangle \cup \langle G_2 \cup I\rangle \cup \langle G_3 \cup I\rangle \cup \langle G_4 \cup I\rangle, *_1, *_2, *_3, *_4)$ be a neutrosophic 4 group where

$\langle G_1 \cup I\rangle$ = {1, 2, 3, 4, I, 2I, 3I, 4I} neutrosophic group under multiplication modulo 5.

$\langle G_2 \cup I\rangle$ = {0, 1, 2, 1 + I, 2 + I, 2I + 2, 2I + 1, I, 2I} a neutrosophic group under multiplication modulo 3,

$\langle G_3 \cup I\rangle$ = $\langle Z \cup I\rangle$ a neutrosophic group under addition and

$\langle G_4 \cup I\rangle$ = {(a, b) | a, b ∈ {1, I, 4, 4I} component-wise multiplication modulo 5};

$\langle G \cup I\rangle$ is a strong neutrosophic 4-group.

Now we give an example of a neutrosophic N-group.

***Example 2.3.2:*** Let $(\langle G \cup I\rangle = \langle G_1 \cup I\rangle \cup \langle G_2 \cup I\rangle \cup G_3 \cup G_4, *_1, *_2, *_3, *_4)$ be a neutrosophic 4-group, where

$\langle G_1 \cup I\rangle$ = {1, 2, 3, 4, I, 2I, 3I, 4I} a neutrosophic group under multiplication modulo 5.
$\langle G_2 \cup I\rangle$ = {0, 1, I, 1+I} a neutrosophic group under multiplication modulo 2.
$G_3$ = $S_3$ and
$G_4$ = $A_5$ the alternating group.

$\langle G \cup I\rangle$ is a neutrosophic N-group (N = 4).

Now as in case of other algebraic structures we define the order of a neutrosophic N-group.

**DEFINITION 2.3.2:** *Let $(\langle G \cup I\rangle = \langle G_1 \cup I\rangle \cup \langle G_2 \cup I\rangle \cup \ldots \cup \langle G_N \cup I\rangle, *_1, \ldots, *_N)$ be a strong neutrosophic N-group. By order of $\langle G \cup I\rangle$ we mean the number of distinct elements in $\langle G \cup I\rangle$. If the number of elements in $\langle G \cup I\rangle$ is finite we say $\langle G \cup I\rangle$ is a finite neutrosophic N-group, otherwise infinite and we denote it by $o(\langle G \cup I\rangle)$.*

It is interesting to note that even if one of the $(\langle G_i \cup I\rangle, *_i)$ is infinite then $\langle G \cup I\rangle$ is an infinite neutrosophic N-group. If all the groups $(\langle G_i \cup I\rangle, *_i)$ are finite then the neutrosophic N-group is finite.

Now we proceed on to define neutrosophic sub N-group of a neutrosophic N-group.

**DEFINITION 2.3.3:** *Let $(\langle G \cup I\rangle = \langle G_1 \cup I\rangle \cup \langle G_2 \cup I\rangle \cup \ldots \cup \langle G_N \cup I\rangle, *_1, \ldots, *_N)$ be a neutrosophic N-group. A proper subset*

*$(P, *_1, \ldots, *_N)$ is said to be a neutrosophic sub N-group of $\langle G \cup I\rangle$ if $P = (P_1 \cup \ldots \cup P_N)$ and each $(P_i, *_i)$ is a neutrosophic subgroup (subgroup) of $(G_i, *_i)$, $1 \le i \le N$.*

It is important to note $(P, *_i)$ for no i is a neutrosophic group. We first illustrate this by the following example.

***Example 2.3.3:*** Let $(\langle G \cup I\rangle = \langle G_1 \cup I\rangle \cup \langle G_2 \cup I\rangle \cup \langle G_3 \cup I\rangle, *_1, *_2, *_3)$ be a neutrosophic 3 group, where

$\langle G_1 \cup I\rangle$ = $\langle Q \cup I\rangle$ a neutrosophic group under multiplication,
$\langle G_2 \cup I\rangle$ = {0, 1, 2, 3, 4, I, 2I, 3I, 4I} neutrosophic group under multiplication modulo 5 and
$\langle G_3 \cup I\rangle$ = {0, 1, 2, 1 + I, 2 + I, I, 2I, 1 + 2I, 2 + 2I} a neutrosophic group under multiplication modulo 3.

$\langle G \cup I\rangle$ is a neutrosophic 3 group.
Take

$$P = \left\{\left\{\left\langle \frac{1}{2n}, 2^n, \frac{1}{(2I)^n}, (2I)^n, I, 1\right\rangle\right\}, (1, 4, I, 4I), 1, 2, I, 2I\right\},$$

P is a neutrosophic sub 3-group where

$$P_1 = \left\langle \frac{1}{2n}, 2^n, \frac{1}{(2I)^n}\,(2I)^n, I, 1\right\rangle$$

neutrosophic group under multiplication, (1, 4, I, 4I) a neutrosophic subgroup under multiplication modulo 5 and {1, 2, I, 2I} a neutrosophic subgroup under multiplication modulo 3. P is a strong neutrosophic sub 3-group.

Now consider T = {Q \ {0}, 1 2 3 4, 1, 2}. T is also sub 3-group but T is not a neutrosophic sub 3-group of $\langle G \cup I\rangle$. Also consider X = {Q \ {0}, (I, 2I, 1, 2), (1, 4, I, 4I)}, we see X is a neutrosophic sub 3-group. Thus we see a strong neutrosophic N-group can have 3 types of subgroups viz.

1. Strong neutrosophic sub N-groups.
2. Neutrosophic sub N-groups.
3. Sub N-groups.

Also a neutrosophic N-group can have two types of sub N-groups.

1. Neutrosophic sub N-groups.
2. Sub N-groups.

Now having defined sub N-groups of several types we do not know when the order of the sub N-group will divide the order of the neutrosophic N-group, to this end we make the following definitions:

**DEFINITION 2.3.4:** *Let* $(\langle G \cup I\rangle = \langle G_1 \cup I\rangle \cup \langle G_2 \cup I\rangle \cup \ldots \cup \langle G_N \cup I\rangle, *_1, \ldots, *_N)$ *be a strong neutrosophic N-group of finite order. Suppose* $(P, *_1, \ldots, *_N)$ *is a strong neutrosophic sub N-group of* $\langle G \cup I\rangle$ *such that* $o(P) / o(\langle G \cup I\rangle)$ *then we call P to be a Lagrange strong neutrosophic sub N-group of* $\langle G \cup I\rangle$*. If every strong neutrosophic sub N-group is a Lagrange strong neutrosophic sub N-group then we call* $\langle G \cup I\rangle$ *to be a strong Lagrange strong neutrosophic N-group.*

*If* $\langle G \cup I\rangle$ *has atleast one Lagrange strong neutrosophic sub N-group then we call* $\langle G \cup I\rangle$ *a weakly Lagrange strong neutrosophic N-group. If* $\langle G \cup I\rangle$ *has no Lagrange strong neutrosophic sub N-group then we call* $\langle G \cup I\rangle$ *a Lagrange free strong neutrosophic N-group.*

*If* $\langle G \cup I\rangle$ *is a strong neutrosophic N-group and if* $\langle G \cup I\rangle$ *has a proper subset T such that T is a neutrosophic sub N-group and not a strong neutrosophic sub N-group and* $o(T) / o\langle G \cup I\rangle$ *then we call T a Lagrange neutrosophic sub N-group.*

*If every neutrosophic sub N-group of* $\langle G \cup I\rangle$ *is a Lagrange neutrosophic sub N-group then we call* $\langle G \cup I\rangle$ *a Lagrange neutrosophic N-group. If* $\langle G \cup I\rangle$ *has atleast one Lagrange neutrosophic sub N-group then we call* $\langle G \cup I\rangle$ *a weakly Lagrange neutrosophic N-group.*

*If the strong neutrosophic N-group has no Lagrange neutrosophic sub N-group then we call* $\langle G \cup I\rangle$ *to be a Lagrange free neutrosophic N-group.*

*Similarly we define for a strong neutrosophic N-group* $\langle G \cup I\rangle$*, sub N-groups; we say a proper subset* $V = \{V_1 \cup \ldots \cup V_N, *_1, \ldots, *_N\}$ *of a sub N-group of* $\langle G \cup I\rangle$ *to be a Lagrange sub*

*N-group if o(V) / o(⟨G ∪ I⟩). If ⟨G ∪ I⟩ has atleast one Lagrange sub N-group then we call ⟨G ∪ I⟩ to be a weak Lagrange N-group.*

*If every sub N-group of ⟨G ∪ I⟩ is Lagrange then we call ⟨G ∪ I⟩ to be a Lagrange N-group. If ⟨G ∪ I⟩ has no Lagrange sub N-group then we call ⟨G ∪ I⟩ to be a Lagrange free N-group.*

It is easily verified. All strong neutrosophic N-group of order p, p a prime are Lagrange free strong neutrosophic N-group, Lagrange free neutrosophic N-group and Lagrange free N-group.

Does their exists any relation between Lagrange neutrosophic N-group and Lagrange strong neutrosophic N-group? Give examples of a Lagrange strong neutrosophic N-group which is Lagrange free neutrosophic N-group? Does their exist a strong neutrosophic N-group which is Lagrange free strong neutrosophic N-group and a Lagrange neutrosophic N-group?

On similar lines we can define these substructures in case of neutrosophic N-groups which are not strong neutrosophic N-groups.

We just give some examples.

***Example 2.3.4:*** Let (⟨G ∪ I⟩ = (⟨$G_1$ ∪ I⟩ ∪ ⟨$G_2$ ∪ I⟩ ∪ ⟨$G_3$ ∪ I⟩, $*_1$, $*_2$, $*_3$), where ⟨$G_1$ ∪ I⟩ = {⟨$Z_6$ ∪ I⟩} group under addition modulo 6, $G_2 = A_4$ and $G_3 = \langle g \mid g^{12} = 1\rangle$ a cyclic group of order 12, o(⟨G ∪ I⟩) = 60.

Take P = (⟨$P_1$ ∪ I⟩ ∪ $P_2$ ∪ $P_2$, $*_1$, $*_2$, $*_3$), a neutrosophic sub 3-group where {⟨$P_1$ ∪ I⟩} = {0, 3, 3I, 3 + 3I}, $P_2$ = $\left\{\begin{pmatrix}1&2&3&4\\1&2&3&4\end{pmatrix},\begin{pmatrix}1&2&3&4\\2&1&4&3\end{pmatrix},\begin{pmatrix}1&2&3&4\\4&3&2&1\end{pmatrix},\begin{pmatrix}1&2&3&4\\3&4&1&2\end{pmatrix}\right\}$, $P_3 = \{1, g^6\}$, o(P) = 10, 10 / 60 so P is a Lagrange neutrosophic sub 3-group.

Take T = (⟨$T_1$ ∪ I⟩ ∪ $T_2$ ∪ $T_3$, $*_1$, $*_2$, $*_3$) where ⟨$T_1$ ∪ I⟩ = {0, 3, 3I, 3 + 3I}, $T_2 = P_2$ and $T_3 = \{g^4, g^8, 1\}$, o(T) = 11 and 11 $\not|$ 60 so T is not a Lagrange neutrosophic sub 3-group.

Consider $\{0, 2, 4\} \cup \left\{\begin{pmatrix} 1 & 2 & 3 & 4 \\ 2 & 1 & 4 & 3 \end{pmatrix}, \begin{pmatrix} 1 & 2 & 3 & 4 \\ 1 & 2 & 3 & 4 \end{pmatrix}\right\} \cup \{1,$ $g^3, g^6, g^9\} = W = W_1 \cup W_2 \cup W_3$, o(W) = 9, 9 ∤ 60. Thus this sub 3-group is not a Lagrange sub 3-group. Take V = {0, 2, 4} $\cup \left\{\begin{pmatrix} 1 & 2 & 3 & 4 \\ 1 & 3 & 4 & 2 \end{pmatrix}, \begin{pmatrix} 1 & 2 & 3 & 4 \\ 1 & 4 & 2 & 3 \end{pmatrix}, \begin{pmatrix} 1 & 2 & 3 & 4 \\ 1 & 2 & 3 & 4 \end{pmatrix}\right\} \cup \{1, g^3, g^6,$ $g^9\} = V_1 \cup V_2 \cup V_3$, o(V) = 10, 10 / 60. V is a Lagrange sub 3 group.

**DEFINITION 2.3.5:** *Let $(\langle G \cup I\rangle = \langle G_1 \cup I\rangle \cup \langle G_2 \cup I\rangle \cup \ldots \cup \langle G_N \cup I\rangle, *_1, \ldots, *_N)$ be a strong neutrosophic N-group of finite order. Let p be a prime such that $p^\alpha$ / o $(\langle G \cup I\rangle)$ and $p^{\alpha+1} \nmid$ o $(\langle G \cup I\rangle)$ and if $\langle G \cup I\rangle$ has a strong neutrosophic sub N-group P of order $p^\alpha$ then we call P a p-Sylow strong neutrosophic sub N-group.*

*If for every prime p such that $p^\alpha$ / $o(\langle G \cup I\rangle)$ and $p^{\alpha+1} \nmid$ o $(\langle G \cup I\rangle)$ we have a strong neutrosophic sub N-group then we call $\langle G \cup I\rangle$ a Sylow strong neutrosophic N-group. Now if $\langle G \cup I\rangle$ has for a prime p, $p^\alpha$ / $o(\langle G \cup I\rangle)$ and $p^{\alpha+1} \nmid o(\langle G \cup I\rangle)$ a neutrosophic sub N-group P, of order $p^\alpha$ then we call P, a p-Sylow neutrosophic sub N-group.*

*If for every prime p we have a p-Sylow neutrosophic sub N-group then we call $\langle G \cup I\rangle$ a Sylow neutrosophic N-group. If $\langle G \cup I\rangle$ has atleast one p-Sylow strong neutrosophic sub N-group then we call $\langle G \cup I\rangle$ a weak Sylow strong neutrosophic N-group. If $\langle G \cup I\rangle$ has atleast one p-Sylow neutrosophic sub N-group then we call $\langle G \cup I\rangle$ a weak Sylow neutrosophic N-group. If $\langle G \cup I\rangle$ has p-Sylow strong neutrosophic sub N-group then we call $\langle G \cup I\rangle$ a Sylow free strong neutrosophic N-group. If $\langle G \cup I\rangle$ has no p Sylow neutrosophic sub N-group we call $\langle G \cup I\rangle$ a Sylow free neutrosophic N-group.*

Inter relations connecting them will give many interesting results.

Can one say or prove the existence of a Sylow strong neutrosophic N-group which is not a Sylow neutrosophic N-

group or can one give an example of a Sylow neutrosophic N-group which is not a Sylow strong neutrosophic N-group? Does there exists any special criteria for a strong neutrosophic N-group to be both Sylow strong neutrosophic N-group and Sylow neutrosophic N-group?

We have enlisted a few problems and can be tacked by any interested reader.

Now we give nice example before we proceed on to define Cauchy neutrosophic N-group.

***Example 2.3.5:*** Let $(\langle G \cup I\rangle = (\langle G_1 \cup I\rangle \cup \langle G_2 \cup I\rangle \cup \langle G_3 \cup I\rangle, *_1, *_2, *_3)$ where

| | | |
|---|---|---|
| $\langle G_1 \cup I\rangle$ | = | {1, 2, 3, 4, 5, 6, 7, 8, 9, 10, I, 2I, 3I, 4I 5I, 6I, 7I, 8I, 9I, 10I}, be a neutrosophic group under multiplication modulo 11. |
| $\langle G_2 \cup I\rangle$ | = | {(a, b) \| a, b $\in$ {1, 2, I, 2I}, neutrosophic group under component wise multiplication modulo 3, |
| $\langle G_3 \cup I\rangle$ | = | {0, 1, 2, I, 2I, 1 + I, 1 + 2I, 2 + I, 2 + 2I}, neutrosophic group under multiplication modulo 3, |

$o(G \cup I) = 45$, $5 / 45$ and $5^2 \nmid 45$, $3^2 / 45$and $3^3 \nmid 45$.

Consider $T = (T_1 \cup T_2 \cup T_3, *_1, *_2, *_3)$ a proper subset of $\langle G \cup I\rangle$ where $T_1 = \langle 1, 10, I, 10I\rangle$, $T_2 = \{(1, 1), (I, I)\}$ and $T_3 = (1, 2I, I)$, $o(T) = 9$. So T is a 3-Sylow strong neutrosophic sub 3-group. We see $\langle G \cup I\rangle$ has no 5-Sylow strong neutrosophic sub 3-group, so $\langle G \cup I\rangle$ is only a weakly Sylow strong neutrosophic 3-group.

Take $S = (S_1 \cup S_2 \cup S_3, *_1 *_2 *_3)$ proper subset in $\langle G \cup I\rangle$ where $S_1 = (1, 10)$, $S_2 = (1, 1)$ and $S_3 = (1, I)$. S is a 5 Sylow neutrosophic sub 3-group.

$U = \{U_1 \cup U_2 \cup U_3, *_1, *_2, *_3\}$ where $U_1 = \{1, I, 10, 10I\}$, $U_2 = \{(1,1), (2, 2), (1, 2), (2, 1)\}$ and $U_3 = \{1, I\}$. U is a 3-Sylow neutrosophic sub 3-group. Thus $\langle G \cup I\rangle$ is a Sylow neutrosophic sub 3-group but only a weak Sylow neutrosophic sub 3-group.

Now consider the set $W = \{W_1 \cup W_2 \cup W_3, *_1, *_2, *_3\}$ where $W_1 = \{1, I, 10, 10I\}$, $W_2 = \langle G_2 \cup I\rangle$ and $W_3 = \{0, 2, 1, 2I,$

I}. Clearly $o(W) = 25 = 5^2$. Let $B = \{B_1 \cup B_2 \cup B_3, *_1 *_2 *_3\}$ where $B_1 = \{1, 10\}$, $B_2 = \{G_2 \cup I\}$ and $B_3 = \{G_3 \cup I\}$, $o(B) = 27$ and B is a neutrosophic sub N-group of order $3^3$.

Now these sub N-groups leads us to define some more new concepts.

**DEFINITION 2.3.6:** *Let $(\langle G \cup I\rangle = \langle G_1 \cup I\rangle \cup \langle G_2 \cup I\rangle \cup \ldots \cup \langle G_N \cup I\rangle, *_1, \ldots, *_N)$ be a strong neutrosophic N-group. If $\langle G \cup I\rangle$ is a Sylow strong neutrosophic N-group and if for every prime p such that $p^\alpha / o(\langle G \cup I\rangle)$ and $p^{\alpha+1} \nmid o(\langle G \cup I\rangle)$ we have a strong neutrosophic sub N-group of order $p^{\alpha+t}$ $(t \geq 1)$ then we call $\langle G \cup I\rangle$ a super Sylow strong neutrosophic N-group.*

We define in the same way super Sylow neutrosophic N-group.

**DEFINITION 2.3.7:** *Let $(\langle G \cup I\rangle = \langle G_1 \cup I\rangle \cup \ldots \cup \langle G_N \cup I\rangle, *_1, \ldots, *_N)$ be a strong neutrosophic N-group of finite order. Suppose $\langle G \cup I\rangle$ is a Sylow neutrosophic N-group. If in addition to this for every prime p, $p^\alpha / o(\langle G \cup I\rangle)$ and $p^{\alpha+1} \nmid o(\langle G \cup I\rangle)$ we have a neutrosophic sub N-group of order $p^{\alpha+t}$ $(t \geq 1)$ then we call $\langle G \cup I\rangle$ a super Sylow neutrosophic N-group.*

It is very clear from the definition that every super Sylow strong neutrosophic N-group is always a Sylow strong neutrosophic N-group. However a Sylow strong neutrosophic N-group in general is not always super Sylow strong neutrosophic N-groups. Interested reader can construct examples of these.

Now we can as in case of other structures define the notion of super weakly Sylow strong neutrosophic N-group. Let $(\langle G \cup I\rangle = \langle G_1 \cup I\rangle \cup \ldots \cup \langle G_N \cup I\rangle, *_1, \ldots, *_N)$ be such that $\langle G \cup I\rangle$ is a Sylow strong neutrosophic N-group and it has atleast for one prime p such that $p^\alpha / o(\langle G \cup I\rangle)$ and $p^{\alpha+1} \nmid o(\langle G \cup I\rangle)$ we have a strong neutrosophic sub N-group of order $p^{\alpha+t}$ where $(t \geq 1)$; then we call $\langle G \cup I\rangle$ a super weakly Sylow strong neutrosophic N-group.

It is once again left as an exercise for the interested reader to verify that all weakly super Sylow strong neutrosophic N-groups are not super Sylow strong neutrosophic N-group. Further we see all super Sylow strong neutrosophic N-groups are always weakly super Sylow strong neutrosophic N-group, but the converse in general is not true. Thus we have the following containment relation. Sylow strong neutrosophic N-group $\subseteq$ Super Sylow strong neutrosophic N-group. Clearly the containment relation is strict.

These results can be defined and analyzed / extended in case of neutrosophic N-groups.

Now we proceed on to define the notion of Cauchy elements and Cauchy neutrosophic elements of a strong neutrosophic N-group or just a neutrosophic N-group.

**DEFINITION 2.3.8:** *Consider a strong neutrosophic N-group ($\langle G \cup I\rangle = \langle G_1 \cup I\rangle \cup \langle G_2 \cup I\rangle \cup \ldots \cup \langle G_N \cup I\rangle$, $*_1$, $*_2$,…, $*_N$) where $o(\langle G \cup I\rangle)$ is finite. An element $x \in (\langle G \cup I\rangle)$ such that $x^m = e_i = 1$ – identity element of $\langle G_i \cup I\rangle$) is said to be a Cauchy element of the neutrosophic N-group $\langle G \cup I\rangle$ if $m \,/\, o\,\langle G \cup I\rangle$. If $m \nmid o\,(\langle G \cup I\rangle)$, we say x is not a Cauchy element or x is called as an anti Cauchy element of $\langle G \cup I\rangle$, similarly if $y \in \langle G \cup I\rangle$ is such that $y^n = I$ and if $n \,/\, o(\langle G \cup I\rangle)$ then we call y the Cauchy neutrosophic element of $\langle G \cup I\rangle$. If we have a $y_1 \in \langle G \cup I\rangle$ with $y_1^{n_1} = I$ and $n_1 \nmid o\,(\langle G \cup I\rangle)$ then we call $y_1$ an anti Cauchy neutrosophic element of $\langle G \cup I\rangle$.*

*If $\langle G \cup I\rangle$ has no anti Cauchy elements and no anti Cauchy neutrosophic elements then we call $\langle G \cup I\rangle$ to be a Cauchy strong neutrosophic N-group. If $\langle G \cup I\rangle$ has no anti Cauchy element or ('or' in the mutually exclusive sense) has no anti Cauchy neutrosophic elements then we call $\langle G \cup I\rangle$ a semi Cauchy strong neutrosophic N-group. If $\langle G \cup I\rangle$ has no Cauchy elements and no Cauchy neutrosophic elements then we call $\langle G \cup I\rangle$ a Cauchy free strong neutrosophic N-group.*

*Thus we see all Cauchy strong neutrosophic N-groups are semi Cauchy strong neutrosophic N-groups. However a semi*

*Cauchy strong neutrosophic N-group in general need not be a Cauchy strong neutrosophic N-group.*

We illustrate these situations by the following example:

***Example 2.3.6:*** Let $(\langle G \cup I\rangle = \langle G_1 \cup I\rangle \cup \langle G_2 \cup I\rangle \cup \langle G_3 \cup I\rangle, *_1, *_2, *_3)$ where $\langle G_1 \cup I\rangle = A_4$. $G_2 = \{1, 2, 3, 4, I, 2I, 3I, 4I\}$ and $G_3 = \{g \mid g^{12} = e\}$. $o(\langle G \cup I\rangle) = 32 = 2^5$.

Now consider $x = \begin{pmatrix} 1 & 2 & 3 & 4 \\ 1 & 3 & 4 & 2 \end{pmatrix} \in A_4$. Clearly $x^3 = \begin{pmatrix} 1 & 2 & 3 & 4 \\ 1 & 2 & 3 & 4 \end{pmatrix}$ but $3 \nmid o(\langle G \cup I\rangle)$ i.e. (3, 32) = 1, so x in not a Cauchy element of $\langle G \cup I\rangle$. Take $4I \in \langle G_2 \cup I\rangle$, $(4I)^2 = I$. 2 / $o\langle G \cup I\rangle$ so 4I is a Cauchy neutrosophic element of $\langle G \cup I\rangle$. Consider $(3I) \in \langle G_2 \cup I\rangle$ $(3I)^4 = I$ and 4/o $\langle G \cup I\rangle$ so 3I is also a Cauchy neutrosophic element. Also $2I \in \langle G_2 \cup I\rangle$ and $(2I)^4 = I$ so 2I is also a Cauchy neutrosophic element of $\langle G \cup I\rangle$. Similarly $I \in \langle G_2 \cup I\rangle$ is such that $I^2 = I$ so I is a Cauchy neutrosophic element of $\langle G \cup I\rangle$.

Take $g^4 \in G_3$, $(g^4)^3 = e$ but $3 \nmid o \langle G \cup I\rangle$ so $g^4$ is not a Cauchy element of $\langle G \cup I\rangle$. We see $\langle G \cup I\rangle$ is semi Cauchy neutrosophic N-group for all Cauchy neutrosophic element. However $\langle G \cup I\rangle$ has anti Cauchy elements and no anti Cauchy neutrosophic elements. But $\langle G \cup I\rangle$ also has Cauchy elements.

So we make yet another definition.

**DEFINITION 2.3.9:** *Let $\langle G \cup I\rangle = \{\langle G_1 \cup I\rangle \cup \langle G_2 \cup I\rangle \cup \ldots \cup \langle G_N \cup I\rangle, *_1, \ldots, *_N)$ be a strong neutrosophic N-group of finite order. If $\langle G \cup I\rangle$ has atleast a Cauchy element and a Cauchy neutrosophic element then we call $\langle G \cup I\rangle$ a weakly Cauchy strong neutrosophic N-group.*

However we cannot in general interrelate semi Cauchy neutrosophic N-group and a weakly Cauchy neutrosophic N-group. But we have Cauchy neutrosophic N-groups to be weakly Cauchy neutrosophic N-groups and a weakly Cauchy neutrosophic N-group in general is not a Cauchy neutrosophic N-group. Also all Cauchy neutrosophic N-groups are semi Cauchy neutrosophic N-groups.

But a semi Cauchy neutrosophic N-group in general is not a Cauchy neutrosophic N-group. But it is important to observe that in general a neutrosophic N-group can be both a semi Cauchy neutrosophic N-group (say) with respect to Cauchy neutrosophic elements and a weakly Cauchy neutrosophic N-group with respect to Cauchy elements. It is up to the interested reader to study these properties.

But one result of importance is that a common class of neutrosophic N-groups enjoy it.

**THEOREM 2.3.1:** *Suppose* $(\langle G \cup I\rangle = \langle G_1 \cup I\rangle \cup \ldots \cup \langle G_N \cup I\rangle, *_1, \ldots, *_N)$ *be a neutrosophic N-group of order p, p a prime. Then the following are true.*

i. $\langle G \cup I\rangle$ *is not a Cauchy neutrosophic N-group.*
ii. $\langle G \cup I\rangle$ *is not a semi Cauchy neutrosophic N group.*
iii. $\langle G \cup I\rangle$ *is not a weakly Cauchy neutrosophic N-group.*

*All elements of finite order are anti Cauchy elements and anti Cauchy neutrosophic elements.*

***Example 2.3.7:*** Consider the neutrosophic N-group $(\langle G \cup I\rangle = \langle G_1 \cup I\rangle \cup G_2 \cup G_3, *_1 *_2 *_3)$ where $\langle G_1 \cup I\rangle = \{1, 2, 3, 4, I, 2I, 3I, 4I\}$, $G_2 = S_3$ and $G_3 = \langle g \mid g^5 = e\rangle$. $o(\langle G \cup I\rangle) = 19$.

It is easily verified no element is a Cauchy element or a Cauchy neutrosophic element. Every element is anti Cauchy and anti Cauchy neutrosophic element as $(g^2)^5 = e$, $(g^3)^5 = e$, $(4I)^2 = I$, $(3I)^4 = I$, $\left[\begin{pmatrix} 1 & 2 & 3 \\ 1 & 3 & 2 \end{pmatrix}\right]^2 = \begin{pmatrix} 1 & 2 & 3 \\ 1 & 2 & 3 \end{pmatrix}$.

Hence the claim.

Now we can define homomorphism between neutrosophic N-groups which we call as N-homomorphisms or N-neutrosophic homomorphisms.

**DEFINITION 2.3.10:** *Let* $\{\langle G \cup I\rangle = \{\langle G_1 \cup I\rangle \cup \langle G_2 \cup I\rangle \cup \ldots \cup \langle G_N \cup I\rangle, *_1, \ldots, *_N\}$ *and* $\{\langle H \cup I\rangle = \langle H_1 \cup I\rangle \cup \langle H_2 \cup I\rangle \cup \ldots \cup \langle H_N \cup I\rangle, *_1, *_2, *_3, \ldots, *_N\}$ *be any two neutrosophic N-group such that if* $(\langle G_i \cup I\rangle, *_i)$ *is a neutrosophic group then* $(\langle H_i \cup I\rangle, *_i)$ *is also a neutrosophic group. If* $(G_t, *_t)$ *is a group then* $(H_t, *_t)$ *is a group.*

*A map* $\phi : \langle G \cup I\rangle$ *to* $\langle H \cup I\rangle$ *satisfying* $\phi(I) = I$ *is defined to be a N homomorphism if* $\phi_i = \phi \mid \langle G_i \cup I\rangle$ *(or* $\phi \mid G_i$*) then each* $\phi_i$ *is either a group homomorphism or a neutrosophic group homomorphism, we denote the N-homomorphism by* $\phi = \phi_1 \cup \phi_2 \cup \ldots \cup \phi_N : \langle G \cup I\rangle \to \langle H \cup I\rangle$.

One can define $\phi$ to be an isomorphism if each $\phi_i$ is an isomorphism i.e. $\phi_i$ is one to one and onto. Similarly one can also define the concept of N-automorphisms.

Now as in case of usual homomorphism kernel $\phi_i = \{x_i \in \langle G_i \cup I\rangle \mid \phi_i(x_i) = 0\}$, i = 1, 2, …, N.

So ker $\phi$ = ker$\phi_1 \cup \ldots \cup$ ker $\phi_N$. Clearly ker $\phi$ is never a neutrosophic sub N-group only a sub N-group as $\phi_i(I) = I$. The study of kernel $\phi$ is more interesting for ker $\phi$ is not a normal sub N-group of $\langle G \cup I\rangle$.

Does their exist neutrosophic N-groups so that ker $\phi$ = normal sub N-group?

Now we proceed on to define the notion of $(p_1, p_2, \ldots, p_N)$-Sylow neutrosophic sub N-group of a neutrosophic N-group.

**DEFINITION 2.3.11:** *Let* $\langle G \cup I\rangle = \{\langle G_1 \cup I\rangle \cup \langle G_2 \cup I\rangle \cup \ldots \cup \langle G_N \cup I\rangle, *_1, \ldots, *_N\}$ *be a neutrosophic N-group. Let* $H = \{H_1 \cup H_2 \cup \ldots \cup H_N\}$ *be a neutrosophic sub N-group of* $\langle G \cup I\rangle$*. We say H is a* $(p_1, p_2, \ldots, p_N)$ *Sylow neutrosophic sub N-group of* $\langle G \cup I\rangle$ *if* $H_i$ *is a* $p_i$ *- Sylow neutrosophic subgroup of* $G_i$*. If none of the* $H_i$*'s are neutrosophic subgroups of* $G_i$ *we call H a* $(p_1, p_2, \ldots, p_N)$ *Sylow free sub N-group.*

Now we illustrate this by the following example:

***Example 2.3.8:*** Let $\langle G \cup I \rangle = \{A_4 \cup D_{2.7} \cup \{1, 2, 3, 4, I, 2I, 3I, 4I\}\}$ be a neutrosophic 3-group H =

$$\left[\begin{pmatrix}1&2&3&4\\1&2&3&4\end{pmatrix},\begin{pmatrix}1&2&3&4\\2&1&4&3\end{pmatrix},\begin{pmatrix}1&2&3&4\\4&3&2&1\end{pmatrix},\begin{pmatrix}1&2&3&4\\3&4&1&2\end{pmatrix}\right]$$

$\cup$ {b, $b^2$, …, $b^6$, $b^7$ = e} $\cup$ {1, I, 4, 4I}] is a (2, 7, 2) - Sylow neutrosophic sub 3-group.

Also K = $\left[\begin{pmatrix}1&2&3&4\\1&3&4&2\end{pmatrix},\begin{pmatrix}1&2&3&4\\1&4&2&3\end{pmatrix},\begin{pmatrix}1&2&3&4\\1&2&3&4\end{pmatrix}\right]$

$\cup$ {1, I, 4, 4I}} is a (3, 2, 2) Sylow neutrosophic 3-group. The order of K is 3 + 2 + 4 = 9. The 3-order of H is 4 + 7 + 4 = 15.

Now we proceed on to define the notion of conjugate neutrosophic sub N-groups and sub N-groups of a neutrosophic N-group.

**DEFINITION 2.3.12:** *Let* $(\langle G \cup I\rangle = \langle G_1 \cup I\rangle \cup \langle G_2 \cup I\rangle \cup \ldots \cup \langle G_N \cup I\rangle, *_1, \ldots, *_N)$ *be a strong neutrosophic N-group. Suppose* $H = \{H_1 \cup H_2 \cup \ldots \cup H_N, *_1, \ldots, *_N\}$ *and* $K = \{K_1 \cup K_2 \cup \ldots \cup K_N, *_1, \ldots, *_N\}$ *are two neutrosophic sub N-groups of* $\langle G \cup I\rangle$, *we say K is a strong conjugate to H or H is conjugate to K if each* $H_i$ *is conjugate to* $K_i$ *(i = 1, 2,…,N) as subgroups of* $G_i$.

In the same way we can define conjugate sub N-groups when $\langle G \cup I\rangle$ is just a neutrosophic N-group.

**Chapter Three**

# NEUTROSOPHIC SEMIGROUPS AND THEIR GENERALIZATIONS

The study of classical theorems like Lagrange, Sylow and Cauchy were always associated only with groups, and groups happen to be a prefect structure with no shortcomings. So Lagrange's theorem for finite groups worked well. So also the Sylow and Cauchy theorems. But we are at a loss to know why such theorems were not adopted for semigroups. In this chapter we initiate to adapt to neutrosophic semigroups, neutrosophic bisemigroups and Neutrosophic N-semigroups. We call semigroups, which satisfy Lagrange theorem as Lagrange semigroups and so on. The chapter has three sections. In section 1 we introduce neutrosophic semigroups and study their special properties. In section two neutrosophic bisemigroups are introduced and analyzed. Section three newly defines the concept of neutrosophic N-semigroups and gives some of their properties.

## 3.1 Neutrosophic Semigroups

In this section for the first time we define the notion of neutrosophic semigroups. The notion of neutrosophic subsemigroups, neutrosophic ideals, neutrosophic Lagrange semigroups etc. are introduced for the first time and analyzed. We illustrate them with examples and give some of its properties

**DEFINITION 3.1.1:** *Let S be a semigroup, the semigroup generated by S and I i.e. $S \cup I$ denoted by $\langle S \cup I \rangle$ is defined to be a neutrosophic semigroup.*

It is interesting to note that all neutrosophic semigroups contain a proper subset which is a semigroup.

***Example 3.1.1:*** Let $Z_{12} = \{0, 1, 2, \ldots, 11\}$ be a semigroup under multiplication modulo 12. Let $N(S) = \langle Z_{12} \cup I \rangle$ be the neutrosophic semigroup. Clearly $Z_{12} \subset \langle Z_{12} \cup I \rangle$ and $Z_{12}$ is a semigroup under multiplication modulo 12.

***Example 3.1.2:*** Let Z = {the set of positive and negative integers with zero}, Z is only a semigroup under multiplication. Let $N(S) = \{\langle Z \cup I \rangle\}$ be the neutrosophic semigroup under multiplication. Clearly $Z \subset N(S)$ is a semigroup.

Now we proceed on to define the notion of the order of a neutrosophic semigroup.

**DEFINITION 3.1.2:** *Let N(S) be a neutrosophic semigroup. The number of distinct elements in N(S) is called the order of N(S), denoted by o(N(S)).*

If number of elements in N(S) is finite we call the neutrosophic semigroup to be finite otherwise infinite. The neutrosophic semigroup given in example 3.1.1 is finite where as the neutrosophic semigroup given in example 3.1.2 is of infinite order.

Now we proceed on to define the notion of neutrosophic subsemigroup of a neutrosophic semigroup N(S).

**DEFINITION 3.1.3:** *Let N(S) be a neutrosophic semigroup. A proper subset P of N(S) is said to be a neutrosophic subsemigroup, if P is a neutrosophic semigroup under the operations of N (S). A neutrosophic semigroup N(S) is said to have a subsemigroup if N(S) has a proper subset which is a semigroup under the operations of N(S).*

It is interesting to note a neutrosophic semigroup may or may not have a neutrosophic subsemigroup but it will always have a subsemigroup.

Now we proceed on to illustrate these by the following examples.

***Example 3.1.3:*** Let $Z^+ \cup \{0\}$ denote the set of positive integers together with zero. $\{Z^+ \cup \{0\}, +\}$ is a semigroup under the binary operation '+'. Now let $N(S) = \langle Z^+ \cup \{0\}^+ \cup \{I\}\rangle$. N(S) is a neutrosophic semigroup under '+'. Consider $\langle 2Z^+ \cup I\rangle = P$, P is a neutrosophic subsemigroup of N(S). Take $R = \langle 3Z^+ \cup I\rangle$; R is also a neutrosophic subsemigroup of N(S).

Now we have the following interesting theorem.

**THEOREM 3.1.1:** *Let N(S) be a neutrosophic semigroup. Suppose $P_1$ and $P_2$ be any two neutrosophic subsemigroups of N(S) then $P_1 \cup P_2$ (i.e. the union) the union of two neutrosophic subsemigroups in general need not be a neutrosophic subsemigroup.*

*Proof:* We prove this by using the following example. Let $Z^+$ be the set of positive integers; $Z^+$ under '+' is a semigroup.

Let $N(S) = \langle Z^+ \cup I\rangle$ be the neutrosophic semigroup under '+'. Take $P_1 = \{\langle 2Z \cup I\rangle\}$ and $P_2 = \{\langle 5Z \cup I \rangle\}$ to be any two neutrosophic subsemigroups of N(S). Consider $P_1 \cup P_2$ we see $P_1 \cup P_2$ is only a subset of N (S) for $P_1 \cup P_2$ is not closed under the binary operation '+'. For take $2 + 4I \in P_1$ and $5 + 5I \in P_2$. Clearly $(2 + 5) + (4I + 5I) = 7 + 9I \notin P_1 \cup P_2$. Hence the claim.

We proved the theorem 3.1.1 mainly to show that we can give a nice algebraic structure to P1 $\cup$ P2 viz. neutrosophic bisemigroups defined in section 3.2. Now we proceed on to define the notion of neutrosophic monoid.

**DEFINITION 3.1.4:** *A neutrosophic semigroup N(S) which has an element e in N(S) such that $e * s = s * e = s$ for all $s \in N(S)$, is called as a neutrosophic monoid.*

It is interesting to note that in general all neutrosophic semigroups need not be neutrosophic monoids.

We illustrate this by an example.

***Example 3.1.4:*** Let N (S) = $\langle Z^+ \cup I\rangle$ be the neutrosophic semigroup under '+'. Clearly N(S) contains no e such that s + e = e + s = s for all s $\in$ N (S). So N (S) is just a neutrosophic semigroup and not a neutrosophic monoid.

Now we give an example of a neutrosophic monoid.

***Example 3.1.5:*** Let N (S) = $\langle Z^+ \cup I\rangle$ be a neutrosophic semigroup generated under '×'. Clearly 1 in N (S) is such that 1 × s = s for all s $\in$ N (S). So N (S) is a neutrosophic monoid.

It is still interesting to note the following:

1. From the examples 3.1.3 and 3.1.4 we have taken the same set $\langle Z^+ \cup I\rangle$ with respect the binary operation '+', $\langle Z^+ \cup I\rangle$ is only a neutrosophic semigroup but $\langle Z^+ \cup I\rangle$ under the binary operation × is a neutrosophic monoid.
2. In general all neutrosophic monoids need not have its neutrosophic subsemigroups to be neutrosophic submonoids. First to this end we define the notion of neutrosophic submonoid.

**DEFINITION 3.1.5:** *Let N(S) be a neutrosophic monoid under the binary operation *. Suppose e is the identity in N(S), that is s * e = e * s = s for all s $\in$ N(S). We call a proper subset P of N(S) to be a neutrosophic submonoid if*

*i. P is a neutrosophic semigroup under '*'.*
*ii. e $\in$ P, i.e., P is a monoid under '*'.*

***Example 3.1.6:*** Let N(S) = $\langle Z \cup I\rangle$ be a neutrosophic semigroup under '+'. N(S) is a monoid. P = $\langle 2Z^+ \cup I\rangle$ is just a neutrosophic subsemigroup whereas T = $\langle 2Z \cup I\rangle$ is a neutrosophic submonoid. Thus a neutrosophic monoid can have both neutrosophic subsemigroups which are different from the neutrosophic submonoids.

Now we proceed on to define the notion of neutrosophic ideals of a neutrosophic semigroup.

**DEFINITION 3.1.6**: *Let N(S) be a neutrosophic semigroup under a binary operation *. P be a proper subset of N(S). P is said to be a neutrosophic ideal of N(S) if the following conditions are satisfied.*

*i. P is a neutrosophic semigroup.*

*ii. for all $p \in P$ and for all $s \in N(S)$ we have $p * s$ and $s * p$ are in P.*

*Note:* One can as in case of semigroups define the notion of neutrosophic right ideal and neutrosophic left ideal. A neutrosophic ideal is one which is both a neutrosophic right ideal and a neutrosophic left ideal. In general a neutrosophic right ideal need not be a neutrosophic left ideal.

Now we proceed on to give example to illustrate these notions.

***Example 3.1.7:*** Let $N(S) = \langle Z \cup I \rangle$ be the neutrosophic semigroup under multiplication.

Take P to be a proper subset of N(S) where $P = \langle 2Z \cup I \rangle$. Clearly P is a neutrosophic ideal of N(S). Since N(S) is a commutative neutrosophic semigroup we have P to be a neutrosophic ideal.

*Note:* A neutrosophic semigroup N(S) under the binary operation * is said to be a neutrosophic commutative semigroup if $a * b = b * a$ for all $a, b \in N(S)$.

***Example 3.1.8:*** Let $N(S) = \left\{ \begin{pmatrix} a & b \\ c & d \end{pmatrix} / a, b, c, d, \in \langle Z \cup I \rangle \right\}$ be a neutrosophic semigroup under matrix multiplication.

Take $P = \left\{\begin{pmatrix} x & y \\ 0 & 0 \end{pmatrix} / x, y \in \langle Z \cup I \rangle\right\}$.

Clearly $\begin{pmatrix} a & b \\ c & d \end{pmatrix}\begin{pmatrix} x & y \\ 0 & 0 \end{pmatrix} \notin P$ ; but $\begin{pmatrix} x & y \\ 0 & 0 \end{pmatrix}\begin{pmatrix} a & b \\ c & d \end{pmatrix} \in P$.

Thus P is only a neutrosophic right ideal and not a neutrosophic left ideal of N (S).

Now we proceed on to define the notion of neutrosophic maximal ideal and neutrosophic minimal ideal of a neutrosophic semigroup N(S).

**DEFINITION 3.1.7:** *Let N(S) be a neutrosophic semigroup. A neutrosophic ideal P of N(S) is said to be maximal if $P \subset J \subset N(S)$, J a neutrosophic ideal then either J = P or J = N(S). A neutrosophic ideal M of N(S) is said to be minimal if $\phi \neq T \subseteq M \subseteq N(S)$ then T = M or $T = \phi$.*

We cannot always define the notion of neutrosophic cyclic semigroup but we can always define the notion of neutrosophic cyclic ideal of a neutrosophic semigroup N(S).

**DEFINITION 3.1.8:** *Let N(S) be a neutrosophic semigroup. P be a neutrosophic ideal of N (S), P is said to be a neutrosophic cyclic ideal or neutrosophic principal ideal if P can be generated by a single element.*

We proceed on to define the notion of neutrosophic symmetric semigroup.

**DEFINITION 3.1.9:** *Let S(N) be the neutrosophic semigroup. If S(N) contains a subsemigroup isomorphic to S(n) i.e. the semigroup of all mappings of the set (1, 2, 3, ..., n) to itself under the composition of mappings, for a suitable n then we call S (N) the neutrosophic symmetric semigroup.*

**Remark:** We cannot demand the subsemigroup to be neutrosophic, it is only a subsemigroup.

**DEFINITION 3.1.10:** *Let N(S) be a neutrosophic semigroup. N(S) is said to be a neutrosophic idempotent semigroup if every element in N (S) is an idempotent.*

***Example 3.1.9:*** Consider the neutrosophic semigroup under multiplication modulo 2, where N (S) = {0, 1, I, 1 + I}. We see every element is an idempotent so N (S) is a neutrosophic idempotent semigroup.

Next we proceed on to define the notion of weakly neutrosophic idempotent semigroup.

**DEFINITION 3.1.11:** *Let N(S) be a neutrosophic semigroup. If N(S) has a proper subset P where P is a neutrosophic subsemigroup in which every element is an idempotent then we call P a neutrosophic idempotent subsemigroup.*

*If N(S) has at least one neutrosophic idempotent subsemigroup then we call N(S) a weakly neutrosophic idempotent semigroup.*

We illustrate this by the following example:

***Example 3.1.10:*** Let N(S) = {0, 2, 1, I, 2I, 1 + I, 2 + 2 I, 1 + 2I, 2 + I} be the neutrosophic semigroup under multiplication modulo 3.

Take P = {1, I, 1 + 2I, 0}; P is a neutrosophic idempotent subsemigroup. So N(S) is only a weakly neutrosophic idempotent semigroup. Clearly N(S) is not a neutrosophic idempotent semigroup as $(1 + I)^2 = 1$ which is not an idempotent of N(S).

**DEFINITION 3.1.12:** *Let N(S) be a neutrosophic semigroup (monoid). An element $x \in N(S)$ is called an element of finite order if $x^m = e$ where e is the identity element in N(S) i.e. (se = es = s for all $s \in S$ ) (m the smallest such integer).*

**DEFINITION 3.1.13:** *Let N (S) be a neutrosophic semigroup (monoid) with zero. An element $0 \neq x \in N(S)$ of a neutrosophic semigroup is said to be a zero divisor if there exist $0 \neq y \in N(S)$ with $x \cdot y = 0$. An element $x \in N(S)$ is said to be invertible if there exist $y \in N(S)$ such that $xy = yx = e$ ($e \in N(S)$, is such that $se = es = s$ for all $s \in N(S)$).*

***Example 3.1.11:*** Let N (S) = {0, 1, 2, I, 2I, 1 + I, 2 + I, 1 + 2I, 2 + 2I} be a neutrosophic semigroup under multiplication modulo 3. Clearly (1 + I) ∈ N (S) is invertible for (1 + I) (1 + I) = 1 (mod 3). (2 + 2I) is invertible for $(2 + 2I)^2 = 1$ (mod 3). N(S) also has zero divisors for (2 + I) I = 2I + I = 0(mod 3). Also (2 + I ) 2 I = 0 (mod 3) is a zero divisor.

Thus this neutrosophic semigroup has idempotents, units and zero divisors. One can define several other properties of semigroups to neutrosophic semigroups as a matter of routine.

## 3.2 Neutrosophic Bisemigroups and their Properties

In this section Neutrosophic bisemigroups are defined analogous to bisemigroups by taking atleast one of the semigroups to be neutrosophic. A stronger version of this viz. strong neutrosophic bisemigroups are defined. Substructures like neutrosophic subbisemigroups neutrosophic biideals are introduced. Condition for a finite neutrosophic bisemigroup to be Lagrange free is also obtained.

Now we proceed on to define the notion of neutrosophic bisemigroup.

**DEFINITION 3.2.1:** *Let (BN(S), *, o) be a nonempty set with two binary operations * and o. (BN(S), *, o) is said to be a neutrosophic bisemigroup if $BN(S) = P_1 \cup P_2$ where atleast one of ($P_1$, *) or ($P_2$, o) is a neutrosophic semigroup and other is just a semigroup. $P_1$ and $P_2$ are proper subsets of BN(S), i.e. $P_1 \not\subseteq P_2$.*

***Example 3.2.1:*** Let (BN(S), *, o) = {0, 1, 2, 3, I, 2I, 3I, S(3), *, o} = ($P_1$, *) ∪ ($P_2$, o) where ($P_1$, *) = {0, 1, 2, 3, I, 2I, 3I, *} and ($P_2$, o) = (S(3), o). Clearly ($P_1$, *) is a neutrosophic semigroup under multiplication modulo 4. ($P_2$, o) is just a semigroup. Thus (BN(S), *, o) is a neutrosophic bisemigroup.

If both ($P_1$, *) and ($P_2$, o) in the above definition are neutrosophic semigroups then we call (BN (S), *, o) a strong neutrosophic bisemigroup. All strong neutrosophic bisemigroups are trivially neutrosophic bisemigroups.

We now give an example of a strong neutrosophic bisemigroup.

***Example 3.2.2:*** Let (BN (S), *, o) be a nonempty set such that BN(S) = {0,1, I, 1 + I, ⟨Z ∪ I⟩, *, o} where $P_1$ = {0, 1, I, 1 + I} and $P_2$ = {⟨ Z ∪ I ⟩, o}; BN(S) is a strong neutrosophic bisemigroup.

We now proceed on to define the notion of neutrosophic subbisemigroup and neutrosophic strong subbisemigroup of a neutrosophic bisemigroup and neutrosophic strong bisemigroup.

**DEFINITION 3.2.2:** *Let (BN (S) = $P_1$ ∪ $P_2$; o, *) be a neutrosophic bisemigroup. A proper subset (T, o, *) is said to be a neutrosophic subbisemigroup of BN (S) if*

i. *$T = T_1 \cup T_2$ where $T_1 = P_1 \cap T$ and $T_2 = P_2 \cap T$ and*
ii. *At least one of ($T_1$, o) or ($T_2$, *) is a neutrosophic semigroup.*

*Note:* We can define for a neutrosophic bisemigroup just a subbisemigroup which need not be a neutrosophic subbisemigroup. But a neutrosophic bisemigroup cannot have a proper neutrosophic strong subbisemigroup.

Now we proceed on to define substructures of the strong neutrosophic bisemigroup.

**DEFINITION 3.2.3:** *Let (BN(S) = $P_1$ ∪ $P_2$, o, *) be a neutrosophic strong bisemigroup. A proper subset T of BN (S) is called the strong neutrosophic subbisemigroup if $T = T_1 \cup T_2$* *with $T_1 = P_1 \cap T$ and $T_2 = P_2 \cap T_2$ and if both $(T_1, *)$ and $(T_2, o)$ are neutrosophic subsemigroups of $(P_1, *)$ and $(P_2, o)$ respectively. We call $T = T_1 \cup T_2$ to be a neutrosophic strong subbisemigroup, if atleast one of $(T_1, *)$ or $(T_2, o)$ is a semigroup then $T = T_1 \cup T_2$ is only a neutrosophic subsemigroup.*

We illustrate this by the following:

***Example 3.2.3:*** Let BN(S) = {0, 1, 2, I, 2I, ⟨Z ∪ I⟩, ×, +} be a neutrosophic strong bisemigroup. Take T = {0, I, 2I, ⟨2Z ∪ I⟩, ×, +} ⊂ BN(S), T is a neutrosophic strong subbisemigroup. Now consider P = {0, 1, 2, ⟨5Z ∪ I⟩, ×, +} = $P_1$ = {0, 1, 2, ×} ∪ $P_2$ = {⟨5Z ∪ I⟩, +} is only a neutrosophic subbisemigroup of BN(S). If we let L = {0, 1, 2, Z, ×, +} = $L_1$ = {0, 1, 2, ×} ∪ (Z, +) = $L_2$ = L, is only just a subbisemigroup.

**THEOREM 3.2.1:** *Let (BN(S), *, o) be a neutrosophic bisemigroup, then B(N(S)) cannot have a neutrosophic strong subbisemigroup.*

*Proof:* Now we are given B(N(S)) = $P_1 \cup P_2$ where only one of $(P_1, *)$ or $(P_2, o)$ is a neutrosophic semigroup so if they have subsemigroup only one of them will be neutrosophic semigroup so it is impossible to have both of them to be neutrosophic semigroups. Hence a neutrosophic bisemigroup cannot have strong neutrosophic subbisemigroup. Hence the claim.

Now we proceed on to define the notion of neutrosophic strong biideal and neutrosophic biideal of a neutrosophic strong bisemigroups and neutrosophic bisemigroups respectively.

**DEFINITION 3.2.4:** *Let (BN (S), *, o) be a strong neutrosophic bisemigroup where $BN(S) = P_1 \cup P_2$ with $(P_1, *)$ and $(P_2, o)$ be any two neutrosophic semigroups. Let J be a proper subset of BN(S) where $I = I_1 \cup I_2$ with $I_1 = J \cap P_1$ and $I_2 = J \cap P_2$ are neutrosophic ideals of the neutrosophic semigroups $P_1$ and $P_2$ respectively. Then I is called or defined as the strong neutrosophic biideal of B(N(S)).*

**DEFINITION 3.2.5:** *Let (BN(S) = $P_1 \cup P_2$ *, o) be any neutrosophic bisemigroup. Let J be a proper subset of B(NS) such that $J_1 = J \cap P_1$ and $J_2 = J \cap P_2$ are ideals of $P_1$ and $P_2$ respectively. Then J is called the neutrosophic biideal of BN(S).*

It is important and interesting to note that as in case of neutrosophic subbisemigroups we cannot in general have a biideal for these neutrosophic bisemigroups. Also a neutrosophic strong bisemigroup can never have a nontrivial neutrosophic biideal or a neutrosophic bisemigroup cannot have a biideal Thus this property distinguishes both subbisemigroup and biideal structures.

The concept of maximality, minimality, quasi maximality and quasi minimality will be now introduced.

**DEFINITION 3.2.6:** *Let (BN(S) $P_1 \cup P_2$ *, o) be a neutrosophic strong bisemigroup. Suppose I is a neutrosophic strong biideal of BN(S) i.e. $I_1 = P_1 \cap I$ and $I_2 = P_2 \cap I$ we say I is a neutrosophic strong maximal biideal of B (N(S)) if $I_1$ is the maximal ideal of $P_1$ and $I_2$ is the maximal ideal of $P_2$.*
*If only one of $I_1$ or $I_2$ alone is maximal then we call I to be a neutrosophic strong quasi maximal biideal of BN (S).*

We now illustrate this by the following example.

***Example 3.2.4:*** Let BN(S) = ({⟨Z ∪ I⟩, 0, 1, 2, I, 2I}, +, × (× under multiplication modulo 3)) be a neutrosophic strong bisemigroup. Take T = {⟨2Z ∪ I⟩, 0, I, 1, 2I, +, ×} = $I_1 \cup I_2$ where $I_1$ = {⟨2Z ∪ I⟩, +} and $I_2$ = {0, I, 1, 2I, ×, multiplication modulo 3}, T is a neutrosophic strong biideal which is maximal. Take J = {⟨8Z ∪ I⟩, {0, 1, I, 2I}, + ×} = $J_1 \cup J_2$ where $J_1$ = {⟨8Z ∪ I⟩, +} and $J_2$ = {0, 1, I, 2I, ×}. Clearly J, is not a maximal ideal of $P_1$ only $J_2$ is a maximal ideal. So J is a neutrosophic strong quasi maximal biideal of BN (S).

Now we proceed on to define these concepts in the case of neutrosophic bisemigroup.

**DEFINITION 3.2.7:** *Let $(BN(S) = P_1 \cup P_2$ *, o) be a neutrosophic semigroup. A neutrosophic biideal $I = I_1 \cup I_2$ of BN(S) is said to be a neutrosophic maximal biideal of B(N(S)) if $I_1 = I \cap P_1$ is a maximal ideal of $P_1$ and $I_2 = I \cap P_2$ is a maximal ideal of $P_2$.*

*Now if only one of $I_1$ or $I_2$ alone is maximal and other not maximal then we call I to be only a neutrosophic quasi maximal biideal of B N(S).*

Now we illustrate this by the following example.

***Example 3.2.5:*** Let B N(S) = {⟨Z ∪ I⟩, 0, 1, 2, 3, 4, 5, 6, …, 11, +, ×} = $P_1 \cup P_2$ = {⟨Z ∪ I⟩} ∪ {0, 1, 2, …, 11} where $P_1$ is a semigroup under '+' and $P_2$ is just a semigroup under multiplication modulo 12. Take U = {⟨2I ∪ I⟩, +} ∪ {0, 2, 4, 6, 8, 10, ×} = $I_1 \cup I_2 \subset P_1 \cup P_2$, U is neutrosophic maximal biideal of BN(S), for $I_1$ and $I_2$ are maximal ideals of $P_1$ and $P_2$ respectively. Now consider U = $T_1 \cup T_2$ = {⟨3I ∪ I⟩, +} ∪ {0, 6, ×} U is only a neutrosophic quasi maximal biideal of BN(S) as {0, 6, ×} is not a maximal ideal of $P_2$.

The definition of minimal biideal is left as an exercise for the reader.

*Note:* Union of any two neutrosophic biideals in general is not a neutrosophic biideal. This is true of neutrosophic strong biideals.

**DEFINITION 3.2.8:** *Let (B, + o) be a non empty set with two binary operations we call B a strong neutrosophic bimonoid if the following conditions are satisfied*

i. *$B = B_1 \cup B_2$ where $B_1$ and $B_2$ are proper subsets of B.*
ii. *$(B_1, +)$ is a neutrosophic monoid.*
iii. *$(B_2, o)$ is a neutrosophic monoid.*

***Example 3.2.6:*** Let (B = $B_1 \cup B_2$ , ×, o) be a non empty set with two binary operations, where $B_1$ = {$Z_6$, semigroup under multiplication modulo 6}, $B_2$ = {0, 1, 2, …, 7, I, 2I, 3I, …, 7I} be a neutrosophic semigroup under multiplication modulo 8. B is a neutrosophic bisemigroup.

***Example 3.2.7:*** Let B = $\{B_1 \cup B_2, \times, +\}$ be a strong neutrosophic bisemigroup where $B_1 = \langle Z \cup I\rangle$ semigroup under multiplication and $B_2 = \{0, 1, 2, 3, I, 2I, 3I\}$ a neutrosophic semigroup under multiplication modulo 4.

*Remark:* All strong neutrosophic bisemigroups are neutrosophic bisemigroups, but clearly from the above example the converse is not true.

***Example 3.2.8:*** Let $(B = B_1 \cup B_2, *, o)$ be a strong neutrosophic bisemigroup where $B_1 = \{\langle Z_4 \cup I\rangle = \{0, 1, 2, 3, I, 2I, 3I\}\}$ a neutrosophic semigroup under multiplication modulo 4 and $B_2 = \langle Z \cup I\rangle$ a neutrosophic semigroup under multiplication. $P = P_1 \cup P_2$ where $P_1 = \{0, 1, 2, I, 2I\}$ is a neutrosophic subsemigroup of $B_1$ and $P_2 = \langle 3Z \cup I\rangle$ is a neutrosophic subsemigroup of $B_2$. Hence B is a strong neutrosophic subbisemigroup of B.
(2) Let $T = T_1 \cup T_2$ where $T_1 = \{0, 1, 2, 3\}$ semigroup under multiplication modulo 4 and $T_2 = \langle 2Z \cup I\rangle$ neutrosophic semigroup under multiplication. T is a neutrosophic sub bisemigroup of B.
(3) Let $R = R_1 \cup R_2$ where $R_1 = \{0, 1, 2\}$ semigroup under multiplication modulo 4 and $R_2 = 2Z$ semigroup under multiplication, R is a subbisemigroup of B.

Thus we see in general a strong neutrosophic bisemigroup can have all 3 types of substructures but a neutrosophic bisemigroup can have only a neutrosophic sub bisemigroup and subbisemigroup. Now we proceed on to define ideals in neutrosophic bisemigroups.

Now we illustrate these by the following examples.

***Example 3.2.9:*** Let $(B = B_1 \cup B_2, o, *)$ be a strong neutrosophic bisemigroup where $B_1 = \langle Z \cup I\rangle$ neutrosophic semigroup under multiplication and $B_2 = \{0, 1, 2, 3, 4, 5, I, 2I, 3I, 4I, 5I\}$ be a neutrosophic semigroup under multiplication modulo 6. $P = P_1 \cup P_2$ where $P_1 = \langle 2Z \cup I \rangle$ is a neutrosophic ideal of $B_1$ and $P_2 = \{0, 2, 2I, 4, 4I\}$ is a neutrosophic ideal of $B_2$. Thus P is a strong neutrosophic biideal of B.

***Example 3.2.10:*** Let $(B = B_1 \cup B_2, *, o)$ be a neutrosophic bisemigroup where $B_1 = Z$, a semigroup under multiplication and $B_2 = \{0, 1, 2, 3, I, 2I, 3I\}$ a neutrosophic semigroup under multiplication modulo 4. Take $J = J_1 \cup J_2$ where $J_1 = \{3Z\}$ is an ideal of $B_1$ and $J_2 = \{0, 2, 2I\}$ is a neutrosophic ideal of $B_2$. Thus J is a neutrosophic biideal.

We can define several other properties we restrain our selves to only a few of them.

**DEFINITION 3.2.9:** *Let $B = (B_1 \cup B_2, *, o)$ be a neutrosophic bisemigroup. The number of distinct elements of B is called the order of B denoted by o(B). If o(B) is finite we call B a finite neutrosophic bisemigroup. If the order of B is infinite we define B to be a infinite neutrosophic bisemigroup.*

The neutrosophic bisemigroup given in example is an infinite neutrosophic bisemigroup where as the one given below is a finite neutrosophic bisemigroup.

***Example 3.2.11:*** Let $B = (B_1 \cup B_2, o, *)$ be a neutrosophic bisemigroup. Here $B_1 = \{0, 1, 2, 3, 4, 5\}$ is a semigroup under multiplication modulo 6. $B_2 = \{0, 1, 2, 3, I, 2I, 3I\}$ is a neutrosophic semigroup under multiplication modulo 4. Thus B is a finite neutrosophic bisemigroup. o(B) = 13.

Now proceed on to define the notion of neutrosophic Lagrange bisemigroup.

**DEFINITION 3.2.10:** *Let $(B = B_1 \cup B_2, *, o)$ be a neutrosophic bisemigroup of finite order. If the order of every neutrosophic subbisemigroups $P = P_1 \cup P_2$ divides the order of B then we call B to be a Lagrange neutrosophic bisemigroup. If at least one neutrosophic subbisemigroup exists such that its order* *divides the order of B we call B to be a weak Lagrange neutrosophic bisemigroup. We call the neutrosophic subbisemigroup whose order divides the order of B to be a Lagrange neutrosophic subbisemigroup.*

*If B has proper neutrosophic subbisemigroups but the order of none of them divide the order of B then we call B to be a Lagrange free neutrosophic bisemigroup.*

We illustrate these by the following example.

***Example 3.2.12:*** Let B = ($B_1 \cup B_2$, o, *) be a finite neutrosophic bisemigroup, where $B_1$ = {$Z_{12}$, the semigroup under multiplication modulo 12} and $B_2$ = {0, 1, 2, 3, 4, I, 2I, 2I, 3I, 4I} a neutrosophic semigroup under multiplication modulo 5. o(B) = 21.

Take P = $P_1 \cup P_2$ where $P_1$ = {0, 6} and $P_2$ = {0, 1, 4, I, 4I}. Clearly o(P) = 7. P is a neutrosophic bisemigroup and 7 / 21. Consider T = $T_1 \cup T_2$ where $T_1$ = {0, 2, 4, 6, 8, 10} and $T_2$ = {0, 1, 4, I, 4I}. T is a neutrosophic bisemigroup and o(T) = 11 but o(T) $\not|$ o(B). Thus B is only a weakly Lagrange neutrosophic bisemigroup.

We give an example of a Lagrange free neutrosophic bisemigroup.

***Example 3.2.13:*** Let (B = $B_1 \cup B_2$, *, o) be a neutrosophic bisemigroup where $B_1$ = {0, 1, 2, 3, 4, 5, I, 2I, 3I, 4I, 5I} a neutrosophic semigroup under multiplication modulo 6. Take $B_2$ = {$Z_{12}$} semigroup under multiplication modulo 12. o(B) = 23 a prime.

Take T = $T_1 \cup T_2$ where $T_1$ = {0, 1, 5, I, 5I} and $T_2$ = {0, 2, 4, 6, 8, 10}. T is a neutrosophic subbisemigroup and o(T) = 11, o(T) $\not|$ o(B) i.e. 11 $\not|$ 23. In fact since the order of B is a prime the order of none of the neutrosophic subbisemigroup will divide the order of B.

In view of this we have the following result.

**THEOREM 3.2.2:** *Let B = ($B_1 \cup B_2$, *, o) be a neutrosophic bisemigroup of finite order say p, p a prime. Then B is a Lagrange free bisemigroup.*

*Proof:* Given B is a finite neutrosophic bisemigroup of order p, p a prime. If P = $P_1 \cup P_2$ is any neutrosophic subbisemigroup then clearly (o (P), p) = 1. Hence B is a Lagrange free neutrosophic bisemigroup.

Now we proceed on to define Cauchy neutrosophic bisemigroup.

**DEFINITION 3.2.11:** *Let B = $B_1 \cup B_2$ be a neutrosophic bimonoid of finite order n. If for all those $x_i \in B$ such that $x_i^t$ = $e_i$ we have t / n (i = 1 or 2). Then we call B to be a Cauchy neutrosophic bisemigroup. (It is to be noted we can have x in B such that $x^2 = x$ or $x^r = 0$ that is why we take only those $x_i$ in $B_i$ where $e_i$ is the identity in $B_i$). If there is atleast one element $x_i$ such that $x_i^r = e_i$ and r / n then we call B to be a weakly Cauchy neutrosophic bisemigroup. If no x exist satisfying this condition we call B to be a Cauchy free neutrosophic bisemigroup.*

Now we proceed to define the concept of Sylow neutrosophic bisemigroups.

**DEFINITION 3.2.12:** *Let (B = $B_1 \cup B_2$, o, *) be a neutrosophic bisemigroup of finite order n. If for every prime p such that $p^\alpha/n$ and $p^{\alpha+1} \nmid n$ we have a neutrosophic subbisemigroup of order $p^\alpha$ we call B a Sylow neutrosophic bisemigroup. If we have atleast one p such that $p^\alpha$ / n and $p^{\alpha+1} \nmid$ n and B has a neutrosophic subbisemigroup of order $p^\alpha$ we call B a weak Sylow neutrosophic bisemigroup.*

*Suppose for no prime p such that $p^\alpha$ / n and $p^{\alpha+1} \nmid$ n we have no neutrosophic subsemigroup of order $p^\alpha$ we call B a Sylow free neutrosophic bisemigroup. We call the neutrosophic subbisemigroup P of order $p^\alpha$ where $p^\alpha$ / o(B) and $p^{\alpha+1} \nmid$ o(B) to be the p-Sylow neutrosophic subbisemigroup.*

***Example 3.2.14:*** Let B = ($B_1 \cup B_2$, *, o) where $B_1$ = {$Z_{18}$, a semigroup under multiplication modulo 18} and $B_2$ = {0, 1, 2, …, 7, I, 2I, …, 7I} semigroup under multiplication modulo 8. B is a neutrosophic bisemigroup of order 33. 3/33 but $3^2 \nmid 33$ also 11/33 and $11^2 \nmid 33$.

Clearly B cannot have a nontrivial neutrosophic bisemigroup of order 3. Let P = $P_1 \cup P_2$ where $P_1$ = {0, 3, 6, 9, 12, 15} and $P_2$ = {0, 4I, 4, 1, I}, P is a neutrosophic sub bisemigroup of order 6 + 5 = 11. Thus 11/33 but $11^2 \nmid 33$. So B is only a weak Sylow neutrosophic bisemigroup.

Now we see the set T = $T_1 \cup T_2$ where $T_1$ = {0, 3, 6, 9, 12, 15} and $T_2$ = {0, 1, I} is a neutrosophic subsemigroup of order 9. 3/33 and $3^2 \nmid 33$ but B has a neutrosophic subsemigroup of order 9.

Next we construct the following example.

***Example 3.2.15:*** Let B = $B_1 \cup B_2$ where $B_1$ = {0, 1, 2, 3, 4, 5, 6, I, 2I, 3I, 4I, 5I, 6I} a semigroup under multiplication modulo 7. Take $B_2$ = {0, 1, 2, 3, 4, 5, I, 2I, 3I, 4I, 5I} a neutrosophic semigroup under multiplication modulo 6. o(B) = 13 + 11 = 24, 2/24, $2^2$/24, $2^3$/24 and $2^4 \nmid 24$; 3/24 and $3^2 \nmid 24$.

Clearly B has no neutrosophic subsemigroup of order 3. Now take V = $V_1 \cup V_2$ where $V_1$ = {0, 1, I} and $V_2$ = {0, 1, I, 4, 4I}. o(V) = 8. Thus B is only a weak Sylow neutrosophic bisemigroup.

***Example 3.2.16:*** Let B = $B_1 \cup B_2$ where $B_1$ = {0, 1, 2} semigroup under multiplication modulo 3. $B_2$ = {0, 1, 2, I, 2I, 1+I, 2+I, 2I+1, 2 + 2I} a semigroup under multiplication modulo 4, o(B) = 12. Clearly 2/12, $2^2$/ 12 but $2^3 \nmid 12$. Also 3/12 and $3^2 \nmid 12$. B has no neutrosophic subbisemigroup of order 3 or 4 so B is Sylow free neutrosophic bisemigroup.

We now proceed on to define when are two neutrosophic subbisemigroups conjugate.

**DEFINITION 3.2.13:** *Let $B = (B_1 \cup B_2, *, o)$ be a neutrosophic bisemigroup. $T = T_1 \cup T_2$ and $S = S_1 \cup S_2$ be two neutrosophic subbisemigroups. We say T is conjugate to S if we can find $x_i$ and $y_j$ $(1 \leq i, j \leq 2)$ such $T_1 x_1 = x_1 S_1$ or $S_1 x_1$) and $T_2 x_2 = x_2 S_2$ (or $S_2 x_2$) We cannot say about the order of $S_i$ and $T_i$, $i = 1, 2$.*

***Example 3.2.17:*** Let $B = B_1 \cup B_2$ be a neutrosophic bisemigroup where $B_1 = Z_{12}$, the semigroup under multiplication modulo 12. $B_2 = \{0, 1, 2, \ldots, 5, I, 2I, \ldots, 5I\}$ a neutrosophic semigroup under multiplication modulo 6. Take $T = T_1 \cup T_2$ and $S = S_1 \cup S_2$ any two neutrosophic subbisemigroups, where $T_1 = \{0, 2, 4, 6, 8, 10\}$ and $T_2 = \{0, 2, 4, 2I, 4I\}$ and $S_1 = \{0, 3, 6, 9\}$ and $S_2 = \{0, 3, 3I\}$, $3T_1 = 2S_1 = \{0, 6\}$ and $6T_2 = 2S_2 = \{0\}$.

Thus the neutrosophic subbisemigroups T and S are conjugate.

*Note:* For two neutrosophic subbisemigroups T and S to be conjugate we do not demand o(T) = o(S).

On similar lines as in case of neutrosophic semigroups we can also define when are two elements in a neutrosophic bisemigroup conjugate. These definitions about substructures and other properties given for neutrosophic bisemigroups can be defined with appropriate modifications in case of strong neutrosophic bisemigroups.

## 3.3 Neutrosophic N-Semigroup

In this section for the first time the notion of neutrosophic N-semigroups are introduced. Substructures like strong neutrosophic sub N-semigroup, neutrosophic sub N semigroup, strong neutrosophic N-ideals, neutrosophic N-ideals, strong Lagrange neutrosophic sub N semigroup, p-Sylow neutrosophic sub N semigroup are introduced and examples are given.

Also some special elements like neutrosophic elements viz. neutrosophic N-ary idempotents, neutrosophic N-ary zero divisors, and neutrosophic Cauchy elements, Cauchy elements are introduced in neutrosophic N-semigroups and studied.

Now we proceed on to define the notion of neutrosophic N-semigroups and their generalizations and particularizations.

**DEFINITION 3.3.1:** *Let $\{S(N), *_1, \ldots, *_N\}$ be a non empty set with N-binary operations defined on it. We call S(N) a neutrosophic N-semigroup (N a positive integer) if the following conditions are satisfied.*

*i.* $S(N) = S_1 \cup \ldots \cup S_N$ *where each* $S_i$ *is a proper subset of S(N) i.e.* $S_i \not\subseteq S_j$ *or* $S_j \not\subseteq S_i$ *if* $i \neq j$.

*ii.* $(S_i, *_i)$ *is either a neutrosophic semigroup or a semigroup for* $i = 1, 2, \ldots, N$.

*Note:* When N = 2, we call S(N) to be a neutrosophic bisemigroup. If all the N-semigroups $(S_i, *_i)$ are neutrosophic semigroups (i.e. for i = 1, 2, …, N) then we call S(N) to be a neutrosophic strong N-semigroup.

Now we will give examples of both neutrosophic N-semigroups and neutrosophic strong N-semigroups.

***Example 3.3.1:*** Let $S(N) = \{S_1 \cup S_2 \cup S_3 \cup S_4 \cup S_5, *_1, \ldots, *_5\}$ be a neutrosophic 5-semigroup where

$S_1$ = $\langle Z \cup I\rangle$ under '+', is a neutrosophic semigroup,

$S_2$ = $\left\{\begin{pmatrix} a & b \\ c & d \end{pmatrix} \mid a, b, c, d \in \langle Q \cup I\rangle\right\}$ under matrix multiplication is a neutrosophic semigroup.

$S_3$ = {0, 1, 2, I, 2I} a neutrosophic semigroup under multiplication modulo 3,

$S_4$ = S(3), the set of all mappings of the set (a, b, c) to itself under the composition of maps and

$S_5$ = $\{Z_{15}$, the semigroup under multiplication modulo 15}.

S(N) is only a neutrosophic 5-semigroup and not a neutrosophic strong 5 semigroup.

Now we proceed on to give an example of a neutrosophic strong N-semigroup.

***Example 3.3.2:*** Let S(N) = $\{S_1 \cup S_2 \cup S_3 \cup S_4, *_1, *_2, *_3, *_4\}$ be such that

$S_1$ = $\{\langle Z \cup I\rangle$, semigroup under multiplication$\}$,
$S_2$ = $\{\langle Q^+ \cup I\rangle$, semigroup under '+', $Q^+$ is the set of all positive rationals$\}$,
$S_3$ = {0, 1, 2, 3, I, 2I, 3I; semigroup under multiplication modulo 4} and
$S_4$ = $\left\{\begin{pmatrix} a & b \\ c & d \end{pmatrix} / a, b, c, d \in \langle Q \cup I\rangle\right\}$ semigroup under matrix multiplication.

It is clearly seen that all the four semigroups are neutrosophic so S(N) is a neutrosophic strong 4-semigroup.

Now we proceed on to define substructures in these two types of neutrosophic N-semigroups.

**DEFINITION 3.3.2:** *Let S(N) = $\{S_1 \cup S_2 \cup \ldots \cup S_N, *_1, \ldots, *_N\}$ be a neutrosophic N-semigroup. A proper subset P = $\{P_1 \cup P_2 \cup \ldots \cup P_N, *_1, *_2, \ldots, *_N\}$ of S(N) is said to be a neutrosophic N-subsemigroup if (1) $P_i = P \cap S$, i = 1, 2,…, N are subsemigroups of $S_i$ in which atleast some of the subsemigroups are neutrosophic subsemigroups.*

***Example 3.3.3:*** Let S(N) = $\{S_1 \cup S_2 \cup S_3 \cup S_4, *_1, *_2, *_3, *_4\}$ be a neutrosophic 4-semigroup where

$S_1$ = $\{Z_{12}$, semigroup under multiplication modulo 12$\}$,
$S_2$ = {0, 1, 2, 3, I, 2I, 3I, semigroup under multiplication modulo 4}, a neutrosophic semigroup.
$S_3$ = $\left\{\begin{pmatrix} a & b \\ c & d \end{pmatrix} / a, b, c, d \in \langle Q \cup I\rangle\right\}$, neutrosophic semigroup under matrix multiplication and
$S_4$ = $\langle Z \cup I\rangle$, neutrosophic semigroup under multiplication.

S(N) is only a neutrosophic 4-semigroup and not a neutrosophic strong 4-semigroup.

Take T = $\{T_1 \cup T_2 \cup T_3 \cup T_4, *_1, *_2, *_3, *_4\}$ where

$T_1$ = $\{0, 2, 4, 6, 8, 10\} \subseteq Z_{12}$,

$T_2$ = $\{0, I, 2I, 3I\} \subset S_2$,

$T_3$ = $\left\{\begin{pmatrix} a & b \\ c & d \end{pmatrix} / a, b, c, d \in \langle Z \cup I \rangle \right\}$ and

$T_4$ = $\{\langle 5Z \cup I \rangle\}$ the neutrosophic semigroup under multiplication.

Clearly T is a neutrosophic sub 4-semigroup of S(N) and not a neutrosophic strong sub 4-semigroup for a neutrosophic N-semigroup cannot have neutrosophic strong sub N-semigroup.

Now we proceed on to illustrate by examples and define neutrosophic strong sub N-semigroup.

**DEFINITION 3.3.3:** *Let S(N) = $\{S_1 \cup S_2 \cup \ldots \cup S_N, *_1, \ldots, *_N\}$ be a neutrosophic strong N-semigroup. A proper subset T = $\{T_1 \cup T_2 \cup \ldots \cup T_N, *_1, \ldots, *_N\}$ of S(N) is said to be a neutrosophic strong sub N-semigroup if each $(T_i, *_i)$ is a neutrosophic subsemigroup of $(S_i, *_i)$ for i = 1, 2,…, N where $T_i = T \cap S_i$.*

If only a few of the $(T_i, *_i)$ in T are just subsemigroups of $(S_i, *_i)$ (i.e. $(T_i, *_i)$ are not neutrosophic subsemigroups then we call T to be a sub N-semigroup of S(N). a neutrosophic strong N-semigroup can have all the 3 types of subsemigroups which we illustrate by the following example.

***Example 3.3.4:*** Let S(N) = $\{S_1 \cup S_2 \cup S_3 \cup S_4 \cup S_5, *_1, *_2, \ldots, *_5\}$ be a neutrosophic strong 5-semigroup where

$S_1$ = $\left\{\begin{pmatrix} a & b \\ c & d \end{pmatrix} | a, b, c, d \in \langle Q \cup I \rangle \right\}$ neutrosophic semigroup under matrix multiplication,

$S_2$ = {0, 1, 2, 3, 4, I, 2I, 3I, 4I} be the neutrosophic semigroup under multiplication modulo 5,

$S_3$ = $\langle Z \cup I\rangle$ the neutrosophic semigroup under multiplication,

$S_4$ = $\langle Q^+ \cup I \rangle$ the neutrosophic semigroup under addition and

$S_5$ = $\left\{\begin{pmatrix} a & b & c \\ d & e & f \end{pmatrix} \mid a, b, c, d, e, f \text{ are in } \langle Q \cup I\rangle\right\}$ neutrosophic semigroup under matrix addition.

S(N) is a neutrosophic strong 5-semigroup.

Now consider the proper subset T = $\{T_1 \cup T_2 \cup T_3 \cup T_4, \cup T_5, *_1, *_2, \ldots, *_5\}$ of S (N) where

$T_1$ = $\left\{\begin{pmatrix} a & b \\ c & d \end{pmatrix} \mid a, b, c, d \in \langle Z \cup I\rangle\right\} \subset S_1$ is a neutrosophic subsemigroup of $S_1$,

$T_2$ = $\{0, 1, I, 2I, 3I, 4I\} \subset S_2$ is a neutrosophic subsemigroup under multiplication modulo 5;

$T_3$ = $\langle 3Z \cup I\rangle \subset S_3$ is a neutrosophic subsemigroup under multiplication,

$T_4$ = $\langle Z^+ \cup I\rangle$ neutrosophic subsemigroup of $S_4$ under addition and

$T_5$ = $\left\{\begin{pmatrix} a & b & c \\ d & e & f \end{pmatrix} \mid a, b, c, d, e, f \in \langle Z \cup I\rangle\right\}$ is a neutrosophic subsemigroup under matrix addition of $S_5$.

T is a neutrosophic strong 5-subsemigroup of S (N). Now consider the proper subset B = $\{B_1 \cup B_2 \cup \ldots \cup B_5, *_1, \ldots, *_5\}$ of S(N) where

$B_1$ = $\left\{\begin{pmatrix} a & b \\ c & d \end{pmatrix} \mid a, b, c, d \in Q\right\}$ is just a subsemigroup of $S_1$ under matrix multiplication,

$B_2$ = {0, I, 2I, 3I, 4I}, is a neutrosophic subsemigroup of $S_2$,

$B_3$ = $\langle 2Z \cup I\rangle$, the neutrosophic subsemigroup under multiplication

$B_4$ = $\{Q^+\}$ subsemigroup under addition and

$B_5$ = $\left\{\begin{pmatrix} a & b & c \\ d & e & f \end{pmatrix} \mid a, b, c, d, e, f \in \langle Z \cup I\rangle\right\}$ is a neutrosophic subsemigroup under matrix addition.

Thus B is only a neutrosophic 5-subsemigroup of S(N).

Now consider A = $\{A_1 \cup A_2 \cup A_3 \cup A_4 \cup A_5, *_1, \ldots, *_5\}$ a proper subset of S(N) which is not a neutrosophic strong 5-sub semigroup or a neutrosophic 5-subsemigroup. For consider

$A_1$ = $\left\{\begin{pmatrix} a & b \\ c & d \end{pmatrix} \mid a, b, c, d \in Q\right\}$ is a subsemigroup under matrix multiplication of $S_1$,

$A_2$ = {0, 1, 2, 3, 4} is a subsemigroup under multiplication modulo 5 of $S_2$,

$A_3$ = $\langle 2Z\rangle$ is the subsemigroup of $S_3$ under multiplication,

$A_4$ = $\{Z^+\}$ is the subsemigroup of $S_4$ under addition and

$A_5$ = $\left\{\begin{pmatrix} a & b & c \\ 0 & 0 & 0 \end{pmatrix} \mid a, b, c \in Q\right\}$ is a subsemigroup of $S_5$ under matrix addition.

Thus A is a 5-subsemigroup which is not neutrosophic. Hence the claim.

Now we proceed on to define the notion neutrosophic strong N-ideal and neutrosophic N-ideal of a neutrosophic strong N-semigroup and neutrosophic N-semigroup respectively.

**DEFINITION 3.3.4:** *Let S(N) = $\{S_1 \cup S_2 \cup \ldots \cup S_N, *_1, \ldots, *_N\}$ be a neutrosophic strong N-semigroup. A proper subset J = $\{I_1 \cup I_2 \cup \ldots \cup I_N\}$ where $I_t = J \cap S_t$ for t = 1, 2, …, N is said to be a neutrosophic strong N-ideal of S(N) if the following conditions are satisfied.*

i. *Each $I_t$ is a neutrosophic subsemigroup of $S_t$, t = 1, 2, …, N i.e. $I_t$ is a neutrosophic strong N-subsemigroup of S(N).*

ii. *Each $I_t$ is a two sided ideal of $S_t$ for t = 1, 2, …, N. Similarly one can define neutrosophic strong N-left ideal or neutrosophic strong right ideal of S(N).*

*A neutrosophic strong N-ideal is one which is both a neutrosophic strong N-left ideal and N-right ideal of S(N).*

***Note:*** It is important and interesting to note that a neutrosophic strong N-semigroup can never have a neutrosophic N-ideal or just a N-ideal it can only have a neutrosophic strong N-ideal.

**DEFINITION 3.3.5:** *Let $S(N) = \{S_1 \cup S_2 \cup \ldots \cup S_N, *_1, \ldots, *_N\}$ be a neutrosophic N-semigroup. A proper subset $P = \{P_1 \cup P_2 \cup \ldots \cup P_N, *_1, \ldots, *_N\}$ of S (N) is said to be a neutrosophic N-subsemigroup, if the following conditions are true*

i. *P is a neutrosophic sub N-semigroup of S(N).*
ii. *Each $P_i = P \cap S_i$, i =1, 2, …, N is an ideal of $S_i$.*

*Then P is called or defined as the neutrosophic N ideal of the neutrosophic N-semigroup S(N).*

Now we can as in case of bisemigroups or neutrosophic bisemigroups define the notion of maximal ideal.

**DEFINITION 3.3.6:** *Let $S(N) = \{S_1 \cup S_2 \cup \ldots \cup S_N, *_1, \ldots, *_N\}$ be a neutrosophic strong N-semigroup. Let $J = \{I_1 \cup I_2 \cup \ldots \cup I_N, *_1, \ldots, *_N\}$ be a proper subset of S(N) which is a neutrosophic strong N-ideal of S(N). J is said to be a neutrosophic strong maximal N-ideal of S(N) if each $I_t \subset S_t$, (t = 1, 2, …, N) is a maximal ideal of $S_t$.*

*It may so happen that at times only some of the ideals $I_t$ in $S_t$ may be maximal and some may not be in that case we call the ideal J to be a neutrosophic quasi maximal N-ideal of S(N). Suppose $S(N) = \{S_1 \cup S_2 \cup \ldots \cup S_N, *_1, \ldots, *_N\}$ is a neutrosophic strong N-semigroup, $J' = \{J_1 \cup J_2 \cup \ldots \cup J_N, *_1, \ldots, *_N\}$ be a neutrosophic strong N-ideal of S(N).*

*J' is said to be a neutrosophic strong minimal N-ideal of S (N) if each $J_i \subset S_i$ is a minimal ideal of $S_i$ for i = 1, 2, …, N. It may so happen that some of the ideals $J_i \subset S_i$ be minimal and some may not be minimal in this case we call J' the neutrosophic strong quasi minimal N-ideal of S(N).*

Now we proceed on to define the notion of neutrosophic maximal N-ideal, neutrosophic minimal N-ideal and their related quasi structures.

*Suppose $S(N) = \{S_1 \cup S_2 \cup \ldots \cup S_N, *_1, \ldots, *_N\}$ be a neutrosophic N-semigroup. $B = \{B_1 \cup B_2 \cup \ldots \cup B_N, *_1, \ldots, *_N\}$ is called the neutrosophic N-ideal of S(N) if*

*(1) Each $B_i = B \cap S_i$ is an ideal of $S_i$, i = 1, 2, …, N.*

*If each of these ideals $B_i$ are maximal ideals of $S_i$ we call B the neutrosophic maximal N-ideal of S(N).*

*If each of these ideals $B_i$ are minimal ideals of $S_i$ then we call B the neutrosophic minimal N-ideal of S(N).*

*If some of the ideals $B_i$ of $S_i$ are maximal and some ideals are not we call B as the neutrosophic quasi maximal N-ideal. If some of the ideals $B_j$ of $S_j$ are minimal then B is called the neutrosophic minimal N-ideal.*

We just give an example.

***Example 3.3.5:*** Let $S(N) = \{S_1 \cup S_2 \cup S_3 \cup S_4, *_1, *_2, *_3, *_4\}$ where

$S_1$ = $\langle Z \cup I \rangle$ the neutrosophic semigroup under multiplication,

$S_2$ = {0, 1, 2, I, 2I} the neutrosophic semigroup under multiplication modulo 3,

$S_3$ = $\{(a, b) \mid a \in Z \cup I \text{ and } b \in Q\}$ component wise multiplication is a neutrosophic semigroup and

$S_4$ = {0, 1, 2, …, 11} semigroup under multiplication modulo 12.

S(N) is a neutrosophic 5 semigroup.

Take T = $\{T_1 \cup T_2 \cup T_3 \cup T_4, *_1, *_2, *_3, *_4\}$ a proper subset of S(N) where

$T_1$ = $\{\langle 5Z \cup I\rangle\}$ neutrosophic subsemigroup under multiplication.
$T_2$ = $\{0, 1, I, 2I\} \subset S_2$ .
$T_3$ = $\{(0, x) \mid x \in Q\}$ is a semigroup under component wise multiplication and
$T_4$ = $\{0, 6\}$ is a semigroup under multiplication modulo 12.

T is a neutrosophic sub 4-semigroup of S(N).

Having just studied the substructure properties now we proceed on to study the order of the neutrosophic N-semigroup.

Let S(N) = $\{S_1 \cup S_2 \cup \ldots \cup S_N, *_1, \ldots, *_N\}$ be a neutrosophic N-semigroup. The order of the neutrosophic N-semigroup S(N) is the number of distinct elements in S(N). If the number of elements in S(N) is finite we call S(N) a finite neutrosophic N-semigroup; if otherwise S(N) is called as an infinite neutrosophic N-semigroup.

***Example 3.3.6:*** Let S(N) = $\{S_1 \cup S_2 \cup S_3 \cup S_4, *_1, *_2, *_3, *_4\}$ be a neutrosophic 4-semigroup where

$S_1$ = $\langle Z \cup I\rangle$, neutrosophic semigroup under multiplication,
$S_2$ = $\{0, 1, I, 1 + I\}$, neutrosophic semigroup under multiplication modulo 2,
$S_3$ = $\left\{\begin{pmatrix} a & b \\ c & d \end{pmatrix} \mid a, b, c, d \in \langle Q \cup I\rangle\right\}$ is a neutrosophic semigroup under matrix multiplication and
$S_4$ = $\{0, 1, 2, 3, 4, \ldots, 9\}$ is a semigroup under multiplication modulo 10.

Clearly S(N) has infinite number of elements so S(N) is an infinite neutrosophic N-semigroup or S(N) is of infinite order or $o(S(N)) = \infty$, or $|S(N)| = \infty$.

***Example 3.3.7:*** Let S(N) = $\{S_1 \cup S_2 \cup S_3 , *_1, *_2, *_3\}$ where

$S_1$ = {0, 1, 2, I, 2I}, a neutrosophic semigroup under multiplication modulo 3.
$S_2$ = S(2) the semigroup of mappings of the set (1, 2) to itself under the composition of mappings and
$S_3$ = {0, 1, 2, 3, 4} semigroup under multiplication modulo 5.

o(S(N)) = 5 + 4 + 5 = 14. Thus S(N) is a neutrosophic 3-semigroup of order 14.

It is important to note that in general even if S(N) is a neutrosophic N-semigroup of finite order still the order of neutrosophic sub N semigroups need not in general divide the order of S(N).

Suppose S(N) is of order n, n a prime then S(N) can have proper neutrosophic sub N-semigroups.

We illustrate this by the following example.

***Example 3.3.8***: Let S(N) = $\{S_1 \cup S_2 \cup S_3 \cup S_4, *_1, *_2, *_3, *_4\}$ where

$S_1$ = {0, 1, 2, …, 11} a semigroup under multiplication modulo 12.
$S_2$ = {0, 1, 2, I, 2I} neutrosophic semigroup under multiplication modulo 5.
$S_3$ = {(S(2)}, semigroup of order 4 and
$S_4$ = {0, 1, 2, …, 7} semigroup under multiplication modulo 8.

S(N) is a neutrosophic 4-semigroup of finite order. o(S(N)) = 12 + 5 + 4 + 8 = 29.

T = $\{T_1 \cup T_2 \cup T_3 \cup T_4, *_1, *_2, *_3, *_4\}$ is a neutrosophic sub 4-semigroup of S(N) where

$T_1$ = $\{0, 2, 4, 6, 8, 10\} \subset Z_{12}$,
$T_2$ = $\{0, 1, I, 2I\} \subset S_2$,
$T_3$ = $\begin{pmatrix} 1 & 2 \\ 1 & 2 \end{pmatrix} \begin{pmatrix} 1 & 2 \\ 2 & 1 \end{pmatrix} \subset S_3$, and
$T_4$ = $\{0, 4\} \subset S_4$.

T is a finite and o(T) = 6 + 4 + 2 + 2 = 14. (14, 29) = 1. Thus we see several interesting things happen in case of neutrosophic N-semigroup. One such is illustrated by the above example.

We can now define some special elements in case of neutrosophic N-semigroups.

**DEFINITION 3.3.7:** *Let $S(N) = \{S_1 \cup S_2 \cup \ldots \cup S_N, *_1, \ldots, *_N\}$ be a neutrosophic N semigroup. An element $x \in S(N)$ is an idempotent if $x \in S_i$ and $x^2 = x$. An element $0 \neq x \in S(N)$ i.e. $x \in S_i$ is said to be a zero divisor if there exists $0 \neq y \in S_i$ with $x y = 0$.*

*Note:* If in a neutrosophic N-semigroup $S(N) = \{S_1 \cup S_2 \cup \ldots \cup S_N, *_1, \ldots, *_N\}$ every $(S_i, *_i)$ is a monoid for i = 1, 2, …, N, then we define S(N) to be a neutrosophic N-monoid.

In general every neutrosophic N-semigroup need not be a neutrosophic N-monoid but every neutrosophic N-monoid is a neutrosophic N-semigroup. Now we proceed on to define the notion of N-ary idempotents, N-ary zero divisors and N-units of a neutrosophic N-semigroup N(S).

**DEFINITION 3.3.8:** *Let $N(S) = \{S_1 \cup S_2 \cup \ldots \cup S_N, *_1, \ldots, *_N\}$ be a neutrosophic N-semigroup. An element $X = (x_1, x_2, \ldots, x_N) \in S(N)$ where each $x_i \in S_i$ is called a N-ary idempotent if $X^2 = \left(x_1^2, x_2^2, \ldots, x_N^2\right) = (x_1, x_2, \ldots, x_N) = X$. i.e. each $x_i \in S_i$ is an idempotent of $S_i$ if one or many of these $x_i$ are neutrosophic elements then X is called as the neutrosophic N-ary-idempotent of S(N).*

First we illustrate this by the following example.

***Example 3.3.9:*** Let $S(N) = \{S_1 \cup S_2 \cup S_3 \cup S_4, *_1, \ldots, *_4\}$ be a neutrosophic 4-semigroup. where

$S_1$ = {0, 1, 2, 3, I, 2I, 3I}, neutrosophic semigroup under multiplication modulo 4.

$S_2$ = {0, 1, 2, 3, …, 11} semigroup under multiplication modulo 12,
$S_3$ = {0, 1, 2, 1 + I, 2I, 2 + I, 1 + 2I, I, 2 + 2I} neutrosophic semigroup under multiplication modulo 3 and
$S_4$ = {(a, b) | a, b ∈ {0, 1, I, 1 + I}} neutrosophic semigroup under component-wise multiplication modulo 2.

Take X = (I, 4, 1+2I, 1+I), X is a neutrosophic 4-ary idempotent of S(N).

Now we proceed on to define N-ary zero divisors and N-ary units.

**DEFINITION 3.3.9:** *Let $S(N) = \{S_1 \cup S_2 \cup \ldots \cup S_N, *_1, \ldots, *_N\}$ be a neutrosophic N-semigroup such that each semigroup has the zero i.e. $x_i\, 0 = 0 x_i = 0$ for all $x_i \in S_i$, $i = 1,2, \ldots, N$. An element $X = (x_1, x_2, \ldots, x_N), \neq (0\ 0\ 0\ 0 \ldots 0)$, $x_i \in S_i$ $(i = 1, 2, \ldots, N)$ in S (N) is said to be a N-ary zero divisor if there exists an element $Y = (y_1, y_2, \ldots, y_N) \neq (0, 0, \ldots, 0)$ in S(N); $y_i \in S_i$ such that*

$$\begin{aligned} XY &= (x_1 y_1, x_2 y_2, \ldots, x_N y_N) \\ &= (0, 0, \ldots, 0) \\ &= Y.X. \end{aligned}$$

*If in the N-ary of X one or more of elements $x_i$ are neutrosophic then we call X to be a neutrosophic N-ary zero divisor of S(N).*

*Note:* It is important to note that if in $S(N) = \{S_1 \cup \ldots \cup S_N, *_1, \ldots, *_N\}$ each of the $S_i$ has 0 such that $x_i *_i 0 = 0 *_i x_i = 0$ we would not be in a position to define N-ary zero divisors.

***Example 3.3.10:*** Let us consider the neutrosophic N-semigroup, $S(N) = \{S_1 \cup S_2 \cup S_3\ *_1, *_2, *_3\}$ where

$S_1$ = $\langle Z^+ \cup I \rangle$, is the neutrosophic semigroup under multiplication;
$S_2$ = $Z_{12}$, semigroup under multiplication modulo 12 and
$S_3$ = S(3) semigroup of all mappings of the set (1 2 3) to itself.

Clearly S(N) cannot have 3-ary zero divisors for $S_1$, and $S_3$ have no zero divisors more so S(N) cannot have 3 neutrosophic 3-ary zero divisors. Thus the condition every semigroup must have 0 is essential for us to define N-ary zero divisors or neutrosophic N-ary zero divisors.

***Example 3.3.11:*** Let $S(N) = \{S_1 \cup S_2 \cup S_3 \cup S_4, *_1, *_2, *_3, *_4\}$ be a neutrosophic 4-semigroup where

$S_1$ = {0, 1, 2, 3, I, 2I, 3I} neutrosophic semigroup under multiplication modulo 4,

$S_2$ = $\{Z_{12}$, the semigroup under multiplication modulo 12};

$S_3$ = $\{(a, b) \mid a, b \in \{0, 1, I, 2I\}\}$ neutrosophic semigroup under component wise multiplication and

$S_4$ = $\left\{\begin{pmatrix} a & b \\ c & d \end{pmatrix} / a, b, c, d, \in Q\right\}$.

Clearly S(N) is a neutrosophic 4-semigroup.

Consider $X = (x_1, x_2, x_3\ x_4)$, $(x_i \in S_i)$ in S(N), i.e., $X = (2I, 6, (0, I), \begin{pmatrix} 1 & 1 \\ 1 & 1 \end{pmatrix}) \neq (0, 0, 0, 0)$.

Let $Y = (2I, 6, (2I, 0), \begin{pmatrix} a & -a \\ -a & a \end{pmatrix})$ in S(N). Clearly XY = (0, 0, 0, 0). We can see YX = (0, 0, 0, 0). X is called as the neutrosophic 4-zero divisor.

*Note:* It may happen we have only XY = 0 and YX ≠ 0 in such cases we say X is a N-ary of left or right zero divisor. For in the same example if we take

$$X = (2I, 6, (0, I), \begin{pmatrix} 1 & 1 \\ 1 & 1 \end{pmatrix}) \text{ but}$$

$Y^1 = (2I, 6, (2I, 0), \begin{pmatrix} a & -a \\ -a & b \end{pmatrix})$

$XY^1 = (0, 0, 0, 0)$ but $Y^1 X \neq 0$ for

$Y^1X = (0, 0, 0, \begin{pmatrix} a-b & -a+b \\ a-b & -a+b \end{pmatrix}) \neq (0, 0, 0, 0)$

Hence the claim.
Next we need to define N-ary of invertible elements in S (N).

**DEFINITION 3.3.10:** *Let $S(N) = (S_1 \cup S_2 \cup \ldots \cup S_N, *_1, \ldots, *_N)$ be a neutrosophic N-monoid. An element $X = (x_1, \ldots, x_N)$ of S(N) where $x_i \in S_i$ is said to be N-ary invertible if there exists a $Y = (y_1, \ldots, y_N)$ in S(N) such that $XY = YX = (e_1, \ldots, e_N)$ where each $e_i \in S_i$ is such that $e_i x_i = x_i e_i = x_i$ for all $x_i \in S_i$. If in X or Y we have neutrosophic elements then we call X to be a neutrosophic N-ary invertible element of S(N).*

Now we illustrate this by the following example.

***Example 3.3.12:*** Let $S(N) = \{S_1 \cup S_2 \cup S_3 \cup S_4, *_1, *_2, *_3, *_4\}$ be a neutrosophic 4-monoid where

$S_1$ = {0, 1, 2, 1+I, 2+I, 2I + 1, 2I + 2, I, 2I}, neutrosophic semigroup under multiplication modulo 3.
$S_2$ = $\{Q^+$, semigroup under multiplication},
$S_3$ = {0, 1, 2, …, 10, semigroup under multiplication modulo 11} and
$S_4$ = {S(3), semigroup}.

Take $X = (2I + 2, 5, 10, \begin{pmatrix} 1 & 2 & 3 \\ 2 & 3 & 1 \end{pmatrix})$ in S(N), now $Y = (2I + 2, \frac{1}{5}, 10, \begin{pmatrix} 1 & 2 & 3 \\ 3 & 1 & 2 \end{pmatrix})$. Now $XY = (1, 1, 1, \begin{pmatrix} 1 & 2 & 3 \\ 1 & 2 & 3 \end{pmatrix}) = YX$.
Clearly X is a neutrosophic 4-ary unit.

Several other properties regarding semigroups can be defined for N-semigroups and neutrosophic N-semigroups.

We just give a brief description of conjugate neutrosophic sub N-semigroups.

**DEFINITION 3.3.11:** *Let $S(N) = (S_1 \cup S_2 \cup \ldots \cup S_N, *_1, \ldots, *_N)$ be a neutrosophic N-semigroup. Let $P = \{ P_1 \cup P_2 \cup \ldots \cup P_N, *_1, \ldots, *_N\}$ and $T(N) = (T_1 \cup T_2 \cup \ldots \cup T_N, *_1, *_2, \ldots, *_N)$ be any two neutrosophic sub N-semigroups. We say P and T are conjugate neutrosophic sub N-semigroups if we have for each pair $(P_i, T_i)$, $x_i P_i = T_i\ x_1$ ; $x_1' \in S_i$ ($P_i$ and $T_i$ are subsemigroups of $S_i$, $i = 1, 2, 3, \ldots, N$).*

***Example 3.3.13:*** Let $S(N) = (S_1 \cup S_2 \cup S_3, *_1, *_2, *_3)$ be a neutrosophic 3-semigroup where

$S_1$ = $\{Z_{12}$, semigroup under multiplication modulo 12$\}$,
$S_2$ = $\{0, 1, 2, 3, 4, 5, I, 2I, 3I, 4I, 5I\}$ neutrosophic semigroup under multiplication modulo 6 and
$S_3$ = $\{Z_4 \times Z_4 \mid ((a, b) \mid a, b \in Z_4)$, semigroup under component wise multiplication modulo 4$\}$.

Take $T = \{T_1 \cup T_2 \cup T_3, *_1, *_2, *_3\}$, a neutrosophic sub 3-semigroup where $T_1 = \{0, 2, 4, 6, 8, 10\} \subset Z_{12}$, $T_2 = \{0, 2, 4, 2I, 4I\}$ is a neutrosophic semigroup. $T_3 = \{(0, 1) (1, 0) (0, 0)\}$. $P = (P_1 \cup P_2 \cup P_3, *_1, *_2, *_3)$ be a neutrosophic sub 3-semigroup where $P_1 = \{0, 3, 6, 9\}$, $P_2 = \{0, 3, 3I\}$ is a neutrosophic semigroup and $P_3 = \{(0, 2) (2, 0) (0, 0)\}$ is a subsemigroup of $S_3$. $3\{0, 2, 6, 4, 8, 10\} = \{0, 6\}$ and $2\{0, 3, 6, 9\} = \{0, 6\}$ i.e. $3\ T_1 = 2P_1 = \{0, 6\}$. $3\{0, 2, 4, 2I, 4I\} = \{0\}$; $\{0, 3, 3I\}2 = \{0\}$. So $3T_2 = 2P_2 = \{0\}$. $2T_3 = P_3\ 1$.

Thus we see T and P are neutrosophic sub 3-semigroups. The following observations are very interesting.

It P and T are conjugate neutrosophic sub N-semigroups say of finite order.

(1) $o(P) \neq o(T)$ in general.
(2) $o(x_i P_i) = o(T_i x_i) \neq o(P_i)$ or $\neq o(T_i)$ in general.

**Chapter Four**

# NEUTROSOPHIC LOOPS AND THEIR GENERALIZATIONS

In this chapter for the first time we introduce the notion of neutrosophic loops and their generalization. A new class of neutrosophic loops of order 4t is introduced. Several interesting properties about them are derived. In this class of neutrosophic loops we have only one neutrosophic loop to be commutative for a given t. Likewise only one neutrosophic loop to be right or left alternative for a given t and no alternative neutrosophic loop. Also these class of neutrosophic loops are simple for they do not have normal neutrosophic subloop. Some of the neutrosophic loops are WIP-loops.

This chapter has three sections. In section 1 we derive and define the neutrosophic loops and give some of its properties. Section 2 defines neutrosophic biloops and section three defines neutrosophic N-loops and gives several interesting properties about them.

## 4.1 Neutrosophic loops and their properties

In this section we introduce the notion of neutrosophic loop. We define several interesting properties about them illustrate them with examples. In this section the new class of neutrosophic loops of order 2(n + 1); n odd is introduced and analyzed.

**DEFINITION 4.1.1:** *A neutrosophic loop is generated by a loop L and I denoted by* $\langle L \cup I \rangle$*. A neutrosophic loop in general need* *not be a loop for $I^2 = I$ and I may not have an inverse but every element in a loop has an inverse.*

Further a neutrosophic loop will always contain a loop. Throughout this book we will denote a neutrosophic loop by $\langle L \cup I\rangle$.

***Example 4.1.1:*** Let $L_5$ (3) be the loop; {e, 1, 2, 3, 4, 5} a loop of order 6. $\{L_5 (3) \cup I\}$ = {e, 1, 2, 3, 4, 5, eI, 1I, 2I, 3I, 4I, 5I} under the table; 2.I 2I = eI, rI. rI = eI, r = 1, 2, 3, 4, 5.

The table of the loop $L_5$ (3) is given below.

| * | e | 1 | 2 | 3 | 4 | 5 |
|---|---|---|---|---|---|---|
| e | e | 1 | 2 | 3 | 4 | 5 |
| 1 | 1 | e | 4 | 2 | 5 | 3 |
| 2 | 2 | 4 | e | 5 | 3 | 1 |
| 3 | 3 | 2 | 5 | e | 1 | 4 |
| 4 | 4 | 5 | 3 | 1 | e | 2 |
| 5 | 5 | 3 | 1 | 4 | 2 | e |

$\langle L_5 (3) \cup I\rangle$ is a neutrosophic loop of order 12.

We as in case of other neutrosophic structures define order of a neutrosophic loop. The number of distinct elements in $\langle L \cup I\rangle$ is called the order of $\langle L \cup I\rangle$. If the number of elements is finite we call $\langle L \cup I\rangle$ a finite loop. If the number of elements in $\langle L \cup I\rangle$ is infinite then $\langle L \cup I\rangle$ is an infinite neutrosophic loop.

Now we proceed on to define the neutrosophic subloop.

**DEFINITION 4.1.2:** *Let $\langle L \cup I\rangle$ be a neutrosophic loop. A proper subset $\langle P \cup I\rangle$ of $\langle L \cup I\rangle$ is called the neutrosophic subloop, if $\langle P \cup I\rangle$ is itself a neutrosophic loop under the operations of $\langle L \cup I\rangle$.*

We now illustrate by an example.

***Example 4.1.2:*** Let $\langle L \cup I\rangle = \langle L_7(4) \cup I\rangle$ be a neutrosophic loop where $L_7(4)$ is a loop. $\langle$e, eI, 2, 2I$\rangle$ is a neutrosophic subloop where $\langle L_7(4) \cup I\rangle$ = {e, 1, 2, 3, 4, 5, 6, 7, eI, 1I, 2I, 3I, 4I, 5I, 6I, 7I}. o$\langle L_7(4) \cup I\rangle$ = 16 and $L_7(4)$ is given by the following table.

| * | e | 1 | 2 | 3 | 4 | 5 | 6 | 7 |
|---|---|---|---|---|---|---|---|---|
| e | e | 1 | 2 | 3 | 4 | 5 | 6 | 7 |
| 1 | 1 | e | 5 | 2 | 6 | 3 | 7 | 4 |
| 2 | 2 | 5 | e | 6 | 3 | 7 | 4 | 1 |
| 3 | 3 | 2 | 6 | e | 7 | 4 | 1 | 5 |
| 4 | 4 | 6 | 3 | 7 | e | 1 | 5 | 2 |
| 5 | 5 | 3 | 7 | 4 | 1 | e | 2 | 6 |
| 6 | 6 | 7 | 4 | 1 | 5 | 2 | e | 3 |
| 7 | 7 | 4 | 1 | 5 | 2 | 6 | 3 | e |

and eI . eI = eI, 3I . 3I = eI, 2I . 2 . I = eI = 4I . 4I = 5I . 5I = 6I . 6I = 7I . 7I = 1I . 1I = eI.

We now proceed on to define a new class of neutrosophic loops. These loops are also even order built using {e, 1, 2, …, n | n an odd number} and the number of elements in them is 2 (n + 1); (n > 3).

**DEFINITION 4.1.3:** *Let $\langle L_n(m) \cup I\rangle$ = {e, 1, 2, ..., n, e.I, 1I, ..., nI}, where n > 3, n is odd and m is a positive integer such that (m, n) = 1 and (m – 1, n) = 1 with m < n. Define on $\langle L_n(m) \cup I\rangle$ a binary operation ' . ' as follows.*

i. *e.i. = i.e. = i for all i $\in L_n(m)$.*
ii. *$i^2$ = e for all i $\in L_n(m)$.*
iii. *iI. iI = e I for all i $\in L_n(m)$.*
iv. *i. j = t where t = (mj – (m – 1)i) (mod n) for all i, j $\in L_n(m)$, i ≠ j, i ≠ e and j ≠ e.*
v. *iI. jI = tI where t = (mj – (m – 1) i) (mod n) for all i I, jI $\in \langle L_n(m) \cup I\rangle$. $\langle L_n(m) \cup I\rangle$ is a neutrosophic loop of order 2 (n + 1).*

For varying m we get different neutrosophic loops, which we denote by $\langle L_n \cup I\rangle$. This new class of neutrosophic loops are of order 4t; t a positive integer. There are only 3 neutrosophic loops of order 12 in $\langle L_5 \cup I\rangle$.

***Example 4.1.3:*** We give the table for $\langle L_5(2) \cup I\rangle$.

| • | e | 1 | 2 | 3 | 4 | 5 | eI | 1I | 2I | 3I | 4I | 5I |
|---|---|---|---|---|---|---|---|---|---|---|---|---|
| e | e | 1 | 2 | 3 | 4 | 5 | eI | 1I | 2I | 3I | 4I | 5I |
| 1 | 1 | e | 3 | 5 | 2 | 4 | 1I | eI | 3I | 2I | 2I | 4I |
| 2 | 2 | 5 | e | 4 | 1 | 3 | 2I | 5I | eI | 4I | 1I | 3I |
| 3 | 3 | 4 | 1 | e | 5 | 2 | 3I | 4I | 1I | eI | 5I | 2I |
| 4 | 4 | 3 | 5 | 2 | e | 1 | 4I | 3I | 5I | 2I | eI | 1I |
| 5 | 5 | 2 | 4 | 1 | 3 | e | 5I | 2I | 4I | 1I | 3I | eI |
| eI | eI | 1I | 2I | 3I | 4I | 5I | eI | 1I | 2I | 3I | 4I | 5I |
| 1I | I | eI | 3I | 5I | 2I | 4I | 1I | eI | 3I | 5I | 2I | 4I |
| 2I | 2I | 5I | eI | 4I | 1I | 3I | 2I | 5I | eI | 4I | 1I | 3I |
| 3I | 3I | 4I | 1I | eI | 5I | 2I | 3I | 4I | 1I | eI | 5I | 2I |
| 4I | 4I | 3I | 5I | 2I | eI | 1I | 4I | 3I | 5I | 2I | eI | 1I |
| 5I | 5I | 2I | 4I | 1I | 3I | eI | 5I | 2I | 4I | 1I | 3I | eI |

This loop is a non commutative and non associative neutrosophic loop of order 12. Thus for all our examples to have non abstract neutrosophic loops we take loops from the new class of neutrosophic loop which are also or order 2(n + 1) or 4t, $t = 2\left(\frac{n+1}{2}\right)$ as n + 1 is also even as n is odd.

Now we define when is a neutrosophic loop commutative.

**DEFINITION 4.1.4:** *Let ($\langle L \cup I\rangle$, o) be a neutrosophic loop, we say $\langle L \cup I\rangle$ is commutative if a.b = b.a for all a, b $\in \langle L \cup I\rangle$. In the new class of neutrosophic loops only the loop* $\langle L_n\left(\frac{n+1}{2}\right) \cup I\rangle$ *is a commutative neutrosophic loop.*

Now we define a notion called strictly non commutative loop.

**DEFINITION 4.1.5:** *A neutrosophic loop* $(\langle L \cup I\rangle, o)$ *is strictly non commutative if* $x\, o y \neq y\, o x$ *for any* $x, y \in \langle L \cup I\rangle$ $x \neq y$, $x \neq e$, $y \neq e$*)*.

The loop $\langle L_5(2) \cup I\rangle$ given in example 4.1.3 is a strictly non commutative neutrosophic loop.

Now we, as in case of other algebraic structures define the notion of Lagrange neutrosophic subloop and their generalizations.

**DEFINITION 4.1.6:** *Let* $(\langle L \cup I\rangle, o)$ *be a neutrosophic loop of finite order. A proper subset P of* $\langle L \cup I\rangle$ *is said to be Lagrange neutrosophic subloop, if P is a neutrosophic subloop under the operation 'o' and* $o(P) / o\langle L \cup I\rangle$*. If every neutrosophic subloop of* $\langle L \cup I\rangle$ *is Lagrange then we call* $\langle L \cup I\rangle$ *to be a Lagrange neutrosophic loop.*

*If* $\langle L \cup I\rangle$ *has no Lagrange neutrosophic subloop then we call* $\langle L \cup I\rangle$ *to be a Lagrange free neutrosophic loop. If* $\langle L \cup I\rangle$ *has atleast one Lagrange neutrosophic subloop then we call* $\langle L \cup I\rangle$ *a weakly Lagrange neutrosophic loop.*

Now we will illustrate these by the following examples.

***Example 4.1.4:*** Let $\langle L_n (m) \cup I\rangle$ be a new class of neutrosophic loops of order 2 (n + 1) where n is a prime.

Then we see all neutrosophic subloops are Lagrange i.e. $\langle L_n(m) \cup I\rangle$ is a Lagrange neutrosophic loop.

***Example 4.1.5:*** Consider the neutrosophic loop $\langle L_{15}(2) \cup I\rangle$ = {e 1, 2, 3, 4, …, 15, eI, 1I, 2I, …, 14I, 15I} of order 32. It is easily verified P = {e, 2, 5, 8, 11, 14, eI, 2I, 5I, 8I, 11I, 14I} is a neutrosophic subloop and, order of P is 12 and 12 $\not|$ 32. Hence the claim. Thus P is not a Lagrange neutrosophic subloop of $\langle L_{15}(2) \cup I\rangle$.

Consider T = {eI, 3I, e, 3} $\subset \langle L_{15}(2) \cup I\rangle$; T is a Lagrange neutrosophic subloop of $\langle L_{15} (2) \cup I\rangle$ as 4 / 32. Thus $\langle L_{15}(2) \cup I\rangle$ is a weakly Lagrange neutrosophic loop.

Now we can define Cauchy element and Cauchy neutrosophic element of a neutrosophic loop $\langle L \cup I \rangle$.

**DEFINITION 4.1.7:** *Let $(\langle L \cup I \rangle, o)$ be a neutrosophic loop of finite order. An element $x \in \langle L \cup I \rangle$ is said to be a Cauchy element if $x^r = e$ and $r / o\langle L \cup I \rangle$. An element $x^{r^1} = e$ which is not a Cauchy element is called as an anti Cauchy element of $\langle L \cup I \rangle$.*

*A Cauchy neutrosophic element y of $\langle L \cup I \rangle$ is one such that $y^t = eI$ and $t / o\langle L \cup I \rangle$. If $\langle L \cup I \rangle$ has its elements to be either Cauchy element or Cauchy neutrosophic element (i.e. $\langle L \cup I \rangle$ has no anti Cauchy element or anti Cauchy neutrosophic element) then we call $\langle L \cup I \rangle$ to be a Cauchy neutrosophic loop. If all elements are anti Cauchy elements then we call $\langle L \cup I \rangle$ a Cauchy free neutrosophic loop.*

*If $\langle L \cup I \rangle$ has atleast one Cauchy element and one Cauchy neutrosophic element then we call $\langle L \cup I \rangle$ to be a weakly Cauchy neutrosophic loop.*

We illustrate these by the following example.

***Example 4.1.6:*** Let $\langle L_{11}(3) \cup I \rangle$ be a neutrosophic loop of order 24 = 4.6 = 4.2.3. In $\langle L_n(3) \cup I \rangle$ every element is such that $x^2 = e$ or $(Ix)^2 = eI$ so $\langle L \cup I \rangle$ is a Cauchy neutrosophic loop. Now we leave it as an exercise for the reader to prove all neutrosophic loops $\langle L_n(m) \cup I \rangle$ are Cauchy neutrosophic loops.

Next we speak about the Sylow theorems for neutrosophic loops. First we give our observations about the new class of neutrosophic loops $\langle L_n \cup I \rangle$.

**DEFINITION 4.1.8:** *A neutrosophic subloop H of the neutrosophic loop $\langle L \cup I \rangle$ is called a p-Sylow neutrosophic subloop of order $p^k$, if $p^k / o\langle L \cup I \rangle$ but $p^{k+1} \not| \; o\langle L \cup I \rangle$ for some prime p.*

Now before we prove the existence of p-Sylow neutrosophic subloop we prove the following for the new class of neutrosophic loop.

The new class of neutrosophic loop $\langle L_n \cup I\rangle$ will have all its loops to be of order 2 (n + 1) where as n is odd, n + 1 is even. We prove the following result.

**THEOREM 4.1.1:** *Let $\langle L_n(m) \cup I\rangle$ be the neutrosophic loop from $\langle L_n \cup I\rangle$. For every t / n there exists t subloops of order 2 (k + 1) where k / t.*

*Proof:* Let $\langle L_n(m) \cup I\rangle = \{e, 1, 2, \ldots, n, eI, 1I, \ldots, nI\}$ For $i < t$ consider the subset $\langle H_i(t) \cup I\rangle = \{e, i, i + t, i + 2t, \ldots, i + t (k – 1), eI, iI, (i + t) I, \ldots, [i + t (k – i)]I\}$ of $\langle L_n(m) \cup I\rangle$. Clearly $e \in \langle H_i(t) \cup I\rangle$ and so $\langle H_i(t) \cup I\rangle$ is itself a neutrosophic loop as $H_i(t)$ is a proper subloop of $L_n(m)$. So $\langle H_i(t) \cup I\rangle$ is a neutrosophic subloop of $\langle L_n(m) \cup I\rangle$; for if we take (i + rt), (i + st), (i + rt) I, (i + st) I$\in \langle H_i(t) \cup I\rangle$. If (i + rt)(i + st) = q then (i + rt)I(i + st)I = qI and q and qI satisfy q = [m (i + st) – (m – 1) (i + rt)] (mod n) and qI = [m (i + st)I – (m –1) (i + rt)I] (mod n) or q = i + ut or qI = (i + ut) I, where

u = (ms – (m – 1)r) (mod k) and
uI = (msI – (m –1) rI) (mod k).

Since q and qI is of the form, i + ut and (i + ut) I, respectively, q, qI $\in \langle H_i(t) \cup I\rangle$. Thus $\langle H_i(t) \cup I\rangle$ is a neutrosophic subloop of $\langle L_n(m) \cup I\rangle$ of order 2 (k + 1) As i can vary from 1 to t there exists t such neutrosophic subloops.

Also we have the following interesting corollary.

**COROLLARY 4.1.1:** *Let $\langle H_i(t) \cup I\rangle$ and $\langle H_j(t) \cup I\rangle$ be two neutrosophic subloops of $\langle L_n(m) \cup I\rangle$ then $\langle H_i(t) \cup I\rangle \cap \langle H_j(t) \cup I\rangle = \{e, eI\}$ where $i \neq j$.*

*Proof:* We shall assume that suppose $\langle H_i(t) \cup I\rangle \cap \langle H_j(t) \cup I\rangle \neq \{e, eI\}$, then there exists a $x \in \langle H_i(t) \cup I\rangle \cap \langle H_j(t) \cup I\rangle$ that is $x \in \langle H_i(t) \cup I\rangle$ and $x \in \langle H_j(t) \cup I\rangle$, $x \neq e$ or e I so

$$x = i + r_1 t$$
$$= j + r_2 t \text{ for some}$$

$1 \leq r_1, r_2 < n / t$, so we must have $i + r_1 t = j + r_2 t$ which implies $i - j = (r_2 - r_1) t$; i.e. $t / (i - j)$. Hence $i = j$ as i and $j \leq t$.

**COROLLARY 4.1.2:** *Let the neutrosophic subloops ($\langle H_i(t) \cup I\rangle$) of $\langle L_n(m) \cup I\rangle$ be as in theorem. Then* $\bigcup_{i=1}^{t} \langle H_i(t) \cup I\rangle = \langle L_n(m) \cup I\rangle$ *for every t dividing n.*

*Proof:* Since $\langle H_i(t) \cup I\rangle$ are neutrosophic subloops of $\langle L_n(m) \cup I\rangle$ we have $\bigcup_{i=1}^{t} \langle H_i(t) \cup I\rangle \subseteq \langle L_n(m) \cup I\rangle$. To prove the equality we have to show $\langle L_n(m) \cup I\rangle \subseteq \bigcup_{i=1}^{t} \langle H_i(t) \cup I\rangle$ that is to show that for $x \in \langle L_n(m) \cup I\rangle$, $x \in \bigcup_{i=1}^{t} \langle H_i(t) \cup I\rangle$. If x = e or eI then nothing to prove. Let x ≠ e or eI. Then for this x and given t we can find integers r and s such that x = rt + s (if xI then xI = (rt + s)I). Then clearly $x \in \langle H_i(t) \cup I\rangle$ that is $x \in \bigcup_{i=1}^{t} \langle H_i(t) \cup I\rangle$. Hence the claim

**COROLLARY 4.1.3:** *The neutrosophic subloops $\langle H_i(t) \cup I\rangle$ and $\langle H_j(t) \cup I\rangle$ in the above theorem are always isomorphic for every t dividing n.*

*Proof:* Given $\langle L_n (m) \cup I\rangle$ is the neutrosophic loop having $\langle H_i(t) \cup I\rangle$ and $\langle H_j(t) \cup I\rangle$ as given distinct neutrosophic subloops. Define mapping f from $\langle H_i(t) \cup I\rangle$ to $\langle H_j(t) \cup I\rangle$ as follows f (xI) = f (x) I for f (I) = I i.e. I is left invariant under the map. Clearly f(e) = e, f(eI) = eI and f(i + rt) = j + rt for $1 \leq r \leq n / t$. Thus f is a natural neutrosophic isomorphism between $\langle H_i(t) \cup I\rangle$ and $\langle H_j(t) \cup I\rangle$. Hence the result.

*Notation:* Let $\langle L_n(m) \cup I\rangle \in \langle L_n \cup I\rangle$. For t / n (say n = kt) we have $\langle H_i(t) \cup I\rangle$ = {e, i, i + t, i + 2t, …, i + (k – 1)t, eI, iI, (i + t)I, (i + 2t)I, …, (i + (k–1)t)I} for i = 1, 2, …, t which denotes the neutrosophic subloop of $\langle L_n(m) \cup I\rangle$ and its order is 2 [(n/t) + 1].

Now we proceed on to give yet another interesting property.

**THEOREM 4.1.2:** *Let $\langle L_n(m) \cup I\rangle \in \langle L_n \cup I\rangle$. If $\langle H \cup I\rangle$ is a neutrosophic subloop of $\langle L_n(m) \cup I\rangle$ of order 2 (t + 1) then t /n.*

*Proof:* Let i , j ∈ $\langle H \cup I\rangle$, (i ≠ j, i ≠ e, j ≠ e, i ≠ eI, j ≠ eI). Then i.j = k where k is given by k = [mj – (m – 1)i ] (mod n). Now ($m_j$ – (m – 1)i ) = j + (m – 1) (j – 1) that is k – j = (m – 1) (j – i). Clearly k – j = (m – 1) (j – i) which implies by basic number theory difference between k and j is a multiple of the difference between i and j.

Since $\langle H \cup I\rangle$ is a neutrosophic subloop, $\langle H \cup I\rangle$ is closed, hence the difference between any two elements will also be a multiple of some number (say d). So $\langle H \cup I\rangle$ contains {e, s, s + d, s + 2d, …, s + [n / (d – 1)]d, eI, sI, (s + d)I (s + 2d) I,.., (s + [n / (d – 1)] d) I} (e and eI belongs to $\langle H \cup I\rangle$) as $\langle H \cup I\rangle$ is a neutrosophic subloop.

This is true for some s such that (1 < s < d). But the set of these elements is nothing but $\langle H_s(d) \cup I\rangle$ whose order is 2 (d + 1) and d/n. Hence the claim.

Now we give a nice characterization theorem about these new class of neutrosophic loops in $\langle L_n \cup I\rangle$.

**THEOREM 4.1.3:** *$\langle L_n(m) \cup I\rangle \in \langle L_n \cup I\rangle$ contains a neutrosophic subloop of order 2 (k + 1) if and only if k/n.*

*Proof:* Let k/n say n = kt to show that there exists a neutrosophic subloop of order k + 1. Since n = kt so t/n, by the theorem just proved there exists a neutrosophic subloop of order 2 (k + 1).

Conversely if there exists a neutrosophic subloop of order 2 (k + 1) then k/n.

**THEOREM 4.1.4:** *For this new class of neutrosophic loops the Lagrange theorem for groups is satisfied by every neutrosophic subloop of $\langle L_n(m) \cup I\rangle$ if and only if n is an odd prime.*

*Proof:* Let n be an odd prime say p. There exists neutrosophic subloops of order 4 and 2 (p + 1) only by earlier theorem. Since, if the order of the neutrosophic subloop is 2 (p + 1) it is trivially $\langle L_n(m) \cup I\rangle$ as order of $\langle L_n(m) \cup I\rangle$ is 2 (p + 1). Now clearly for the neutrosophic subloop of order 4 we have 4 / 2 (p + 1) since p is odd. Hence the Lagrange theorem for groups is satisfied.

Conversely let n be not a prime number say n = rs, $1 \leq r, s \leq n$. To show that Lagrange theorem for groups is not satisfied by all neutrosophic subloops of $\langle L_n(m) \cup I\rangle$ we have to show that there always exists a neutrosophic subloop of $\langle L_n(m) \cup I\rangle$ which does not satisfy Lagrange theorem for groups. We have n = rs, $(1 < r, s < n)$.

Now for this integer r we can have either

(i) $r^2 \geq n$ or

(ii) $r^2 < n$.

We make use of this two mutually exclusive conditions to prove the result.

**Case i:** If $r^2 \geq n$, consider the neutrosophic subloop $\langle H_i(s) \cup I\rangle$ (as s/n). Clearly $o(\langle H_i(s) \cup I\rangle) = 2\ (r + 1)$. Now if the Lagrange's Theorem for groups is satisfied by $\langle H_i(s) \cup I\rangle$ then it implies r + 1 / n + 1. Since r / n implies r + 1 / n + 1. We must have $r^2 < n$, which is a contradiction to our assumption. So $\langle H_i(s) \cup I\rangle$ does not satisfy the Lagrange's theorem for finite groups.

**Case ii:** If $r^2 < n$ then $s^2 > n$ (as n = rs). So the neutrosophic subloop $\langle H_j(r) \cup I\rangle$ for any $j \in \{1, 2, \ldots, r\}$ does not satisfy the Lagrange theorem. Hence the claim.

We shall prove the new class of neutrosophic loops, $\langle L_n \cup I\rangle$ has only 2-Sylow neutrosophic subloops of minimal order 4.

**THEOREM 4.1.5:** *Let $\langle L_n(m) \cup I\rangle \in \langle L_n \cup I\rangle$ (order of $\langle L_n(m) \cup I\rangle$ is 2 (n + 1)). Let $n + 1 = p^k r$ where (p, r) = 1, p a prime. Then there exists a p-Sylow neutrosophic subloop of order 2 $(p^k)$ if and only if $p^k - 1 / r - 1$.*

*Proof:* Suppose there exists a p-Sylow neutrosophic subloop of order $p^k$ then $p^k - 1 / (p^k r - 1)$. But this implies $p^k - 1/ r - 1$.

Conversely if $p^k - 1/r - 1$ then $p^k - 1 / p^k r - 1$ using the theorem there exist a neutrosophic subloop of $2p^k$.

Now we see that our above theorem is perfectly valid for we prove for any odd prime p we cannot have p-Sylow neutrosophic subloops.

**THEOREM 4.1.6:** *For any neutrosophic loop in the class of neutrosophic loops $\langle L_n \cup I\rangle$ i.e. for $\langle L_n(m) \cup I\rangle$ of $\langle L_n \cup I\rangle$ there exists only 2-Sylow neutrosophic subloops.*

*Proof:* Just above we have proved if there exists a p-Sylow neutrosophic subloop of order $2p^k$ then $(p^k - 1) / (r - 1)$ where r is given by $p^k r = n + 1$ and (p, r) = 1. Suppose $p \neq 2$ i.e. p is an odd prime. Since (n + 1) is even r is even as $p^k r = n + 1$. Since $(p^k - 1)$ is even and r – 1 is odd as r is even so $(p^k - 1) \nmid (r - 1)$. Thus if p is an odd prime there cannot exists any p-Sylow neutrosophic subloops. Hence the result.

***Remark**:* It is intentionally used in the theorem a p-Sylow neutrosophic subloop of order $2p^k$. This is mainly to show later on p is an even prime.

Now we can give the p-Sylow neutrosophic subloops of $\langle L_5(3) \cup I\rangle$ = {e, 1, 2, 3, 4, 5, eI, I, 2I, 3I, 4I}. Clearly $P_1$ = {e, eI, 1, 1I}, $P_2$ = {e, eI, 2, 2I}, $P_3$ = {e, eI, 3, 3I} $P_4$ = {e, eI, 4, 4I}, $P_5$ = {e, eI, 5, 5I} are 2-Sylow neutrosophic subloops of $L_5$ (3).

Now we proceed on to define the notion of normal neutrosophic subloop of a neutrosophic loop.

**DEFINITION 4.1.9:** *Let $\langle L \cup I\rangle$ be a neutrosophic loop. A neutrosophic subloop $\langle H \cup I\rangle$ of $\langle L \cup I\rangle$ is said to be a normal neutrosophic subloop of $\langle L \cup I\rangle$ if*

*i.* $\langle H \cup I\rangle x = x \langle H \cup I\rangle$
*ii.* $(\langle H \cup I\rangle x) y = \langle H \cup I\rangle (x y)$
*iii.* $y (x (\langle H \cup I\rangle)) = (y x) (\langle H \cup I\rangle)$ *for all* $x, y \in \langle L \cup I\rangle$.

*A neutrosophic loop $\langle L \cup I\rangle$ is simple if it does not contain any non trivial normal neutrosophic subloop.*

Now we prove some interesting results about the new class of neutrosophic loops $\langle L_n \cup I\rangle$.

**THEOREM 4.1.7:** *Let $\langle L_n(m) \cup I\rangle \in \langle L_n \cup I\rangle$. Then $\langle L_n(m) \cup I\rangle$ does not contain any non trivial normal neutrosophic subloop.*

*Proof:* Let $\langle H_j(t) \cup I\rangle$ be a neutrosophic subloop of $\langle L_n(m) \cup I\rangle \in \langle L_n \cup I\rangle$.

***Case i:*** If t = n then $\{H_i(t) \cup I\} = \{e, i , eI, i I\}$ be a neutrosophic subloop. For this i we can find j, k $\in \langle L_n(m) \cup I\rangle$ such that (i, j) k $\neq$ i (j. k). Then $(\langle H_i(t) \cup I\rangle$ j) k $\neq$ $(\langle H_i(t) \cup I\rangle)$(j k). So $\langle H_i(t) \cup I\rangle$ is not a normal neutrosophic subloop.

***Case ii:*** If t $\neq$ n that is t < n then $\langle H_i(t) \cup I\rangle$ = {e, i, i + t, i + 2t,…, i + (n / t – 1) t, eI, i .I (i + t) I, (i + 2t) I, …, [i + (n/t – 1) t] I}. Take j $\notin \langle H_i(t) \cup I\rangle$ then $(\langle H_i(t) \cup I\rangle)$. j = $(\langle H_r(t) \cup I\rangle)$ \ ({e, eI} $\cup$ {j, jI}) where r is given by r = (mj – (m –1) i ) mod t. Now take k $\in \langle H_r(t) \cup I\rangle$ then $((\langle H_i(t) \cup I\rangle)$. j), k = $(\langle H_r(t) \cup I\rangle)$ \ (k, kI) $\cup$ (j, jI, k, kI) = $\langle A \cup I\rangle$ say.

Now if j, k, jI, kI $\in \langle H_i(t) \cup I\rangle$ then $(\langle H_i(t) \cup I\rangle)$ (j k) = $\langle H_i(t) \cup I\rangle \neq \langle A \cup I\rangle$ and if j, k, jI, kI $\notin \langle H_i(t) \cup I\rangle$, then e, eI $\notin$ $(\langle H_i(t) \cup I\rangle)$ (j. k) and so $\langle H_j(t) \cup I\rangle$ (j.k) $\neq \langle A \cup I\rangle$, so in

case ii also we get that $\langle H_i(t) \cup I\rangle$ is not a normal neutrosophic subloop.

Then we have very important result regarding the new class of neutrosophic loops $\langle L_n \cup I\rangle$.

**THEOREM 4.1.8:** *Each neutrosophic loop in $\langle L_n \cup I\rangle$ is simple.*

Now we proceed on to define Moufang Bruck, Bol, WIP, right alternative, left alternative neutrosophic loops.

**DEFINITION 4.1.10:** *A neutrosophic loop $\langle L \cup I\rangle$ where L is a loop is said to be neutrosophic Moufang loop if its satisfies anyone of the following identities.*

- *i.* $(x\,y)\,(zx) = (x\,(y\,z)\,x)$.
- *ii.* $((x\,y)\,z)\,y = x\,(y\,(z\,y))$.
- *iii.* $x\,(y\,(xz)) = ((xy)\,x)z$ *for all* $x, y, z \in \langle L \cup I\rangle$.

It is left as an exercise for the reader to prove that $\langle L_n \cup I\rangle$ does not contain any neutrosophic Moufang loops.

*Hint:* Use the property Moufang loops are diassociative.

Next we define neutrosophic Bruck loop.

**DEFINITION 4.1.11:** *Let $\langle L \cup I\rangle$ be a neutrosophic loop $\langle L \cup I\rangle$ is said to be a neutrosophic. Bruck loop if*

- *i.* $(x\,(y\,x))\,z = x\,(y\,(x\,z))$ *and*
- *ii.* $(xy)^{-1} = x^{-1}\,y^{-1}$ *for all* $x, y, z \in \langle L \cup I\rangle$ *whenever a neutrosophic element has no inverse we do expect* $x, y, \in \langle L \cup I\rangle$ *to satisfy condition (ii).*

It is important to note that none of the neutrosophic loops in the new class of loops $\langle L_n \cup I\rangle$ are Bruck neutrosophic loops.

Now we proceed on to define neutrosophic Bol loops.

**DEFINITION 4.1.12:** *A neutrosophic loop $\langle L \cup I\rangle$ is called a Bol neutrosophic loop if ((x y) z) y = x ((y z) y) for all x, y, z, $\in \langle L \cup I\rangle$.*

**DEFINITION 4.1.13:** *A neutrosophic loop $\langle L \cup I\rangle$ is said to be right alternative if (x y) y = x (y y) for all x, y $\in \langle L \cup I\rangle$ and left alternative if (x x) y = x (x y) for all x, y $\in \langle L \cup I\rangle$ and is an alternative neutrosophic loop if both it is a right and left neutrosophic alternative loop.*

**DEFINITION 4.1.14:** *A neutrosophic loop $\langle L \cup I\rangle$ is called a weak inverse property loop (WIP – neutrosophic loop) if (x y) z = e imply x (y z) = e and (xI yI) zI = eI and xI(yI zI) = eI for all x, y, z $\in \langle L \cup I\rangle$ (e the identity element of L).*

We mainly prove which are the identities satisfied by the new class of neutrosophic loops.

**THEOREM 4.1.9:** *The class of neutrosophic loops $\langle L_n \cup I\rangle$ contains exactly one left alternative neutrosophic loop and one right alternative neutrosophic loop and does not contain any alternative neutrosophic loop.*

*Proof:* Suppose $\langle L_n(m) \cup I\rangle \in \langle L_n \cup I\rangle$ be right alternative then (i. j). j = i (j.j), (iI jI)jI = iI (jI jI) for all i, j, iI, jI in $\langle L_n(m) \cup I\rangle$. If i = j or iI = jI or e, eI $\in$ {i, j} then above equality holds good trivially so take i ≠ j, iI ≠ jI, e ≠ j, eI ≠ jI, e ≠ i, eI ≠ iI. Now (i. j.) j = t where t = (mj – (m – 1) (mj – (m – 1)i) mod n and i (jj) = i (as jj = e ). So we must have t ≡ i (mod n) or (mj – (m – 1) (mj – (m-1)i ) = i(mod n) or $(m^2 – 2m)$ (i – j) ≡ 0 (mod n). If $\langle L_n(m) \cup I\rangle$ is to be right alternative we have this equation for all i, j (iI, jI) $\in \langle L_n(m) \cup I\rangle$; i, j $\in$ {1, 2,…, n} i ≠ j. Hence $(m^2 – 2m)$ = 0 or = m = 2 is the only solution. Thus $\langle L_m(2) \cup I\rangle$ is the only neutrosophic loop in $\langle L_n \cup I\rangle$ which is right alternative.

Suppose $\langle L_n(m) \cup I\rangle$ is left alternative then i (ij) = (i.i) j for all i, j $\in \langle L_n(m) \cup I\rangle$. Take i, j $\in \langle L_n(m) \cup I\rangle$ (i ≠ j, j ≠ e and i ≠ e) j ≠ eI, i ≠ eI then i (i.j) = t where t is given by, t (m(mj – (m – 1)i) – (m – 1)i)(mod n) and (i.i)j = j(i.i = e). Thus we have t ≡ j (mod n) that is $(m^2 - 1)(i - j) \equiv 0 \pmod n$ this must be true for al i and j. Hence $m^2 - 1 \equiv 0 \pmod n$. Now $m = n - 1$ is the unique solution for this equation since, $m \neq n$ and $(m - 1, n) = 1$) Thus $\langle L_n(n-1) \cup I\rangle$ is the only left alternative neutrosophic loop in $\langle L_n \cup I\rangle$.

If the neutrosophic loop is to be alternative we must have the loop to be both left and right alternative hence both conditions must be satisfied i.e., a common solutions to both the equations $(m^2 - 1) \equiv 0 \pmod n$ and $(m^2 - 2m) \equiv 0 \pmod n$ simultaneously is impossible an $n > 3$. Thus the new class of neutrosophic loops has no alternative neutrosophic loop.

It is important to observe that right alternative and left alternative loops in $\langle L_n \cup I\rangle$ are not commutative.

Next we prove none of the loops in the class $\langle L_n \cup I\rangle$ are Moufang neutrosophic loops.

**THEOREM 4.1.10:** *The new class of neutrosophic loops $\langle L_n \cup I\rangle$ does not contain any Moufang neutrosophic loop.*

*Proof:* Let $\langle L_n \cup I\rangle$ be the new class of neutrosophic loops. We prove the result under two cases

i. $\langle L_n(m) \cup I\rangle \in \langle L_n \cup I\rangle$ is commutative
ii. $\langle L_n(m) \cup I\rangle$ the neutrosophic loop of $\langle L_n \cup I\rangle$ is non commutative.

**Case i:** Let $\langle L_n(m) \cup I\rangle$ be a commutative neutrosophic loop. To prove $\langle L_n(m) \cup I\rangle$ is not Moufang, it is enough to show that $(x y)(z x) \neq x((y z) x)$ for atleast one triple $x, y, z \in \langle L_n(m) \cup I\rangle$ (Note what we prove for x, y, z holds good verbatim for xI, yI, zI so we discuss only for x, y, z to avoid repetition).

Take $z = x$ ($x \neq e$ or eI and $y \neq e$, or eI). Then $(xy)(zx) = (xy)(xx) = xy$ (as $xx = e$ in $\langle L_n(m) \cup I\rangle$) and $x((y z) x) = (x (yx)) x$.

If $\langle L_n(m) \cup I\rangle$ is Moufang we must have $(xy)(zx) = x((y z) x)$ for all $x, y, z \in \langle L_n(m) \cup I\rangle$. So $x y = x((yx)x)$. That is $y = ye = y(x x) = (y x) x$ (as $xx = e$ for all $x \in \langle L_n(m) \cup I\rangle$ hence $\langle L_n(m) \cup I\rangle$ must be right alternative which is not possible. Hence if the loop is commutative then it is not a Moufang loop.

**Case ii:** Let the neutrosophic loop $\langle L_n(m) \cup I\rangle$ be non commutative that is we have atleast a pair (x, y) of distinct elements different from identity such that $yx \neq xy$. So $(x y) (y x) \neq e$. Putting $z = y$ in the Moufang identity we get $(xy) (zx) = (xy) (yx) \neq e$ but $x ((yz) x) = x ((yy)x) = x.x. = e$. Hence claim.

Now we prove the new class of neutrosophic loop $\langle L_n \cup I\rangle$ does not contain any Bol loop.

**THEOREM 4.1.11:** *The new class of neutrosophic loops $\langle L_n \cup I\rangle$ does not contain any Bol loop.*

*Proof:* Recall that a neutrosophic loop $\langle L \cup I\rangle$ is Bol if it satisfies $((x y)z) y = x ((y z)y)$ for all x, y, z, $\in \langle L \cup I\rangle$.
We prove the results under 2 cases

i. When $L_n(m)$ is not right alternative.
ii. When $L_n(m)$ is right alternative.

**Case i:** If $\langle L_n(m) \cup I\rangle$ is not right alternative then there exists x, $y \in \langle L_n(m) \cup I\rangle$ ($x \neq y$, $x \neq e$ or e I and $y \neq e$ or e I) such that $x(yy) \neq (xy)y$.
Take $z = y$ in

$$\begin{aligned} ((xy)z)y &= x((yz)y) \text{ then} \\ ((xy)y)y &= x((yy)y) \\ ((xy)y)y &\neq xy. \end{aligned}$$

But

$$\begin{aligned} x((yz)y) &= x((yy)y) \\ &= xy. \end{aligned}$$

So when $\langle L_n(m) \cup I\rangle$ is not right alternative, it is not a Bol loop.

**Case ii:** Let $\langle L_n(m) \cup I\rangle$ be a neutrosophic right alternative loop. We know $\langle L_n(2) \cup I\rangle$ is the only right alternative loop in $\langle L_n \cup$

I⟩. If we take x = y = k + 1 and z = k, k < n then we get ((x y) z) y = zy = t where t = (k + 2) (mod n). Now x (yz) y) = r (say) where r = k + 5 (mod n).

If $\langle L_n(2) \cup I\rangle$ is a neutrosophic Bol loop then k + 5 = k + 2 (mod 5) that is n/3 which is not possible as n > 3. So even when $\langle L_n(2) \cup I\rangle$ is right alternative the neutrosophic loop $\langle L_n(2) \cup I\rangle$ is not a Bol loop.

Now we prove the new class of neutrosophic loops do not contain any neutrosophic Bruck loop.

**THEOREM 4.1.12:** *The new class of neutrosophic loops* $\langle L_n \cup I\rangle$ *does not contain a neutrosophic Bruck loop.*

*Proof:* Recall that a neutrosophic loop L is Bruck if it satisfies (x (y x)) z = x (y (xz)) (1) for x, y, z ∈ ⟨L ∪ I⟩. $(xy)^{-1} = x^{-1} y^{-1}$ (2) for x, y ∈ L whenever xy have inverse. We will show that for any neutrosophic loop $\langle L_n(m) \cup I\rangle \in \langle L_n \cup I\rangle$ there exist x, y, z ∈ $\langle L_n(m) \cup I\rangle$ for which (1) is not true. Let $\langle L_n(m) \cup I\rangle$ = {e, 1, 2, …, n, eI, 1I, 2I, …, nI} ∈ $\langle L_n \cup I\rangle$.

**Case i:** If $\langle L_n(m) \cup I\rangle$ is not right alternative take x = z = y + 1 with y < n, then x (y (x z)) = x (y ( x x)) = xy = t where t is given by t = (my – (m – 1)x) (mod n) and (x (y x)) z = (x (y x)) x = r where r is given by r ≡ $((-m^3 + 2m^2 - m + 1)\, x + (m^3 - 2m^2 + m)\, y)$ (mod n). If $\langle L_n(m) \cup I\rangle$ is to be neutrosophic Bruck loop we must have t ≡ r (mod n) so we get $(m^3 - 2m^2)(x - y) \equiv 0$ (mod n) or (m – 2) ≡ 0 (mod n) as ($(m^2, n) = 1$ and x = y + 1) or m = 2 which is a contradiction as $\langle L_n(m) \cup I\rangle$ is not right alternative. So $\langle L_n(m) \cup I\rangle$ is not a neutrosophic Bruck loop in this case.

**Case ii:** Let $\langle L_n(m) \cup I\rangle$ be a right alternative neutrosophic loop. Then it is not left alternative. So there exists x, z ∈ $\langle L_n(m) \cup I\rangle$ (x ≠ e, x ≠ z and z ≠ e) such that (x x) z ≠ x (x z) that is z ≠ x(xz). Now take y = x in (x (y x)) z = x (y (xz)). Then (x (y x))z = (x (xx)) z = xz but x (y (xz)) = x (x (xz)) ≠ xz. So $\langle L_n(m) \cup I\rangle$ is not a neutrosophic Bruck loop.

Now we obtain a necessary and sufficient condition for a neutrosophic loop in the class of loop $\langle L_n \cup I\rangle$ to be a WIP-neutrosophic loop.

**THEOREM 4.1.13:** *Let $\langle L_n(m) \cup I\rangle \in \langle L_n \cup I\rangle$. Then the neutrosophic loop $\langle L_n(m) \cup I\rangle$ is a weak inverse property (WIP) loop if and only if $(m^2 - m + 1) \equiv 0 \pmod n$.*

*Proof:* We know $\langle L \cup I\rangle$ is a WIP neutrosophic loop with identity e; if (x y) z = e implies x (y z) = e then obviously (xI yI) zI = eI and xI (yI zI) = eI where x, y, z, xI, yI, zI $\in \langle L \cup I\rangle$. We show the working for x, y, z on similar lines the result holds good for xI, yI, zI.

Suppose $\langle L_n(m) \cup I\rangle \in \langle L_n \cup I\rangle$ is a WIP loop. Choose x, y, z $\in \langle L_n(m) \cup I\rangle$ such that (x – y, n) = 1 and z = xy. Now z = xy implies z = [my – (m – 1) x] mod n. Since $\langle L_n(m) \cup I\rangle$ is a neutrosophic WIP loop and (xy) z = e, we must have x (yz) = e or yz = x. That is x = [mz – (m – 1)y] mod n.

Putting the value of z from z = (my – (m – 1)x) (mod n) in x = (mz – (m – 1)y) mod n we get $(m^2 - m + 1)(x - y) \equiv 0$ (mod n) or $m^2 - m + 1 \equiv 0$ (mod n).

Conversely if $(m^2 - m + 1) \equiv 0$ (mod n) then it is easy to see that x = (mz – (m – 1) y) (mod n) holds good whenever z = (my – (m – 1) x) (mod n) hold good i.e. (x y) z = e implies x (y z) = e and x, y and z are distinct elements of $\{\langle L_n(m) \cup I\rangle \setminus (e, eI)\}$. However if any one x or y or z is equal to e or x = y then WIP identity holds trivially.

Hence $\langle L_n(m) \cup I\rangle$ is a WIP neutrosophic loop. It is interesting to not left or right alternative neutrosophic loop of $\langle L_n \cup I\rangle$ is not a WIP loop.

We leave it for the reader to check $\langle L_7(3) \cup I\rangle \in \langle L_7 \cup I\rangle$ is a WIP neutrosophic loop. It is still interesting to note that the new class of neutrosophic loops $\langle L_n \cup I\rangle$ does not contain any associative neutrosophic loop.

**THEOREM 4.1.14:** *The new class of neutrosophic loops $\langle L_n \cup I\rangle$ does not contain any associative neutrosophic loop i.e. a neutrosophic group.*

*Proof:* Let $\langle L_n(m) \cup I\rangle \in \langle L_n \cup I\rangle$. To prove $\langle L_n(m) \cup I\rangle$ is not an associative neutrosophic loop it is sufficient to find a triple (x, y, z) in $\langle L_n(m) \cup I\rangle$ such that $(x\ y)\ z \neq x\ (y\ z)$. Take three distinct non identity elements i, j, k in $\langle L_n(m) \cup I\rangle$ such that $(i.j) \neq k$ and $i \neq (j.\ k)$. Now $(i.j)\ k = r$ where r is given by $r = [mk - (m - 1)(mj - (m - 1)i]$ [mod n] and $i(jk) = t$ where t is given by $t = (m(mk - (m - 1)) - (m - 1)\ i]$ (mod n).

If $(i\ j)\ k = i\ (jk)$ then we must have $r \equiv t$ (mod n ) or $(m^2 - m)(k - i) \equiv 0$ (mod n) but $(m^2 - m, n) = 1$. So the above equation gives $k - i \equiv 0$ (mod n) or $k = i$ a contradiction to our assumption. Hence the result.

Having seen all these properties for this new class of neutrosophic loops $\langle L_n(m) \cup I\rangle$ we now define the results in general for neutrosophic loops. It has become pertinent to mention here that we do not have a class of loop which is naturally obtained so we are in a difficult position to see these loops or their neutrosophic analogue.

**DEFINITION 4.1.15:** *Let $\langle L \cup I\rangle$ be a neutrosophic loop of finite order. If P is a neutrosophic subloop of $\langle L \cup I\rangle$ and if $o(P) / o(\langle L \cup I\rangle)$ then we call P a Lagrange neutrosophic subloop. If every neutrosophic subloop is a Lagrange neutrosophic subloop then we call the neutrosophic loop $\langle L \cup I\rangle$ to be a Lagrange neutrosophic loop. If $\langle L \cup I\rangle$ has no Lagrange neutrosophic subloop we call $\langle L \cup I\rangle$ to be a Lagrange free neutrosophic loop. If $\langle L \cup I\rangle$ has atleast one Lagrange neutrosophic subloop then we call $\langle L \cup I\rangle$ to be a weakly Lagrange neutrosophic loop.*

Next we define the notion of p-Sylow neutrosophic subloops.

**DEFINITION 4.1.16:** *Let $\langle L \cup I\rangle$ be a neutrosophic loop of finite order. Let p be a prime such that $p^{\alpha} / o(\langle L \cup I\rangle)$ and $p^{\alpha+1} \nmid o(\langle L \cup I\rangle)$ and if $\langle L \cup I\rangle$ has a proper neutrosophic subloop P of order $p^{\alpha}$ then we call P a p-Sylow neutrosophic subloop. If*

*$\langle L \cup I\rangle$ has atleast one p-Sylow neutrosophic subloop then we call $\langle L \cup I\rangle$ a weakly Sylow neutrosophic loop. If $\langle L \cup I\rangle$ has no p-Sylow neutrosophic subloop then we call $\langle L \cup I\rangle$ a Sylow free neutrosophic loop.*

*If for every prime p such that $p^{\alpha} / o(\langle L \cup I\rangle)$ and $p^{\alpha+1} \nmid o(\langle L \cup I\rangle)$ we have a p-Sylow neutrosophic subloop then we call $\langle L \cup I\rangle$ a Sylow neutrosophic loop. If in addition $\langle L \cup I\rangle$ being a Sylow neutrosophic loop we have for every prime p, $p^{\alpha} / o(\langle L \cup I\rangle)$ and $p^{\alpha+1} \nmid o(\langle L \cup I\rangle)$ we have a neutrosophic subloop of order $p^{\alpha+t}$ $(t \geq 1)$ then we call $\langle L \cup I\rangle$ a super Sylow neutrosophic loop.*

It is important and interesting to make note of the following

(1) Every super Sylow neutrosophic loop is a Sylow neutrosophic loop but a Sylow neutrosophic loop in general is not a super Sylow neutrosophic loop. The reader is requested to construct a non abstract example for the same.
(2) Every Sylow neutrosophic loop is a weakly Sylow neutrosophic loop, however a weakly Sylow neutrosophic loop in general is not a Sylow neutrosophic loop.

Here also a innovative reader can construct a non abstract example of the same. Now we can say only one thing. All neutrosophic loops of order $2^n$ $(n \geq 1)$ cannot have p-Sylow neutrosophic subloop p a prime.

We had talked about Cauchy elements in the new class of neutrosophic loops now we proceed on to study Cauchy element in general neutrosophic loops.

**DEFINITION 4.1.17:** *Let $\langle L \cup I\rangle$ be a neutrosophic loop of finite order. An element $x \in \langle L \cup I\rangle$ is said to be a Cauchy element of $\langle L \cup I\rangle$ if $x^n = 1$ and $n / o(\langle L \cup I\rangle)$; $y \in \langle L \cup I\rangle$ is said to be a Cauchy neutrosophic element of $\langle L \cup I\rangle$ if $y^m = I$ and $m / o(\langle L \cup I\rangle)$. If $y_1^{m_1} = I$ but $m_1 \nmid o(\langle L \cup I\rangle)$ then we call $y_1$ a anti Cauchy neutrosophic element of $\langle L \cup I\rangle$. If $x_1^{n_1} = 1$ and $n_1 \nmid o(\langle L \cup I\rangle)$,*

*$x_1$ is called the anti Cauchy element of $\langle L \cup I\rangle$. If $\langle L \cup I\rangle$ has no anti Cauchy element and no anti neutrosophic Cauchy element then we call $\langle L \cup I\rangle$ to be a Cauchy neutrosophic loop.*

If $\langle L \cup I\rangle$ has atleast one Cauchy neutrosophic element and one Cauchy element then we call $\langle L \cup I\rangle$ to be a weakly Cauchy neutrosophic loop. If $\langle L \cup I\rangle$ has no anti Cauchy element (or no anti Cauchy neutrosophic element) then $\langle L \cup I\rangle$ is called semi Cauchy neutrosophic loop 'or' is used in the strictly mutually exclusive sense.

*It is interesting to see that all Cauchy neutrosophic loops are weakly Cauchy neutrosophic loops but clearly Cauchy neutrosophic loops are in general not a weakly Cauchy neutrosophic loop.*

Also we see every Cauchy neutrosophic loop is semi Cauchy neutrosophic loop and semi Cauchy neutrosophic loop is never a Cauchy neutrosophic loop. However we see we do not have any relation between weakly Cauchy neutrosophic loop or semi Cauchy neutrosophic loop. Interested reader is expected construct examples and counter examples for these types of Cauchy neutrosophic loops. Several more properties can be derived about these neutrosophic loops.

## 4.2 Neutrosophic Biloops

Now we proceed on to define the notion of neutrosophic biloops. The study of biloop have been carried out in [48]. Here we define the several interesting properties about neutrosophic biloops.

**DEFINITION 4.2.1:** *Let $(\langle B \cup I\rangle, *_1, *_2)$ be a non empty neutrosophic set with two binary operations $*_1$, $*_2$, $\langle B \cup I\rangle$ is a neutrosophic biloop if the following conditions are satisfied.*

i. *$\langle B \cup I\rangle = P_1 \cup P_2$ where $P_1$ and $P_2$ are proper subsets of $\langle B \cup I\rangle$.*
ii. *$(P_1, *_1)$ is a neutrosophic loop.*
iii. *$(P_2, *_2)$ is a group or a loop.*

*We now illustrate this by the following example.*

***Example 4.2.1:*** Let (⟨B ∪ I⟩, $*_1$, $*_2$) = ({e, 1, 2, 3, 4, 5, eI, 1I, 2I, 3I, 4I, 5I} ∪ {g | $g^5$ = e'}, $*_1$, $*_2$). Clearly ⟨B ∪ I⟩ = $P_1 \cup P_2$ where $P_1$ = ({e, 1, 2, 3, 4, 5, eI, 1I, 2I, 3I, 4I, 5I}) ∪ {g | $g^5$ = e'} = $P_2$. ⟨B ∪ I⟩ is a neutrosophic biloop. All biloops are not in general neutrosophic biloops.

We can say the order of the neutrosophic biloop is the number of distinct elements of ⟨B ∪ I⟩. If the number of elements in ⟨B ∪ I⟩ is finite we call the neutrosophic biloop to be finite. If the number of elements in ⟨B ∪ I⟩ is infinite we call the neutrosophic biloop to be infinite. The neutrosophic biloop given in the example 4.2.1 is finite. Infact order of a neutrosophic biloop is denoted by o(⟨B ∪ I⟩) and o(⟨B ∪ I⟩) in the above example is 17.

Now we give yet another example of a neutrosophic biloop.

***Example 4.2.2:*** Let ⟨B ∪ I⟩ = {⟨$L_7(3)$ ∪ I⟩ ∪ {Z, group under addition}}. ⟨B ∪ I⟩ is an infinite neutrosophic biloop.

Now we proceed in to define the notion of neutrosophic subbiloops.

**DEFINITION 4.2.2:** *Let (⟨B ∪ I⟩, $*_1$, $*_2$) be a neutrosophic biloop. A proper subset P of ⟨B ∪ I⟩ is said to be a neutrosophic subbiloop of ⟨B ∪ I⟩ if (P, $*_1$, $*_2$) is itself a neutrosophic biloop under the operations of ⟨B ∪ I⟩.*

We make the following observations.

**THEOREM 4.2.1:** *(P, $*_1$) or (P, $*_2$) taken from definition 4.2.2 is not a neutrosophic loop or group.*

*Proof:* We can prove this by example for take (⟨B ∪ I⟩ $*_1$, $*_2$) = {⟨$L_5(3)$ ∪ I⟩, 1, g, $g^2$ $g^3$}. P = {e, 1, $g^2$, eI, 1I, I} is not a loop or group under the operation $*_1$ or $*_2$ but $P = P_1 \cup P_2 = \{e, 1, eI, 1I\} \cup \{1, g^2\}$ is a neutrosophic subbiloop of P.

We prove the following Theorem.

**THEOREM 4.2.2:** *Let $(\langle B \cup I\rangle = B_1 \cup B_2, *_1, *_2)$ be a neutrosophic biloop. $(P, *_1, *_2)$ a proper subset of $\langle B \cup I\rangle$ is a neutrosophic biloop if and only if $P_i = P \cap B_i$, $i = 1, 2$ and $P_1$, is a neutrosophic subloop and $P_2$ is a subgroup.*

*Proof:* If $(P, *_1, *_2)$ is a neutrosophic biloop then we can say $P = (P_1 \cup P_2, *_1, *_2)$ where $(P_1, *_1)$ is a neutrosophic loop and $(P_2, *_2)$ is a group. Thus $P_1 = P \cap B_1$ is a neutrosophic subloop of $B_1$ and $P_2 = P \cap B_2$ is a subgroup of $B_2$. Hence the claim.

Suppose $P = P_1 \cup P_2$ with $P_1 = P \cap B_1$ and $P_2 = P \cap B_2$ are neutrosophic subloop and a subgroup, obviously P is a neutrosophic subbiloop.

Now we illustrate this by the following example.

***Example 4.2.3:*** Let $(B = B_1 \cup B_2, *_1, *_2)$ be a neutrosophic biloop where $B_1 = \{\langle L_3(2) \cup I\rangle\}$ and $B_2 = D_{2.6}$. Take $P = (\{e, eI, 4, 4I\} \cup \{b, b^2, b^3, b^4, b^5, b^6 = 1\} *_1, *_2)$, P is a neutrosophic biloop of B but P as a set is neither a neutrosophic loop under $*_1$ nor a group under the binary operation $*_2$.

Now we see o(B) = 24 and o(P) = 10 but 10 $\not|$ 24 so we see the Lagrange theorem for finite groups is not satisfied. Take $V = (V_1 \cup V_2, *_1, *_2)$ where $V_1 = \{e, eI, 2, 2I\}$ and $V_2 = \{g^3, 1\}$. V is a neutrosophic subbiloop and o(V) / o(B) i.e. 6/24.

Now we make some more definitions in this direction.

**DEFINITION 4.2.3:** *Let $(B = B_1 \cup B_2, *_1, *_2)$ be a finite neutrosophic biloop. Let $P = (P_1 \cup P_2, *_1, *_2)$ be a neutrosophic biloop. If o(P) / o(B) then we call P a Lagrange neutrosophic subbiloop of B.*

*If every neutrosophic subbiloop of B is Lagrange then we call B to be a Lagrange neutrosophic biloop. If B has atleast one Lagrange neutrosophic subbiloop then we call B to be a*

*weakly Lagrange neutrosophic biloop. If B has no Lagrange neutrosophic subbiloops then we call B to be a Lagrange free neutrosophic biloop.*

It is easy to verify all Lagrange neutrosophic biloops are weakly Lagrange neutrosophic biloops.

However the converse is not true. It is left as an exercise for the reader to find or characterize those neutrosophic biloops which are (1) Lagrange (2) weakly Lagrange.

Now we proceed on to prove some interesting theorem.

**THEOREM 4.2.3:** *All neutrosophic biloops of prime order are Lagrange free.*

*Proof:* Let $B = (B_1 \cup B_2, *_1, *_2)$ be a finite neutrosophic biloop of order p, p a prime. Suppose P is any neutrosophic subbiloop of B then clearly $o(P) \nmid o(B)$. Thus $(B = B_1 \cup B_2, *_1, *_2)$ is a Lagrange free neutrosophic biloop.

Now we give another example.

***Example 4.2.4:*** Let $(B = (B_1 \cup B_2), *_1, *_2)$ where $B_1 = \{\langle L_7(3) \cup I\rangle *_1\}$ and $B_2 = \{g \mid g^7 = 1\}$; B is a neutrosophic biloop of order 23.

Take $P = P_1 \cup P_2 = \{e, eI, 3I, 3\} \cup \{g \mid g^7 = 1\}$, P is a neutrosophic subbiloop of B and o(P) = 11 and 11/ 23. Thus what ever be the neutrosophic subbiloops we see B is a Lagrange free neutrosophic biloop.

Yet we give another example.

***Example 4.2.5:*** Let $B = (B_1 \cup B_2, *_1, *_2)$ where $B_1 = \{\langle L_5(3) \cup I\rangle, *_1\}$ and $B_2 = \{g \mid g^8 = 1\}$, o(B) = 20. Take $(P = P_1 \cup P_2, *_1, *_2)$ where $P_1 = \{e, eI, 3, 3I\} \subset B_1$ and $P_2 = \{1, g^4\} \subset B_2$. P is a neutrosophic subbiloop and o(P) = 6 and $6 \nmid 20$.

Take $T = T_1 \cup T_2$ where $T_1 = \{e, eI, 2, 2I\} \subset B_1$ and $T_2 = \{g^2, g^4, g^6, 1\}$; o(T) = 8 and $8 \nmid 20$. So T is not a Lagrange neutrosophic subbiloop. P is also not a Lagrange neutrosophic subbiloop.

Take S = $\{S_1 \cup S_2, *_1, *_2\}$ where $S_1$ = {e, eI, 3, 3I} and $S_2$ = $\{1\} \subset B_2$ we see S is a neutrosophic subbiloop of order 5 and 5 / 20. So S is a Lagrange neutrosophic subbiloop. We see the neutrosophic loop $\langle L_5(3) \cup I\rangle$, can have neutrosophic biloops of order only 4 and $\langle g \mid g^8 = 1\rangle$ can have subloops of order 1, 2 and 4. Thus we see B is a weakly Lagrange neutrosophic biloop.

We see in the above example the neutrosophic loop and the group satisfy the Lagrange theorem separately but as a neutrosophic biloop, they do not satisfy the Lagrange theorem for finite groups.

Now we as in case of other neutrosophic bistructures here we define Cauchy element and Cauchy neutrosophic element of a neutrosophic bigroup.

**DEFINITION 4.2.4:** *Let B = $(B_1 \cup B_2, *_1, *_2)$ be a neutrosophic biloop of finite order. An element $x \in B$ such that $x^n = 1$ is called a Cauchy element if n / o(B) otherwise x is an anti Cauchy element of B.*

*We call an element y of B with $y^m = I$ to be a Cauchy neutrosophic element if m / o(B); otherwise y is a anti Cauchy neutrosophic element of B. If every element in B is either a Cauchy element or a Cauchy neutrosophic element then we call B to be a Cauchy neutrosophic biloop.*

We now illustrate this with some examples.

***Example 4.2.6:*** Let B = $(B_1 \cup B_2, *_1, *_2)$ where $B_1 = \langle L_7(3) \cup I\rangle$, a neutrosophic loop and $B_2 = \{g \mid g^{16}= 1\}$ cyclic group of order 16. Clearly o(B) = 32.

It is easily verified that B is a Cauchy neutrosophic biloop.

We give get another example of a Cauchy neutrosophic biloop.

***Example 4.2.7:*** Let $(B = B_1 \cup B_2, *_1, *_2)$ be a neutrosophic biloop where $B_1 = \langle L_5(3) \cup I\rangle$ and $B_2 = \langle g \mid g^4 = 1\rangle$. B is a Cauchy neutrosophic biloop.

We give another example.

***Example 4.2.8:*** Let ($B = B_1 \cup B_2, *_1, *_2$) be a neutrosophic biloop when $B_1 = \langle L_5(3) \cup I\rangle$ and $B_2 = \langle g \mid g^7 = 1\rangle$, o(B) = 19. No element in B is a Cauchy element. So B is a Cauchy free neutrosophic biloop.

Now we can give a characterization theorem for a neutrosophic biloop for which the neutrosophic loop is taken from the new class of neutrosophic loops to be Cauchy neutrosophic biloops by choosing a proper group.

**THEOREM 4.2.4:** *Let $B = (B_1 \cup B_2, *_1, *_2)$ where $B_1$ is a new class neutrosophic loop of order 2 (n + 1) and choose group, $B_2$ so that $2 (n + 1) + 2^u = 2^t$ where $2^u = o(B_2)$. B is a Cauchy neutrosophic loop, only when $n + 1 = 2^w$ for some w.*

*Proof:* Let $\langle L_n(m) \cup I\rangle$ be a neutrosophic loop from the new class of neutrosophic loops. Order of $\langle L_n(m) \cup I\rangle = 4t = 2 (n + 1)$. It $t = 2^s$ then choose $B_2$ to be a group of order $2^u$ then $(\langle L_n(m) \cup I\rangle \cup B_2)$ will be a Cauchy neutrosophic biloop.

It t is not of the form $2^s$ then choose the order of $B_2$ to a number m such that $4t + m = 2^r$. This is always possible but we cannot in general say the elements in $B_2$, will be Cauchy with respect to $B = (\langle L_n(m) \cup I\rangle \cup B_2)$.

Thus is very clear from the following example. For choose $B = (\langle L_n(m) \cup I\rangle \cup B_2)$ a neutrosophic biloop.

The order of $\langle L_{25}(m) \cup I\rangle = 52$, choose order of $B_2$ to be 12 a cyclic group of order 12 or $A_4$ or $D_{26}$. Now $o(B) = 64 = 2^6$. But we have $x \in B$ such that $x^3 = 1$, so x is not Cauchy. Hence our claim.

One of the most important question is that if ($B = B_1 \cup B_2, *_1, *_2$) is a finite neutrosophic biloop say order of B is n. Can we say for every t / n as in case of finite groups say B has a Cauchy element x such that $x^t = 1$ and Cauchy neutrosophic element such that $y^{t_1} = I$. We say this is impossible in general in case of bistructures so even in case of neutrosophic biloops.

We illustrate this by the following example.

***Example 4.2.9:*** Let B = (⟨$L_{17}(3) \cup I$⟩ ∪ $B_2$) where $B_2$ = ⟨g | $g^6$ = 1⟩ be a finite neutrosophic biloop of order 42. Clearly 7/42 but B has no element x such that $x^7 = 1$ or no element y such that $y^7 = I$, for every element in ⟨$L_{17}(3) \cup I$⟩ is such that $x^2 = 1$ or $(xI)^2 = I$ and of course $B_2$ is a cyclic group of order 6 so every element in $B_2$ is of order 2 or 3 only.

Hence the claim. That is we redo the definition and call a finite neutrosophic biloop to be Cauchy as Cauchy condition is not true in case of biloop.

Now we have already talked about Sylow property for the new class of neutrosophic loops. We will define and discuss about Sylow property in case of finite neutrosophic biloop. Before this we will define a new class of neutrosophic loops.

**DEFINITION 4.2.5:** *Let B = ($B_1 \cup B_2$, $*_1$, $*_2$) be a neutrosophic biloop. We say B is a new class of neutrosophic biloops if $B_1$ ∈ ⟨$L_n \cup I$⟩ and $B_2$ is any group. i.e. B = (⟨$L_n(m) \cup I$⟩ ∪ $B_2$, $*_1$, $*_2$).*

Now we proceed on to define Sylow structure on the neutrosophic biloops.

**DEFINITION 4.2.6:** *Let (B = $B_1 \cup B_2$, $*_1$, $*_2$) be a neutrosophic biloop of finite order. Let p be a prime such that $p^\alpha$ / o(⟨B⟩) and $p^{\alpha+1} \nmid$ o(B). If B = $B_1 \cup B_2$ has a neutrosophic subbiloop P of order $p^\alpha$ then we call P the p-Sylow neutrosophic subbiloop of B.*

*If (B = $B_1 \cup B_2$, $*_1$, $*_2$) has atleast one p-Sylow neutrosophic subbiloop then we call B a weakly Sylow neutrosophic biloop. If B has no p-Sylow neutrosophic subbiloop then we call B a Sylow free neutrosophic biloop.*

*If for every prime p such that $p^\alpha$ / o($B_1 \cup B_2$) and $p^{\alpha+1} \nmid$ o($B_1 \cup B_2$) we have a p-Sylow neutrosophic subbiloop then we call (B = $B_1 \cup B_2$, $*_1$, $*_2$) to be a Sylow-neutrosophic biloop.*

We first give some examples before we proceed to prove some of its properties.

***Example 4.2.10:*** Let $(B = B_1 \cup B_2, *_1, *_2)$ be a neutrosophic biloop where $B_1 = \langle L_5(4) \cup I \rangle$ and $B_2 = \langle g \mid g^8 = 1 \rangle$. Clearly o(B) = 20. Now $2^2 / 20$ and $2^3 \not| \; 20$. Also $5 / 20$ and $5^2 \not| \; 20$.

This neutrosophic biloop has no 2-Sylow neutrosophic subbiloop but has a 5-Sylow neutrosophic subbiloop given by $P = (P_1 \cup P_2, *_1, *_2)$ where $P_1 = \{eI, e, 3, 3I\}$; $P_2 = \{1\}$.
Thus B is only a weakly Sylow neutrosophic biloop.

We have the following interesting observations. We see all Sylow neutrosophic biloops are weakly Sylow neutrosophic biloop. Clearly a weakly Sylow neutrosophic biloop is not a Sylow neutrosophic biloop.

Now we also prove we have a non trivial class of Sylow free neutrosophic biloops.

***Example 4.2.11:*** Let $(\langle B \cup I \rangle = B_1 \cup B_2, *_1, *_2)$ be a finite neutrosophic biloop where $B_1 = \langle L_9(8) \cup I \rangle$ and $B_2 = \{S_3\}$, the symmetric group of degree 3. o(B) = 26, $2 / 26$ and $2^2 \not| \; 26$ also $13/26$ but $13^2 \not| \; 26$. B has no 2-Sylow neutrosophic subloop or 13-Sylow neutrosophic subloop. Thus B is a Sylow free neutrosophic biloop.

We define neutrosophic Moufang biloop, neutrosophic Bol biloop, neutrosophic Bruck biloop, neutrosophic WIP-biloop and so on.

**DEFINITION 4.2.7:** *Let $(B = B_1 \cup B_2, *_1, *_2)$ be a neutrosophic biloop. We say B is a neutrosophic Moufang biloop if B for all proper subsets $(P = P_1 \cup P_2, *_1, *_2)$ where P is a neutrosophic subbiloop of B in which $P_1$ is a proper neutrosophic Moufang subloop of $B_1$.*

Thus we do not demand the whole neutrosophic loop to satisfy the Moufang identity but every neutrosophic subloop satisfies the Moufang identity for we can have several such neutrosophic biloop by varying the subgroups. We have a class of neutrosophic biloops to satisfy the Moufang identity.

***Example 4.2.12:*** Let $(B = (B_1 \cup B_2), *_1, *_2)$ be a neutrosophic biloop where $B_1 = \langle L_5(3) \cup I\rangle$ and $B_2 = A_4$, B is a neutrosophic Moufang biloop.

For take any neutrosophic subloop of $\langle L_5(3) \cup I\rangle$. It will be of the form (e, eI, 3I, 3). It is easily verified (e, eI, 3, 3I) satisfies the Moufang identify.

Thus $(P = P_1 \cup P_2, *_1, *_2)$ where

$$P_1 = \{e, eI, 3I, 3\} \text{ and}$$

$$P_2 = \left\{\begin{pmatrix} 1 & 2 & 3 & 4 \\ 1 & 2 & 3 & 4 \end{pmatrix}, \begin{pmatrix} 1 & 2 & 3 & 4 \\ 2 & 1 & 4 & 3 \end{pmatrix}\right\}$$

is clearly a neutrosophic subbiloop.

So B is a Moufang neutrosophic biloop. Thus we have the following interesting result. We saw the new class of neutrosophic loops was not Moufang but by defining Moufang neutrosophic biloop we see this new class of neutrosophic biloops $\langle L_n \cup I\rangle$ are Moufang provided n is a prime.

**THEOREM 4.2.5:** *Let $(B = (B_1 \cup B_2), *_1, *_2)$ be the new class neutrosophic biloops B is a neutrosophic Moufang biloop if $B_1 \in \langle L_n \cup I\rangle$ where n is a prime and $B_2$ any group.*

*Proof:* Given $B = (B_1 \cup B_2), *_1, *_2)$ is a neutrosophic biloop from the new class of neutrosophic biloops i.e. $B_1 \in \langle L_n \cup I\rangle$. To show if n is a prime B is a Moufang neutrosophic biloop.

We know if n is a prime then every loop in the class of loops $\langle L_n \cup I\rangle$ has only 2 types of neutrosophic subbiloop (1) the whole neutrosophic loop $L_n(m)$ (2) {e, eI, tI, t} where $1 \leq t \leq n$ a neutrosophic subloop of order 4 and it has no other subloop.

We have just proved all neutrosophic loops given by the following table.

| * | e | eI | t | tI |
|---|---|---|---|---|
| e | e | eI | t | tI |
| eI | eI | eI | tI | tI |
| t | t | tI | e | eI |
| tI | tI | tI | eI | eI |

satisfies the Moufang identity. Hence $\{B = (B_1 \cup B_2), *_1, *_2\}$ where $B_1 = \langle L_n(m) \cup I\rangle \in \langle L_n \cup I\rangle$ is a neutrosophic Moufang biloop whatever be $B_2$, we show n is a prime is essential for n is a prime we show the neutrosophic biloop is not Moufang for the neutrosophic loop $\langle L_n(m) \cup I\rangle$ (n-not a prime) may have neutrosophic subloops which do not satisfy the Moufang identity.

We illustrate this by the following example.

***Example 4.2.13:*** Take $(B = (B_1 \cup B_2), *_1, *_2)$ where $B_1 = \langle L_{15}(2) \cup I\rangle$ and $B_2 = \langle A_4\rangle$. B is a neutrosophic biloop. Take $P = (P_1 \cup P_2, *_1, *_2)$ where $P_1$ = {e, 2, 5, 8, 11, 14, eI, 2I, 5I, 8I, 11I, 14I}; $P_1$ is a not Moufang subloop given by the table.

| | e | 2 | 5 | 8 | 11 | 14 | eI | 2I | 5I | 8I | 11I | 14I |
|---|---|---|---|---|---|---|---|---|---|---|---|---|
| e | e | 2 | 5 | 8 | 11 | 14 | eI | 2I | 5I | 8I | 11I | 14I |
| 2 | 2 | e | 8 | 14 | 5 | 11 | 2I | eI | 8I | 14I | 5I | 11I |
| 5 | 5 | 14 | e | 11 | 2 | 8 | 5I | 14I | eI | 11I | 2I | 8I |
| 8 | 8 | 11 | 2 | e | 14 | 5 | 8I | 11I | 2I | eI | 14I | 5I |
| 11 | 11 | 8 | 14 | 5 | e | 2 | 11I | 8I | 14I | 5I | eI | 2I |
| 14 | 14 | 5 | 11 | 2 | 8 | I | 14I | 5I | 11I | 2I | 8I | eI |
| eI | eI | 2I | 5I | 8I | 11I | 14I | eI | 2I | 5I | 8I | 11I | 14I |
| 2I | 2I | eI | 8I | 14I | 5I | 11I | 2I | eI | 8I | 14I | 5I | 11I |
| 5I | 5I | 14I | eI | 11I | 2I | 8I | 5I | 14I | eI | 11I | 2I | 8I |
| 8I | 8I | 11I | 2I | eI | 14I | 5I | 8I | 11I | 2I | eI | 14I | 5I |
| 11I | 11I | 8I | 14I | 5I | eI | 2I | 11I | 8I | 14I | 5I | eI | 2I |
| 14I | 14I | 5I | 11I | 2I | 8I | eI | 14I | 5I | 11I | 2I | 8I | eI |

For take in the identity (xy) (zx) = (x(yz)) x
x = 2  y = 5  and z = 8

Now (xy) (zx) = (25) (82) = 8.11 = 14.

$$\begin{aligned}(x(yz))x &= (2(58))\,2 \\ &= (2.11)\,2 \\ &= 5.2 = 14\end{aligned}$$

satisfies Moufang identity.

Now in the identity
(xy) (zx) = (z (yz)) x
take x = 8, y = 14 and z = 2.
Now

$$\begin{aligned}(xy)\,(zx) &= (8.14)\,(2.8) \\ &= 5.14 = 8. \\ (x(yz))x &= [8\,(14.2)]\,.8 \\ &= (8.5).8 \\ &= 2.8 \\ &= 14.\end{aligned}$$

Thus (xy)(zx) ≠ (x(yz)) x, when x = 8, y = 14 and z = 2. So this neutrosophic subloop is not a neutrosophic Moufang subloop. Hence by taking $\{T = T_1 \cup T_2, *_1, *_2\}$ where $T_1 = P_1$ and $P_2$ = $\left\{\begin{pmatrix}1 & 2 & 3 & 4\\ 2 & 1 & 4 & 3\end{pmatrix}, \begin{pmatrix}1 & 2 & 3 & 4\\ 1 & 2 & 3 & 4\end{pmatrix}\right\}$ we see T is not a neutrosophic Moufang subbiloop. So $\{\langle L_{15}(2) \cup I\rangle \cup A_4\}$ is not a neutrosophic Moufang biloop as it has neutrosophic subbiloop which are not Moufang.

Now we proceed on to define the notion of neutrosophic Bol biloop.

**DEFINITION 4.2.8:** *Let $(B = B_1 \cup B_2, *_1, *_2)$ be a neutrosophic biloop, we say B is a neutrosophic Bol biloop if every proper subset $P = (P_1 \cup P_2, *_1, *_2)$ of B which is a neutrosophic subbiloop of B is such that $(P_1, *_1)$ is a Bol loop. i.e. every proper neutrosophic subloop of the neutrosophic loop $(B_1, *_1)$ must be a Bol loop.*

***Example 4.2.14:*** Let B = $\{B_1 \cup B_2, *_1, *_2\}$ be a neutrosophic biloop where $B_1 = \langle L_5(3) \cup I\rangle$ and $B_2 = A_3$. $B_1$ is a neutrosophic Bol biloop. For the proper neutrosophic subloops of $\langle L_5(3) \cup I\rangle$ are only of the form {e, eI, t, tI} given by the table.

| | e | eI | t | tI |
|---|---|---|---|---|
| e | e | eI | t | tI |
| eI | eI | eI | tI | tI |
| t | t | tI | e | eI |
| tI | tI | tI | eI | eI |

It is easily verified the neutrosophic subloop satisfies the Bol identity. Thus we have the following interesting theorem.

**THEOREM 4.2.6:** *Let $(B = B_1 \cup B_2, *_1, *_2)$ be a neutrosophic biloop from the new class of neutrosophic biloops i.e. $(B_1 \in \langle L_n \cup I\rangle)$, B is a Bol neutrosophic biloop when n is a prime.*

*Proof:* Given $(B = B_1 \cup B_2, *_1, *_2)$ is a neutrosophic biloop from the new class of neutrosophic biloops. So $B_1 \in \langle L_n \cup I\rangle)$ and $B_2$ is any group.

It is given n is a prime. Let $B_1 = \langle L_n(m) \cup I\rangle$ where n is a prime say p. When n is a prime the only neutrosophic subloops of $\langle L_n(m) \cup I\rangle$ are $\langle L_n(m) \cup I\rangle$ and neutrosophic subloops of order 4 given by $P_1$ = {e, eI, t, tI}, $1 \leq t \leq n$; $P_1$ given by the table.

| | e | eI | t | tI |
|---|---|---|---|---|
| e | e | eI | t | tI |
| eI | eI | eI | tI | tI |
| t | t | tI | e | eI |
| tI | tI | tI | eI | eI |

Now we have to verify the neutrosophic subloops {e, eI, t, tI} satisfies the Bol identity ((xy)z)y = (x ((yz)y) for x, y, z, ∈ {e, eI, t, tI}. Clearly Bol identity is satisfied by {e, eI, t, tI}. So when n is a prime $(\langle L_n(m) \cup I\rangle \cup B_2)$ is a neutrosophic Bol biloop when n is not a prime we will show $\langle L_n(m) \cup I\rangle$ need not have all of its neutrosophic subloops to satisfy Bol identity.

We illustrate this by the following example.

***Example 4.2.15:*** Let B = ($B_1 \cup B_2$, $*_1$, $*_2$) be a neutrosophic biloop where $B_1 = \langle L_{15}(2) \cup I\rangle$, n is not a prime n = 15. To show all its neutrosophic subloops are not Bol neutrosophic loops.

Consider a neutrosophic subloop given by $P_1$ = {e, 1, 4, 7, 10, 13, eI, 1.I, 4I, 7I, 10I, 13I}. Let $P_2$ be a subgroup of $B_2$. To prove B is not a Bol neutrosophic biloop, it is sufficient if we show $B_1 = \langle L_{15}(2) \cup I\rangle$ has a neutrosophic subloop which does not satisfy the Bol identity ((xy)z) y = (z ((yz)y) Let x = 7, y = 10 and z = 13.

$$
\begin{aligned}
[(xy)z]y &= [(7.10)13]\ 10 \\
&= \{[2 \times 10 - 1.7]\ (\text{mod } 15)]\ 13\}\ 10 \\
&= (13.13)\ .10 = e.10 = 10.
\end{aligned}
$$

$$
\begin{aligned}
x[(yz)y] &= 7\ [(10.13).10] \\
&= 7\{[(2 \times 13 - 1.10) \text{ mod } 15].\ 10\} \\
&= 7\{1.0\} = 7\ [\{2 \times 10 - 1\}\ (\text{mod } 15)] \\
&= 7.4\ (4 \times 2 - 17)\ (\text{mod } 15) = 1.
\end{aligned}
$$

Since the neutrosophic subloop is not Bol we see (P = $P_1 \cup P_2$, $*_1$, $*_2$) does not satisfy the Bol identity so (B = $B_1 \cup B_2$, $*_1$, $*_2$) = ($\langle L_{15}(2) \cup I\rangle \cup B_2$, $*_1$, $*_2$) is not a neutrosophic Bol loop.

Now we proceed on to define the notion of neutrosophic alternative biloops.

**DEFINITION 4.2.9:** *Let B = ($B_1 \cup B_2$, $*_1$, $*_2$) be a neutrosophic biloop. B is said to be a neutrosophic right (left) alternative biloop if every proper neutrosophic subbiloop of B is neutrosophic right (left) alternative biloop of B.*

***Example 4.2.16:*** Let (B = $B_1 \cup B_2$, $*_1$, $*_2$) be a neutrosophic biloop where $B_1 = \langle L_5(3) \cup I\rangle$ and $B_2$ any group. Any neutrosophic subloop of $B_1 = \langle L_5(3) \cup I\rangle$ satisfies the right (left) alternative identity. For the only type of neutrosophic subloops are $P_1 = \{e, eI, t, tI\}$; $1 \leq t \leq 5$. Clearly $P_1$ satisfies the right (left) alternative identity. So $B = (B_1 \cup B_2, *_1, *_2)$ is a neutrosophic right (left) alternative biloop.

We prove the following theorem.

**THEOREM 4.2.7:** *Let $(B = B_1 \cup B_2, *_1, *_2)$ be a neutrosophic biloop where $B_1 \in \langle L_n \cup I\rangle$. B is a neutrosophic right (left) alternative biloop only when n is a prime.*

*Proof:* Given $(B = B_1 \cup B_2, *_1, *_2)$ is a neutrosophic biloop taken from the new class of neutrosophic biloops i.e. $B_1 \in \langle L_n \cup I\rangle$. It is assumed from the theorem that n is a prime. To show $B_1$ is a neutrosophic right (left) alternative loop. Since when n is prime the only subloops of $\langle L_n(m) \cup I\rangle \in \langle L_n \cup I\rangle$ are the neutrosophic subloops of order 4 and the loop $\langle L_n(m) \cup I\rangle$. The subloops take only the form {e, eI, t, tI} so every neutrosophic subbiloop of $B = (B_1 \cup B_2)$ are (left) right alternative. Hence the claim.

It is important mention if n is not a prime $B = (B_1 \cup B_2)$ where $B_1 = \langle L_n(m) \cup I\rangle$, when m = 2 we get the right alternative neutrosophic biloop $\langle L_n(2) \cup I\rangle$,. If $m = n - 1$ then we get the left alternative neutrosophic biloop as the loops of $\langle L_n(2) \cup I\rangle$ and $\langle L_n(n - 1) \cup I\rangle$,are right and left alternative respectively.

We prove this by the following example.

***Example 4.2.17:*** Let $(B = B_1 \cup B_2, *_1, *_2)$ be a neutrosophic biloop where $B_1 = \langle L_{45}(8) \cup I\rangle$; n = 45 a non prime and $B_2$ is any group. We will show a neutrosophic subbiloop of B, which does not satisfy right (left) alternative identity.

Take $P = (P_1 \cup P_2, *_1, *_2)$ where $P_1$ = {e, 1, 6, 11, 16, 21, 26, 31, 36, 41, eI, 1.I, 6I, 11I, 16I, 21I, 26I, 31I, 36I and 41I} and $P_2$ any subgroup of $B_2$. Clearly $P_1$ is a neutrosophic subloop of $L_{45}(8)$.

Consider the identity (xy) y = x (yy). Take x = 16 and y = 41

$$
\begin{aligned}
(x\,y)y &= [16.41]\,41 \\
&= [(41 \times 8 - 16 \times 7) \ (\text{mod } 45)]\,41 \\
&= 36.41 = (41 \times 8 - 36 \times 7) \ (\text{mod } 45) \\
&= 31. \\
x\,(yy) &= 16\,(41.41) = 16e \\
&= 16.
\end{aligned}
$$

So $x(yy) \neq (xy)y$ hence $(P = P_1 \cup P_2)$ does not satisfy the right alternative identity. Thus $(\langle L_n(m) \cup I\rangle \cup B_2, *_1, *_2)$ is not a neutrosophic right alternative biloop. Take the left alternative identity.

$$
\begin{aligned}
(xx)y &= x\,(xy) \\
x &= 16 \\
y &= 41 \\
(xx)y &= (16.16)\,41 \\
&= 41. \\
x(xy) &= 16\,(16.41) \\
&= 16\,[[41 \times 8 - 16 \times 7]\ (\text{mod } 45)] \\
&= 16\,[36] \\
&= 36 \times 8 - 16 \times 7 \\
&= 41.
\end{aligned}
$$

So this set satisfies the left alternative identity. Take $x = 36$ and $y = 6$

$$
\begin{aligned}
(xx)y &= (36.36).6 = 6(\text{mod } 45) \\
x\,(xy) &= 36\,(36.6)\ \ (\text{mod } 45) \\
&= 36\,[6 \times 8 - 36 \times 7]\ (\text{mod } 45) \\
&= 36.21\ (\text{mod } 45) \\
&= 21 \times 8 - 36 \times 7\ (\text{mod } 45) \\
&= 6\ (\text{mod } 45).
\end{aligned}
$$

satisfies left alternative identity.

Thus the neutrosophic biloop $\{\langle L_{45}(8) \cup I\rangle \cup B_2\}$ is not a neutrosophic right alternative biloop. The reader is advised to check whether the biloop is left alternative.

How ever we show that we have a class of neutrosophic left alternative biloops and neutrosophic right alternative biloops even when n is not a prime.

**THEOREM 4.2.8:** *Let* $(B = B_1 \cup B_2, *_1, *_2)$ *be a neutrosophic biloop. If* $B_1 = \langle L_n(n-1) \cup I \rangle$ *and* $B_2$ *any group* ($n > 3$ *any odd number*) *B is always a neutrosophic left alternative biloop.*

*Proof:* $B_1 = \langle L_n(n-1) \cup I \rangle$ ($n > 3$, n any odd) is always a left alternative neutrosophic loop is the new class of neutrosophic loop of $\langle L_n \cup I \rangle$. So all neutrosophic subloops of $L_n(n-1)$ are left alternative. Hence every proper neutrosophic subbiloop of $(\langle L_n(n-1) \cup I \rangle \cup B_2, *_1, *_2)$ for any group $B_2 = G$ is left alternative. Thus $(\langle L_n(n-1) \cup I \rangle \cup B_2)$ is a neutrosophic left alternative biloop. Hence the claim.

**THEOREM 4.2.9:** *Let* $(B = B_1 \cup B_2, *_1, *_2)$ *be a neutrosophic biloop where* $B_1 = \langle L_n(2) \cup I \rangle \in \langle L_n \cup I \rangle$ *is always a neutrosophic right alternative biloop for all n* ($n > 3$ *and n an odd number*) *and* $B_2$ *any group.*

*Proof:* Given $(B = B_1 \cup B_2, *_1, *_2)$ is a neutrosophic biloop with $B_1 = \langle L_n(2) \cup I \rangle$ ($n > 3$ and n an odd number). We know $\langle L_n(2) \cup I \rangle$ is a neutrosophic right alternative loop. So all its neutrosophic subloops are also right alternative. Thus every proper neutrosophic subbiloop of $(\langle L_n(2) \cup I \rangle \cup B_2)$ is a neutrosophic right alternative subbiloop. Hence $(\langle L_n(2) \cup I \rangle \cup B_2)$ is a neutrosophic right alternative biloop.

Now we proceed on to define neutrosophic normal subbiloop of a neutrosophic biloop.

**DEFINITION 4.2.10:** *Let* $(B = B_1 \cup B_2, *_1, *_2)$ *be a neutrosophic biloop. A neutrosophic subbiloop* $H = (H_1 \cup H_2, *_1, *_2)$ *is said to be a neutrosophic normal subbiloop of B if*

i. $x H_1 = H_1 x$
ii. $(H_1 x) y = H_1 (xy)$

*iii.* $y (x H_1) = (y x) H_1$ *for all* $x, y \in B_1$
*iv.* $H_2$ *is a normal subgroup of* $B_2$.

*We call a neutrosophic biloop to be a simple neutrosophic biloop if it has no nontrivial neutrosophic normal subbiloops.*

We now give an example of a neutrosophic biloop.

***Example 4.2.18:*** Let $(B = B_1 \cup B_2, *_1, *_2)$ be a neutrosophic biloop where $B_1 = \langle L_{11}(2) \cup I\rangle$ and $B_2 = A_5$. (The alternating subgroup of the symmetric group of $S_5$). B has no neutrosophic normal subloop. So $B = B_1 \cup B_2$ is a simple neutrosophic biloop.

Now it may happen that one of the neutrosophic loop or the group has normal subloop or normal subgroup. How to differentiate this for this we define the notion of seminormal and semisimple.

**DEFINITION 4.2.11:** *Let* $(B = B_1 \cup B_2, *_1, *_2)$ *be a neutrosophic biloop if only one of the neutrosophic loop or the group is simple then we call the neutrosophic biloop to be a semi-simple neutrosophic biloop.*

Here it is important to note that the term 'or' is used in the mutually exclusive sense.

Now we give an example of the same.

***Example 4.2.19:*** Let $B = (B_1 \cup B_2, *_1, *_2)$ be a neutrosophic biloop. Suppose $B_1 = \langle L_7(3) \cup I\rangle$ is simple and $B_2 = S_n$ $(n \geq 5)$ the symmetric group of degree n then B is a neutrosophic semisimple biloop and the neutrosophic subbiloop $P = P_1 \cup P_2$ where $P_1 = \{e, eI, t, tI\} \cup \{A_n\}$ is a neutrosophic seminormal biloop of B $(n \geq 5)$.

Now proceed on to define the notion of strong neutrosophic biloop.

**DEFINITION 4.2.12:** *Let ($B = B_1 \cup B_2$, $*_1$, $*_2$) be non empty set with two binary operations. B is said to be strong neutrosophic biloop.*

i. *$B = B_1 \cup B_2$ is a union of proper subsets of B.*
ii. *($B_1$, $*_1$) is a neutrosophic loop.*
iii. *($B_2$, $*_2$) is a neutrosophic group.*

Now we just illustrate this by the following example.

***Example 4.2.20:*** Let ($B = B_1 \cup B_2$, $*_1$, $*_2$) where $B_1 = \langle L_5(2) \cup I \rangle$ is a neutrosophic loop and $B_2$ = {1, 2, 3, 4, I, 2I, 3I, 4I} under multiplication modulo 5 is a neutrosophic group. Thus B is a strong neutrosophic biloop.

All strong neutrosophic biloops are neutrosophic biloops but clearly a neutrosophic biloop in general is not a strong neutrosophic biloop.

All results derived and definitions given in case of neutrosophic biloops can be derived with appropriate modification. Thus we don't derive or define any of them but leave it as an exercise for the reader to prove.

Now we define a neutrosophic biloop of type II as follows.

**DEFINITION 4.2.13:** *Let B = ($B_1 \cup B_2$, $*_1$, $*_2$) be a proper set on which is defined two binary operations $*_1$ and $*_2$. We call B a neutrosophic biloop of type II if the following conditions are satisfied*

.

i. *$B = B_1 \cup B_2$ where $B_1$ and $B_2$ are proper subsets of B.*
ii. *($B_1$, $*_1$) is a neutrosophic loop.*
iii. *($B_2$, $*_2$) is a loop.*

*(Clearly we call a neutrosophic biloop of type I is one in which $B_1$ is a neutrosophic loop and $B_2$ is a group).*

Now type II neutrosophic biloop enjoy some more extra properties. First we give some examples.

***Examples 4.2.21:*** Let B = ($B_1 \cup B_2$, $*_1$, $*_2$) where $B_1 = \langle L_7(3) \cup I \rangle$ and $B_2 = L_5(2)$ then B is a neutrosophic biloop of type II. All properties defined in case of neutrosophic biloops of type I can be easily extended to type II neutrosophic biloops.

Now we define a new class of neutrosophic biloops of type II as follows.

**DEFINITION 4.2.14:** *Let B = ($B_1 \cup B_2$, $*_1$, $*_2$) be a neutrosophic biloop of type II. We say B is a new class of neutrosophic biloop of type II if ($B_1$, $*_1$) = $\langle L_n(m) \cup I\rangle$ i.e. $B_1 \in \langle L_n \cup I\rangle$ and ($B_2$, $*_2$) = $L_{n_1}(m_1) \in L_{n_1}$.*

*Note:* It is very important to note that $n \neq n_1$ for if $n = n_1$ then the sets $B_1$ and $B_2$ will not be proper subsets of B. These new class of neutrosophic biloops of type II enjoy several properties.

So we mention some of the properties enjoyed only by this new class of neutrosophic biloops of type II. The definitions given for type I neutrosophic biloop hold good with simple modifications.

**THEOREM 4.2.10:** *The new class of neutrosophic biloops of type II are simple neutrosophic biloops.*

*Proof:* Given (B = $B_1 \cup B_2$, $*_1$, $*_2$) is a new class of neutrosophic biloop of type II i.e. $B_1 \in \langle L_n \cup I\rangle$ and $B_2 \in L_{n_1}$. ($n \neq n_1$). We know all loops in $L_{n_1}$ are simple. Also we have proved all neutrosophic loops in $\langle L_n \cup I\rangle$ are simple i.e. they have no nontrivial normal subloops. Thus B is a simple neutrosophic biloop of type II.

**THEOREM 4.2.11:** *The order of all neutrosophic biloops of type II are of even order.*

*Proof:* Given (B = $B_1 \cup B_2$, $*_1$, $*_2$) is a neutrosophic biloop of type II so $B_1 = \langle L_n(m) \cup I\rangle) = 2 (n +1)$ (where n + 1 is even) and o($L_{n_1}$) = $n_1 + 1$ and $n_1 + 1$ is even. Thus o(B) = 2(n + 1) + $n_1$ +

1 = even number + even number = a even number. Hence the claim.

**THEOREM 4.2.12:** *Let $(B = B_1 \cup B_2, *_1, *_2)$ be a new class of neutrosophic biloops of type II. If $B_1 = \langle L_n(2) \cup I\rangle$ and $B_2 = L_{n_1}(2)$ then B is a right alternative neutrosophic biloop of type II $(n \neq n_1)$.*

*Proof:* Clearly $B_1 = \langle L_n(2) \cup I\rangle$ is a right alternative neutrosophic loop and $B_2 = L_{n_1}(2)$ is a right alternative loop so B is a neutrosophic biloop of type II $(n \neq n_1)$.

**THEOREM 4.2.13:** *Let $(B = B_1 \cup B_2, *_1, *_2)$ where $B_1 = \langle L_n(n - 1) \cup I\rangle$ and $B_2 = \langle L_n(n - 1)\rangle$ are left alternative neutrosophic loop and left alternative loop. Then B is a left alternative neutrosophic biloop of type II.*

*Proof:* Clearly $B = \langle L_n(n - 1) \cup I\rangle \cup L_{n_1}(n_1 - 1), *_1, *_2)$ is a left alternative neutrosophic biloop of type II by the fact both $\langle L_n(n - 1) \cup I\rangle$ and $L_{n_1}(n_1 - 1)$ satisfy left alternative identity $(n \neq n_1)$.

Several other properties enjoyed by these classes of neutrosophic biloops can be derived. This task is left for the reader.

## 4.3 Neutrosophic N-loop

Now we proceed onto define the notion of N-loop where N > 2 when N = 2 we get the biloop. In this section we introduce several new definitions and new properties, enjoyed by the neutrosophic N-loop.

**DEFINITION 4.3.1:** *Let $S(B) = \{S(B_1) \cup \ldots \cup S(B_N), *_1, \ldots, *_N\}$ be a non empty neutrosophic set with N binary operations. S(B) is a neutrosophic N-loop if $S(B) = S(B_1) \cup \ldots \cup S(B_N)$, $S(B_i)$ are* *proper subsets of S(B) for $1 \leq i \leq N$) and some of $S(B_i)$ are neutrosophic loops and some of the $S(B_j)$ are groups.*

*The order of the neutrosophic N-loop is the number of elements in S(B). If the number of elements in S(B) is in finite then we say S(B) is an infinite neutrosophic N-loop of infinite order. Thus even if one of the $S(B_i)$ is infinite we see the neutrosophic N-loop S(B) is infinite.*

***Example 4.3.1:*** Let $S(B) = \{S(B_1) \cup S(B_2) \cup S(B_3) \cup S(B_4), *_1, *_2, *_3, *_4\}$ be a neutrosophic 4-loop, where $S(B_1) = \{\langle L_5(3) \cup I\rangle\}$, $S(B_2) = \{S_4\}$, $S(B_3) = \{\langle L_{17}(8) \cup I\rangle\}$ and $S(B_4) = \{Q \setminus \{0\}$; group under multiplication$\}$. Clearly S(B) is a neutrosophic 4 loop of infinite order since $o(S(B_4)) = \infty$.

***Example 4.3.2:*** Let $S(B) = \{S(B_1) \cup S(B_2) \cup S(B_3), *_1, *_2, *_3\}$ where $S(B_1) = \{\langle L_5(3) \cup I\rangle\}$, $S(B_2) = \{g \mid g^{12} = 1\}$ and $S(B_3) = S_3$. S(B) is a neutrosophic 3-loop and $o(S(B)) = 30$.

Now we just define the substructure of the neutrosophic N-loop.

**DEFINITION 4.3.2:** *Let $S(B) = \{S(B_1) \cup S(B_2) \cup \ldots \cup S(B_N), *_1, \ldots, *_N\}$ be a neutrosophic N-loop. A proper subset $(P, *_1, \ldots, *_N)$ of S(B) is said to be a neutrosophic sub N loop of S (B) if P itself is a neutrosophic N-loop under the operations of S (B).*

*Note:* It is interesting and important to note that $(P, *_i)$ need not have any structure. This is evident from the following example.

***Example 4.3.3:*** Let $S(B) = \{S(B_1 \cup S(B_2) \cup S(B_3), *_1, *_2, *_3\}$ where $S(B_1) = \{\langle L_5(3) \cup I\rangle\}$, $S(B_2) = \langle g \mid g^{12} = 1\rangle$ and $S(B_3) = S_3$, is a neutrosophic 3-loop. Take

$$P = \{e, eI, 2, 2I, 1, g^6, \begin{pmatrix} 1 & 2 & 3 \\ 1 & 2 & 3 \end{pmatrix}, \begin{pmatrix} 1 & 2 & 3 \\ 2 & 1 & 3 \end{pmatrix}\}.$$

P is a neutrosophic sub 3-loop but $(P, *_1)$ or $(P, *_2)$ or $(P, *_3)$ do not have any algebraic structure. They are not even closed under these binary operations. But $P = (P_1 \cup P_2 \cup P_3, *_1, *_2, *_3)$ where $P_1 = \{e, eI, 2, 2I\}$, $P_2 = \{1, g^6\}$ and

$$P_3 = \{\begin{pmatrix} 1 & 2 & 3 \\ 1 & 2 & 3 \end{pmatrix}, \begin{pmatrix} 1 & 2 & 3 \\ 2 & 1 & 3 \end{pmatrix}\}$$

is a neutrosophic sub 3-loop of S(B).

In view of this we give the following theorem.

**THEOREM 4.3.1:** *Let $S(B) = \{S(B_1) \cup S(B_2) \cup \ldots \cup S(B_N), *_1, *_2, *_3, \ldots, *_N\}$ be a neutrosophic N-loop. Let $(P, *_1, \ldots, *_N)$ be a proper subset of S(B), P is a neutrosophic sub N-loop of S(B) if and only if $P_i = P \cap S(B_i)$, i =1, 2, …, N and $P = P_1 \cup P_2 \cup \ldots \cup P_N$ and $(P_i, *_i)$ is a substructure of $(S(B_i), *_i)$ (i = 1, 2, …, N).*

*Proof:* If $(P_i, *_i)$ is a substructure of $(S(B_i), *_i)$, clearly $(P = P_1 \cup P_2 \cup \ldots \cup P_N, *_1, \ldots, *_N)$ is a neutrosophic sub N-loop of S(B). Conversely if given $(P, *_1, \ldots, *_N)$ is a sub-N-loop of S(B) then take $P_i = P \cap S(B_i)$, we have $P = (P_1 \cup P_2 \cup \ldots \cup P_N, *_1, \ldots, *_N)$ is given to be a neutrosophic sub N-loop, so each of $(P_i, *_i)$ is either a neutrosophic loop or a group. Hence the claim.

Next we can define neutrosophic N-loops of level II.

**DEFINITION 4.3.3:** *Let $S(L) = \{L_1 \cup L_2 \cup \ldots \cup L_N, *_1, \ldots, *_N\}$ where*

i. *$L = L_1 \cup \ldots \cup L_N$ is such that each $L_i$ is a proper subset of L, $1 \leq i \leq N$.*
ii. *Some of $(L_i, *_i)$ are a Neutrosophic loops.*
iii. *Some of $(L_j, *_j)$ are just loops.*
iv. *Some of $(L_K, *_K)$ are groups and rest of*
v. *$(L_t, *_t)$ are neutrosophic groups.*

*Then we call $L = L_1 \cup L_2 \cup \ldots \cup L_N, *_1, \ldots, *_N\}$ to be a neutrosophic N-loop of level II.*

The interested reader can analyze the relation between the neutrosophic N-loops and neutrosophic N-loops of level II.

Now what all we define for neutrosophic N-loops will be carried out to neutrosophic N-loops of level II with appropriate modifications.

**DEFINITION 4.3.4:** *Let $(L = L_1 \cup L_2 \cup \ldots \cup L_N, *_1, *_2, \ldots, *_N\}$ be a neutrosophic N-loop of finite order. Suppose P is a proper subset of L, which is a neutrosophic sub N-loop. If o(P) / o(L) then we call P a Lagrange neutrosophic sub N-loop. If every neutrosophic sub N-loop is Lagrange then we call L to be a Lagrange neutrosophic N-loop. If L has atleast one Lagrange neutrosophic sub N-loop then we call L to be a weakly Lagrange neutrosophic N-loop. If L has no Lagrange neutrosophic sub N-loop then we call L to be a Lagrange free neutrosophic N-loop.*

Now we know as in case of other algebraic structures every Lagrange neutrosophic N-loop is weakly Lagrange neutrosophic N-loop.

Now we proceed on to give some examples of them.

***Example 4.3.4:*** Let $L = (L_1 \cup L_2 \cup L_3, *_1, *_2, *_3)$ be a neutrosophic 3 loop. $L_1 = \{\langle L_5(3) \cup I\rangle, *_1\}$, $L_2 = \langle g \mid g^4 = e\rangle$ and $L_3 = \{D_{2.6}$ i.e. a, b $\in D_{26}$ with $a^2 = b^6 = 1$ and b a b = a}. o(L) = 28. Take $\{P = P_1 \cup P_2 \cup P_3, *_1, *_2, *_3\}$ where

$$P_1 = \{e, eI, 2, 2I\}, P_2 = \{g^2, 1\} \text{ and}$$

$P_3 = \{b, b^2, b^3, b^4, b^5, b^6 = 1\}$ is a neutrosophic sub-3 loop of L, o(P) = 12, o(P) $\not|$ o(L) i.e. 12 $\not|$ 28 so P is not a Lagrange neutrosophic sub 3-loop of L.

Take $\{T = T_1 \cup T_2 \cup T_3, *_1, *_2, *_3\}$ where $T_1 = \{e, eI, 3, 3I\}$, $T_2 = \{e, g, g^2, g^3\}$ and $T_3 = \{1, b, b^2, b^3, b^4, b^5\}$ is a neutrosophic sub 3-loop. o(T) = 14, o(T) / o(L) so T is a Lagrange neutrosophic sub 3-loop of L. Thus L is only a weakly Lagrange neutrosophic 3-loop.

***Example 4.3.5:*** Let $L = (L_1 \cup L_2 \cup L_3, \cup L_4, *_1, *_2, *_3, *_4)$ be a neutrosophic 4-loop of finite order; where $L_1 = \{\langle L_5 (2), \cup I \rangle\}$, $L_2 = \{D_{2.4} \mid a^2 = b^4 = 1$, b a b = a}, $L_3 = \{S_3\}$ and $L_4 = \{g \mid g^{10} =$

e}, o(L) = 36. Take P = {$P_1 \cup P_2 \cup P_3 \cup P_4$, $*_1$, $*_2$, $*_3$, $*_4$} where $P_1$ = {e, eI, 4, 4I}, $P_2$ = {1 a},

$$P_3 = \left\{\begin{pmatrix}1 & 2 & 3\\1 & 2 & 3\end{pmatrix}, \begin{pmatrix}1 & 2 & 3\\2 & 3 & 1\end{pmatrix}, \begin{pmatrix}1 & 2 & 3\\3 & 1 & 2\end{pmatrix}\right\} \text{ and}$$

$P_4$ = {$g^2$, $g^4$, $g^6$, $g^8$, 1}. o(P) = 14. o(P) $\not|$ o(L). So P is not a Lagrange neutrosophic sub 4-loop.

Take T = {$T_1 \cup T_2 \cup T_3 \cup T_4$ $*_1$, $*_2$, $*_3$, $*_4$}, a proper subset of L where $T_1$ = {e, eI, 3, 3I},

$$T_3 = \left\{\begin{pmatrix}1 & 2 & 3\\2 & 1 & 3\end{pmatrix}, \begin{pmatrix}1 & 2 & 3\\1 & 2 & 3\end{pmatrix}\right\},$$

$T_2$ = {1, b, $b^2$, $b^3$} and $T_4$ = {e, $g^5$}. T is a Lagrange neutrosophic sub N-loop, o(T) / o(L) i.e. 12/36.

Thus L is only a weakly Lagrange neutrosophic N-loop. The following theorem establishes that there exists a non trivial class of Lagrange free neutrosophic N-loop.

**THEOREM 4.3.2:** *Let L = {$L_1 \cup L_2 \cup \ldots \cup L_N$, $*_1$, $*_2$, …, $*_N$} to be a neutrosophic N-loop of finite order n, where n is a prime. Then L is a Lagrange free neutrosophic N-loop.*

*Proof:* Given L = {$L_1 \cup L_2 \cup \ldots \cup L_N$, $*_1$, $*_2$, …,$*_N$} be a neutrosophic N-loop of order n, n a prime. Now whatever neutrosophic sub N-loop P exist we see o(P) $\not|$ o(L) for any neutrosophic sub N-loop. Thus if o(L) is a prime; L is a Lagrange free neutrosophic N-loop.

We illustrate this by an example.

***Example 4.3.6:*** Let L = ($L_1 \cup L_2 \cup L_3$, $\cup L_4$, $*_1$, $*_2$, $*_3$, $*_4$), be a neutrosophic 4-loop where $L_1$ = {⟨$L_5$ (3), $\cup$ I⟩}, $L_2$ = {g | $g^3$ = e}, $L_3$ = {$D_{2.5}$} and $L_4$ = $A_4$, o(L) = 37 a prime. Now we see P = {$P_1 \cup P_2 \cup P_3 \cup P_4$; $*_1$, $*_2$, $*_3$, $*_4$} where $P_1$ = {e, eI, 3, 3I}, $P_2$ = $L_2$, $P_3$ = {1, a} and

$$P_4 = \left\{ \begin{pmatrix} 1 & 2 & 3 & 4 \\ 1 & 2 & 3 & 4 \end{pmatrix}, \begin{pmatrix} 1 & 2 & 3 & 4 \\ 1 & 3 & 4 & 2 \end{pmatrix}, \begin{pmatrix} 1 & 2 & 3 & 4 \\ 1 & 4 & 2 & 3 \end{pmatrix} \right\}$$

is a neutrosophic sub 4-loop of L. o(P) = 12, (12, 37) = 1. Thus even if neutrosophic sub N-loop P exists clearly o(P) $\not|$ o(L) as o(L) is a prime. Hence the claim.

Now we can proceed on to define the notion of p-Sylow neutrosophic sub N-loop.

**DEFINITION 4.3.5:** *Let $L = \{L_1 \cup L_2 \cup .... \cup L_N, *_1, *_2, ..., *_N\}$ be a neutrosophic N-loop of finite order. Let p be a prime if $p^{\alpha}$ / o(L) and $p^{\alpha+1} \not|$ o(L). If $P = \{P_1 \cup P_2 \cup ... \cup P_N, *_1, ..., *_N\}$ be a neutrosophic sub N-loop of order $p^{\alpha}$ then we call P a p-Sylow neutrosophic sub N-loop of L. If for every prime p such that $p^{\alpha}$ / o(L) and $p^{\alpha+1} \not|$ o(L) we have a p-Sylow neutrosophic sub N-loop then we call L a Sylow neutrosophic N-loop. If L has atleast one p-Sylow neutrosophic sub N-loop then we call L a weakly Sylow neutrosophic N-loop. If L has no p-Sylow neutrosophic sub N-loop then we call L to be a Sylow free neutrosophic N-loop.*

Clearly every Sylow neutrosophic N-loop is a weakly Sylow neutrosophic N-loop. It is left as an exercise for the reader to prove that all neutrosophic N-loops of order p, p a prime are Sylow free neutrosophic N-loops.

We know when G is any finite group of order n we have for every divisor t of n we have an element x in G such that $x^t = 1$.

**DEFINITION 4.3.6:** *Let $\{L = L_1 \cup L_2 \cup ... \cup L_N, *_1, *_2, ..., *_N\}$ be a neutrosophic N-loop of finite order n. An element $x \in L$ is called a Cauchy element if $x^t = 1$ and t/n. An element $y \in L$ is called a Cauchy neutrosophic element if $y^m = I$ and m/n. If t $\not|$ n, x is called as anti Cauchy element.*

*If m $\not|$ n, y is called as an anti Cauchy neutrosophic element of L. We call a neutrosophic N-loop to be Cauchy if*

*every element is either a Cauchy element of L or a Cauchy neutrosophic element of L.*

If L has no Cauchy neutrosophic element and no Cauchy element then L is a Cauchy free neutrosophic N-loop.

***Example 4.3.7:*** Let L = ($L_1 \cup L_2 \cup L_3, \cup L_4, *_1, *_2, *_3, *_4$) where $L_1 = \{\langle L_5 (3 \cup I) \rangle\}$; $L_2 = \langle g \mid g^4 = 1\rangle$, $L_3 = S_3$ and $L_4 = \{D_{2.4} = a, b, a^2 = b^4 = 1, bab = a\}$ is a neutrosophic 4-loop of order 30, 2 / 30, 4 $\not|$ 30, 5 / 30, 6 / 30 and 3 / 30. Clearly L has Cauchy neutrosophic elements and Cauchy elements of order 2. i.e. $x \in L_1$ is such that $x^2 = e$ and $xI \in L_1$ is such that $(xI)^2 = I$. $b \in L_4$ is such that $b^4 = 1$ and $\begin{pmatrix} 1 & 2 & 3 \\ 2 & 3 & 1 \end{pmatrix} \in L_3$ is such that

$$\left[\begin{pmatrix} 1 & 2 & 3 \\ 2 & 3 & 1 \end{pmatrix}\right]^3 = \begin{pmatrix} 1 & 2 & 3 \\ 1 & 2 & 3 \end{pmatrix}.$$

We have no Cauchy element or Cauchy neutrosophic element of order 5. Still we see L is not a Cauchy neutrosophic 4-loop for $b^4 = 1$ and 4 $\not|$ 30. Thus it is interesting to note that the Cauchy neutrosophic N-loops are different from usual Cauchy theorem for groups of finite order. Now we give an example of Cauchy free neutrosophic N-loop.

***Example 4.3.8:*** Let L = ($L_1 \cup L_2 \cup L_3, \cup L_4, *_1, *_2, *_3, *_4$) where $L_1 = \{\langle L_5 (3 \cup I)\rangle\}$, $L_2 = \langle g \mid g^4 = 1\rangle$, $L_3 = \{S_3\}$ and $L_4 = \{D_{24} \mid a, b \in D_{24}; a^2 = b^4 = 1, b\, a\, b = a\}$ is a neutrosophic 4-loop of finite order, o(L) = 29.

Now L has $x \in L_1$ such that $x^2 = e$ and $xI \in L_1$ with $(xI)^2 = I$. Also $g \in L_2$ is such that $g^3 = e'$, $b \in D_{24}$ is such that $b^4 = 1$. Thus no element is a Cauchy element or a Cauchy neutrosophic element of L. Thus L is a Cauchy free neutrosophic 4-loop. Now we see all properties and results defined in case of neutrosophic N-groups can be easily extended to the case of neutrosophic N-loop.

We make a fast recollection of the definitions of neutrosophic Moufang N-loops, neutrosophic Bol N-loops, neutrosophic

Bruck N-loops, neutrosophic WIP N-loops and neutrosophic alternative N-loops.

**DEFINITION 4.3.7:** *Let $L = \{L_1 \cup L_2 \cup \ldots \cup L_N, *_1, *_2, \ldots, *_N\}$ be a neutrosophic N-loop. L is said to be a Moufang neutrosophic N-loop if L satisfies any one of three identities*

i. $(x\, y)\, (z\, x) = (x\, (y\, z))\, x$
ii. $((x\, y)\, z)\, y = x\, (y\, (zy))$
iii. $x\, (y\, (xz)) = ((xy)\, x)\, z$

*for all x, y, z in L. A neutrosophic N-loop L is called a Bruck neutrosophic N-loop if $(x\, (y\, x))\, z = x\, (y\, (xz))$ and $(x\, y)^{-1} = x^{-1}\, y^{-1}$ for all $x, y, z \in L$ whenever $x^{-1}$, $y^{-1}$ exist.*

*A neutrosophic N-loop L is called a Bol loop if $((xy)z)y = x((yz)y)$ for all $z, x, y \in L$. A neutrosophic N-loop is right alternative if $(xy)y = x(yy)$ for all $x, y \in L$ and left alternative if $(xx)\, y = x\, (xy)$ for all $x,y \in L$. L is said to be a neutrosophic alternative N-loop if L is both a neutrosophic left alternative N-loop and neutrosophic right alternative N-loop.*

A neutrosophic N-loop is said to be a weak inverse property WIP-N-loop if $(xy)\, z = e$ implies $x\, (y\, z) = e$ for all $x, y, z \in L$; e is the identity element of L.

We just give examples of neutrosophic left alternative N-loops, neutrosophic right alternative N-loops and neutrosophic WIP-N-loops.

***Example 4.3.9:*** Let $L = (L_1 \cup L_2 \cup L_3, \cup L_4, *_1, *_2, *_3, *_4)$ where $L_1 = \{\langle L_5\, (4) \cup I\rangle\}$, $L_2 = \{S_3\}$, $L_3 = D_{2.6}$ and $L_4 = \{g \mid g^{12} = e^1\}$. Clearly L is a neutrosophic left alternative N-loop.

***Example 4.3.10:*** Let $L = (L_1 \cup L_2 \cup L_3, \cup L_4 \cup L_5, *_1, \ldots, *_5\}$ where $L_1 = \langle L_7\, (2) \cup I\rangle$, $L_2 = \langle L_{13}\, (2) \cup I\rangle$, $L_3 = A_4$, $L_4 = D_{2.3}$ and $L_5 = S_3$. L is a neutrosophic 5-loop which is a neutrosophic right alternative 5-loop.

The interested reader can take the task of proving L in the above example is a neutrosophic right alternative 5-loop.

***Example 4.3.11:*** Consider the neutrosophic N-loop, $L = (L_1 \cup L_2 \cup L_3, \cup L_4, *_1, *_2, *_3, *_4)$ with $L_1 = \langle L_7(3) \cup I\rangle$ where $L_7(3)$ is given by the following table,

$L_7 (3) = \{e, 1, 2, 3, 4, 5, 6, 7\}$;

| $*_1$ | e | 1 | 2 | 3 | 4 | 5 | 6 | 7 |
|---|---|---|---|---|---|---|---|---|
| e | e | 1 | 2 | 3 | 4 | 5 | 6 | 7 |
| 1 | 1 | e | 4 | 7 | 3 | 6 | 2 | 5 |
| 2 | 2 | 6 | e | 5 | 1 | 4 | 7 | 3 |
| 3 | 3 | 4 | 7 | e | 6 | 2 | 5 | 1 |
| 4 | 4 | 2 | 5 | 1 | e | 7 | 3 | 6 |
| 5 | 5 | 7 | 3 | 6 | 2 | e | 1 | 4 |
| 6 | 6 | 5 | 1 | 4 | 7 | 3 | e | 2 |
| 7 | 7 | 3 | 6 | 2 | 5 | 1 | 4 | e |

Now operation of elements in $\langle L_7(3) \cup I\rangle = \{e, 1, 2, 3, 4, 5, 6, 7, eI, 1I, 2I, 3I, 4I, 5I, 6I, 7I\}$ is as follows. Now eI. eI = eI 1I 1I = e.I, kI. kI = eI, $1 \leq k \leq 7$ for $r \neq s$, rI. sI = (r.s) I, r.s $\in L_7$ (3) and r.s. $\in \{1, 2, 3, 3, 4, 5, 6, 7\}$. So $\langle L_7 (3) \cup I\rangle$ is a WIP loop.

Take $L_2 = S_3$, $L_3 = A_4$ and $L_4 = \{D_{2.7} \mid a, b, a^2 = b^7 = 1, b\ a\ b = a\}$. Clearly L is a WIP neutrosophic 4-loop of finite order.

Now we proceed on to define neutrosophic normal sub N-loop of a neutrosophic N-loop.

**DEFINITION 4.3.8:** *Let $L = \{L_1 \cup L_2 \cup \ldots \cup L_N, *_1, *_2, \ldots, *_N\}$ be a neutrosophic N-loop. A proper subset H of L is said to be a neutrosophic normal sub N-loop of L if the following conditions are satisfied.*

*i.* *H is a neutrosophic sub N-loop of L*

*ii.* *xH = H x*

*(H x) y = H (xy)*

*y (x H) = (y x) H*

*for all x, y $\in$ L.*

*If the neutrosophic N-loop L has no trivial normal sub-N-loop, we call L to be a simple neutrosophic N-loop.*

We first give example of them.

***Example 4.3.12:*** Let L = ($L_1 \cup L_2 \cup L_3, \cup L_4, *_1, *_2, *_3, *_4$) be a neutrosophic 4-loop, where $L_1 = \langle L_7(3) \cup I \rangle$, $L_2 = S_3$, $L_3 = \{g \mid g^7 = 1\}$ and $L_4 = \{A_4\}$. L is a simple neutrosophic 4-loop, for $L_1$ is a simple neutrosophic loop as it has no normal subloops.

Several other properties enjoyed by loops and N-loops can be by appropriate methods transformed and studied for neutrosophic N-loops.

Now we proceed on to define the notion of strong neutrosophic N-loop.

**DEFINITION 4.3.9:** *Let $\{\langle L \cup I\rangle = L_1 \cup L_2 \cup \ldots \cup L_N, *_1, *_2, \ldots, *_N\}$, be a nonempty set with N-binary operations where*

i. *$\langle L \cup I\rangle = L_1 \cup L_2 \cup \ldots \cup L_N$ where each $L_i$ is a proper subset of $\langle L \cup I\rangle$; $1 \leq i \leq N$*
ii. *($L_i$, $*_i$) is a neutrosophic loop, at least for some i.*
iii. *($L_j$, $*_j$) is a neutrosophic group.*

*Then we call $\{\langle L \cup I\rangle, *_1, \ldots, *_N\}$ to be a strong neutrosophic N-loop.*

Now we give some examples to illustrate our definition.

***Example 4.3.13:*** Let $\{\langle L \cup I\rangle = L_1 \cup L_2 \cup L_3, *_1, *_2, *_3\}$ where $L_1 = \langle L_5(2) \cup I\rangle$, $L_2 = \langle L_7(2) \cup I\rangle$ and $L_3 = \{1, 2, I, 2I\}$. $\{\langle L \cup I\rangle\}$ is a strong neutrosophic 3-loop.

***Example 4.3.14:*** Let $\{\langle L \cup I\rangle = L_1 \cup L_2 \cup L_3, *_1, *_2, *_3\}$ be a neutrosophic 3-loop where $L_1 = \langle L_5(3) \cup I\rangle$, $L_2 = \langle L_7(2) \cup I\rangle$ and $L_3 = \langle Z \cup I\rangle$. $\langle L \cup I\rangle$ is a strong neutrosophic 3-loop.

We say a strong neutrosophic N-loop is commutative if each $(L_i, *_i)$ is a commutative structure. A strong neutrosophic N-loop is said to be finite if the number of elements in $\langle L \cup I \rangle = \{L_1 \cup L_2 \cup \ldots \cup L_N, *_1, \ldots, *_N\}$ has only a finite number of distinct elements, otherwise $\langle L \cup I \rangle$ is said to be an infinite neutrosophic N-loop. The above examples give both finite strong neutrosophic N-loops and an infinite strong neutrosophic N-loop.

We now give an example of a commutative strong neutrosophic N-loop.

***Example 4.3.15:*** Let $\{\langle L \cup I \rangle = L_1 \cup L_2 \cup L_3, *_1, *_2, *_3\}$ where $L_1 = \langle L_5 (3) \cup I \rangle$, $L_2 = \langle 1, 2, 3, I, 2I, 3I \rangle$ and $L_3 = \langle L_7(4) \cup I \rangle$. Clearly $\langle L \cup I \rangle$ is a commutative neutrosophic 3-loop of finite order.

Now we proceed on to define neutrosophic strong or strong neutrosophic sub N-loop.

**DEFINITION 4.3.10:** *Let $\{\langle L \cup I \rangle = \{ L_1 \cup L_2 \cup \ldots \cup L_3, *_1, \ldots, *_N\}$ be a strong neutrosophic N-loop. Let P be a proper subset of $\langle L \cup I \rangle$. We say P is a strong neutrosophic sub N-loop of $\langle L \cup I \rangle$ if $(P, *_1, \ldots, *_N)$ is a strong neutrosophic N-loop under the binary operations $*_1, \ldots, *_N$.*
*The following are as a matter of routine.*

i. *$(P, *_1, \ldots, *_N)$ is not a loop or group under any of the operations $*_1, \ldots, *_N$.*
ii. *$P = P_1 \cup P_2 \cup \ldots \cup P_N$ where each $P_i$ is a proper subset of P; $1 \leq i \leq N$ and $P_i = P \cap L_i$, $1 \leq i \leq N$. (each $P_i$ is a neutrosophic subloop or a neutrosophic subgroup)*

*Otherwise P cannot in general be a strong neutrosophic sub N-loop of $\langle L \cup I \rangle$.*

We illustrate this by the following example.

***Example 4.3.16:*** Let $\{\langle L \cup I\rangle = L_1 \cup L_2 \cup L_3, *_1, *_2, *_3\}$ be a strong neutrosophic 3 loop where $L_1 = \langle L_5(3) \cup I\rangle$, $L_2 = \langle L_{15}(7) \cup I\rangle$, and $L_3 = \{1, 2, 3, 4, I, 2I, 3I, 4I\}$. Take P = {e, eI, 3, 3I, 14, 14I, 1I, 4, 4I}. Clearly $P \subset \langle L \cup I\rangle$ but P is not closed under any of the 3 binary operations.

Take $P_1 = P \cap L_1 = \{e, eI, 3, 3I\}$, $P_2 = P \cap L_2 = \{e, eI, 14, 14I\}$ and $P_3 = P \cap L_3 = \{1, I, 4, 4I\}$. $(P = P_1 \cup P_2 \cup P_3, *_1, *_2, *_3)$ is a neutrosophic 3-subloop. $o(\langle L \cup I\rangle) = 30$, $o(P) = 12$, $o(P) \not| \; o(\langle L \cup I\rangle)$.

Thus we are interested in characterizing those neutrosophic N-subloops whose order divides the order of the strong neutrosophic N-loop.

To this end we make the following definition. Throughout this section we assume $\langle L \cup I\rangle$ to be a strong neutrosophic N-loop.

**DEFINITION 4.3.11:** *Let $\{\langle L \cup I\rangle = L_1 \cup L_2 \cup \ldots \cup L_N, *_1, *_2, \ldots, *_N\}$ be a strong neutrosophic N-loop of finite order. Suppose $(P = P_1 \cup P_2 \cup \ldots \cup P_N, *_1, \ldots, *_N)$ be any proper subset of $\langle L \cup I\rangle$ such that P is a neutrosophic strong sub N-loop of $\langle L \cup I\rangle$. If $o(P) / o(\langle L \cup I\rangle)$ then we call P a Lagrange strong neutrosophic sub N-loop of $\langle L \cup I\rangle$. If every neutrosophic strong sub N-loop P of $\langle L \cup I\rangle$ divides the $o(\langle L \cup I\rangle)$ i.e. $o(P) / o(\langle L \cup I\rangle)$ then we call $\langle L \cup I\rangle$ itself to be a Lagrange strong neutrosophic N-loop.*

*If $\langle L \cup I\rangle$ has atleast one Lagrange strong neutrosophic sub N-loop we call $\langle L \cup I\rangle$ to be a weak Lagrange strong neutrosophic N-loop. If $\langle L \cup I\rangle$ has no Lagrange strong neutrosophic sub N-loop then we call $\langle L \cup I\rangle$ to be a Lagrange free strong neutrosophic N-loop.*

*Interrelations existing between these definitions can be obtained as a matter of routine. Now we illustrate this with some examples.*

***Example 4.3.17:*** Let $\{\langle L \cup I\rangle = L_1 \cup L_2 \cup L_3, *_1, *_2, *_3\}$ be a finite neutrosophic 3-loop, where $L_1 = \langle L_5(3) \cup I\rangle$, $L_2 = \langle L_{13}(3)$ $\cup$ I⟩ and $L_3$ = {⟨$Z_6 \cup$ I⟩ group under '+' modulo 6}. Is ⟨L ∪ I⟩ a Lagrange strong neutrosophic loop?

**THEOREM 4.3.3:** *Let {⟨L ∪ I⟩ = $L_1 \cup L_2 \cup \ldots \cup L_N$, $*_1, \ldots, *_N$} be a strong neutrosophic N-loop of order n, n a prime. Then ⟨L ∪ I⟩ is a Lagrange free strong neutrosophic N-loop.*

*Proof:* Given o(⟨L ∪ I⟩) = n, n a prime. So any proper subset P which is a neutrosophic strong N-subloop is such that (o(P), o(⟨L ∪ I⟩) = 1 So no strong neutrosophic sub N-loop is Lagrange hence ⟨L ∪ I⟩ is a Lagrange free strong neutrosophic N-loop.

It is easily verified whatever be the number of neutrosophic sub N-loop in ⟨L ∪ I⟩ still none of them can be a Lagrange neutrosophic sub N-loop. Hence the claim.

On similar line's we define Sylow neutrosophic N-loop and Cauchy neutrosophic N-loop.

**DEFINITION 4.3.12:** *Let {⟨L ∪ I⟩ = $L_1 \cup L_2 \cup \ldots \cup L_N$, $*_1, \ldots, *_N$} be a strong neutrosophic N-loop of finite order. Let p be a prime such that $p^\alpha$/ o(⟨L ∪ I⟩) and $p^{\alpha+1} \nmid$ o(⟨L ∪ I⟩).*

*Suppose ⟨L ∪ I⟩ has a neutrosophic strong N-subloop P of order $p^\alpha$ then we call P to be a p-Sylow strong neutrosophic N-subloop. If for every prime p, such that $p^\alpha$/ o(⟨L ∪ I⟩) and $p^{\alpha+1} \nmid$ o(⟨L ∪ I⟩), we have a p-Sylow neutrosophic strong N-subloop then we call ⟨L ∪ I⟩ to be a Sylow strong neutrosophic N-loop.*

*If ⟨L ∪ I⟩ has atleast one p-Sylow strong neutrosophic N-subloop then we call ⟨L ∪ I⟩ to be a weak Sylow strong neutrosophic N-loop. If ⟨L ∪ I⟩ has no p-Sylow strong neutrosophic sub N-loop then we call ⟨L ∪ I⟩ to be Sylow free strong neutrosophic N-loop.*

*We have a class of Sylow free neutrosophic N-loop, for all neutrosophic N-loops of order n, n a prime are Sylow free neutrosophic N-loop.*

Thus it may so happen that we may have a strong neutrosophic N-loop. $\langle L \cup I\rangle$ is of finite order it may still happen that $\langle L \cup I\rangle$ is a Sylow free strong neutrosophic N-loop. But if for every prime p, $p^{\alpha} / o(\langle L \cup I\rangle)$ and $p^{\alpha+1} \not| \ o(\langle L \cup I\rangle)$, we may have proper strong neutrosophic N-subloop of order $p^{\alpha+t}$ $(t \geq 1)$.

To this end we define a new concept called pseudo Sylow neutrosophic N-loops.

**DEFINITION 4.3.13:** *Let $\{\langle L \cup I\rangle = L_1 \cup L_2 \cup \ldots \cup L_N, *_1, \ldots, *_N\}$ be a neutrosophic strong N-loop of finite order. Let $\langle L \cup I\rangle$ be a Sylow free strong neutrosophic N-loop. We call a proper subset $T = \{T_1 \cup T_2 \cup \ldots \cup T_N, *_1, \ldots, *_N\}$ of $\langle L \cup I\rangle$ to be a pseudo p-Sylow strong neutrosophic sub N-loop of $\langle L \cup I\rangle$ if for a prime p, $p^{\alpha} / o(\langle L \cup I\rangle)$ and $p^{\alpha+1} \not| \ o(\langle L \cup I\rangle)$ we have a strong neutrosophic sub N-loop of order $p^{\alpha+t}$ $(t \geq 1)$.*

*If for every prime p we have a pseudo p-Sylow strong neutrosophic sub N-loop then we say $\langle L \cup I\rangle$ is a pseudo Sylow strong neutrosophic N-loop. If in a Sylow free strong neutrosophic N-loop $\langle L \cup I\rangle$, we have at least one pseudo p-Sylow strong neutrosophic sub N-loop then we call $\langle L \cup I\rangle$ a weak pseudo Sylow strong neutrosophic N-loop.*

Now we proceed on to define the notion of Cauchy element and Cauchy neutrosophic element of a neutrosophic N-loop $\langle L \cup I\rangle$.

**DEFINITION 4.3.14:** *Let $\{\langle L \cup I\rangle = L_1 \cup L_2 \cup \ldots \cup L_N, *_1, \ldots, *_N\}$ be a neutrosophic strong N-loop of finite order. If for $x \in \langle L \cup I\rangle$ we have $x^t = 1$ and if $t / o(\langle L \cup I\rangle)$ then we call x to be a Cauchy element of $\langle L \cup I\rangle$. If every element x in $\langle L \cup I\rangle$ is a Cauchy element then we call $\langle L \cup I\rangle$ a Cauchy strong neutrosophic N-loop.*

*Let $y \in \langle L \cup I\rangle$ if $y^r = I$ and $r / o(\langle L \cup I\rangle)$ then we call y to be a neutrosophic Cauchy element of $\langle L \cup I\rangle$. If every neutrosophic element $y \in \langle L \cup I\rangle$ is a neutrosophic Cauchy element of $\langle L \cup I\rangle$ then we call $\langle L \cup I\rangle$ to be a neutrosophic Cauchy strong neutrosophic N-loop. If $\langle L \cup I\rangle$ is both a Cauchy*

*strong neutrosophic N-loop and a neutrosophic Cauchy neutrosophic N-loop then we call $\langle L \cup I \rangle$ to be strong Cauchy strong neutrosophic N-loop.*

The reader is requested to illustrate this by an example.

Now several interesting results can be derived, we make a mention that as classical Cauchy Theorem for finite groups is not true so only we have to define Cauchy element Cauchy neutrosophic elements of a neutrosophic N-loop $\langle L \cup I \rangle$.

Now we just define when a neutrosophic N-loop is Moufang. On similar lines as a matter of routine one can define neutrosophic Bol N-loop, neutrosophic Bruck N-loop, neutrosophic WIP- N-loop and so on.

**DEFINITION 4.3.15:** *Let $\{\langle L \cup I \rangle = L_1 \cup L_2 \cup \ldots \cup L_N, *_1, \ldots, *_N\}$ be a neutrosophic strong N-loop. We say $\{\langle L \cup I \rangle = L_1 \cup L_2 \cup \ldots \cup L_N, *_1, \ldots, *_N\}$ is a Moufang neutrosophic N-loop if every proper subset P of $\langle L \cup I \rangle$ which is a neutrosophic strong sub N-loop satisfies the Moufang identity. We do not demand the totality of the neutrosophic strong N-loop to satisfy the Moufang identity.*

*We just illustrate this by the following example.*

***Example 4.3.18:*** Let $\{\langle L \cup I \rangle = L_1 \cup L_2 \cup L_3, *_1, *_2, *_3\}$ be a neutrosophic N-loop where $L_1 = \langle L_7(3) \cup I \rangle$, $L_2 = \langle L_5(2) \cup I \rangle$ and $L_3 = \{1, 2, 3, 4, I, 2I, 3I, 4I\}$, multiplication modulo 5. Now take any neutrosophic sub N-loop P of $\langle L \cup I \rangle$, P satisfies the Moufang identity. Thus $\langle L \cup I \rangle$ is a neutrosophic Moufang N-loop. Clearly every set of elements of $\langle L \cup I \rangle$ does not satisfy the Moufang identity.

In the same way other neutrosophic N-loops which satisfy special identities are defined.

**DEFINITION 4.3.16:** *Let $\{\langle L \cup I \rangle = L_1 \cup L_2 \cup \ldots \cup L_N, *_1, \ldots, *_N\}$ be a neutrosophic N-loop. $\langle L \cup I \rangle$ is said to be neutrosophic strong Moufang N-loop if $\langle L \cup I \rangle$ satisfies any one of the following identities.*

*i.* $(x\ y)\ (z\ x) = (x\ (y\ z))\ x$

*ii.* $((x\ y\ ))\ y = x\ (y\ (z\ y))$

*iii.* $x\ (y\ (x\ z\ )) = ((\ x\ y)\ x\ )\ z$ *for all* $x, y, z \in \langle L \cup I \rangle$.

*It is easily seen all neutrosophic strong Moufang N-loops are neutrosophic Moufang N-loops however the converse is not true.*

*This is evident from the example given.*

Now on similar lines one can define neutrosophic strong WIP N-loop, neutrosophic strong left alternative N-loop, neutrosophic strong Bol N-loop and so on. We illustrate some of these by examples.

***Example 4.3.19:*** Let $\{\langle L \cup I \rangle = L_1 \cup L_2 \cup \ldots \cup L_N, *_1, \ldots, *_N\}$ be a neutrosophic N-loop. Take N = 3 and $L_1 = \{\langle L_n(m) \cup I \rangle\}$, $m^2 - m + 1 \equiv 0 \pmod n\}$ and $L_2 = \{0, 1, I, 1 + I\}$ and $L_3 = \{1, 2, I, 2I\}$. Clearly $\langle L \cup I \rangle$ is a neutrosophic strong WIP-N-loop. $\langle L \cup I \rangle$ is only a neutrosophic Bol N-loop. $\langle L \cup I \rangle$ is only a neutrosophic alternative N-loop.

***Example 4.3.20:*** Let $\langle L \cup I \rangle = \{L_1 \cup L_2 \cup L_3, *_1, *_2, *_3\}$ where $L_1 = \langle L_n(2) \cup I \rangle$, $L_2 = \{1, 2, I, 2I\}$ and $L_3 = \langle L_{21}(2) \cup I \rangle$. $\langle L \cup I \rangle$ is a neutrosophic strong right alternative N-loop. Clearly $\langle L \cup I \rangle$ is only a neutrosophic left alternative N-loop. Also $\langle L \cup I \rangle$ is a neutrosophic Moufang N-loop.

Now we proceed on to give yet another example of a neutrosophic strong left alternative N-loop.

***Example 4.3.21:*** Let $\{\langle L \cup I \rangle = L_1 \cup L_2 \cup L_3, *_1, *_2, *_3\}$ where $L_1 = \{\langle L_7 (6) \cup I \rangle\}$, $L_2 = \{\langle L_{15} (14) \cup I \rangle\}$ and $L_3 = \{0, 1, 2, 3, 4, I, 2I, 3I, 4I\}$. Clearly $\langle L \cup I \rangle$ is a strong neutrosophic 3-loop which is a neutrosophic strong left alternative 3-loop. Is this a neutrosophic right alternative 3-loop?

Now we define deficit sub N-loops.

**DEFINITION 4.3.17:** *Let $\{\langle L \cup I\rangle = L_1 \cup L_2 \cup ... \cup L_N, *_1, ..., *_N\}$, be a neutrosophic N-loop. Let P be a proper subset of $\langle L \cup I\rangle$ such that*

$$P = \left\{P_{L_1} \cup ... \cup P_{L_t}, *_{i_1}, ..., *_{i_t} \,|1\leq i_1, ..., i_t \leq N \text{ and } 1<t<N\right\}.$$

*If P is a neutrosophic t-loop ($*_{i_p} = *_j$, $1 \leq j \leq N$; $1 \leq i_p \leq t$ and $P_{i_p} = P \cap L_j$) then we call P a neutrosophic (N-t) deficit N-subloop of $\langle L \cup I\rangle$.*

***Example 4.3.22:*** Let $\langle L \cup I\rangle = \{L_1 \cup L_2 \cup L_3, \cup L_4 *_1, *_2, *_3, *_4\}$ be a neutrosophic 4-loop where $L_1 = \langle L_5(3) \cup I\rangle$, $L_2 = S_3$, $L_3 = A_4$ and $L_4 = \langle g \mid g^8 = 1\rangle$.
Let

$$P = \left\{e, eI, 2, 2I, \begin{pmatrix}1&2&3\\1&2&3\end{pmatrix}, \begin{pmatrix}1&2&3\\2&1&3\end{pmatrix}, \begin{pmatrix}1&2&3&4\\1&2&3&4\end{pmatrix}, \begin{pmatrix}1&2&3&4\\2&1&4&3\end{pmatrix}\right\}$$

P is a neutrosophic (N – 2) deficit sub 4-loop (N = 4). Take

$$T = \left\{e, eI, 4, 4I, \begin{pmatrix}1&2&3&4\\1&2&3&4\end{pmatrix}, \begin{pmatrix}1&2&3&4\\2&1&4&3\end{pmatrix}, \begin{pmatrix}1&2&3&4\\3&4&1&2\end{pmatrix}, \begin{pmatrix}1&2&3&4\\4&3&2&1\end{pmatrix}\right\}$$

T is a neutrosophic (N – 2) (N = 4) deficit sub 4-loop of $\langle L \cup I\rangle$.

Having defined neutrosophic (N – r) deficit sub N-loops we define some more interesting properties. The (N – r) deficit sub N-loops can be defined for strong neutrosophic N-loops also.

**DEFINITION 4.3.18:** *Let $\{\langle L \cup I\rangle = L_1 \cup L_2 \cup ... \cup L_N, *_1, ..., *_N\}$ be a neutrosophic N-loop of finite order. A neutrosophic (N – t) deficit sub N-loop P of $\langle L \cup I\rangle$ is said to be Lagrange neutrosophic (N – t) deficit sub N-loop if $o(P) / o(\langle L \cup I\rangle)$. If every neutrosophic (N – t) deficit sub N-loop of $\langle L \cup I\rangle$ is Lagrange neutrosophic then we call $\langle L \cup I\rangle$ a Lagrange (N – t) deficit neutrosophic N-loop. If $\langle L \cup I\rangle$ has atleast one Lagrange*

*(N – t) deficit neutrosophic sub N-loop then we call ⟨L ∪ I⟩ a weak Lagrange (N – t) deficit neutrosophic N-loop.*

**DEFINITION 4.3.19:** *Let {⟨L ∪ I⟩ = $L_1 \cup L_2 \cup \ldots \cup L_N$, $*_1, \ldots, *_N$} be a finite neutrosophic N loop. Let P = {$P_{i_1} \cup \ldots \cup P_{i_t}$ | $t < N, 1 \leq i_1, i_2, \ldots, i_t \leq N$ and $P_{i_r} = P \cap L_q$, $1 \leq q \leq N$, $1 \leq r \leq t$} be a Neutrosophic (N – t) deficit N-subloop of ⟨L ∪ I⟩.*

*If $o(L_{i_1} \cup L_{i_2} \cup \ldots \cup L_{i_t}) = m$ and if p is a prime such that $p^{\alpha} / m$ and $p^{\alpha+1} \not| \; m$ and if $o(P) = p^{\alpha}$ then we call P a p-Sylow neutrosophic (N – t) deficit N-subloop. (We do not have any relation with o⟨L ∪ I⟩. (p, o(⟨L ∪ I⟩)) = 1 or even $p^{\alpha+1} \not| \; o(⟨L \cup I⟩)$).*

*If for every prime p related with a (N – t) deficit sub N-loop we have an associated p-Sylow neutrosophic (N – t) deficit N-subloop we call the neutrosophic N-loop ⟨L ∪ I⟩ to be a Sylow (N – t) deficit neutrosophic N-loop. We define ⟨L ∪ I⟩ to be a weak Sylow (N – t) deficit neutrosophic N-loop if ⟨L ∪ I⟩ has atleast one p-Sylow (N – t) deficit neutrosophic subloop for $1 < t < N$.*

The above 2 definitions can be defined in case of strong neutrosophic N-loop with appropriate changes. We illustrate this by the following example.

***Example 4.3.23:*** Let {⟨L ∪ I⟩ = $L_1 \cup L_2 \cup L_3 \cup L_4$ $*_1, *_2, *_3, *_4$} be a neutrosophic 3-loop of finite order, where $L_1$ = ⟨$L_5(3)$ ∪ I⟩, $L_2 = A_4$ and $L_3 = G = ⟨g \mid g^{12} = 1⟩$. Clearly o(⟨L ∪ I⟩) = 36. This can have only 2 deficit neutrosophic subloop.

They can be 2-subloop of the form $L_1 \cup L_2 \cup \phi$ or $L_1 \cup \phi \cup L_3$. $o(L_1 \cup L_2 \cup \phi) = 24 = (L_1 \cup \phi \cup L_3)$; $2^3 / 24$ and $2^4 \not| \; 24$, $3 / 24$; $3^2 \not| \; 24$.

Clearly ⟨L ∪ I⟩ cannot be a Sylow (3 – 1) deficit neutrosophic 3-loop as ⟨L ∪ I⟩ can have a 3-Sylow deficit neutrosophic sub 3-loop. To find whether ⟨L ∪ I⟩ has a 2-Sylow deficit neutrosophic sub 3-loop.

Take B = {e, eI, 2, 2I, 1, $g^3$, $g^6$, $g^9$) ⊂ $L_1 \cup \phi \cup L_3$. B is a 2-Sylow deficit neutrosophic sub 3-loop.
Take

$$A = \left\{ e, eI, 2, 2I, \begin{pmatrix} 1 & 2 & 3 & 4 \\ 1 & 2 & 3 & 4 \end{pmatrix}, \begin{pmatrix} 1 & 2 & 3 & 4 \\ 2 & 1 & 4 & 3 \end{pmatrix}, \begin{pmatrix} 1 & 2 & 3 & 4 \\ 4 & 3 & 2 & 1 \end{pmatrix}, \begin{pmatrix} 1 & 2 & 3 & 4 \\ 3 & 4 & 1 & 2 \end{pmatrix} \right\}$$

A is also a 2-Sylow deficit neutrosophic sub 3-loop.

We give yet another example.

***Example 4.3.24:*** Take $\langle L \cup I \rangle = \{L_1 \cup L_2 \cup L_3, \cup L_4 *_1, *_2, *_3, *_4\}$ where $L_1 = \langle L_5 (3) \cup I \rangle$, $L_2 = G = A_4$. $L_3 = \langle L_7 (2) \cup I \rangle$ and $L_4 = \langle g \mid g^6 = 1 \rangle$. The (N-t) values are 3 or 2 i.e. the possible combinations are 8.

$$\begin{aligned}
C_1 &= L_1 \cup L_2 \cup L_3 \cup \phi \\
C_2 &= L_1 \cup L_2 \cup \phi \cup L_4 \\
C_3 &= L_1 \cup \phi \cup L_3 \cup L_4 \\
C_4 &= \phi \cup L_2 \cup L_3 \cup L_4 \\
C_5 &= L_1 \cup \phi \cup \phi \cup L_4 \\
C_6 &= L_1 \cup L_2 \cup \phi \cup \phi \\
C_7 &= \phi \cup L_2 \cup L_3 \cup \phi \text{ and} \\
C_8 &= \phi \cup \phi \cup L_3 \cup L_4
\end{aligned}$$

$$\begin{aligned}
o(C_1) &= 40 & \quad o(C_5) &= 18 \\
o(C_2) &= 30 & \quad o(C_6) &= 24 \\
o(C_3) &= 34 & \quad o(C_7) &= 28 \\
o(C_4) &= 34 & \quad o(C_8) &= 22
\end{aligned}$$

$C_8$ has 11 Sylow but has no 2-Sylow (4 – 2) deficit sub 4-loop. Take P = {I, e, eI, 2, 2I, 1, g, $g^2$, $g^3$, $g^4$, $g^5$}; o(P) = 11 and 11 / $o(C_8)$.

The reader is expected to work for the p-Sylow (4 – t) deficit neutrosophic sub 4-loops.

Also give an example of a strong neutrosophic 5 loop and find (5 – 2) deficit strong neutrosophic sub 5-loop.

**Chapter Five**

# NEUTROSOPHIC GROUPOIDS AND THEIR GENERALIZATIONS

This chapter very briefly introduces in two sections the notion of neutrosophic groupoids and neutrosophic bigroupoids and neutrosophic N-groupoids. First section defines neutrosophic groupoids and enumerates some of its properties. Section two first introduces the notion of neutrosophic bigroupoids and suggests the reader to define several notions analogous to those done in case of neutrosophic biloops and neutrosophic bisemigroups. The later part of the section introduces the notion of neutrosophic N-groupoids and enumerates the analogous definitions to be defined by the reader. As groupoids are generalizations of both semigroups on one side and loops on other side it would not be difficult to construct new definitions as both these concepts are dealt elaborately in this book.

## 5.1 Neutrosophic Groupoids

In this section we first introduce the notion of neutrosophic groupoids. We also define some of their properties. It is left as a work of researcher to develop more results and properties. Problems related with them are proposed in the final chapter of this book. New notions like centre, direct product, conjugate pair are introduced.

**DEFINITION 5.1.1:** *Let (G, *) be a groupoid. A neutrosophic groupoid is defined as a groupoid generated by $\{\langle G \cup I\rangle\}$ under the operation *.*

We give a few examples of neutrosophic groupoids.

***Example 5.1.1:*** Let G = {a, b, ∈ $Z_3$ such that a * b = a + 2b (mod 3)} be a groupoid. $\langle G \cup I\rangle$ = {0, 1, 2, I, 2I, 1 + I, 2 + I, 1 + 2I, 2 + 2I, *} is a neutrosophic groupoid. For instance for 1 + 2I and 2 + I in $\langle G \cup I\rangle$ we have (1 + 2I) * (2 + I) = 1 + 2I + 4 + 2I (mod 3) = 2 + I (mod 3).

***Example 5.1.2:*** Let ($\langle Z^+ \cup I\rangle$, *) be the set of positive integers with a binary operation * where a * b = 2a + 3b for a, b ∈ $\langle Z^+ \cup I\rangle$.

For consider

$$5 * (4 * 1) = 5 * [8 + 3] = 5 * 11 = 10 + 33 = 43.$$

Now

$$(5 * 4) * 1 = (10 + 12) * 1 = 22 * 1 = 44 + 3 = 47.$$

Clearly 5 * (4 * 1) ≠ (5 * 4) * 1, so ($Z^+$, *) is a groupoid.

Consider $\langle Z^+ \cup I, *\rangle$ is a neutrosophic groupoid. Elements in $Z^+ \cup I$ = {a + bI / a, b ∈ $Z^+$}.

Now if 2 + 5I, 7 + I ∈ $\langle Z^+ \cup I\rangle$;

$$(2 + 5I) * (7 + I) = 2 (2 + 5I) + 3 (7 + I) = 4 + 10 I + 21 + 3I = 25 + 13I.$$

In general a * b ≠ b * a for a , b ∈ $Z^+$. Now

$$(7 + I) * (2 + 5I) = 2 (7 + I) + 3 (2 + 5I) = 14 + 2I + 6 + 10I$$

$= \quad 20 + 12\,I.$

Clearly $(7 + I) * (2 + 5I) \neq (2 + 5I) * (7 + I)$.

Now we define the order of a neutrosophic groupoid and the neutrosophic subgroupoid.

**DEFINITION 5.1.2:** *Let $\{\langle G \cup I\rangle, *\}$ be a neutrosophic groupoid. The number of distinct elements in $\{\langle G \cup I\rangle, *\}$ is called the order of the neutrosophic groupoid. If $\{\langle G \cup I\rangle, *\}$ has infinite number of elements then we say the neutrosophic groupoid is infinite. If $\{\langle G \cup I\rangle, *\}$ has only a finite number of elements then we say $\{\langle G \cup I\rangle, *\}$ is a finite neutrosophic groupoid.*

The neutrosophic groupoid given in example 5.1.1 is finite where as the neutrosophic groupoid given in example 5.1.2 is infinite.

**DEFINITION 5.1.3:** *We say a neutrosophic groupoid $\{\langle G \cup I, *\rangle\}$ is commutative if $a * b = b * a$ for all $a, b \in \langle G \cup I\rangle$.*

**DEFINITION 5.1.4:** *Let $\{\langle G \cup I\rangle, *\}$ be a neutrosophic groupoid. A proper subset P of $\langle G \cup I\rangle$ is said to be a neutrosophic subgroupoid if (P, *) is a neutrosophic groupoid. (L, *) is a subgroupoid if L is a subgroupoid and has no neutrosophic elements in them.*

We illustrate this by the following example.

***Example 5.1.3:*** Let $\langle Z_{10} \cup I\rangle = \{0, 1, 2, 3, \ldots, 9, I, 2I, \ldots, 9I, 1 + I, 2 + I, \ldots, 9 + 9I\}$; define * on $\langle Z_{10} \cup I \rangle$ by $a * b = 3a + 2b \pmod{10}$ for all $a, b \in \langle Z_{10} \cup I\rangle$.

Let

$a = 2 + 5I$ and $b = 7 + 3I$.

$$\begin{aligned} a * b &= [3 (2 + 5I) + 2 (7 + 3I)] \pmod{10} \\ &= (6 + 15\,I + 14 + 6I) \pmod{10} \\ &= I. \end{aligned}$$

$$
\begin{aligned}
b * a &= [3\,(7 + 3I) + 2\,(2 + 5I)] \pmod{10} \\
&= (21 + 9I + 4 + 10I) \pmod{10} \\
&= 4\,I.
\end{aligned}
$$

Thus $\{\langle G \cup I\rangle, *\}$ is not a commutative neutrosophic groupoid but it is a finite neutrosophic groupoid.
Take

$$
\begin{aligned}
P &= \langle (0, 5) \cup I\rangle \\
&= \{5, 5I, 0\ 5 + 5I\} \\
&= \{5, 5I, 0, 5 + 5I\}
\end{aligned}
$$

is a neutrosophic subgroupoid. $L = (Z_{10}, *)$ is just a subgroupoid of $\langle Z_{10} \cup I\rangle$.

Thus in a neutrosophic groupoid we can define two substructures called subgroupoid and neutrosophic subgroupoid. Now we see in general in case of finite neutrosophic groupoids the order of neutrosophic subgroupoid or subgroupoid does not divide the order of the neutrosophic groupoid.

We just record this important fact. If $\{\langle G \cup I\rangle, *\}$ is a neutrosophic groupoid then it always has a proper subset which is a subgroupoid.

We just illustrate this by the following example.

***Example 5.1.4:*** Let $\langle Z_4 \cup I\rangle = [0, 1, 2, 3, I, 2I, 3I, 1 + 2I, 1 + I, 1 + 3I, 2 + I, 2 + 2I, 2 + 3I, 3 + I, 3 + 2I, 3 + 3I]$. $\langle Z_4 \cup I \rangle$ is a neutrosophic groupoid under the operation * where

$$a * b = 2a + b \pmod 4 \text{ i.e.}$$

if $a = 3 + 2I$, $b = 1 + 2I$.

$$
\begin{aligned}
a * b &= [2\,(3 + 2I) + (1 + 2I)] \pmod 4 \\
&= (6 + 4I + 1 + 2I) \pmod 4 \\
&= 3 + 2I.
\end{aligned}
$$

$o(\langle Z_4 \cup I \rangle) = 16$.

Let $P = \{0, 2, 2I, 2 + 2I\}$, P is a neutrosophic subgroupoid and $o(P)\ /\ o(\langle Z_4 \cup I \rangle)$. $\{0, 2, 2 + 2I\} = T$ is a neutrosophic subgroupoid $o\,(T) \nmid o\,(\langle Z_4 \cup I \rangle)$.

Thus based of these observations we make the following definitions.

**DEFINITION 5.1.5:** *Let $\{\langle G \cup I\rangle, *\}$ be a finite neutrosophic groupoid. If (P, *) is a neutrosophic subgroupoid (subgroupoid) such that $o(P) / o(\langle G \cup I\rangle)$ then we call P to be a Lagrange neutrosophic subgroupoid (subgroupoid). If every neutrosophic subgroupoid or subgroupoid is Lagrange then we call $\{\langle G \cup I\rangle, *\}$ to be a Lagrange neutrosophic groupoid.*

*If $\{\langle G \cup I\rangle, *\}$ has atleast one Lagrange subgroupoid or Lagrange neutrosophic subgroupoid then we call $\{\langle G \cup I\rangle, *\}$ to be a weakly Lagrange neutrosophic groupoid.*

If $\{\langle G \cup I\rangle, *\}$ has no Lagrange subgroupoid or Lagrange neutrosophic subgroupoid then we call $\{\langle G \cup I\rangle, *\}$ to be a Lagrange free neutrosophic groupoid.

We define the notion of Sylow and Cauchy neutrosophic groupoids.

**DEFINITION 5.1.6:** *Let $\{\langle G \cup I\rangle, *\}$ be a finite neutrosophic groupoid. An element $x \in \langle G \cup I\rangle$ with $x^n = 1$ is said to be a Cauchy element if $n / o(\langle G \cup I\rangle)$. An element $y \in \langle G \cup I\rangle$ with $y^m = I$ is said to be a neutrosophic Cauchy element if $m / o(\langle G \cup I\rangle)$.*

*In a neutrosophic groupoid if every element which is such that $x^n = 1$ is a Cauchy element and if $y^m = I$, $y \in \langle G \cup I\rangle$ is a Cauchy neutrosophic element then we call $\langle G \cup I\rangle$ to be a Cauchy Neutrosophic groupoid. If $\langle G \cup I\rangle$ has atleast some Cauchy neutrosophic element or Cauchy element then $\langle G \cup I\rangle$ is defined as a weak Cauchy Neutrosophic groupoid. If $\langle G \cup I\rangle$ has no Cauchy element or Cauchy neutrosophic element then we call $\langle G \cup I\rangle$ to be a Cauchy free neutrosophic groupoid.*

(It is important to note we may have elements x in $\langle G \cup I\rangle$ with $x^2 = x$ or $x^n = 0$ we do not list them in our definition)

Interested reader can construct examples of them for that is not a difficult task. Now we proceed onto define the notion of p-Sylow neutrosophic subgroupoid and Sylow neutrosophic groupoid.

**DEFINITION 5.1.7:** *Let {⟨G ∪ I⟩, *} be a finite neutrosophic groupoid of order n, suppose P is a prime such that $p^{\alpha}/n$ and $p^{\alpha+1} \not| n$.*

*If we have a neutrosophic subgroupoid T of ⟨G ∪ I⟩ such that T is order $p^{\alpha}$ then we call T to be a p-Sylow neutrosophic subgroupoid of {⟨G ∪ I⟩, *}. If for every prime p such that $p^{\alpha}/o$ (⟨G ∪ I⟩) and $p^{\alpha+1} \not|$ o(⟨G ∪ I⟩) we have a neutrosophic subgroupoid of order $p^{\alpha}$ then we call {⟨G ∪ I⟩, *} to be a Sylow neutrosophic groupoid.*

*If {⟨G ∪ I⟩, *} has at least one p-Sylow neutrosophic subgroupoid then we call {⟨G ∪ I⟩, *} to be a weak Sylow neutrosophic groupoid. If for any prime p such that $p^{\alpha}$ / o⟨G ∪ I⟩ and $p^{\alpha+1} \not|$ o⟨G ∪ I⟩ we don't have an associated p-Sylow neutrosophic subgroupoid of order $p^{\alpha}$ then we call {⟨G ∪ I⟩, *} to be a Sylow free neutrosophic groupoid.*

We just give some illustrative examples before we proceed on to define the notion of neutrosophic Bol, Bruck and P-groupoids.

***Example 5.1.5:*** Let ⟨G ∪ I⟩ = {0, 1, 2, 3, I, 2I, 3I, 1 + I, 2 + I, 1 + 3I, 2 + 3I, 2I + 1, 2I + 2, 3 + I, 3 + 2I, 3 + 3I}, define a binary operation * on ⟨G ∪ I⟩ by a * b = a + 2b (mod 4) for a, b ∈ ⟨G ∪ I⟩.

Clearly o(⟨G ∪ I⟩) = 16. ⟨G ∪ I⟩ has a neutrosophic subgroupoid of order 3 given by P = {0, 2, 2I}. It has also neutrosophic subgroupoid of order 4. Thus {⟨G ∪ I⟩, *} can only be a weakly Lagrange neutrosophic groupoid.

We give yet another example.

***Example 5.1.6:*** Let ⟨G ∪ I⟩ = {0, 2, 2I, 1 + I, 2 + 2I, 3 + 3I} be a neutrosophic groupoid under the binary operation * where a * b = a + 2b (mod 4). Now o(⟨G ∪ I⟩) = 6, 3 / 6 and $3^2 \not| 6$, 2 / 6 and $2^2 \not| 6$. We have P = {0, 2I} a neutrosophic subgroupoid of order 2; for the related table is

| * | 0 | 2I |
|---|---|---|
| 0 | 0 | 0 |
| 2I | 2I | 2I |

Now consider the subset T = {0, 2, 2I}. T is a neutrosophic subgroupoid of order 3 given by the following table.

| * | 0 | 2 | 2I |
|---|---|---|---|
| 0 | 0 | 0 | 0 |
| 2 | 2 | 2 | 2 |
| 2I | 2I | 2I | 2I |

*Thus the neutrosophic groupoid given in this example is a Sylow neutrosophic groupoid.*

*But this neutrosophic groupoid is not Lagrange for consider the set V = {0, 2, 2I, 2 + 2I} given by the following table.*

| * | 0 | 2 | 2I | 2+2I |
|---|---|---|---|---|
| 0 | 0 | 0 | 0 | 0 |
| 2 | 2 | 2 | 2 | 2 |
| 2I | 2I | 2I | 2I | 2I |
| 2+2I | 2+2I | 2+2I | 2+2I | 2+2I |

o(V) = 4, $4 \nmid 6$.

Now we proceed on to define a new notion called super Sylow neutrosophic groupoids.

**DEFINITION 5.1.8:** *Let $\langle G \cup I\rangle$ be a neutrosophic groupoid of finite order. Suppose $\{\langle G \cup I\rangle, *\}$ is a Sylow neutrosophic groupoid and if in addition for all relevant prime p such that $p^{\alpha}/o(\langle G \cup I\rangle)$ and $p^{\alpha+1} \nmid o(\langle G \cup I\rangle)$ we have a neutrosophic subgroupoid of order $p^{\alpha+t}$ ($t \geq 1$ and $p^{\alpha+t} < 0(\langle G \cup I\rangle)$ for atleast one 't' then we call $\langle G \cup I\rangle$ to be a super Sylow neutrosophic groupoid.*

It is interesting to note that the Sylow neutrosophic groupoid given in example 5.1.7 is a super Sylow neutrosophic groupoid. For $o(\langle G \cup I \rangle) = 6$, $3^2 = 9 > 6$ but $2^2 \nmid 6$ by $2^2 < 6$ and we have a neutrosophic subgroupoid of order 4. Hence the claim. It is still interesting to note that every super Sylow neutrosophic groupoid is a Sylow neutrosophic groupoid but a Sylow neutrosophic groupoid in general need not always be a super Sylow neutrosophic groupoid.

Now we proceed onto define the notion of Bruck, Bol, Moufang, alternative and P-groupoids.

**DEFINITION 5.1.9:** *A neutrosophic groupoid* $\{\langle G \cup I \rangle, *\}$ *is said to be a Moufang neutrosophic groupoid if it satisfies the Moufang identity* $(x * y) * (z * x) = (x * (y * z) * x$ *for all* $x, y, z, \in \langle G \cup I \rangle$.

*A neutrosophic groupoid* $\{\langle G \cup I \rangle, *\}$ *is said to be a Bol neutrosophic groupoid if* $\{\langle G \cup I \rangle, *\}$ *satisfies the Bol identity i.e.* $((x * y) * z) * y = x * ((y * z) * y)$ *for all* $x, y, z \in \langle G \cup I \rangle$.

*A neutrosophic groupoid* $\{\langle G \cup I \rangle, *\}$ *is said to be P-neutrosophic groupoid if* $(x * y) * x = (x * (y * x)$ *for all* $x, y \in \langle G \cup I \rangle$.

*A neutrosophic groupoid* $\langle G \cup I \rangle$ *is said to a right alternative neutrosophic groupoid if it satisfies the identity* $(x * y) * y = x * (y * y)$ *for* $x, y \in \langle G \cup I \rangle$. $\langle G \cup I \rangle$. *is said to be left alterative neutrosophic groupoid if* $(x * x) * y = x * (x * y)$ *for all* $x, y \in \langle G \cup I \rangle$. *A neutrosophic groupoid is alternative if it is both right and left alternative simultaneously.*

Now we proceed on to define neutrosophic left ideal of a neutrosophic groupoid $\langle G \cup I \rangle$.

**DEFINITION 5.1.10:** *Let* $\{\langle G \cup I \rangle, *\}$ *be a neutrosophic groupoid. A proper subset H of* $\langle G \cup I \rangle$ *is said to be a neutrosophic subgroupoid of* $\langle G \cup I \rangle$ *if* $(H, *)$ *itself is a neutrosophic groupoid.*

*A non empty subset P of the neutrosophic groupoid $\langle G \cup I\rangle$ is said to be a left neutrosophic ideal of the neutrosophic groupoid $\langle G \cup I\rangle$ if*

i. *P is a neutrosophic subgroupoid.*
ii. *For all $x \in \langle G \cup I\rangle$ and $a \in P$, $x * a \in P$.*

*One can similarly define right neutrosophic ideal of a neutrosophic groupoid $\langle G \cup I\rangle$. We say P is a neutrosophic ideal of the neutrosophic groupoid $\langle G \cup I\rangle$ if P is simultaneously a left and a right neutrosophic ideal of $\langle G \cup I\rangle$.*

Now we proceed on to define the notion of neutrosophic normal subgroupoid of a neutrosophic groupoid $\langle G \cup I\rangle$.

**DEFINITION 5.1.11:** *Let $\langle G \cup I\rangle$ be a neutrosophic groupoid. A neutrosophic subgroupoid V of $\langle G \cup I\rangle$ is said to be a neutrosophic normal subgroupoid of $\langle G \cup I\rangle$ if*

i. $aV = Va$
ii. $(Vx)\, y = V\,(xy)$
iii. $y\,(xV) = (yx)\,V$

*for all x, y, $a \in \langle G \cup I\rangle$.*

*A neutrosophic groupoid is said to be neutrosophic simple if it has no nontrivial neutrosophic normal subgroupoids.*

Now we define yet another new notion called neutrosophic normal groupoids.

**DEFINITION 5.1.12:** *Let $\langle G \cup I\rangle$ be a neutrosophic groupoid. Let H and K be two proper neutrosophic subgroupoids of G with $H \cap K = \phi$ we say H is neutrosophic conjugate with K if there exists a $x \in H$ such that $H = xK$ (or $Kx$) (or in the mutually exclusive sense).*

We can define direct product of neutrosophic groupoids which is also a neutrosophic groupoid.

**DEFINITION 5.1.13:** *Let $\{\langle G \cup I\rangle, *_1\}$, $\{\langle G \cup I\rangle, *_2\}$, ..., $\{\langle G \cup I\rangle, *_n\}$ be n neutrosophic groupoids $*_i$ binary operations defined on each $\langle G_i \cup I\rangle$, i = 1, 2,..., n. The direct product of $\langle G_1 \cup I\rangle$, $\langle G_2 \cup I\rangle$, ..., $\langle G_n \cup I\rangle$ denoted by $\langle G \cup I\rangle = \langle G_1 \cup I\rangle \times \langle G_2 \cup I\rangle \times \ldots \times \langle G_n \cup I\rangle = \{(g_1, g_2, \ldots, g_n) \mid g_i \in \langle G_i \cup I\rangle\ 1 \leq i \leq n\}$, component wise multiplication of $G_i$ makes $\langle G \cup I\rangle$ a neutrosophic groupoid.*

*For if $g = (g_i, \ldots, g_n)$ and $h = (h_1, h_2, \ldots, h_n)$ in $\langle G \cup I\rangle$ then $g * h = (g_1 *_1 h_1, g_2 *_2 h_2, \ldots, g_n *_n h_n)$. Clearly $g * h \in \{\langle G \cup I\rangle, *\}$. Thus $\{\langle G \cup I\rangle, *\}$ is a neutrosophic groupoid.*

The notion of direct product helps in finding neutrosophic subgroupoids, neutrosophic normal subgroupoids, neutrosophic ideals and also helps in construction of more and more neutrosophic groupoids satisfying the above conditions. In fact one can also relax in the definition of the direct product of neutrosophic groupoids we can also take some groupoids instead of neutrosophic groupoids.

***Example 5.1.7:*** Let $\{\langle G_1 \cup I\rangle, *_1\}$, $(G_2, *_2)$ and $\{\langle G_3 \cup I\rangle, *_3\}$, be any three groupoids where

$\langle G_1 \cup I\rangle$ = {0, 1, 2, I, 2I, 1 + I, 1 + 2I, I + 2, 2I + 2} with $*_1$ as $a *_1 b = 2a + 1b \pmod 3$).

$(G_2, *_2)$ = {$Z_{12}$, with $a * b = a + 5b \pmod{12}$} and

$\{\langle G_3 \cup I\rangle, *_3\}$ = {(0, 1, 2, 3, I, 2I, 3I, 1 + I, 1 + 2I, 1 + 3I, 2 + I, 2 + 2I, 2 + 3I 3 + I, 3 + 2I, 3 + 3I), $a *_3 b = 2a + b \pmod 4$}.

Let

$\{\langle G \cup I\rangle, *\} = \langle G_1 \cup I\rangle \times \{G_2\} \times \langle G_3 \cup I\rangle = \{(g_1, g_2, g_3) \mid g_i \in G_i$, or $\langle G_i \cup I\rangle, 1 \leq i \leq 3\}$.

Let $X = (x_1, x_2, x_3)$ and $Y = (y_1, y_2, y_3) \in \{\langle G \cup I\rangle, *\}$.

Define

$X * Y = (x_1, x_2, x_3) * (y_1, y_2, y_3)$
$= \{(x_1 *_1 y_1), (x_2 *_2 y_2), (x_3 *_3 y_3)\}$
$= \{2x_1 + y_1 \pmod 3, x_2 + 5y_2 \pmod{12}, 2x_3 + y_3 \pmod 4\}.$

i.e. if x = (I, 5, 1 + 3I) and y = (2 + 2I, 7, 3I)

$x * y = (I *_1, 2 + 2I, 5 *_2 7, 1 + 3I *_3 3I)$
$= [(2I + 2 + 2I) \bmod 3, (5 + 5.7 \pmod{12})$
$= (2 (1 + 3I) + 3I) \pmod 4)]$
$= (I + 2, 4, 2 + I) \in (\langle G \cup I \rangle).$

In neutrosophic groupoids we can have either left inverse or left identity, left zero divisor and so on [Likewise right inverse, right identity, right zero divisor and so on]. To define inverse right (left) we need the notion of neutrosophic groupoids with right (left) identity.

**DEFINITION 5.1.14:** *Let {⟨G ∪ I⟩, *} be a neutrosophic groupoid, we say an element e ∈ ⟨G ∪ I⟩ is a left identity if e * a = a for all a ∈ ⟨G ∪ I⟩, similarly right identity of a neutrosophic groupoid can be defined. If e happens to be simultaneously both right and left identity we say the neutrosophic groupoid has an identity.*

*Similarly we can say an element 0 ≠ a ∈ ⟨G ∪ I⟩ has a right zero divisor if a * b = 0 for some b ≠ 0 in ⟨G ∪ I⟩ and $a_1$ in ⟨G ∪ I⟩ has left zero divisor if $b_1 * a_1 = 0$ (both $a_1$ and $b_1$ are different from zero). We say ⟨G ∪ I⟩ has zero divisor if a * b = 0 and b * a = 0 for a, b ∈ ⟨G ∪ I⟩ \ {0}.*

Now we proceed on to define the notion of centre of the neutrosophic groupoid ⟨G ∪ I⟩.

**DEFINITION 5.1.15:** *Let {⟨G ∪ I⟩, *} be a neutrosophic groupoid, the neutrosophic centre of the groupoid ⟨G ∪ I⟩ is C(⟨G ∪ I⟩) = {a ∈ ⟨G ∪ I⟩} \ a * x = x * a for all x ∈ ⟨G ∪ I⟩}.*

Now we proceed on to define the notion of conjugate pair in a neutrosophic groupoid.

**DEFINITION 5.1.16:** *Let $\{\langle G \cup I\rangle, *\}$ be a neutrosophic groupoid of order n. ($n < \infty$). We say a, b $\in \langle G \cup I\rangle$ is a conjugate pair if a = b * x (or x * b for some x $\in \langle G \cup I\rangle$) and b = a * y (or y * a for some y $\in \langle G \cup I\rangle$). An element a in $\langle G \cup I\rangle$ is said to be right conjugate with b in $\langle G \cup I\rangle$ if we can find x, y $\in \langle G \cup I\rangle$ such that a * x = b and b * y = a (x * a = b and y * b = a).*

## 5.2 Neutrosophic Bigroupoids and their generalizations

In this section we proceed on to define the new notion neutrosophic bigroupoids, neutrosophic N-groupoids and analyze some of its properties. All semigroups are groupoids i.e., the class of semigroups are contained in the class of groupoids. Likewise we can say the class of bisemigroups is contained in the class bigroupoids we just define and indicate how the definitions and other results can be extended in case of bigroupoid from bisemigroups.

**DEFINITION 5.2.1:** *Let (BN(G), *, o) be a non empty set with two binary operations * and o. (BN(G), *, o) is said to be a neutrosophic bigroupoid if*

*$BN(G) = G_1 \cup G_2$ where at least one of ($G_1$, *) or ($G_2$, o) is a neutrosophic groupoid and other is just a groupoid. $G_1$ and $G_2$ are proper subsets of BN(G); i.e., $G_1 \not\subseteq G_2$.*

Now we illustrate this by an example.

***Example 5.2.1:*** Let (BN(G), *, o) be a neutrosophic bigroupoid with $BN(G) = G_1 \cup G_2$ where

$G_1$ = $\{\langle Z_{10} \cup I\rangle \mid a * b = 2a + 3b \pmod{10};\ a, b \in \langle Z_{10} \cup I\rangle\}$ and

$G_2$ = $\{Z_{13} / a * b = 3a + 10b \pmod{13};\ a, b \in Z_{13}\}$.

BN(G) is a neutrosophic bigroupoid. If both $(G_1, *)$ and $(G_2, *)$ are neutrosophic groupoids in the above definition then we call BN(G) a strong neutrosophic bigroupoid.

It is easily verified that all neutrosophic strong bigroupoids are neutrosophic bigroupoids but not conversely.

Now we proceed to give an example of a strong neutrosophic bigroupoid.

***Example 5.2.2:*** Let (BN(G), *, o) be a non-empty set such that $BN(G) = \{\langle Z \cup I\rangle \cup \langle Z_{12} \cup I\rangle = G_1 \cup G_2, *, o$ where $(\langle Z \cup I\rangle, *)$ is a neutrosophic groupoid defined by $a * b = 5a + 2b$ for all $a, b \in \langle Z \cup I\rangle\}$ and $\{\langle Z_{12} \cup I\rangle$, o is a neutrosophic groupoid given by $a \text{ o } b = 8a + 4b \pmod{12}$ for all $a, b \in \langle Z_{12} \cup I\rangle\}$. B(N(G)) is a strong neutrosophic bigroupoid.

As in case of neutrosophic bisemigroups we can define the notion of neutrosophic sub-bigroupoid and sub-bigroupoid in case of neutrosophic bigroupoids. In case of strong neutrosophic bigroupoids, we can define 3 sub-structures viz.

1. Strong neutrosophic sub-bigroupoids.
2. Neutrosophic sub-bigroupoids and
3. Sub-bigroupoids

This is a simple exercise left for the reader.

*Now we can define as in the case of neutrosophic bisemigroups biideals in neutrosophic bigroupoids. Here it has become important to mention that in case of strong neutrosophic bigroupoid we can have only strong neutrosophic biideal. Also, for any neutrosophic bigroupoid we can have only neutrosophic biideal. Thus a strong neutrosophic bigroupoid cannot have neutrosophic biideal or just biideal. Likewise a neutrosophic bigroupoid cannot have strong neutrosophic biideal or a biideal. This is the marked difference between the biideals and subbigroupoids in neutrosophic bigroupoids.*

*We just illustrate these situations by the following example.*

***Example 5.2.3:*** Let BN(G) = $\{G_1 \cup G_2, *_1, *_2\}$ where

$G_1$ = $\{\langle Z_{12} \cup I\rangle$ a neutrosophic groupoid defined by $*_1$ as a $*_1$ b = 8a + 4b (mod 12) for all a, b $\in \langle Z_{12} \cup I\rangle$ and

$G_2$ = $\{(a, b)/ a, b \in \langle Z_4 \cup I\rangle\}$. Component wise multiplication i.e., (a, b) $*_2$ (a', b') = (3a + a'(mod 4), (3b + b') mod 4) $\{(a *_2 a', b *_2 b')$ where a $*_2$ a' = 3a + a' (mod 4).

Let

$P_1$ = $\{0, 6, 6I\} \subset G_1$ is a neutrosophic ideal of $G_1$.

$P_2$ = $\{\langle 2, 2I\rangle\} \subset G_2$ be the neutrosophic ideal generated by $\langle(2, 2I)\rangle$.

Clearly $P_1 \cup P_2$ is a strong neutrosophic biideal of BN(G).

The notion of strong neutrosophic maximal bi-ideal, strong neutrosophic minimal biideal, strong neutrosophic quasi maximal biideal, strong neutrosophic quasi minimal biideal can be defined in case of strong neutrosophic bigroupoids. Likewise, neutrosophic maximal biideal, neutrosophic minimal biideal, neutrosophic quasi maximal biideal and neutrosophic quasi minimal biideal can be defined in case of neutrosophic bigroupoids.

Now we proceed on to define the notion of neutrosophic N-groupoids.

**DEFINITION 5.2.2:** *Let N(G) = $\{G_1 \cup G_2 \cup \ldots \cup G_N, *_1, \ldots, *_N\}$ be a non-empty set with N-binary operations, N(G) is called a neutrosophic N-groupoid if some of the $G_i$'s are neutrosophic groupoids and some of them are neutrosophic semigroups and N(G) = $G_1 \cup G_2 \cup \ldots \cup G_N$ is the union of the proper subsets of N(G).*

*It is important to note that a $G_i$ is either a neutrosophic groupoid or a neutrosophic semigroup. We call a neutrosophic N-groupoid to be a weak neutrosophic N-groupoid if in the* *union* $N(G) = G_1 \cup G_2 \cup \ldots \cup G_N$ *some of the* $G_i$*'s are neutrosophic groupoids, some of the* $G_j$*'s are neutrosophic semigroups and some of the* $G_k$*'s are groupoids or semigroups 'or' not used in the mutually exclusive sense. The order of the neutrosophic N-groupoids are defined as that of N-groupoids. Further we call a neutrosophic N-groupoid to be commutative if each* $(G_i, *_i)$ *is commutative for* $i = 1, 2, \ldots, N$.

Let $\langle G \cup I\rangle = \{G_1 \cup G_2 \cup G_3 \cup \ldots \cup G_N, *_1, \ldots, *_N\}$ be a neutrosophic N-groupoid a proper subset P of $\langle G \cup I\rangle$ is called a neutrosophic sub-N-groupoid if P itself under the N-operations of $\langle G \cup I\rangle$ is a neutrosophic N-groupoid. We can as in the case of N-groupoids define other sub N-structures like Lagrange neutrosophic sub N-groupoids, p-Sylow neutrosophic sub N-groupoid, neutrosophic normal sub N-groupoid, and conditions when are two neutrosophic sub N-groupoids N-conjugate and so on. We now just define the new notion of N-quasi loop.

**DEFINITION 5.2.3:** *Let* $\langle G \cup I\rangle = \{G_1 \cup G_2 \cup G_3 \cup \ldots \cup G_N, *_1, \ldots, *_N\}$ *be a non-empty set with N-binary operations. We call* $\langle G \cup I\rangle$ *a neutrosophic N-quasi loop if at least one of the* $(G_j, *_j)$ *are neutrosophic loops. So a neutrosophic sub-N-quasiloop will demand one of the subsets* $P_i$ *contained* $G_i$ *to be a neutrosophic subloop.*

All properties pertaining to the sub-structures can be derived as in the case of neutrosophic N-groupoids. We define N-quasi semigroups.

**DEFINITION 5.2.4:** *Let* $\langle G \cup I\rangle = \{G_1 \cup G_2 \cup G_3 \cup \ldots \cup G_N, *_1, \ldots, *_N\}$ *be a non-empty set with N-binary operations with each* $G_i$ *a proper subset of* $\langle G \cup I\rangle$, $i = 1, 2, \ldots, N$. $\langle G \cup I\rangle$ *is a neutrosophic N-quasi semigroup if some of the* $(G_i, *_i)$ *are neutrosophic loops and the rest are neutrosophic semigroups.*

*Note:* We do not have in the collection any neutrosophic groupoid or groupoid. Likewise we define neutrosophic N-quasi groupoid as a non-empty set with N-binary operations $*_1, \ldots, *_N$ on $G_1, G_2, \ldots, G_N$ where $\langle G \cup I\rangle = \{G_1 \cup G_2 \cup G_3 \cup \ldots \cup G_N,$

$*_1, \ldots, *_N\}$ where $(G_i, *_i)$ are either neutrosophic groups or neutrosophic groupoid or used in the mutually exclusive sense.

Further each $G_i$ is a proper subset of $\langle G \cup I\rangle$, i = 1, 2, …, N. Now all notions defined for neutrosophic N-loops can be easily extended to the class of neutrosophic groupoids.

Interested reader can work in this direction. To help the reader in chapter 7 if this book 25 problems are suggested only on neutrosophic groupoids and neutrosophic N-groupoids. All identities studied in case of neutrosophic loops and neutrosophic N-loops can also be defined and analyzed in case of neutrosophic groupoids.

## Chapter Six

# MIXED NEUTROSOPHIC STRUCTURES

Here in this chapter we define the notion of mixed neutrosophic structures and their dual and illustrate it with examples. We give only hints for definition for it can be done as a matter of routine.

Further these mixed structures will have applications when the domain values are taken from different algebraic structures. We give examples of them so that it makes the reader understand the concept easily. This chapter has only one section.

**DEFINITION 6.1.1:** *Let* $\{\langle M \cup I\rangle = M_1 \cup M_2 \cup \ldots \cup M_N, *_1, \ldots, *_N\}$, $(N \geq 5)$ *we call* $\langle M \cup I\rangle$ *a mixed neutrosophic N-structure if*

i. $\langle M \cup I\rangle = M_1 \cup M_2 \cup \ldots \cup M_N$, *each* $M_i$ *is a proper subset of* $\langle M \cup I\rangle$.
ii. *Some of* $(M_i, *_i)$ *are neutrosophic groups.*
iii. *Some of* $(M_j, *_j)$ *are neutrosophic loops.*
iv. *Some of* $(M_k, *_k)$ *are neutrosophic groupoids.*
v. *Some of* $(M_r, *_r)$ *are neutrosophic semigroups.*
vi. *Rest of* $(M_t, *_t)$ *can be loops or groups or semigroups or groupoids. ('or' not used in the mutually exclusive sense*

*(From this the assumption* $N \geq 5$ *is clear).*

***Example 6.1.1:*** Let $\langle M \cup I\rangle = \{M_1 \cup M_2 \cup M_3 \cup M_4 \cup M_5 \cup M_6, *_1, *_2, *_3, *_4, *_5, *_6\}$ where

$(M_1, *_1)$ = $\langle L_5 (3) \cup I\rangle$ a neutrosophic loop.
$(M_2, *_2)$ = $\langle Z \cup I\rangle$ under addition, a neutrosophic group.
$(M_3, *_3)$ = $\langle Z_{10} \cup I\rangle$, $Z_{10}$ semigroup under multiplication modulo 10.
$(M_4, *_4)$ = $\langle Z_6 \cup I\rangle = \{0, 1, 2, 3, 4, 5, I, 2I, 3I, 4I, 5I / a *_4 b = (2a + 3b) \pmod 6\}$.
$M_5$ = $S_4$, a group.
$M_6$ = $Z_8$, semigroup under multiplication modulo 8.

$\langle M \cup I\rangle$ is a mixed neutrosophic 6-structure.

Now we define the mixed dual neutrosophic N-structure.

**DEFINITION 6.1.2:** *$\{\langle D \cup I\rangle = D_1 \cup D_2 \cup \ldots \cup D_N, *_1, *_2, \ldots, *_N\}$, $N \geq 5$ be a non empty set on which is defined N-binary operations. We say $\langle D \cup I\rangle$ is a mixed dual neutrosophic N-structure if the following conditions are satisfied.*

- i. *$\langle D \cup I\rangle = D_1 \cup D_2 \cup \ldots \cup D_N$ where each $D_i$ is a proper subset of $\langle D \cup I\rangle$*
- ii. *For some i, $(D_i, *_i)$ are groups*
- iii. *For some j, $(D_j, *_j)$ are loops*
- iv. *For some k, $(D_k, *_k)$ are semigroups*
- v. *For some t, $(D_t, *_t)$ are groupoids.*
- vi. *The rest of $(D_m, *_m)$ are neutrosophic groupoids or neutrosophic groups or neutrosophic loops or neutrosophic semigroup 'or' not used in the mutually exclusive sense.*

We illustrate this by the following example.

***Example 6.1.2:*** Let $\{\langle D \cup I\rangle = D_1 \cup D_2 \cup D_3 \cup D_4 \cup D_5, *_1, *_2, *_3, *_4, *_5\}$ where $D_1 = L_5 (2)$, $D_2 = A_4$, $D_3 = S(3)$, $D_4 = \{Z_{10}$ such that $a *_4 b = 2a + 3b \pmod{10}\}$ and $D_5 = \langle L_7 (3) \cup I\rangle$. $\langle D \cup I\rangle$ is a mixed dual neutrosophic 5-structure.

Now we proceed on to define weak mixed neutrosophic N-structure.

**DEFINITION 6.1.3:** *Let $\langle W \cup I\rangle = \{W_1 \cup W_2 \cup \ldots \cup W_N, *_1, *_2, \ldots, *_N\}$ be a non empty set with N-binary operations $*_1, \ldots, *_N$. $\langle W \cup I\rangle$ is said to be a weak mixed neutrosophic structure, if the following conditions are true.*

i. *$\langle W \cup I\rangle = W_1 \cup W_2 \cup \ldots \cup W_N$ is such that each $W_i$ is a proper subset of $\langle W \cup I\rangle$.*
ii. *Some of $(W_i, *_i)$ are neutrosophic groups or neutrosophic loops.*
iii. *Some of $(W_j, *_j)$ are neutrosophic groupoids or neutrosophic semigroups*
iv. *Rest of $(W_k, *_k)$ are groups or loops or groupoids or semigroups. i.e. In the collection $\{W_i, *_i\}$ all the 4 algebraic neutrosophic structures may not be present.*

*At most 3 algebraic neutrosophic structures are present and atleast 2 algebraic neutrosophic structures are present. Rest being non neutrosophic algebraic structures.*

*We just illustrate this by the following example.*

***Example 6.1.3:*** Let $\langle W \cup I\rangle = \{W_1 \cup W_2 \cup W_3 \cup W_4 \cup W_5, *_1, *_2, *_3, *_4, *_5\}$ where

$W_1 = \langle L_5(3) \cup I\rangle$,
$W_2 = \{\langle Z_{12} \cup I\rangle$ semigroup under multiplication modulo 12},
$W_3 = S_3$,
$W_4 = S(5)$ and
$W_5 = \{Z_6 \mid a *_5 b = 2a + 4b \pmod 6\}$.

$\langle W \cup I\rangle$ is a weakly mixed neutrosophic 5-structure.

We can define dual of weakly mixed neutrosophic N-structure.

**DEFINITION 6.1.4:** *Let $\{\langle V \cup I\rangle = V_1 \cup V_2 \cup \ldots \cup V_N, *_1, \ldots, *_N\}$ be a non empty set with N-binary operations. We say $\langle V \cup I\rangle$ is a weak mixed dual neutrosophic N-structure if the following conditions are true.*

i. *$\langle V \cup I\rangle = V_1 \cup V_2 \cup \ldots \cup V_N$ is such that each $V_i$ is a proper subset of $\langle V \cup I\rangle$*
ii. *Some of $(V_i, *_i)$ are loops or groups*
iii. *Some of $(V_j, *_j)$ are groupoids or semigroups*
iv. *Rest of the $(V_k, *_k)$ are neutrosophic loops or neutrosophic groups or neutrosophic groupoids or neutrosophic semigroups.*

***Example 6.1.4:*** Let $\{\langle V \cup I\rangle = V_1 \cup V_2 \cup V_3 \cup V_4, *_1, *_2, *_3, *_4\}$ where $V_1 = L_7(3)$, $V_2 = S_3$. $V_3 = S(5)$ and $V_4 = \{\langle L_{15}(8) \cup I\rangle\}$. Clearly $\langle V \cup I\rangle$ is a weakly mixed neutrosophic 4-algebraic structure.

We define order of the mixed neutrosophic N-structure $\langle M \cup I\rangle$ as the number of distinct elements in $\langle M \cup I\rangle$.

Now we define sub N-structure for one class clearly it can be done to all other structures with appropriate modifications.

**DEFINITION 6.1.5:** *Let $\{\langle M \cup I\rangle = M_1 \cup M_2 \cup \ldots \cup M_N, *_1, \ldots, *_N\}$ where $\langle M \cup I\rangle$ is a mixed neutrosophic N-algebraic structure. We say a proper subset $\{\langle P \cup I\rangle = P_1 \cup P_2 \cup \ldots \cup P_N, *_1, \ldots, *_N\}$ is a mixed neutrosophic sub N-structure if $\langle P \cup I\rangle$ under the operations of $\langle M \cup I\rangle$ is a mixed neutrosophic N-algebraic structure.*

We illustrate them by the following examples.

***Example 6.1.5:*** Let $\{\langle M \cup I\rangle = M_1 \cup M_2 \cup M_3 \cup M_4 \cup M_5 \cup M_6, *_1, *_2, *_3, *_4, *_5, *_6\}$ be a mixed neutrosophic 6-structure, where

$M_1$ = $\langle L_5(3) \cup I\rangle$,
$M_2$ = {I, 2I, 1, 2 multiplication modulo 3},
$M_3$ = $\{\langle Z_6 \cup I\rangle$; neutrosophic semigroup under multiplication modulo 6},
$M_4$ = {0, 1, 2, 3, I, 2I, 3I, a * b (2a + b) mod 4},
$M_5$ = $S_3$ and

$M_6$ = {$Z_{10}$, semigroup under multiplication modulo 10}.

Take $\langle W \cup I\rangle$ = {$W_1 \cup W_2 \cup W_3 \cup W_4 \cup W_5 \cup W_6$, $*_1$, $*_2$, $*_3$, $*_4$, $*_5$, $*_6$} where

$W_1$ = {eI, 2I, e, 2},
$W_2$ = {I, 1},
$W_3$ = {0, 3, 3I},
$W_4$ = {0, 2, 2I},
$W_5$ = {$A_3$} and
$W_6$ = {0, 2, 4, 6, 8}.

Clearly W is a mixed neutrosophic sub N-structure.

It is important to note that a mixed neutrosophic N-structure can have weak mixed neutrosophic sub N-structure. But a weak mixed neutrosophic sub N-structure cannot in general have a mixed neutrosophic sub N structure.

***Example 6.1.6:*** Let {$\langle W \cup I\rangle = W_1 \cup W_2 \cup W_3 \cup W_4 \cup W_5 \cup W_6$, $*_1$, …, $*_N$} be a mixed neutrosophic 6-structure, where

$W_1$ = $\langle L_5(3) \cup I\rangle$
$W_2$ = {1, 2, I, 2I},
$W_3$ = {$\langle Z_6 \cup I\rangle$; $Z_6 \cup I$ semigroup under multiplication modulo 6},
$W_4$ = {($\langle Z_8 \cup I\rangle$, a $*_4$ b = 2 a + 6b (mod 8)}
$W_5$ = {g | $g^4 = 1$} and
$W_6$ = {$S_3$}.

$\langle W \cup I\rangle$ is a mixed neutrosophic 6-structure.

Take {$\langle T \cup I\rangle = T_1 \cup T_2 \cup T_3 \cup T_4 \cup T_5 \cup T_6$} where $T_1$ = {e, eI, 3, 3I}, $T_2$ = {1, 2}, $T_3$ = {0, 3, 3I, I}, $T_4$ = {$Z_8$}, $T_5$ = {1, $g^2$} and $T_6$ = {$A_3$}, $\langle T \cup I\rangle$ is a weak mixed neutrosophic sub 6-structure which is not a mixed neutrosophic substructure.

Now we proceed on to define weak mixed deficit neutrosophic sub N-structures.

**DEFINITION 6.1.6:** *Let* $\langle W \cup I\rangle = \{W_1 \cup W_2 \cup \ldots \cup W_N, *_1, *_2, \ldots, *_N\}$ *be a mixed neutrosophic N-structure. We call a finite non empty subset P of* $\langle W\cup I\rangle$, *to be a weak mixed deficit neutrosophic sub N-structure if* $P = \{P_1 \cup P_2 \cup \ldots \cup P_t, *_1, \ldots, *_t\}$, $1 < t < N$ *with* $P_i = P \cap L_k$, $1 \le i \le t$ *and* $1 \le k \le N$ *and some* $P_i$*'s are neutrosophic groups or neutrosophic loops some of the* $P_j$*'s are neutrosophic groupoids or neutrosophic semigroups and rest of the* $P_k$*'s are groups or loops or groupoids or semigroups.*

***Example 6.1.7:*** Let $\langle M \cup I\rangle = \{M_1 \cup M_2 \cup \ldots\cup M_6, *_1,\ldots, *_6\}$ where

$M_1$ = $\{\langle L_5(3) \cup I\rangle\}$,
$M_2$ = $\{\langle Z_6 \cup I\rangle$, semigroup under multiplication modulo 6},
$M_3$ = {1, 2, 3, 4, I, 2I, 3I, 4I, neutrosophic group under multiplication modulo 5},
$M_4$ = {0, 1, 2, 3, I, 2I, 3I neutrosophic groupoid with binary operation $*_4$ such that a $*_4$ b (3a + 2b) (mod 4)},
$M_5$ = $Z_{12}$, group under '+' and
$M_6$ = $\{Z_4$, semigroup under multiplication modulo 6}.

Take $P = P_1 \cup P_2 \cup P_3 \cup P_4 \cup P_5$ where $P_1$ = {e, eI, 2, 2I}, $P_2$ = {3, 3I}, $P_3$ = {1, 2, 3, 4}, $P_4$ = {0, 2, 2I}, $P_5 = \{0, 2\} \subset M_{6.1.}$ Clearly P is a weak deficit mixed neutrosophic sub 6-stucture of $M \cup I$.

Now we just hint how to define Lagrange substructures.

**DEFINITION 6.1.7:** *Let* $\langle M \cup I\rangle = \{M_1 \cup M_2 \cup \ldots \cup M_N, *_1, \ldots, *_N\}$ *be a mixed neutrosophic N-structure of finite order. A proper subset P of* $\langle M \cup I\rangle$ *which is neutrosophic sub N-structure is said to be Lagrange if* $o(P) / o(\langle M \cup I\rangle)$. *If every mixed neutrosophic sub N-structure is Lagrange then we call* $\langle M \cup I\rangle$ *to be a Lagrange mixed neutrosophic N-structure.*

If $\langle M \cup I\rangle$ has no Lagrange mixed neutrosophic structure then we call $\langle M \cup I\rangle$ to be a free Lagrange neutrosophic N-structure.

Now on similar lines we define Lagrange weak deficit mixed neutrosophic sub N structure and Lagrange weak deficit mixed neutrosophic N-structure.

**DEFINITION 6.1.8:** *Let* $\{\langle M \cup I\rangle = M_1 \cup M_2 \cup \ldots \cup M_N, *_1, \ldots, *_N\}$ *be a mixed neutrosophic N-structure of finite order. We call a proper subset* $P = \{P_1 \cup P_2 \cup \ldots \cup P_t, 1 < t < N\}$ *which is weak mixed neutrosophic sub N-structure to be Lagrange if* $o(P) / o\langle M \cup I\rangle$.

*If every proper subset of P which is a weak mixed neutrosophic sub N-structure is Lagrange then we call* $M \cup I$ *to be a Lagrange mixed weak neutrosophic N-structure. If* $\langle M \cup I\rangle$ *has atleast one Lagrange mixed weak neutrosophic sub N-structure then we call* $\langle M \cup I\rangle$ *to be a weak Lagrange mixed weak neutrosophic N-structure. If* $\langle M \cup I\rangle$ *has no Lagrange mixed weak neutrosophic sub N-structure then we call* $\langle M \cup I\rangle$ *to be a Lagrange free mixed weak neutrosophic N-structure.*

Now on similar lines we define the notion of Lagrange, weak Lagrange and Lagrange free in case of weak mixed deficit neutrosophic N-structure. Thus this work is left as an exercise for the reader! One can easily construct examples of these. Now we define Sylow structure for the mixed neutrosophic N-structure as follows.

**DEFINITION 6.1.9:** *Let* $\{\langle M \cup I\rangle = M_1 \cup M_2 \cup \ldots \cup M_N, *_1, \ldots, *_N\}$ *be a mixed neutrosophic N-structure of finite order. Let p be a prime such that* $p^\alpha / o(\langle M \cup I\rangle)$ *and* $p^{\alpha+1} \nmid o(\langle M \cup I\rangle)$. *If* $\langle M \cup I\rangle$ *has a subset P which is a mixed neutrosophic substructure of order* $p^\alpha$ *then we call P a p-Sylow mixed neutrosophic sub N-structure.*

*If for every prime p, with* $p^\alpha / o(\langle M \cup I\rangle)$ *and* $p^{\alpha+1} \nmid o(\langle M \cup I\rangle)$ *we have a p-Sylow mixed neutrosophic sub N-structure then we call* $\langle M \cup I\rangle$ *to be Sylow mixed neutrosophic*

*N-structure. If $\langle M \cup I\rangle$ has atleast one p-Sylow mixed neutrosophic sub N-structure then we call $\langle M \cup I\rangle$ to be a weak Sylow mixed neutrosophic N-structure. If $\langle M \cup I\rangle$ has no p-Sylow mixed neutrosophic sub N-structure then we call $\langle M \cup I\rangle$ to be a Sylow free mixed neutrosophic N-structure.*

On similar lines we can define Sylow mixed weak neutrosophic N-structure, weak Sylow mixed weak neutrosophic N-structure and Sylow free mixed weak neutrosophic N-structure. In the same way Sylow deficit mixed neutrosophic N-structure and so on can be defined.

We just define the notion of Cauchy neutrosophic element and Cauchy element of a mixed neutrosophic N-structure.

**DEFINITION 6.1.10:** *Let $\langle M \cup I\rangle = \{M_1 \cup M_2 \cup \ldots \cup M_N, *_1, \ldots, *_N\}$ be a mixed neutrosophic N-structure of finite order. We say an element $x \in \langle M \cup I\rangle$ is a Cauchy element if $x^n = 1$ and $n / o(\langle M \cup I\rangle)$. If $y \in \langle M \cup I\rangle$ with $y^m = I$ and $m / o(\langle M \cup I\rangle)$ then we call y a Cauchy neutrosophic element of $\langle M \cup I\rangle$.*

Several interesting properties as in case of other neutrosophic N-structures can be derived for mixed neutrosophic N-structure also.

## Chapter Seven

# PROBLEMS

In this chapter some problems about the neutrosophic structures and their neutrosophic N-structures ($N \geq 2$) are given. It has become essential to mention here that in this book lots about neutrosophic semigroup and neutrosophic loops and their generalizations have been dealt with; we have restrained ourselves from elaborately dealing with neutrosophic groupoids. Further as the class of groupoids contains the class of semigroups an associative structure and also it contains the class of loops we have given a very few definitions about neutrosophic groupoids and their generalization. As our main motivation is to make the reader do problems about neutrosophic groupoids and their generalizations, we have given nearly 25 problems.

We wish to state here throughout the text we have given problems then and there in the text for the reader. Most of them are simple exercises.

1. Can one extend biorder to N-order in case of a neutrosophic N-group?

2. Define N-centre of a neutrosophic N-order and illustrate with examples.

3. Construct a neutrosophic 4-group isomorphism $\phi$ from $\langle G \cup I\rangle$ to itself where $\langle G \cup I\rangle = \{A_4 \cup \{1, 2, 1 + I, 1 + 2I, 2 + I, 2I + 2, I, 2I, 0\} \cup \langle g / g^9 = e\rangle \cup \{1, 2, 3, 4, I, 2I, 3I, 4I\}\}$ for which Ker $\phi$ is a non trivial sub-4-group.

4.
   i. Does the neutrosophic 4-group given in problem (3) have anti Cauchy elements?
   ii. Is it a semi Cauchy neutrosophic 4-group?
   iii. Is it a weakly Cauchy neutrosophic 4-group?
   iv. Does it have Cauchy neutrosophic elements?
   v. Does it have Lagrange sub 4 group?
   vi. Does it have p-Sylow sub 4-group?

5. Is $\langle G \cup I \rangle$ given in problem 3, page 195 a Sylow neutrosophic 4-group? Justify your claim.

6. Can $\langle G \cup I \rangle = \{0, 1, 2, 1 + I, 2 + I, 2I + 1, 2I + 2, I, 2I\} \cup \{S_4\} \cup \{D_{2.5}\}\}$ have normal neutrosophic 3-groups? Find its neutrosophic sub-3 groups? Is this 3-group super Sylow? Justify your claim?

7. Give an example of a neutrosophic 3-group having a (2, 7, 5) Sylow neutrosophic sub-3-group.

8. Give an example of a neutrosophic 4-group having a (3, 5, 7, 11)-Sylow sub-4-group.

9. Find the means to find the number of $(p_1, \ldots, p_N)$- Sylow neutrosophic sub N-group of a neutrosophic N-group $(\langle G \cup I \rangle = \langle G_1 \cup I \rangle \cup \langle G_2 \cup I \rangle \cup \ldots \cup \langle G_N \cup I \rangle, *_1, \ldots, *_N)$.

10. Prove or disprove $(p_1, \ldots, p_N)$-Sylow neutrosophic sub N-groups are conjugate?

11. Find a necessary and sufficient condition for any two $(p_1, \ldots, p_N)$- Sylow neutrosophic sub N-groups to be conjugate?

12. Does their exists a (2, 2, 2, 2, 2) - Sylow sub neutrosophic 5-group?

13. When does the order of every neutrosophic subgroup P of N(G) divide the order of the finite neutrosophic group N(G)? Characteristic them!

14. When does the order of every pseudo neutrosophic (sub) group L of N(G) divide the order of the finite neutrosophic group N(G)? Characterize them.

15. Can conditions be put on N(G) or its neutrosophic subgroup (or pseudo neutrosophic subgroup)so that the partition of $\langle G \cup I \rangle$ is possible by right or left cosets?

16. Give examples of neutrosophic Moufang biloops.

17. Does their exist a Cauchy neutrosophic biloop which is Bruck?

18. Define the following:

    i. S-neutrosophic Bruck biloop
    ii. S-neutrosophic Bol biloop
    iii. S-neutrosophic WIP-biloop
    iv. S-neutrosophic Alternative biloop,

    by giving examples of each.

19. Define a neutrosophic normal sub N-loop. Give examples of them.

20. Define neutrosophic simple N-loop. Give examples.

21. Is $\langle L \cup I \rangle = \{L_1 \cup L_2 \cup L_3, *_1, *_2, *_3\}$ where $L_1 = \langle L_5(3) \cup I \rangle$, $L_2 = S_3$ and $L_3 = \langle g \mid g^6 = 1 \rangle$; a neutrosophic simple 3-loop? Justify your answer.

22. Define Lagrange (N – t) deficit neutrosophic sub-N-group; illustrate this with examples.

23. Define Sylow (N – t) deficit neutrosophic sub N-group and give examples.

24. Define Lagrange (N – t) deficit neutrosophic sub N-groupoid. Give examples of them.

25. Define Lagrange (N – t) deficit neutrosophic sub N-semigroup. Give some examples.

26. Find all the (5 – t), t < 5, Lagrange deficit neutrosophic 5 subloops of the neutrosophic 5-loop. $\langle L \cup I \rangle = (L_1 \cup L_2 \cup L_3 \cup L_4 \cup L_5, *_1, *_2, *_3, *_4)$ where $L_1 = \langle L_5(3) \cup I \rangle$, $L_2 = \langle L_7(4) \cup I \rangle$, $L_3 = S_3$, $L_4 = A_4$ and $L_5 = D_{2.7}$.

27. For the problem (26) find all Sylow (5 – t), (t < 5) deficit neutrosophic 5-sub-loops.

28. Define deficit Cauchy and deficit Cauchy neutrosophic elements of a finite neutrosophic N-Loop.

29. Does their exist a class of neutrosophic groupoids which are inner commutative?

30. Define neutrosophic Moufang groupoids?

31. Is the groupoid $\{\langle Z_7 \cup I \rangle \mid a * b = 3a + 4b \pmod 7,\ a, b \in \langle Z_7 \cup I \rangle\}$ a neutrosophic Moufang groupoid?

32. Define the following and illustrate them with examples.

    i. neutrosophic P-groupoid
    ii. neutrosophic Bol groupoid
    iii. neutrosophic WIP-groupoid
    iv. neutrosophic alternative groupoid.

33. Does the neutrosophic groupoid $N(G) = \{\langle Z_{12} \cup I \rangle \mid a * b = 2a + 4b \pmod{12}$ for $a, b \in \langle Z_{12} \cup I \rangle\}$ fall under any one of the following four classes given in the above problem?

34. Give an example of neutrosophic bigroupoid which is Moufang.

35. Define a left (right) alternative neutrosophic bigroupoid and give an example of neutrosophic bigroupoid not alternative but only left (or right) alternative.

36. Does their exists a neutrosophic bigroupoid which is both WIP and Moufang? Justify your claim!

37. Give an example of Lagrange neutrosophic bigroupoid.

38. Give an example of Sylow neutrosophic bigroupoid.

39. Does their exist a Sylow neutrosophic bigroupoid which is not a Lagrange neutrosophic bigroupoid. Justify your answer!

40. Can one say their can be a relation between Sylow neutrosophic bigroupoid and Lagrange neutrosophic bigroupoid?

41. Give an example of a weak Lagrange neutrosophic bigroupoid.

42. Prove all neutrosophic bigroupoids of prime order are Lagrange free and Sylow free neutrosophic bigroupoids.

43. Define Cauchy elements and Cauchy neutrosophic elements in a finite neutrosophic bigroupoid. Illustrate with examples.

44. Define neutrosophic biideals in a neutrosophic bigroupoid.

45. Does the strong neutrosophic bigroupoid have just neutrosophic biideals; justify your claim.

46. Define neutrosophic N-groupoids which are

    i. Lagrange neutrosophic N-groupoids

ii. Weak Lagrange neutrosophic N-groupoids
iii. Lagrange free neutrosophic N-groupoids

give examples of each.

47. Define neutrosophic normal sub N-groupoid of a neutrosophic N-groupoid. Give an example of each of the

i. neutrosophic N-groupoid having a neutrosophic normal sub N-groupoid.
ii. neutrosophic N-groupoid which has no neutrosophic normal sub N-groupoid. (i.e., simple neutrosophic N-groupoid)

48. Give an example of finite neutrosophic N-groupoid which is Lagrange.

49. Give an example of finite neutrosophic N groupoid of composite order which has no Cauchy element or Cauchy free element.

50. Define neutrosophic N-ideals in a neutrosophic N-groupoid and illustrate with examples.

51. Give an example of a (3, 7, 2, 3, 5) - Sylow neutrosophic sub 5 groupoid of a 5-groupoid.

52. Define (N – t) deficit neutrosophic sub N groupoid and illustrate it with examples.

53. Give an example of (9 – 3) deficit neutrosophic sub 9 groupoid.

# INDEX

## A

## B

## C

**D**

**F**

**G**

**I**

**L**

**M**

**N**

**O**

**P**

**R**

**S**

**W**

**Z**

# ABOUT THE AUTHORS

**Dr.W.B.Vasantha Kandasamy** is an Associate Professor in the Department of Mathematics, Indian Institute of Technology Madras, Chennai, where she lives with her husband Dr.K.Kandasamy and daughters Meena and Kama. Her current interests include Smarandache algebraic structures, fuzzy theory, coding/ communication theory. In the past decade she has guided 11 Ph.D. scholars in the different fields of non-associative algebras, algebraic coding theory, transportation theory, fuzzy groups, and applications of fuzzy theory of the problems faced in chemical industries and cement industries. Currently, four Ph.D. scholars are working under her guidance.

She has to her credit 612 research papers of which 209 are individually authored. Apart from this, she and her students have presented around 329 papers in national and international conferences. She teaches both undergraduate and post-graduate students and has guided over 45 M.Sc. and M.Tech. projects. She has worked in collaboration projects with the Indian Space Research Organization and with the Tamil Nadu State AIDS Control Society. This is her 25th book.

She can be contacted at vasantha@iitm.ac.in
You can visit her work on the web at: http://mat.iitm.ac.in/~wbv

---

**Dr.Florentin Smarandache** is an Associate Professor of Mathematics at the University of New Mexico in USA. He published over 75 books and 100 articles and notes in mathematics, physics, philosophy, psychology, literature, rebus. In mathematics his research is in number theory, non-Euclidean geometry, synthetic geometry, algebraic structures, statistics, neutrosophic logic and set (generalizations of fuzzy logic and set respectively), neutrosophic probability (generalization of classical and imprecise probability). Also, small contributions to nuclear and particle physics, information fusion, neutrosophy (a generalization of dialectics), law of sensations and stimuli, etc.

He can be contacted at smarand@unm.edu